\documentclass[12equipped]{article}
\usepackage{amssymb,amsmath}

\input epsf.tex

\title{Wall-crossing structures in Donaldson-Thomas invariants, integrable systems  and Mirror Symmetry}

\author {Maxim Kontsevich, Yan Soibelman}

\begin{document}
\maketitle

\newcommand{\CC}{{\mathcal C}}
\newcommand{\LL}{{\mathcal L}}
\newcommand{\MM}{{\mathcal M}}
\newcommand{\NN}{{\mathcal N}}
\newcommand{\OO}{{\mathcal O}}
 \newcommand{\ZZ}{{\mathcal Z}}
\renewcommand{\O}{{\mathcal O}}
\newcommand{\Ome}{{\Omega}^{3,0}}
\newcommand{\E}{{\mathcal E}}
\newcommand{\F}{{\mathcal F}}
\newcommand{\TT}{{\mathcal T}}
\newcommand{\g}{{\mathfrak g}}
\newcommand{\GGamma}{{\underline{\Gamma}}}

\renewcommand{\k}{{\bf k}}
\newcommand{\kk}{\overline{\bf k}}

\newcommand{\op}[1]{\operatorname{#1}}

\newtheorem{thm}{Theorem}[subsection]
\newtheorem{defn}[thm]{Definition}
\newtheorem{lmm}[thm]{Lemma}
\newtheorem{rmk}[thm]{Remark}
\newtheorem{prp}[thm]{Proposition}
\newtheorem{conj}[thm]{Conjecture}
\newtheorem{exa}[thm]{Example}
\newtheorem{cor}[thm]{Corollary}
\newtheorem{que}[thm]{Question}
\newtheorem{ack}{Acknowledgements}

\newcommand{\B}{{\bf B}}
\newcommand{\C}{{\bf C}}
\newcommand{\K}{{\bf k}}
\newcommand{\R}{{\bf R}}
\newcommand{\N}{{\bf N}}
\newcommand{\Z}{{\bf Z}}
\newcommand{\Q}{{\bf Q}}
\newcommand{\G}{\Lambda}
\newcommand{\A}{A_{\infty}}
\newcommand{\M}{{\mathsf{M}}}
\newcommand{\epi}{\twoheadrightarrow}
\newcommand{\mono}{\hookrightarrow}
\newcommand\ra{\rightarrow}
\newcommand\uhom{{\underline{Hom}}}
\renewcommand\O{{\cal O}}
\newcommand{\epp}{\varepsilon}





\tableofcontents

\section{Introduction}

\subsection{DT-invariants and Poisson manifolds}

The paper is devoted to the notion of {\it wall-crossing structure} and its constructions and applications in various situations. It is motivated by our previous work on Mirror Symmetry and Donaldson-Thomas invariants (see [KoSo1,2,6,7]) where examples of wall-crossing structures appeared for the first time.

We consider two types of wall-crossing structures in this paper: the one related to the theory of Donaldson-Thomas invariants (DT-invariants for short) and the one related to Mirror Symmetry. Since our main motivation is the former, let discuss it in detail.

It was proposed in [KoSo1], Section 1.5 (see also [KoSo7], Section 7.2) that so-called numerical Donaldson-Thomas invariants counting semistable objects in $3$-dimensional Calabi-Yau categories (DT-invariants of $3CY$ categories for short)  introduced in loc.cit. encode a geometric object, which is a (formal) Poisson manifold. Construction of the Poisson manifold relies on the wall-crossing formulas (WCF for short) introduced in loc.cit. 

Conversely, we suggested there that the DT-invariants  can be recovered from that formal Poisson manifold. Therefore,  collections of DT-invariants satisfying WCF are in one-to-one correspondence with a certain class of formal  Poisson manifolds.

Wall-crossing formulas of this type appear in physics (as samples of the numerous literature on the subject we mention [AlManPerPi], [Den], [DenMo], [GaMoNe 1-4]).

WCF can be naturally realized as identities in the group of (formal) Poisson automorphisms of the Poisson torus naturally associated with the Grothendieck group of the category.  In the categorical framework wall-crossing formulas depend on a choice of stability condition on the category. It was explained in [KoSo1], Section 2 that in fact wall-crossing formulas arise in a more general framework of stability data on graded Lie algebras. The above-discussed case of $3CY$ categories endowed with stability condition corresponds to the case of torus Lie algebras endowed with the stability data given by the DT-invariants of the category.

On the other hand, in our work on Mirror Symmetry (see [KoSo2,6]) a different class of wall-crossing formulas appeared (later it was considered  in [GroSie1,2] in the higher-dimensional case). They can be realized as identities in  the group of (formal) volume-preserving transformations of a complex torus. This type of formulas is related to the counting of pseudo-holomorphic discs and hence does not depend on the stability condition on the relevant Fukaya category.

Similarities between two types of wall-crossing formulas mentioned above lead to a question: is there a structure which makes the formulas similar? The answer is positive. We call it wall-crossing structure (WCS for short) in this paper. 

Besides of the general formalism of WCS we also discuss several new situations in which they appear, most notably, the case of complex integrable systems ``with central charge" (see [KoSo1] and Section 4.2 below). Those include Hitchin integrable systems or Seiberg-Witten integrable systems. To make a link with the above-mentioned categorical version of DT-invariants, we observe that in many cases the (universal cover of the) base of complex integrable system can be thought of as a subspace in the space of stability conditions on some $3CY$ category. In the case of Hitchin integrable system this category is the Fukaya category of the local Calabi-Yau $3$-fold associated with the spectral curve. An interesting fact is that using the geometry of the base of the integrable system we can construct a collection of integers which are similar to DT-invariants (e.g. they satisfy WCF). It is natural to expect that they coincide with the DT-invariants of the above-mentioned Fukaya category.

More precisely, in this paper we are going to discuss  {\it three} different ways to produce collections of integer numbers which enjoy WCF (and subsequently define WCS). 

The first construction of the invariants is based on the count of certain gradient trees on the base of the integrable system (they can be called ``tropical DT-invariants'' because of that).

The second construction extracts DT-invariants from the geometry of the formal neighborhood of a singular curve (``wheel of projective lines'') in a certain algebraic variety. This variety is a compactification of the mirror dual to the total space of the integrable system.

Finally, for a ``good'' non-compact Calabi-Yau $3$-fold one can define DT-invariants of its Fukaya category (we expect they can be defined for any Calabi-Yau $3$-fold). Moreover the moduli space of deformations of the Calabi-Yau $3$-fold is naturally a base of a complex integrable system with central charge. Hypothetically the base can be embedded into the space of stability conditions on the Fukaya category.

We conjecture (see Conjecture 1.2.2) that all three approaches agree in the cases when they all can be applied.

In the case of non-compact Calabi-Yau $3$-folds, the formal Poisson manifold mentioned at the beginning of this subsection is a ``completion at infinity'' of an algebraic Poisson variety. In some examples this variety can be realized as the moduli space of local systems on a punctured curve.

\subsection{Three constructions}
Here we  describe the above-mentioned three approaches (constructions) with more details.

1. For a complex integrable system $\pi:M\to B$ with central charge the following geometry arises on its base $B$.

Denote by $B^0\subset B$ an open dense subset parametrizing non-degenerate fibers of $\pi$.
There is  a local system $\underline{\Gamma}\to B^0$  with fibers $\underline{\Gamma}_b=H_1(\pi^{-1}(b),\Z)$. The local system carries an integer skew-symmetric pairing. In general the pairing can be degenerate, hence $M$ is a Poisson manifold only.

There is a well-defined central charge $Z\in \Gamma(B^0,\underline{\Gamma}^{\vee}\otimes \C)$, where $\underline{\Gamma}^{\vee}$ is the dual local system.
The central charge can be thought of as a local embedding 
$B^0\to \underline{\Gamma}^{\vee}_b\otimes \C$.

Then (under some conditions on $B$) we assign to every generic point $b\in B^{0}$ and every $\gamma\in \underline{\Gamma}_b$  an integer number $\Omega_b^{trop}(\gamma)\in \Z$ which we informally call tropical DT-invariant.
Our construction uses  tropical trees on $B$ with external vertices at the smooth part of $B^{sing}=B-B^0$, as well as wall-crossing formulas  from [KoSo1]. The construction is reminiscent to the attractor flow story in supergravity (see [Den]) recast in mathematical terms in [KoSo7]. Edges of the tropical trees are also gradient lines of the functions on $B$ given by $b\mapsto |Z_b(\gamma)|^2, \gamma\in \GGamma_b$.

2. Here we assume for simplicity that the skew-symmetric form on $\underline{\Gamma}$ is non-degenerate, hence $M$ is a holomorphic symplectic manifold (this assumption is not necessary and will be relaxed in the main body of the paper). Let $\omega^{2,0}$ denote the holomorphic symplectic form. Also assume that the above integrable system is endowed with a holomorphic Lagrangian section $s:B\to M$.
Then we can assign to the above data a filtered associative algebra of finite type over $\C$ (in fact over $\Z$). Roughly speaking, it is the algebra of endomorphisms (in the Fukaya category of $(M,Re(\omega^{2,0})$) of the Lagrangian submanifold $s(B)$ with filtration coming from areas of pseudo-holomorphic discs.  Let $M^{\vee}$ denotes the  affine scheme of finite type, which is the spectrum of this algebra. It can be thought of as a mirror dual to the symplectic manifold $(M,Re(\omega^{2,0}))$. Actual geometric construction goes along the lines of [KoSo2] and uses the corresponding WCF.
Then $M^{\vee}$ is a complex symplectic manifold of the same complex dimension $2n$. Hypothetically, for a reductive group $G$ and the corresponding Hitchin integrable system on a smooth projective curve, the space  $M^{\vee}$ can be identified with the moduli space of $^LG$-local systems on the same curve (Betti realization), where $^LG$ is the Langlands dual group.

We will reconstruct the collection of integer numbers satisfying WCF using the geometry of $M^{\vee}$. Namely, with $M^{\vee}$ one can canonically associate a $\Z PL$-space 
$Sk:=Sk({M^{\vee}})$ called the {\it skeleton} which sits in the Berkovich spectrum of $M^{\vee}$ (see [KoSo2] and Section 6.5. below for the details). Hypothetically, each point $b\in B^0$ gives rise to a piecewise linear embedding $i_b:\C^{\ast}\simeq \R^2-\{0\} \to Sk$.
In the case of Hitchin integrable systems it can be interpreted in terms of the asymptotic behavior of the monodromy of a connection depending on small parameter.

From the point of view of the geometry of $M^{\vee}$ this embedding can be interpreted such as follows.
We have a (partial) Poisson compactification $\overline{M^{\vee}}$ of $M^{\vee}$ by normal crossing divisors and a singular curve $C$ which is a ``wheel" of projective lines $\C{\bf P}^1$ in $\overline{M^{\vee}}-M^{\vee}$, and such that locally near $C$ the space $\overline{M^{\vee}}$ is isomorphic to a toric Poisson variety endowed with a wheel of $1$-dimensional toric strata. The embedding $i_b$ gives rise to an element $Z_b\in H^1(U_{\varepsilon}(C)\cap {M^{\vee}},\Z)\otimes \C$, where $U_{\varepsilon}(C)$ is a small tubular neighborhood of $C$. In other words, we have a linear functional $Z_b:H_1(U_{\varepsilon}(C)\cap {M^{\vee}},\Z)\simeq \Gamma_b\to \C$.

Simple arguments from the deformation theory show that, after a choice of $Z_b$, deformations of the above local toric model are parametrized by collections $(\Omega_b^{MS}(\gamma))_{\gamma\in \Gamma_b-\{0\}}$ satisfying the Support Property from [KoSo1]. Moreover, varying the point $b\in B^0$ we arrive to the collection of numbers satisfying WCF.

3. There is a class of algebraic complex integrable systems with central charge associated with non-compact Calabi-Yau $3$-folds. The base of the integrable system associated with a Calabi-Yau $3$-fold $X$ is (roughly) isomorphic to the moduli space of deformations of $X$. It looks plausible that all Hitchin systems arise in this way (see [DiDoPa] for the $A-D-E$ case). Hypothetically, any point $b\in B^0$ gives rise to a stability condition on the category ${\mathcal F}(X)$, the Fukaya category of $X$. According to the general theory of [KoSo1,5] with  a   stability condition on ${\mathcal F}(X)$ one can associate a collection $\Omega^{cat}_b(\gamma)$ of ``categorical" DT-invariants of ${\mathcal F}(X)$.

\begin{conj} $\Omega_b^{trop}(\gamma)=\Omega_b^{MS}(\gamma)=\Omega_b^{cat}(\gamma)$.

\end{conj}

This conjecture should be the guiding line for the paper.

\begin{rmk} For Hitchin system with the group  $SL(n,\C)$ Gaiotto, Moore and Neitzke proposed  an interpretation of the invariants $\Omega^{cat}_b(\gamma)$ as counting invariants of certain ``networks''  on the spectral curve of the Hitchin system (for $n=2$ they are geodesics of the quadratic differential defined by the point $b\in B^{0}$). The trees in the definition of $\Omega_b^{trop}(\gamma)$ are different from the networks, since they are subsets of the base $B$ rather than of the spectral curve. We do not have a ``counting'' interpretation for the numbers $\Omega_b^{MS}(\gamma)$.

\end{rmk}

The concept of wall-crossing structure (WCS for short) underlying all three constructions  is discussed in the next section, starting with WCS in a vector space. Intuitively, WCS is given by a collection of group elements parametrized by pairs of points outside of codimension one ``walls" and satisfying some consistency conditions.  We observe that WCS can be described in different ways, in particular, as a sheaf of sets.  

Furthermore, WCS is determined by simpler data, which we call {\it initial data} for WCS. In the case of complex integrable systems with central charge discussed later in the paper, the initial data are ultimately related to the behavior of the affine structure on the base near the discriminant set $B^{sing}$.  
The above conjecture can be reformulated as an equivalence of three wall-crossing structures defined in three different ways.

\subsection{Content of the paper}
After the detailed discussion of the Conjecture 1.2.1
let us briefly explain other topics which we discuss in the paper. 

Sections 2,3 are devoted to the concept of wall-crossing structure and examples. We introduce several useful notions like e.g. support of WCS (this concept is related to the Support Property from [KoSo1] which in turn controls the support of DT-invariants). We introduce the notion of attractor flow (the latter goes back to supergravity, see [Den], [DenGrRa], [DenMo]) and define initial data in terms of trees with edges which are trajectories of the attractor flow. 
The initial data can be thought of as a space of ``boundary values" which are assigned to ``free ends'' of attractor trees.

We start Section 4 with a brief discussion of  complex integrable systems from the point of view of Hodge theory. In fact we consider not only polarized integrable systems but semipolarized as well. In the latter case  fibers are semiabelian varieties with polarized quotients. In the case of Hitchin the semipolarized integrable systems appears when the Higgs field has singularities. Then we introduce the notion of a complex integrable system with the central charge. We also explain the construction of WCS and initial data for complex integrable systems with central charge. The initial data are related to the behavior of the integrable system at the discriminant set.

The approach to DT-invariants via wheels of projective lines is the subject of Section 5. The idea which we have already discussed above is that DT-invariants can be interpreted as ``coordinates'' on the moduli space of deformations of the formal neighborhood of a wheel of projective lines in a Poisson toric variety.

The relationship of WCS and SYZ picture of Mirror Symmetry is discussed in Section 6.  Among other things we argue that the mirror dual to the total space of an integrable system with central charge is an affine scheme of finite type over $\Z$. We also stress the role of canonical $B$-field, which is a $2$-torsion. We explain how the set up of Section 5 appears in this framework. Roughly speaking the mirror dual is a log Calabi-Yau whose skeleton is isomorphic to the base of the integrable system. Wheel of lines is related to a choice of central charge.

Deformation theory of non-compact Calabi-Yau threefolds is the subject of Section 7. We discuss there not only the smoothness of the moduli space but also  some plausible assumptions under which the Fukaya category and Donaldson-Thomas invariants are well-defined. Moduli spaces of such deformations serve as bases of the corresponding integrable systems. This gives a generalization of the work [DiDoPa] where the case of local Calabi-Yau $3$-folds associated with ALE spaces was considered. We should also mention the papers [Don], [DonMar] which initiated mathematical works on the relationship between Calabi-Yau $3$-folds and complex integrable systems. In that case the authors considered {\it compact} Calabi-Yau $3$-folds, differently from [DiDoPa]. The corresponding complex integrable systems and their relationship to WCS are discussed in Section 9.

Section 8 is essentially devoted to $GL(r)$ Hitchin integrable systems with possibly irregular singularities (although it also contains other interesting topics like deformation theory of complex Lagrangian manifolds in Section 8.2). In the case of Hitchin integrable systems many structures discussed earlier for general complex integrable systems admit non-trivial interpretations in terms of (irregular) spectral curves. Our approach to the notion of irregular spectral curve is   non-standard (in particular it is quite different from the one in [Bo1]).
Roughly speaking, the spectral curve is defined as an effective divisor in the Poisson surface obtained from the compactified cotangent bundle of the initial curve by a series of blow-ups. The existence of the smooth locus in the base of complex integrable system formed by such spectral curves is related to the existence of a solution to additive Deligne-Simpson problem (see Section 8.3). In that case the general machinery of Section 4 can be applied.

We remark that Section 8.6 contains several interesting conjectures about the mirror dual to the total space of $GL(r)$ Hitchin integrable system which are related to different topics which we do not discuss here. In particular, the conjectures about extension of the (twistor) family over $\C^{\ast}$ of mirror duals to the whole line $\C$  relate those mirror duals to  WKB asymptotics of flat section of connections with a small parameter and to the corresponding theory of resurgent functions.  This relationship has a flavor of ``non-linear Hodge theory of infinite rank" and deserves  further study.

In Section 9 we discuss WCS, attractor flow and DT-invariants in the framework
of compact Calabi-Yau threefolds. As we have already mentioned,  the corresponding complex integrable
systems were studied by Donagi and Markman. They are nonpolarized.
In this framework one still expects the WCS but the initial data 
are determined by the values of DT-invariants not only at conifold points (i.e.
$B^{sing}$ in the above notation) but also by their values at the so-called attractor points. The
value of DT-invariants at a generic conifold point should be equal to $1$. This restriction 
 is not completely clear “from first principles”. The values of DT-invariants at the attractor points are arbitrary integers.

Section 10
is devoted to a version of WCS for the Lie algebra of volume-preserving vector fields on an algebraic torus. This WCS arises naturally
in Mirror Symmetry, in the study of SYZ picture of mirror dual families of collapsing Calabi-Yau manifolds (see [KoSo2], [GroSie1,2]). The formalism of WCS in this case 
 is a bit more complicated
than the one developed in the main body of the paper. Nevertheless, there are several surprising similarities between them. It is natural to suggest that certain (not yet discovered) unified structure is hidden behind.

In  the Appendix we  describe a cocycle with coefficients in $\Z/ 2\Z$ associated with a skew-symmetric form. It is used in the definition of the canonical $B$-field in Section 6.

Finally, we should say that the main goal of our paper is to describe  the general picture of the rich geometry of Wall Crossing Structures. We have tried to formulate (sometimes in the form of conjectures and assumptions) of what {\it should be true}. As a result, besides of proven theorems the paper contains many ideas and new projects. On the other hand, many aspects of the story are not  discussed  (or just touched) in this paper, in particular  quantum versions of the results or the relation to canonical bases in cluster algebras, etc.

\section{Wall-crossing structures}

Wall-crossing formulas presented in [KoSo1] are identities in certain pronilpotent groups of automorphisms. It is convenient to axiomatize the corresponding structure, which appears in a completely different situations. In particular it generalizes the notion of stability data on a graded Lie algebra introduced in the loc.cit.

In this section $\Gamma$ denotes a fixed finitely-generated free abelian group, i.e. $\Gamma\simeq \Z^k$ for some $k\in\Z_{\ge 0}$.
The associated real vector space is $\Gamma_\R:=\Gamma\otimes\R$. We will denote by $\g$  a  fixed $\Gamma$-graded Lie algebra over  $\Q$,
$$\g=\bigoplus_{\gamma\in \Gamma} \g_{\gamma}.$$

\subsection{Wall-crossing structures on a vector space}

\subsubsection{Nilpotent case}
 
Let us assume that the set
$$\op{Supp}{\g}:=\{\gamma\in \Gamma\,|\,\g_\gamma\ne 0\}\subset \Gamma$$
is finite and is contained in an open half-space in $\Gamma_\R$. In particular, all elements of $\op{Supp}{\g}$ are non-zero, i.e. $\g_0=0$.
Under our assumption the Lie algebra $\g$ is nilpotent. Let us denote by $G$ the corresponding nilpotent group. The exponential map $\op{exp}:\g\to G$ is a bijection of sets.

The finite union of hyperplanes $\gamma^\perp \subset \Gamma^*_\R$ (``wall associated with $\gamma$") will be denoted by $\op{Wall}_{\g}$. Its complement has a finite number of connected components which are open convex domains in $\Gamma^*_\R$. These components
are exactly open strata in the natural stratification of $\Gamma^*_\R$ associated with the finite collection of hyperplanes  $\left(\gamma^\perp \right)_{\gamma \in\op{Supp}{\g}}$. 
Notice that different elements $\gamma\in \op{Supp}(\g)$ can give the same hyperplane,
$$\gamma_1^\perp =\gamma_2^\perp\Longleftrightarrow \gamma_1\parallel \gamma_2.$$

\begin{defn} A (global) wall-crossing structure ((global) WCS for  short) for $\g$ is an assignment 
$$(y_1,y_2)\to g_{y_1,y_2}\in G$$ for any $y_1,y_2\in \Gamma_\R^*-\op{Wall}_{\g}$ which 
is locally constant in $y_1,y_2$, satisfies the cocycle condition
 $$g_{y_1,y_2}\cdot g_{y_2,y_3}=g_{y_1,y_3}\,\,\,\forall y_2,y_2,y_3\in \Gamma_\R^*-\op{Wall}_{\g}$$
 and such that in the case when the straight interval connecting $y_1$ and $y_2$ intersects only one of hyperplanes $\gamma^\perp$
 then 
 $$\op{\log}(g_{y_1,y_2})\in \bigoplus_{\gamma':\gamma'\parallel \gamma} \g_{\gamma'}\,.$$

\end{defn}

It follows from the definition that we can associate with any stratum $\tau$ of codimension 1 (which is an open domain in $\gamma^\perp$
for some $\gamma\in\op{Supp}{\g}$) a ``jump'' 
 $$g_\tau:=g_{y_1,y_2}$$
  where points $y_1,y_2$ are such that
 $y_1(\gamma)>0,y_2(\gamma)<0$, and the interval connecting $y_1$ with $y_2$ intersects $\tau$ and no other  strata of codimension $\ge 1$ (hence this interval  does not intersect  other hyperplanes in our collection, except $\gamma^\perp$). 
Obviously, a WCS is uniquely determined by the collection of jumps $\left(g_\tau\right)_{\op{codim}\tau=1}$, satisfying the cocycle condition for each stratum of codimension 2. 
 
\subsubsection{Description in terms of sheaves and groups}

Notice that the complement $\Gamma_\R^*-\op{Wall}_{\g}$ contains two distinguished components $U_+,U_-$ 
(which are different iff $\g\ne 0$) consisting of points $y\in \Gamma_\R^*$ such that
 $y(\gamma)>0$ (resp. $y(\gamma)<0$) for all $\gamma\in \op{Supp}{\g}$.
 Hence with any global WCS $\sigma=(g_{y_1,y_2})$ we can associate an element 
 $$g_{+,-}:=g_{y_+,y_-}\in G\,,\,\,\,y_\pm\in U_\pm.$$
We will prove later in this section (see Theorem 2.1.6) that the map $\sigma\mapsto g_{+,-}$ provides a bijection between the set of wall-crossing structures and $G$ (considered as a set). 

For  any point $y\in \Gamma_\R^*$  we have a decomposition of $\g$ (considered as a vector space) into the
direct sum of three vector spaces
 $$\g=\g_-^{(y)}\oplus \g_0^{(y)}\oplus \g_+^{(y)}$$
corresponding to components $\g_\gamma$ such that $y(\gamma)\in \R$ is negative, zero  or positive respectively.
Obviously all these subspaces are $\Gamma$-graded Lie subalgebras of $\g$. We denote by $G_-^{(y)},G_0^{(y)},G_+^{(y)}$ the corresponding nilpotent subgroups of $G$. Then it is easy to see that the multiplication map
$$G_{-}^{(y)}\times  G_{0}^{(y)}\times G_{+}^{(y)}\to G\,\,\,\,,({g}_-,{g}_0,{g}_+)\mapsto
{g}_-\cdot {g}_0\cdot {g}_+$$
is a bijection. Hence any element ${g}\in G$ can be uniquely decomposed as the product
 $${g}={g}_-^{(y)}{g}_0^{(y)}{g}_+^{(y)}.$$

We denote by $\pi_y:G\to G_0^{(y)}=G_-^{(y)}\backslash G/G_+^{(y)}$ the canonical projection to the double coset. In the above notation we have $\pi_y(\g)=\g_0^{(y)}$. We claim that there exists a sheaf of sets on $\Gamma_{\R}^{\ast}$ with the stalk over $y\in \Gamma_{\R}^{\ast}$ given by $G_0^{(y)}$. 

This is a particular case of the following general construction.
Suppose we are given:

a) a topological space $M$;

b) a set $S$;

c) an assignment to any point $m\in M$ of a set $S_m$ and a surjection $\pi_m: S\to S_m$, such that for any two elements $s_1,s_2\in S$ the set $\{m\in M|\pi_m(s_1)=\pi_m(s_2)\}$ is open in $M$. 

Then the above data give rise to a  sheaf of sets ${\mathcal S}$ on $M$ in the following way. Its \'etal\'e space
 ${\mathcal S}^{\mbox{\'et}}$ consists of pairs $\{(m,s')|m\in M, s'\in S_m\}$. A  base of topology is given by the sets $W_{s,U}=\{(m,s')|m\in U,s'=\pi_m(s)\}$, where $s$ runs through the set $S$ and $U$ runs through the set of open subsets of $M$.

One can easily prove the following result.

\begin{lmm} The projection ${\mathcal S}^{\mbox{\'et}}\to M, (m,s')\mapsto m$ is a local homeomorphism.

\end{lmm}
Using the standard equivalence between sheaves of sets and local homeomorphisms (\'etale maps), we obtain

\begin{cor} The above construction gives rise to a sheaf of sets ${\mathcal S}$ such that the stalk ${\mathcal S}_m$ is equal to $S_m$ for any $m\in M$.

\end{cor}

Let us apply this Lemma to the case $M=\Gamma_{\R}^{\ast}, S=G, S_y=G_0^{(y)}, y\in M$ and the map $\pi_y: G\to G_0^{(y)}=G_{-}^{(y)}\backslash G/G_{+}^{(y)}, g\mapsto g_0^{(y)}$ given by the canonical projection to the double coset. It is easy to see that the openness condition from c) is satisfied.

\begin{defn} The corresponding sheaf of sets is called the sheaf of wall-crossing structures and is denoted by $WCS_{\g}$.

\end{defn}

Sheaf $WCS_\g$ is constructible with respect to the natural stratification of $\Gamma_{\R}^{\ast}$ given by  the finite arrangement of hyperplanes $\gamma^\perp\subset \Gamma_\R^*$, where $\gamma\in \op{Supp}{\g}$.

 Notice that if $y\in \Gamma_{\R}^{\ast}-\op{Wall}_{\g}$ then the stalk at $y$ is $G_0^{(y)}=\{1\}$. If the point $y$ belongs to a stratum of codimension one (i.e. $y$ lies on exactly one wall $\gamma^\perp$) then the Lie algebra of the corresponding stalk is $G_0^{(y)}$ where
$Lie(G_0^{(y)})=\oplus_{\gamma^{\prime}\parallel \gamma}\g_{\gamma^{\prime}}$. Finally, the stalk at $y=0$ is the whole group $G$.

It will be important for the future  to study the space of sections of $WCS_\g$ in the following situation. Let
$l\subset \Gamma^*_\R$ be a straight line intersecting $U_+$ and $U_-$. We endow $l$  with the direction from $U_+$ to $U_-$ and require that it does not intersect strata of codimension bigger or equal than $2$. Let $y_1,...,y_n$ be the ordered along $l$ set of intersection points
with walls. Then the set of sections $\Gamma(l,WCS_\g)$ is just the product $\prod_{i=1}^n G^{(y_i)}_0$.
 The natural map
 $$\bigoplus_{i=1}^n \g^{(y_i)}_0\to \g$$
 is a bijection, hence any element $g\in G$ can be uniquely decomposed into the ordered product of elements of $G^{(y_i)}_0$.
 We conclude that the set of sections $\Gamma(l,WCS_\g)$ can be identified naturally with $G$.

Let us consider the following three sets:

i) $S_1=G$;

ii) $S_2=\Gamma(\Gamma_{\R}^{\ast}, WCS_{\g})$;

iii) $S_3$ being the set of all wall-crossing structures on $\Gamma_{\R}^{\ast}$.

There are three maps
$S_1\to S_2\to S_3\to S_1$ given such as follows:

1) the map $S_1\to S_2$ sends $g\in G$ to the section $s_g$ such that $s_g(y)=g_0^{(y)}$;

2) the map $S_2\to S_3$ sends a section $s\in \Gamma(\Gamma_{\R}^{\ast}, WCS_{\g})$ to the unique WCS such that for any straight interval connecting two points $y_1,y_2$, intersecting a hyperplane $\gamma^\perp$ at one point $y_0$ and not intersecting other hyperplanes, and such that 
$$y(\gamma_1)<0, \,y_2(\gamma)>0$$
the transformation $g_{y_1,y_2}$ coincides with $s(y_0)\in G^{(y_0)}_0\subset G$,

3) the map $S_3\to S_1$ sends a WCS $\sigma$ to the corresponding element $g_{+,-}:=g_{+,-}(\sigma)$. 

It is clear that the maps $S_1\to S_2$ and $S_3\to S_1$ are well-defined. It is not obvious that the map $S_2\to S_3$ indeed takes values in the set of wall-crossing structures.

\begin{prp} The map $S_2\to S_3$ is well-defined.

\end{prp}

{\it Proof.} Our description of the map defines jumps $g_\tau$ for all strata of codimension one.
 We need to check the cocycle condition in codimension $2$.
Let   $\rho$ be any stratum of codimension 2. It lies in a finite collection of $k\ge 2$ different hyperplanes
 $(\gamma_i^\perp)_{i=1,\dots,k}$. Near $\rho$ we have $2k$ open strata and $2k$ strata of codimension 1. Among $2k$ open  strata we have two distinguished ones containing points $y_1,y_2$ respectively, where 
 $$y_1(\gamma_i)<0,y_2(\gamma_i)>0\,\,\,\,\,\forall i=1,\dots,k\,.$$
  There are two paths (up to homotopy) connecting $y_1$ and $y_2$ contained in the union of $2k$ open and $2k$ codimension one strata near $\tau$ and intersecting  each of the hyperplanes $(\gamma_i^\perp)_{i=1,\dots,k}$ at one point. 
    We want to prove that the composition of jumps along one path coincides with the similar composition along another one.
    It follows from the definition of the sheaf $WCS_\g$ that both compositions coincide with 
$g_0^{(y)}$ where $y$ is any point in $\rho$. Hence the cocycle condition is satisfied. $\blacksquare$

\begin{thm} The above three maps are bijections, and their composition is the identity map.

\end{thm}

We see that sets $S_1,S_2,S_3$ are canonically identified with each other.

{\it Proof.} We split the proof into three lemmas.

\begin{lmm} The map $S_1\to S_2$ is a bijection.

\end{lmm}

This map is obviously injective, because the  $s_g(0)=g$ for any $g\in G$. It is surjective because the point $0\in \Gamma_{\R}^{\ast}$ belongs to the closure of any stratum (which is a conical set). Therefore any section is uniquely determined by its value at $0$. This proves  first Lemma.

\begin{lmm} The composition $S_1\to S_2\to S_3\to S_1$ is the identity map.

\end{lmm}

Given $g\in G$ let us choose a line $l\subset \Gamma_{\R}^{\ast}$ such that it intersects both sets $U_+$ and $U_-$ and does not intersect strata of codimension greater or equal than $2$. Let us endow the line with the direction from $U_+$ to $U_-$. Let us denote by $y_1,...,y_m$ the ordered points of intersection of $l$ with  walls. Then $g$ coincides with the ordered product of $g_0^{(y_i)}$
 and hence ${g}=g_{+,-}$. This proves the second Lemma.

\begin{lmm} The map $S_3\to S_1$ is injective.

\end{lmm}

Observe that for any point $y\in \Gamma_{\R}^{\ast}$ which belongs to the stratum of codimension $1$ (which is open in a  wall) there exists a line $l$ as in the previous Lemma and such that $y\in l$. Because of our choice of the line, there exists a unique $i_0, 1\le i_0\le m$ that $y=y_{i_0}$. We know that the element $g_{+,-}$ determines uniquely all elements $g_0^{(y_i)}$, in particular $g_0^{(y_{i_0})}=g_0^{(y)}$.  We conclude that all jumps are determined uniquely by the element $g_{+,-}$. This means that the corresponding WCS is also determined uniquely. This proves the third Lemma. Combined together, the three lemmas give the proof of the Theorem. $\blacksquare$.

\begin{rmk} The set of sections of the sheaf ${WCS}_{\g}$ on any open $U\subset \Gamma_\R^{\ast}$ can be described
 as the set of locally constant maps from the set of connected components of intersections of codimension one strata with $U$ to corresponding subgroups of $G$, which satisfy the cocycle condition near points of strata of codimension two. Also, let $U\subset \Gamma_{\R}^{\ast}$ be an open convex subset. We define cones $C_{\pm}(U)\subset \Gamma_{\R}$ as cones generated by $\gamma\in \Gamma$ such that $\pm y(\gamma)>0$ for all $y\in U$. We denote by $G_{\pm}(U)=exp(\oplus_{\gamma\in C_{\pm}(U)}\g_{\gamma})$ the corresponding nilpotent Lie groups.
Then $\Gamma(U,\op{WCS}_{\g})\simeq G_-(U)\backslash G/G_+(U)$.
There is also a description of  the set $\Gamma(U,\op{WCS}_{\g})$ similar to the Definition 2.1.1. Namely, in the Definition 2.1.1 we consider pairs $y_1,y_2\in U-\op{Wall}_{\g}$.

\end{rmk}

 \subsubsection{Pronilpotent case}
 
 Let  $\g=\oplus_{\gamma \in \Gamma}\g_{\gamma}$ be a graded Lie algebra. We do not impose any restrictions on $\op{Supp}(\g)$. In particular we do not assume that the support belongs to a half-space in $\Gamma_{\R}$ (cf. Section 2.1.1).

 Let $C\subset \Gamma_\R:=\Gamma\otimes \R$ be a  convex cone. We assume that $C$ is {\it strict}, which means that the closure of $C$ does not contain a line, or, equivalently, $C$ is contained in the positive octant (in some coordinates on $\Gamma_\R$). Yet another equivalent condition:
 there exists $\phi\in  \Gamma_{\R}^{\ast}$ such that the restriction of $\phi$ to the cone $C$ is a proper map to $\R_{\ge 0}$.

 In this case we define a pronilpotent Lie algebra $\g_C$ as an infinite product
 $$\g_C:=\prod_{\gamma\in C\cap \Gamma -\{0\}} \g_\gamma$$
 and denote by $G_C$ the corresponding pronilpotent group. The exponential map identifies $\g_C$ and $G_C$.

Lie algebra $\g_C$ is the projective limit of nilpotent Lie algebras
$$\g_{C,\phi}^{(k)}=\bigoplus_{\gamma \in C\cap \Gamma-\{0\},\phi(\gamma)\le k}\g_{\gamma}=\g_{C}/m_{C,\phi}^{(k)},$$
where 
$$m_{C,\phi}^{(k)}=\bigoplus_{\gamma \in C\cap \Gamma,\phi(\gamma)> k}\g_{\gamma}$$
is the Lie ideal in $\g_C$, and $\phi\in  \Gamma_{\R}^{\ast}$  is the above function.

Denote by $G_{C,\phi}^{(k)}=exp(\g_{C,\phi}^{(k)})$ the corresponding nilpotent group and by $pr_{C,\phi}^{(k)}:G_C\to G_{C,\phi}^{(k)}$ the natural epimorphism of groups.

Then the sheaf of sets $WCS_{\g_C}$ is defined as the projective limit of the sheaves $WCS_{\g_{C\phi}^{(k)}}$. It follows from the end of the Remark 2.1.10 that for any open convex subset $U\in \Gamma_{\R}^{\ast}$ the set of sections $WCS_{\g_C}(U)$ admits the following description:

a) For any $y_1,y_2\in U$ which do not belong to $(\cup_{\gamma\in C\cap (\Gamma-\{0\})}\gamma^{\perp})\cap U$ we are given an element $g_{y_1,y_2}\in G_C$ satisfying the cocycle condition.

b)  The projections of these elements to  $G_{C,\phi}^{(k)}(U)$ satisfy the condition from Definition 2.1.1.

The latter condition informally means that the ``jump"   at the generic point of the hyperplane $H=\gamma^{\perp}\subset \Gamma_{\R}^{\ast}$, where $\gamma\in C-\{0\}$ belongs to the subgroup $G_{C,H}:=exp(\prod_{\gamma^{\prime}\in C\cap \Gamma, (\gamma^{\prime})^{\perp}=H}\g_{\gamma^{\prime}})$.

In the next definition we extend this picture to the case when the cone $C$ is not fixed in advance and can depend on a point of $\Gamma_{\R}^{\ast}$.

\begin{defn} Let $\g=\oplus_{\gamma \in \Gamma}\g_{\gamma}$ be the graded Lie algebra, as before. We define the sheaf $WCS_{\g}$ on $\Gamma_{\R}^{\ast}$ such as follows: for any open subset $U$ the set of sections $\Gamma(U,WCS_{\g})$ consists of a family of elements $g(y,\gamma)\in \g_{\gamma}$ such that $y\in U, \gamma\in \Gamma-\{0\}$ and  $y(\gamma)=0$ satisfying the following condition:

For any $y\in U$ there exists a neighborhood $U_y\subset U$ and strict  convex cone $C_{y,U_y}\subset \Gamma_{\R}$ such that for any $y_1\in U_y$ the element $g(y_1,\gamma)\ne 0$ iff $\gamma\ne 0$ and $\gamma\in C_{y,U_y}$.

Furthermore, let us fix $\phi\in  \Gamma_{\R}^{\ast}$ such that its restriction to the closure $\overline{C_{y,U_y}}$ is a proper map to $\R_{\ge 0}$. Then we require that for any $k>0$ the map  
$$y_1\mapsto pr_{C_{y,U_y}}^{(k)}(exp(\sum_{\gamma}g(y_1,\gamma))), y_1\in U_y$$
 is an element of the set of sections  $\Gamma(U_y,WCS_{\g^{(k)}_{C_{y,U_y},\phi}})$.

\end{defn}

Notice that if $\op{Supp}(\g)$ is finite and contained in an open half-space in $\Gamma_{\R}$ then the above definition agrees with the one given in Section 2.1.2.
We also have the following pronilpotent analog of the Theorem 2.1.6.

\begin{prp} Assume that $\op{Supp}(\g)\subset C-\{0\}$, where $C$ is the cone described at the beginning of this subsection. Then the 
set of global sections $\sigma=\left(g(y,\gamma)\right)$ of $ WCS_{\g}$
is in the natural one-to-one correspondence with elements of $g\in {G}_C$.

\end{prp}

{\it Proof.} Follows from the nilpotent case. $\blacksquare$

Having a section $s\in \Gamma(U,WCS_{\g})$ for an open $U\subset \Gamma^{\ast}_{\R}$ we define its support $Supp(s)\subset U\times \Gamma_{\R}$ as a minimal closed, conic in the direction of $\Gamma_{\R}$ set which contains the set of pairs $(y,\gamma), y\in \Gamma_{\R}^{\ast},\gamma\in \Gamma$ such that $y(\gamma)=0$ and $log(g_0^{(y)})_{\gamma}\in \g_{\gamma}-\{0\}$.

\subsection{Wall-crossing structure on a topological space}

Let us consider the following data:

1) A Hausdorff locally connected topological space $M$ (then we will speak about WCS on $M$).

2) A local system of finitely-generated free abelian groups of finite rank $\pi:\underline{\Gamma}\to M$.

3) A local system of $\underline{\Gamma}$-graded Lie algebras $\underline{\g}=\oplus_{\gamma\in \underline{\Gamma}} \underline{\g}_{\gamma}\to M$ over the field $\Q$.

4) A homomorphism of sheaves of abelian groups $Y:\underline{\Gamma}\to \underline{Cont}_{M}$, where $\underline{Cont}_{M}$ is the sheaf of real-valued continuous functions on $M$.

Equivalently we can interpret $Y$ locally as a continuous map from a domain in $M$ to $\Gamma_{\R}^{\ast}$. Then we define the pull-back sheaf $WCS_{\underline{\g},Y}:=Y^{\ast}(WCS_{\g})$, where $WCS_{\g}$ is the sheaf of sets on $\Gamma_{\R}^{\ast}$ constructed in the previous subsection.

\begin{defn} A (global) wall-crossing structure on $M$ is a global section of $WCS_{\g,Y}$. The support of WCS $\sigma$  is a closed subset of $tot(\underline{\Gamma}\otimes \R)$  whose fiber over any point $m\in M$ is described such as follows:
it is a strict convex closed cone $Supp_{m,\sigma}\subset \underline{\Gamma}_m\otimes \R$ which is equal to the support of the germ of
$WCS_{\underline {\g},m}$ at the point $Y(m)\in \underline{\Gamma}_m^{\ast}\otimes \R$ associated with the section $\sigma$.
\end{defn}

This definition makes obvious the  functoriality of the notion of WCS with respect to pullbacks.

\subsection{Examples of wall-crossing structures}

1) Let us fix a free abelian group of finite rank $\Gamma$ together with a $\Gamma$-graded Lie algebra $\g=\oplus_{\gamma \in \Gamma}\g_{\gamma}$ and a homomorphism of abelian groups $Z:\Gamma\to \C$ (central charge). Then we take  $M=\R/2\pi \Z$, and define on $M$ constant local systems with fibers $\Gamma$ and $\g$. We set $Y_{\theta}(\gamma)=Im(e^{-i\theta}Z(\gamma))$, where $\theta\in \R$. Then a  WCS associated with this choice is the same as  stability data on $\g$ in the sense of [KoSo1]. Family of stability data from [KoSo1] parametrized by the topological space $M$ is the same as WCS on $\R/2\pi \Z\times Hom(\Gamma,\C)$ with constant local systems $\underline{\Gamma}, \underline{\g}$ and the above map $Y$.

2) Let us fix $\Gamma, \g, Z\in Hom(\Gamma,\C)$ as in Example 1.
Let $\widehat{M}=\R/2\pi \Z\times Hom(\Gamma,\C)$. We endow  $\widehat{M}$ with constant local systems with fibers $\Gamma$ and $\g$ respectively. 

The  subset $M_Z=\R/2\pi \Z\times \{Z\}\subset \widehat{M}$ is isomorphic to the one from the previous example. Then interpreting WCS on $M_Z$ as the pullback sheaf $(Y_{|M_Z})^{\ast}({WCS_{\g}})$ we conclude that this WCS can be extended to a neighborhood of $Z$. 

Thus we have a WCS for nearby central charges. Using compactness of the circle $\R/2\pi \Z$ we conclude that for any stability data on $\g$ with the central charge $Z_0$ there exists a germ of universal family of stability data with central charges in a neighborhood of $Z_0$.  The above construction gives an alternative to [KoSo1] way to define the notion of continuous family of stability data on a graded Lie algebra. As a byproduct we can interpret wall-crossing formulas from [KoSo1] in terms of WCS.

3) Assume that in the example 1) we have an involution $\eta:\g\to \g$ which maps $\g_{\gamma}\to \g_{-\gamma}$. Let us choose a local system on $M=\R/2\pi\Z$ obtained from the trivial local systems from Example 1) by identification $\theta\mapsto \theta+\pi$ on $\R$, $\gamma\mapsto -\gamma$ on $\Gamma$ and $x\mapsto \eta(x)$ on $\g$. The corresponding WCS can be identified with symmetric stability data on $\g$ from [KoSo1].

4) Assume that $\Gamma$ is endowed with an integer skew-symmetric form $\langle \bullet,\bullet\rangle$. Let us fix central charge $Z$ and set $\g=\oplus_{\gamma\in \Gamma}\Q\cdot e_{\gamma}$, where 
$$[e_{\gamma_1},e_{\gamma_2}]=(-1)^{\langle \gamma_1,\gamma_2\rangle}\langle \gamma_1,\gamma_2\rangle e_{\gamma_1+\gamma_2}.$$
We will call it the {\it torus Lie algebra}.
All previous examples can be specified to this case.

In this case one can encode the WCS as a collection of numbers $\Omega(\gamma)$ which are called DT-invariants for the torus Lie algebras and the central charge $Z$. The relation with Definition 2.1.11 is as follows:
$$g(e^{i\theta},\gamma)=\sum_{k|\gamma, k\ge 1}{\Omega(\gamma/k)\over{k^2}}e_{\gamma},$$
where $Z(\gamma)\in e^{i\theta}\R_{>0}$ 

5) The quantum version of the previous example deals with the Lie algebra
$\g=\oplus_{\gamma\in \Gamma}\Q(q^{1/2})\cdot \widehat{e}_{\gamma}$ where
$$[\widehat{e}_{\gamma_1},\widehat{e}_{\gamma_2}]={q^{\langle \gamma_1,\gamma_2\rangle/2}- q^{-\langle \gamma_1,\gamma_2\rangle/2}\over{q^{1/2}-q^{-1/2}}}\widehat{e}_{\gamma_1+\gamma_2}.$$
Here $\widehat{e}_{\gamma}={\widehat{e}_{\gamma}^{quant}\over{q^{1/2}-q^{-1/2}}}$ are the normalized generators of the quantum torus $\widehat{e}_{\gamma_1}^{quant}\widehat{e}_{\gamma_2}^{quant}=
q^{\langle \gamma_1,\gamma_2\rangle/2}\widehat{e}_{\gamma_1+\gamma_2}^{quant}$.
We will call it the {\it quantum torus Lie algebra}.

\begin{rmk} In the previous two examples each group $G_U$ contains a subgroup $G_U^{adm}$ of admissible (or quantum admissible in the example 4)) series (see [KoSo5] for the definition). Then in the definition of WCS we can require that $g_{m_1,m_2,U}\in G_U^{adm}$. This leads to the integrality  of the corresponding DT-invariants.

\end{rmk}

6) Let $\CC$ be an ind-constructible $3CY$ category with the class map $cl:K_0(\CC)\to \Gamma$ (see [KoSo1] for the terminology and notation).  Let  $G\subset Aut(\CC,cl)$ be a subgroup which preserves  a connected component $Stab_0(\CC,cl)\subset Stab(\CC,cl)$ of the set of constructible stability conditions on $\CC$.  In particular $G$ acts on $\Gamma$. We require that $G$ acts on $Stab_0(\CC,cl)$ freely and properly discontinuously, and also that $G$ contains the group $\Z$ generated by the shift functor $[1]$. Consider $M:=Stab_0(\CC,cl)$.
Notice that there exists a WCS on $M$ with the constant local systems $\Gamma$ and $\g$ given in the Examples 4), 5) and $Y=Im(e^{-i\theta}Z)$. Then this WCS descends to $M/G$. In practice $G$ is a stabilizer of $Stab_0(\CC,cl)$ in $Aut(\CC,cl)$. 

Mirror Symmetry predicts that the moduli space of Calabi-Yau $3$-folds endowed with a square of a holomorphic volume form carries a local system of $3CY$ categories (Fukaya categories) endowed with stability condition. Thus we expect that there exists the corresponding WCS on an appropriate topological space (see Sections 7.3, 9.1 about that). The group $G$ in this case is the fundamental group of the above moduli space.

7) Let $(Q,W)$ be a quiver with polynomial potential (more generally, we can consider a smooth algebra with potential).
Then we constructed in [KoSo5] an invertible series $A$ (called quantum DT-series) in the quantum torus with generators $e_{\gamma}, \gamma\in \Gamma_+\otimes \R\simeq C:=\R_{\ge 0}^n$, where $\Gamma_+\simeq \Z_{\ge 0}^n$ is the cone of dimension vectors. By the Theorem 2.1.6 it defines a WCS on $\Gamma_{\R}^{\ast}$. 

Let $Z:\Gamma_+\to \C$ be a central charge. We assume that $Y=Im(Z)$ is positive on $C-\{0\}$. Consider the straight line $l\subset \Gamma_{\R}^{\ast}$ given by $t\mapsto Re(Z)+tIm(Z)\in \Gamma_{\R}^{\ast}$. Intersections of this line with the walls (i.e. points $t\in \R$ for which there exists $\gamma\in \Gamma_+-\{0\}$ such that $Re(Z(\gamma))+tIm(Z(\gamma))=0$) correspond to rays $\alpha$ in the upper-half plane with  vertex in the origin. The clockwise factorization formula $A=\prod_{\alpha}^{\longrightarrow}A_{\alpha}$ (see [KoSo1,5]) can be interpreted as the previously discussed product formula for the element $g_{+,-}\in G_C$.

8) In many examples it is natural to consider $\Gamma$-graded Lie algebras where $\Gamma$ is a finitely generated abelian group, possibly with torsion.
For example, for the Fukaya category ${\mathcal F}(X)$ of a Calabi-Yau $3$-fold $X$ one should take $\Gamma=H_3(X,\Z)$, which can have a non-zero torsion.

In such cases the considerations of the previous and this sections still work, if we replace $\Gamma$ by $\Gamma^{free}$ which is the quotient of $\Gamma$ by the torsion subgroup $\Gamma^{tors}$. Then having the $\Gamma$-graded Lie algebra $\g$ we define the graded Lie algebra $\g_{\Gamma^{free}}=\oplus_{\mu\in \Gamma^{free}}\g_{\mu}^{free}$ where $\g_{\mu}^{free}=\oplus_{\gamma\,mod\,\Gamma^{tors}=\mu} \g_{\gamma}$.

\section{Initial data and attractor flow}

Definitions and construction of this section, which might look a bit unmotivated, will be used below in Sections 4.5 and 9.3 when we will discuss the set of ``initial values" which determine DT-invariants associated with complex integrable systems.

\subsection{Initial data}

Let $\Gamma$ be a free abelian group of finite rank endowed with a skew-symmetric integer form $\langle\bullet,\bullet\rangle: \bigwedge^2\Gamma\to \Z$. Let $\g:=\gamma_{\Gamma}=\oplus_{\gamma\in \Gamma}\g_{\gamma}$ be a $\Gamma$-graded Lie algebra over a commutative ring  of characteristic zero satisfying the condition $[\g_{\gamma_1},\g_{\gamma_2}]=0$ as soon as $\langle \gamma_1,\gamma_2\rangle=0$. E.g. the Lie algebras from examples 4), 5) in the previous subsection satisfy this condition. We denote by $\Gamma_0$ the kernel of $\langle\bullet,\bullet\rangle$. Then one has a decomposition $\g=\g_{\Gamma_0}\oplus\g_{\Gamma-\Gamma_0}$, where $\g_{\Gamma_0}=\oplus_{\gamma\in \Gamma_0}\g_{\gamma}$ is a central subalgebra and $\g_{\Gamma-\Gamma_0}=\oplus_{\gamma\in \Gamma-\Gamma_0}\g_{\gamma}$ is its complement. The skew-symmetric form $\langle\bullet,\bullet\rangle$ gives rise to a homomorphism of abelian groups $\iota:\Gamma\to \Gamma^{\vee}$. Then $\Gamma_{0}=Ker\,\iota$. Since a WCS in the abelian case of  $\g_{\Gamma_0}$ is something very simple (just a collection of elements of $\g_{\gamma} $ with support in a strict cone, we are going to discuss only the ``non-trivial" part which is the induced WCS for $\g_{\Gamma-\Gamma_0}$. In what follows we assume that if $\gamma\in \Gamma_0$ then $\g_{\gamma}=0$.


Let us consider the ``global" version of the above situation.
Namely, assume we are given a WCS, say, $\sigma$, on a smooth manifold  $B^0$ such that the corresponding local systems $\underline{\Gamma}, \underline{\g}$ satisfy the same properties as $\Gamma, g$ above.
More precisely, $\underline{\Gamma}$ is endowed with an integer skew-symmetric pairing  $\langle\bullet,\bullet\rangle$ such that if for $b\in B^0$ we have $\langle \gamma_1,\gamma_2\rangle=0, \gamma_i\in \underline{\Gamma}_b$ then for the corresponding components of $\g_b$ we have $[\g_{b,\gamma_1},\g_{b,\gamma_2}]=0$. Let $\underline{\Gamma}_0\subset \underline{\Gamma}$ be the kernel of  $\langle\bullet,\bullet\rangle$. Then assume that if $\gamma\in \Gamma_{0,b}$ then $\g_{b,\gamma}=0, b\in B^0$.

Now we will make an additional assumption on $B^0$ which will be justified later in the case of complex integrable systems (strange notation for $B^0$ is also borrowed from there). Namely, we assume 
that the homomorphism of sheaves $Y$ interpreted  locally as a continuous map from $B^0$ to a real vector space, is
a smooth submersion.
Moreover, we assume that $B^0$ is endowed with a foliation such that locally near $b\in B^0$ the leaf $M:=M_b$ containing $b$ is identified via $Y$ with
an affine space   over the vector space $\iota(\underline{\Gamma}_{\R,b})$.

For each leaf $M$ let us define a smooth manifold $M^{\prime}_{\Z}$ as the set of pairs $(m,\gamma), m\in M, \gamma\in \underline{\Gamma}_m$ which satisfy the condition that $\gamma\in {\underline{\Gamma}}_m-\underline{{\Gamma}}_{0,m} $ and such that $ Y(m)(\gamma)=0$.
Notice that $dim\, M^{\prime}_{\Z}=dim\,M-1$. 

We define a bigger set $M^{\prime}\supset M^{\prime}_{\Z}$ as the set of pairs 
$$\{(m,v)\in tot (\underline{\Gamma}_{\R})|Y(m)(v)=0\}.$$
Clearly $dim\,M^{\prime}=dim\,M+rk\,\Gamma-1$.
\begin{defn} The attractor flow on $M^{\prime}$ is defined by the vector field
$\dot{v}=0, \dot{m}=\iota(v), v\in \GGamma_{m,\R}$. It preserves the  manifold $M^{\prime}_{\Z}$ and hence induces the ``integer" attractor flow on it.

\end{defn}

By our assumptions the vector field does not vanish on $M^{\prime}_{\Z}$.

Similarly we define $(B^0)^{\prime}$ and $(B^0)^{\prime}_{\Z}$ as the union of above-defined sets over all leaves $M$. The attractor flow extends to the both bigger manifolds.

Assume that we are given a WCS $\sigma$ on $B^0$. Therefore we have a piecewise constant map $a: (B^0)^{\prime}_{\Z}\to tot(\underline{\Gamma})$ which assigns to a point $(b,\gamma)$ the element $a_b(\gamma)\in \underline{\g}_{b,\gamma}$, which is the $\gamma$-component of the corresponding section of $WCS_{\g,Y}$ (we identify naturally the latter  with the logarithm of the corresponding element of the pronilpotent group, see Proposition 2.1.12).

The discontinuity set $\op{W}_a$ of the map $a$ (this set is an analog of ``walls of first kind" from [KoSo1]) belongs to the set of pairs $(b,\gamma)\in Supp_{\sigma}\cap tot(\underline{\Gamma})-tot(\underline{\Gamma}_0)$ such that $\gamma=\gamma_1+\gamma_2$, with  $\langle\gamma_1,\gamma_2\rangle \ne 0$ (this condition implies that all vectors $\gamma,\gamma_1,\gamma_2$ do not belong to $\underline{\Gamma}_{0,b}$, while $(b,\gamma_i)\in Supp_{\sigma}, i=1,2$). Clearly the discontinuity set of the map $a$ is a locally-finite  hypersurface in  $(B^0)^{\prime}_{\Z}$ which is locally a pull-back of a $\Z PL$ hypersurface in $\underline{\Gamma}_b^{\vee}\otimes \R$.

The set $Supp_{\sigma}$ can be very complicated. For example in the case of $SL_2$ Hitchin integrable system considered in [GaMoNe2] the fibers of $Supp_{\sigma}$ should coincide with cones of invariant measures for singular foliations on surfaces given by real parts of certain quadratic differentials. 

We believe that in general the support of WCS has a ``fractal" structure similar to the following one. Let us consider the set of points $(x,y)\in \R^2$ such that either $x\in \R-\Q, y=0$ or $x={p\over{q}}\in \Q, (p,q)=1$ and $0\le y\le {1\over{q}}$. This is a closed subset in $\R^2$, and fibers of the projection $(x,y)\mapsto x$ are compact convex sets.

Although $Supp_{\sigma}$ could be complicated, it is natural to expect that one can find an ``upper bound" on it, which is a closed subset $C^+\subset (B^0)^{\prime}\subset tot(\underline{\Gamma}_{\R})$ satisfying the following properties:

i) Fibers of $C^+$ under the natural projection  $tot(\underline{\Gamma}_{\R})\to B^0$ are strict convex cones;

ii)  The set $C^+$ is preserved by the ``inverse attractor flow"  $\dot{v}=0, \dot{b}=-\iota(v), v\in \GGamma_{b,\R}, t>0$.



Suppose we know $C^+$ by some a priori (e.g. geometric) reasons. Then it gives us an ``upper bound" for the discontinuity set $\op{W}_a$.  More precisely, let us define  $\op{Wall}^+\subset C^+$ as the set  of pairs $(b,\gamma)$ such that $\gamma=\gamma_1+\gamma_2$ where $\gamma_1,\gamma_2\in C^+\cap (B^0)^{\prime}_{\Z}$ such that $\langle \gamma_1,\gamma_2\rangle \ne 0$ (in particular it follows that $\gamma_1$ and $\gamma_2$ are non-parallel vectors). Thus if $Supp_{\sigma}\subset C^+$ then $\op{W}_a\subset \op{Wall}^+$.

\begin{prp} The attractor flow is transversal to $\op{Wall}^+$ at the interior of the latter set.

\end{prp}
{\it Proof.}  It suffices to check the statement on a fixed leaf $M$. We can work locally and assume that $M$ is an affine space, $\Gamma$ is a fixed lattice endowed with an integer skew-symmetric form. Let fix $\gamma\in \Gamma$ and consider the attractor flow $y\mapsto y+t\iota(\gamma), t>0$. On the discontinuity variety we have $\gamma=\gamma_1+\gamma_2$ and $(y+t\iota(\gamma))(\gamma_m)=0, m=1,2$. Then at the intersection point of the attractor flow with the discontinuity subvariety  we have $y(\gamma_1)+t\langle \gamma,\gamma_1\rangle=0$ for some (maybe many) $t>0$. Since $\langle\gamma,\gamma_1\rangle=-\langle\gamma_1,\gamma_2\rangle\ne 0$ we conclude that $t$ is determined uniquely. This implies the transversality and finishes the proof. $\blacksquare$

Notice that for each point $(b,\gamma)\in (B^0)^{\prime}$ there exists  maximal possible $t_{max}:=t_{max}(b,\gamma)\in (0,+\infty]$ such that the trajectory
$\gamma=const, b\mapsto b+\iota(\gamma)t, t\in [0,t_{max})$ does exist. The above trajectory considered for $t\in [0,t_{max})$ will be called {\it the maximal positive trajectory of the point $(b,\gamma)$}.

Let us define an open subset ${\mathcal T}_{(B^0)^{\prime}_{\Z}}\subset (B^0)^{\prime}_{\Z}$ called the {\it tail set} consisting of points $(b,\gamma)\in C^+$ such  that their maximal positive trajectories with respect to the attractor flow do not intersect the set $\op{Wall}^+$, belong to $C^+$, and moreover the same properties hold for all nearby points $(b^{\prime},\gamma^{\prime})\in (B^0)^{\prime}$. There is a local system $\underline{\g}_{loc}$ over ${\mathcal T}_{(B^0)^{\prime}_{\Z}}$ with the fiber $\underline{\g}_{b,\gamma}$ over $(b,\gamma)\in {\mathcal T}_{(B^0)^{\prime}_{\Z}}$.

\vspace{2mm}
{\bf Tail Assumption} {\it For any open $U\subset  (B^0)^{\prime}$   the subset of points $(b, \gamma)\in U$ such that their maximal positive trajectories intersect the tail set ${\mathcal T}_{(B^0)^{\prime}_{\Z}}$, is dense in $U$.} 

\vspace{2mm}

\begin{rmk} Typically $\underline{\g}_{loc}$ is trivial of rank one. E.g. in Example 4) from the previous subsection the fiber is $\Q$, while in Example 5) it is $\Q(q^{1/2})$.

\end{rmk}

\begin{defn} The initial data of a WCS bounded by $C^+$ is the restriction of the map $a$ to ${\mathcal T}_{(B^0)^{\prime}_{\Z}}$.

\end{defn}

As we will explain below, under some additional conditions the initial data uniquely determine its WCS. Moreover, in some cases one can reconstruct WCS just from the knowledge of initial data. This explains the meaning of this notion.

\subsection{Attractor trees}

In what follows we are going to consider  metrized rooted trees with finitely many edges oriented toward tails. Here  ``metrized tree" understood as a length metric space.
Internal vertices have outcoming valency at least $2$, internal edges have finite length, while tail edges can be infinite.  Our convention is  that the root vertex has valency $1$.

Let us assume the notation of the previous subsection.

\begin{defn} An attractor tree is a metrized rooted tree $T$ endowed with a continuous map $f: T\to M$ to a leaf $M\subset B^0$ and a lift $f^{\prime}: T-\{Vertices\}\to M^{\prime}_{\Z}$.
We assume that $f^{\prime}$ maps  edges of $T$ to trajectories of the attractor flow, and the metric on each edge of $T$ is given by $|dt|$, where $t$ is the time parameter for attractor flow on its lifting. We assume that all tail edges are maximal positive trajectories of the corresponding internal vertices of $T$. 
We assume
the balancing condition $\sum_i \gamma_i^{out}=\gamma^{in}$ is satisfied at each internal vertex $v$. Here $\gamma^{in}$ is the speed of the $f^{\prime}$-lift of the only edge incoming from $v$, and $\gamma_i^{out}$ are speeds of the $f^{\prime}$-lifts of all outcoming edges. We assume that all $\gamma_i^{out}$ are pairwise distinct and there exist $i_1,i_2$ such that $\langle \gamma_{i_1}^{out}, \gamma_{i_2}^{out}\rangle\ne 0$. 
\end{defn}

To an attractor tree $T$ with a root $b$ and root edge $\gamma$ we can assign its  combinatorial type in the following way.
Namely, let us consider an abstract rooted tree ${\mathcal T}$ corresponding to $T$ and a collection of velocities of all its edges, including tails. The velocities can be treated as elements of $\GGamma_b$ via the parallel transport along edges of $T$. Then the {\it combinatorial type of $T$ at $(b,\gamma)$} consists of the above abstract tree and the above subset of velocities in $\GGamma_b$. Varying $(b,\gamma)$ we conclude that combinatorial types form a local system over $(B^0)^{\prime}_{\Z}$ with countable fibers. 

It is easy to see that for any combinatorial type at $(b,\gamma)$ an attractor tree with this combinatorial type is uniquely determined by the collection $\{l_e\}$ of lengths of its inner edges $e$. Moreover the lengths $l_e$ and the vector $Y(b)\in \GGamma_{b,\R}^{\ast}$ satisfy a system of linear equations with integer coefficients arising from the following two conditions:

 a)  $Y(f(v))(f^{\prime}(u))=0$, where $v$ is a vertex of $T$ and $u$ is a point on an edge adjacent to $v$ and sufficiently close to $v$;
 
 b) For any inner edge $e$ connecting vertices $v_1$ and $v_2$ we have $Y(f(v_2))-Y(f(v_1))=\iota(f^{\prime}(u))l_e$, where $u$ is any point of $e$.

For a fixed attractor tree $T$ with the root $b$ and root edge $\gamma$ let us consider the set of attractor trees with sufficiently close roots, combinatorial types and lengths. This set (which can be thought of as a germ of the universal deformation of $T$) can be identified with an open domain in the vector subspace of the vector space $\GGamma_{b,\R}^{\ast}\oplus \R^{\{inner\,\,edges\}}$ defined by the above system of linear equations. In particular, for any vertex $v$ of $T$ the point $Y(f(v))$ runs through an open domain in a vector subspace $H_v\subset \GGamma_{b,\R}^{\ast}$ defined over $\Q$.
In particular, we see that the set of roots of attractor trees which are close to $T$ and have the same combinatorial type is locally an open domain in a vector subspace of $\GGamma_{b,\R}^{\ast}$. 

\begin{defn} We say that attractor tree $T$ is locally planar (this property will depend on its combinatorial type only) if for each internal vertex $v$ of $T$ the corresponding vectors $\gamma^{out}_i$ span a $2$-dimensional vector subspace in $\GGamma_v\otimes \R$.

\end{defn}

\begin{prp} For an attractor tree $T$ with the root $b$ and the root edge $\gamma$
the set of roots of all sufficiently close attractor trees of the same combinatorial type has codimension $\ge 1$ if $T$ is locally planar and has  codimension $\ge 2$ otherwise (this codimension is the same as $codim(H_b)$ in the above notation). Moreover in the former case any sufficiently close attractor tree is uniquely determined by its root.

\end{prp}

{\it Proof.} Let us call a vertex $v$ of $T$ non-planar if the vector $\gamma_i^{out}$ outcoming from $v$ span a vector space of dimension $\ge 3$. Let us prove the second part of the Proposition. For that let us assume that $T$ contains a non-planar vertex $v_0$. 
Then the set $Y(f(v_0))$  belongs to a vector subspace  $H_{v_0}\subset \GGamma_{b,\R}^{\ast}$ of codimension $\ge 3$ defined over $\Q$ because of the conditions $Y(f(v_0))(\gamma_i^{out})=0$.  Consider the shortest path $v_0\leftarrow v_1\leftarrow....\leftarrow v_n=b$ of vertices of $T$ joined by edges. We will prove by induction that for all $0\le i\le n-1$ we have $codim(H_{v_i})\ge 3$. Then $codim(H_b)=codim(H_{v_n})\ge 2$. 

The induction step is given by the following lemma.

\begin{lmm} Consider a germ of the universal deformation of a given attractor tree $T$. For any edge  $e: w_2\to w_1$ connecting two internal vertices of a variable tree $T_s$ we have the following: if $codim(H_{w_1})\ge 3$ then $codim(H_{w_2})\ge 3$.

\end{lmm}

{\it Proof of the Lemma.} Let $\gamma=f^{\prime}(u)\in \GGamma_b$ denote the velocity of the edge $e$. Then $Y(f(w_2))=Y(f(w_1))-l_e\iota(\gamma)$, where $l_e$ is the length of $e$.  It follows that $Y(f(w_2))$ belongs to a vector subspace of $H_{w_1}+\R\cdot \iota(\gamma)$. If the latter subspace has codim $\ge 3$ we are done. Assume that it has codimension $2$. Then there exist two linearly independent vector $\mu_1,\mu_2\in  \GGamma_b$ such that this vector subspace is equal to $\mu_1^{\perp}\cap \mu_2^{\perp}$. Let us denote by $\{\gamma_j^{out}\}$ the set of velocities of edges outgoing from $w_2$. Obviously $\gamma$ belongs to this set. The covector $Y(f(w_2))$ is orthogonal to all vectors $\mu_1,\mu_2, \{\gamma_j^{out}\}$. If the vector subspace in $\GGamma_{b,\R}$ generated by the latter set has dimension $\ge 3$ the we are done. Hence we can assume that it has dimension $2$.   Then $\gamma$ and all $\gamma_j^{out}$ are linear combinations of $\mu_1$ and $\mu_2$.  Since $Y(f(w_2))\in \op{Wall}^+$ we know that there exist $\gamma_{j_1}^{out},\gamma_{j_2}^{out}$ such that $\langle\gamma_{j_1}^{out},\gamma_{j_2}^{out}\rangle\ne 0$. It follows that $\langle\mu_1,\mu_2\rangle\ne 0$. Thus we have a $2$-dimensional vector space generated by linearly independent vectors $\mu_1,\mu_2$ and the vector $\gamma$ in this vector space such that
$\langle\mu_i,\gamma\rangle=0, i=1,2$ (the latter follows from the fact that $\mu_i(\iota(\gamma))=0, i=1,2$). We conclude that $\gamma=0$.
This proves the Lemma and the second part of the Proposition.

In order to prove first part we assume that all vertices are locally planar.  And then we again proceed by induction by the number of vertices. The statement is obvious for the tree which has only one vertex (root vertex) and one edge (the root edge which coincides with the tail edge).  Assume that $T$ is planar and contains at least one internal vertex. Let us choose a vertex $v$ such that the only  edges outcoming from $v$ are tails edges. Let us denote by $\gamma^{in}, (\gamma_i^{out})$ the velocities of edges attached to  $v$. Let us denote by $T^{\prime}$ the tree obtained from $T$ by deleting the vertex $v$ and all outcoming tail edges, and extending the incoming edge by maximal positive trajectory $e$ of the attractor flow with velocity $\gamma^{in}$.  In this way we obtain a map from the germ of the universal deformation of $T$ to the one of $T^{\prime}$.
We claim that this is a local homeomorphism. Indeed, the vertex $Y(f(v))$ belongs to  the codimension $2$ vector subspace of $\Gamma_{\R}^{\ast}$ given by $\cap_i(\gamma_i^{out})^{\perp}$. Hence the point $Y(f(v))-l\iota(\gamma^{in})$ (here $l$ is the length of the incoming edge) varies in the open domain of the hyperplane $(\gamma^{in})^{\perp}$. Therefore the point $Y(f(v))$ is (locally) uniquely determined by the tree $T^{\prime}$ as the intersection point of the trajectory $Y(f(e))\subset (\gamma^{in})^{\perp}$ with $\cap_i(\gamma_i^{out})^{\perp}$.

Notice that $T^{\prime}$ has one less vertex than $T$. Continuing by induction we reduce the problem to the case of one root vertex and one edge. This completes the proof of Proposition.
$\blacksquare$

From now on we assume that we are given an upper bound $C^+$ as in the previous subsection. 

\begin{defn}
 We say that the attractor tree is bound by $C^+$ if the image of $f^{\prime}$ belongs to $C^+$ (then one can easily see that for any internal vertex $v$ the $f^{\prime}$-lift of the only outcoming edge starts on $\op{Wall}^+$).
\end{defn}

It follows from the properties of $C^+$ that the attractor tree is bound by $C^+$ if an only if  its tail edges are bound by $C^+$.

Every attractor tree $T$ has finitely many tail edges which are invariant with respect to the (positive) attractor flow. Let us denote by $T^0$ the tree obtained by deleting all tail edges. Every edge of such a tree joins two vertices.

\vspace{2mm}
{\bf Compactness Assumption} 

{\it There exists an open dense subset 
$(B^0)^{\prime\prime}_{\Z}\subset  (B^0)^{\prime}_{\Z}$ with the following property: 
for every $(b,\gamma)\in (B^0)^{\prime\prime}_{\Z}$ there exists a compact subset $K_{(b,\gamma)}\subset (B^0)^{\prime}$ and an open neighborhood $U$ of $(b,\gamma)$ such that for every attractor tree $T$ with the root and root edge in $U$ the corresponding tree $T^0$ belongs to $K_{(b,\gamma)}$.}

\vspace{2mm}

{\bf Mass Function Assumption} 

{\it There exists a morphism of sheaves $X: \underline{\Gamma}\to \underline{Cont}_{(B^0)}$ such that the pull-back of $X$ to $(B^0)^{\prime}$
considered as a continuous function in $(b,v)\in tot(\underline{\Gamma}_{\R})$  decreases (non-strictly) along the attractor flow $v=const, \dot{b}=i(v)$ as long as $(b,v)\in C^+$ and is strictly positive on the set  $C^+-tot(\underline{\Gamma}_{0,\R})$ and strictly decreasing along the flow on this set.}

\vspace{2mm}

Imposing the above three assumptions (Tail, Compactness and Mass), let us consider the graph $G:=G(b,\gamma)$ obtained as the union of all attractor trees with the root at a fixed $(b,\gamma)\in (B^0)^{\prime\prime}_{\Z}$. Then there are finitely many attractor trees which form the graph and that the obtained graph is acyclic.  Assume that the root $b$ runs through the set of generic points satisfying the condition $Y(b)\in \gamma^{\perp}$. The Proposition 3.2.3 implies that 
$G(b,\gamma)$ is locally planar (with the obvious generalization of the Definition 3.2.2 to graphs). The genericity here means that $Y(b)$ belongs to the complement of the locally finite union of codimension $\ge 2$ subspaces of $\Gamma_{\R}^{\ast}$.

\begin{prp} WCS with fixed $\underline{\g}$ and the support belonging to $C_{\Z}^+:=C^+\cap (B^0)^{\prime}_{\Z}$ is uniquely determined by its initial data.

\end{prp}

{\it Proof.} 

Fix a generic point $(b,\gamma)\in C^+_{\Z} $. Let us consider the maximal acyclic graph $G$ described above. All its tails belong to ${\mathcal T}_{(B^0)^{\prime}_{\Z}}$ (otherwise we can enlarge the graph). Then we reconstruct the value  $a(b,\gamma)$ by induction, starting with the  restriction of the function $a$ to ${\mathcal T}_{(B^0)^{\prime}_{\Z}}$ (initial data) and moving toward the point $(b,\gamma)$ along the edges of $G$. Since $G$ is acyclic, for any  internal vertex $(b^{\prime},\gamma^{\prime})$ we can uniquely compute $a(b^{\prime},\gamma^{\prime})$ from the axioms of WCS in $2$-dimensional case. Finally we compute $a(b,\gamma)$ by induction. The Proposition is proved. $\blacksquare$

\begin{rmk} Since the function $a$ is locally-constant we can reconstruct WCS from the knowledge of $a$ on a dense open subset of $(B^0)^{\prime}_{\Z}$. We do not claim that the procedure given in the proof of Proposition 3.2.6 produces the data $a(m,\gamma)$ which correspond to a WCS. The reason for that is that the procedure in the proof ensures that the cocycle condition is satisfied for some (but possibly not all) strata of codimension $2$ (see the end of Section 2.1.1). We need more geometric conditions in order to be sure that all strata of codimension $2$ are taken into account.

\end{rmk}

\subsection{Initial data for WCS in a vector space}

We assume the set up and the notation of the beginning of the previous subsection, i.e. we have a fixed  lattice $\Gamma$ with a fixed integer skew-symmetric form $\langle\bullet,\bullet\rangle$, a fixed $\Gamma$-graded Lie algebra $\g$, etc. We also fix a closed strict convex cone $C\subset \Gamma_{\R}$. We define $B^0= \Gamma_{\R}^{\ast}$ and $Y=id$. We use the Poisson structure on $\Gamma_{\R}^{\ast}$ induced by $\langle\bullet,\bullet\rangle$. 
We set $C^+=\{(b,v)\in (B^0)^{\prime}|b\in B^0,v\in C, b(v)=0\}$.

Then the Tail Assumption is  satisfied. Indeed for any $\gamma\in C\cap (\Gamma-\Gamma_0)$ there are only finitely many $\gamma_1, \gamma_2\in C\cap (\Gamma-\Gamma_0)$ such that $\gamma=\gamma_1+\gamma_2, \langle\gamma_1,\gamma_2\rangle\ne 0$. Then the set of points$\{b|(b,\gamma)\in \op{Wall}^+\}$ is a finite union of $(rk\,\Gamma-2)$-dimensional hyperplanes in $\gamma^{\perp}$, such that all of them are transversal to the flow $\dot{b}=\iota(\gamma)$.

Therefore for sufficiently large times the attractor flow does not intersect $\op{Wall}^+$. This implies the following

\begin{cor} $\pi_0({\mathcal T}_{(B^0)^{\prime}_{\Z}})\simeq \{\gamma\in C\cap (\Gamma-\Gamma_0)\}$.

\end{cor}

In other words for any $\gamma\in \Gamma-\Gamma_0$ we have a unique connected component of the (integer) tail set, which contains ``sufficiently large'' parts of rays in the direction of $\gamma$.

The Compactness Assumption follows from the finiteness of the set of combinatorial types of attractor trees with given $(b,\gamma)$ and such that velocities of all edges belong to $C$.

The Mass Function Assumption is more tricky. In order to construct ``mass function" $X$ let us choose coordinates $y_1,...,y_{2n}, t_1,...,t_m$ in $\Gamma_{\R}$ such that $(y_i)$ are symplectic coordinates and $(t_i)$ are coordinates on the center $\Gamma_{0,\R}$. Let us denote by $x_1,...,x_{2n}, s_1,..., s_m$ the dual coordinates on $\Gamma_{\R}^{\ast}$. Let us also choose a bounded strictly increasing smooth function $f:\R\to \R$ (e.g. $f(x)=arctan(x)$). 

Then we define 

$$X(b,\gamma)=\sum_{i_1,i_2}f(x_{i_1})y_{i_2}\omega^{i_1i_2}+L(\gamma).$$

Here $b=(x_1,...,x_{2n}, s_1,...,s_m), \gamma=(y_1,...,y_{2n}, t_1,...,t_m)$ and $L\in \Gamma_{\R}^{\ast}$ is a covector independent on $b$, $(\omega^{ij})$ is the symplectic form on $\Gamma^{symp}$.

The condition $\dot{X}(b,\gamma)>0$ if $\dot{b}=\iota(\gamma)$ is satisfied for $\gamma\in \Gamma-\Gamma_0$ and arbitrary $L$. Indeed, 
$\dot{b}$ is given $\dot{x}_{i_1}=\sum_{i_2}y_{i_2}\omega^{i_1i_2},\dot{s}_j=0$. It follows that $\dot{X}(b,\gamma)$ is strictly positive. Moreover for sufficiently large $L\in C^{\vee}$ we will have $X(b,\gamma)>0$ for any $\gamma\in C-\{0\}$.

In the Remark 3.2.7 we warned the reader that for any initial data there exists at most one corresponding WCS, but its existence is not guaranteed in general. It follows from the Proposition 3.3.2 below that in the case of WCS in a vector space there is no problem with the existence.

From the decomposition $\g=\g_{\Gamma_0}\oplus \g_{\Gamma-\Gamma_0}$ we obtain a similar decomposition $\g_C:=\oplus_{\gamma\in C\cap \Gamma}\g_{\gamma}=\g_{\Gamma_0,C}\oplus\g_{\Gamma-\Gamma_0,C}$. As we know a WCS on $B^0=\Gamma_{\R}^{\ast}$ with the support in $C^+$ is the same as an element of the pronilpotent group $G_{\Gamma-\Gamma_0,C}=exp(\g_{\Gamma-\Gamma_0,C})$. The initial data is given by the restriction of the map $a$ to the tail set, which defines an embedding $\psi: G_{\Gamma-\Gamma_0,C}\to \prod_{\gamma\in C\cap (\Gamma-\Gamma_0)}\g_{\gamma}$.

\begin{prp} $\psi$ is a bijection of sets.

\end{prp}

{\it Proof.} Let us choose an additive map $\eta:\Gamma\to \Z$ such that $\eta(C-\{0\})\subset \R_{>0}$.  Then $G_{\Gamma-\Gamma_0,C}=\varprojlim_kG_{C,k}$, where each nilpotent group $G_k:=G_{C,k}$ is defined similarly to $G_{C,\Gamma-\Gamma_0}$ by taking the exponent of the Lie algebra which is the quotient by the ideal generated by $\g_{\gamma}$ with $\eta(\gamma)>k, \gamma\in (\Gamma-\Gamma_0)\cap C$. Then we have a sequence of (compatible with respect to the index $k$) maps $\psi_k:G_k\to \prod_{\gamma\in C\cap(\Gamma-\Gamma_0),\eta(\gamma)\le k}\g_{\gamma}$. We will prove that $\psi$ is a bijection by induction in $k$. For $k=0$ there is nothing to prove. Assume that $\psi_{k-1}$ is a bijection. Notice that the fiber of the natural projection $G_k\to G_{k-1}$ is a torsor over the abelian Lie group 
$exp(\prod_{\gamma\in C\cap(\Gamma-\Gamma_0),\eta(\gamma)=k}\g_{\gamma}$). Similarly the fiber of the natural projection $\prod_{\gamma\in C\cap(\Gamma-\Gamma_0),\eta(\gamma)\le k}\g_{\gamma}\to \prod_{\gamma\in C\cap(\Gamma-\Gamma_0),\eta(\gamma)\le k-1}\g_{\gamma}$ is a torsor over the $\Q$-vector space $\prod_{\gamma\in C\cap(\Gamma-\Gamma_0),\eta(\gamma)=k}\g_{\gamma}$. For a fixed point $g_{k-1}\in G_{k-1}$ we have an isomorphism of torsors.  This implies that $\psi_k$ is a bijection. $\blacksquare$

\begin{rmk} Last part of the proof of the above Proposition is very transparent in the language of the factorization $G=G_-^{(y)}G_0^{(y)}G_+^{(y)}$ from the previous Section: multiplication of the left factor by an element from $Centr(G_0^{(y)})$ is equivalent to the multiplication of the factor from 
$G_0^{(y)}$ by this element.

\end{rmk}

The above considerations imply the following alternative description of the initial data for WCS in a vector space. For a $\gamma\in \Gamma-\Gamma_0$ consider the Lie subalgebra $\g^{\gamma}=\oplus_{\gamma^{\prime}\in C\cap (\Gamma-\Gamma_0),\langle \gamma^{\prime},\gamma\rangle=0}\g_{\gamma^{\prime}}$. Then we have a homomorphism of Lie algebras $pr$ from $\g^{\gamma}$ onto the abelian Lie algebra $\g^{\gamma}_{ab}=\oplus_{\gamma^{\prime}\in C\cap (\Gamma-\Gamma_0),\gamma^{\prime}\parallel \gamma}\g_{\gamma^{\prime}}$. The restriction $a(\gamma)$ of $a(y,\gamma)$ to the tail set is given by the $\gamma$-component of the element $log(\,g_0^{(\iota(\gamma))})$. Indeed the element $g_0^{(y)}$ stabilizes for sufficiently large $t$ as long as we follow the attractor flow $y\mapsto y+t\iota(\gamma)$.

In the framework of Section 2.3, Example 4) let us fix an isomorphism $\Gamma\simeq \Z^I$. Then consider the initial data  given by: $a_{init}(\gamma)=0, \gamma\notin \Z_{\ge 1}e_i$ and $a_{init}(ke_i)={1\over{k^2}}$ otherwise. Here $(e_i)_{i\in I}$ is the standard basis in $\Z^I$. By the above Proposition 3.2.6 we have a unique WCS with these initial data and support in the cone $C=\R_{\ge 0}^I$.

Recall that for a lattice with a basis and an integer skew-symmetric form we can construct a quiver $Q$ with the set of of vertices $I$ and the number of arrows $i\to j$ equal to $\langle e_i,e_j\rangle$. Let us  choose a generic potential $W$ for this quiver.

\begin{conj} The group element $g=1+...$ corresponding to the above WCS coincides with the DT-series from [KoSo1].

\end{conj}

Similarly, if we replace above ${1\over{k^2}}$ by ${q^{1/2}-q^{-1/2}\over{k(q^{k/2}-q^{-k/2})}}$ and take the Lie algebra from Example 5) (quantum torus) then conjecturally we obtain the quantum DT-series from [KoSo1] , Sect. 8 (up to multiplication by a central  series in the generators $\widehat{e}_{\gamma}, \gamma\in C\cap \Gamma_0-\{0\}$, see [KoSo1] for the notation). Recall that it is related to the theory of (quantum) cluster varieties. Finally, we remark that the canonical group element $g$ corresponding to the WCS derived from a $3CY$ category in the formal way by means of the transformation between positive and negative chambers does not have to coincide with the motivic DT-series of that category (but  the latter defines {\it some} WCS and the corresponding canonical group element). The comparison is sometimes possible in case if the category has a ``good" set of generators (e.g. $3$-dimensional spherical generators) which can serve as ``initial data" in the categorical framework.

\section{Geometry of complex integrable systems}

Complex integrable system is usually understood as a holomorphic generically surjective map $\pi:(X,\omega^{2,0})\to B$ of a complex analytic symplectic manifold of dimension $2n$ to a complex analytic manifold of dimension $n$ such that generic fibers are holomorphic Lagrangian submanifolds. A generic fiber is acted locally transitively by an abelian Lie algebra of dimension $n$. In many interesting situations generic fibers are Zariski open subsets in complex abelian varieties. We can compactify these fibers, and obtain a fibration by abelian varieties over an open dense subset $B^0\subset B$. This will be the situation discussed in Section 4.1.1. In Section 4.1.2 we will generalize the story to semiabelian case.
  The behavior of an integrable system near the discriminant locus $B^{sing}=B-B^0$ is more complicated, although in the generic point of the discriminant one can find an explicit local model (see Section 4.6). Our philosophy is that all the information necessary for the construction of a wall-crossing structure (which is our principal goal) is already encoded in the geometry of $B^0$ 
(see e.g. Completeness Assumption in Section 4.4).
Hence in what follows we will use a slightly nonstandard terminology. Complex integrable systems in the usual sense recalled above we will call {\it full complex integrable systems}. Hence a full integrable system can have  e.g. singular fibers. Complex integrable systems in our sense has semiabelian fibers such that  the first integer homology of  fibers form a local system of lattices.

\subsection{Integrable systems and variations of Hodge structure}

\subsubsection{Case of pure Hodge structure}

Let $(X^0,\omega^{2,0})$ be a complex  analytic symplectic manifold of complex dimension $2n$.\footnote{Many of the results below can be generalized to the case of smooth algebraic varieties.}
Assume we are given  an holomorphic map $\pi:X^0\to B^0$ such that for any $b\in B^0$ the fiber $\pi^{-1}(b)$ is a complex Lagrangian submanifold of $X^0$, which is in fact a torsor over an  abelian variety endowed with a covariantly constant integer polarization. We will call such data a {\it (polarized) complex integrable system}.

Let $\GGamma$ be a local system of free abelian groups over $B^0$ with a fiber $\GGamma_b:=H_1(\pi^{-1}(b),\Z), b\in B^0$. Polarization gives rise to a covariantly constant 
skew-symmetric bilinear form $\langle\bullet,\bullet\rangle:\bigwedge^2\GGamma\to \underline{\Z}_{B^0}$ which induces a covariantly constant symplectic form on $\GGamma\otimes \Q$. 
In this case we will speak about local system of symplectic lattices.

Then map $\gamma\mapsto \int_{\gamma}\omega^{2,0}$ gives rise to a morphism of sheaves of abelian groups
$$\alpha:=\int \omega^{2,0}: \GGamma\to \Omega^{1,cl}_{B^0,hol},$$
where $\Omega^{1,cl}_{B^0,hol}$ denotes the sheaf of holomorphic closed $1$-forms on $B^0$.

Let $U\subset B^0$ be a simply connected domain. Let us choose a basis $(\gamma_1,...,\gamma_{2n})$ of $\GGamma(U)$. Then the homomorphism $\alpha$ gives rise to a collection of holomorphic closed $1$-forms $\alpha_i=\int_{\gamma_i}\omega^{2,0}, 1\le i\le 2n$ which can be written on $U$ as $\alpha_i=dz_i, 1\le i\le 2n, z_i\in {\mathcal O}(U)$. The collection of functions $(z_1,...,z_{2n})$ defines a holomorphic map $Z:U\to \C^{2n}$.
Let $\omega_{ij}=\langle\gamma_i,\gamma_j\rangle$ and $(\omega^{ij})_{i,j}\in Mat(2n,\Q)$ be the matrix inverse to $(\omega_{ij})_{i,j}$. It is easy to see that:

1) $\sum_{i,j}\omega^{ij}dz_i\wedge dz_j=\sum_{i,j}\omega^{ij}\alpha_i\wedge\alpha_j=0$.

2) the $(1,1)$-form $\sqrt{-1}\sum_{i,j}\omega^{ij}dz_i\wedge d\overline{z}_j=\sqrt{-1}\sum_{i,j}\omega^{ij}\alpha_i\wedge\overline{\alpha}_j$ is positive (hence it defines a K\"ahler metric on $U$).

It follows from 2) that $Z$ is an immersion. It follows from 1) that $Z(U)$ is a Lagrangian submanifold.

Recall the well-known fact that  near each point $b\in B^0$ the structure of polarized integrable system is determined by the triple $(\GGamma,\langle\bullet,\bullet\rangle,\alpha)$ satisfying 1) and 2).

Indeed, suppose we are given
a symplectic lattice $(\Gamma, \langle\bullet,\bullet\rangle)$ and a holomorphic Lagrangian  embedding of a neighborhood $U$ of of $b$ to $\Gamma^{\vee}\otimes \C$ defined up to a shift, such that for any $b_1\in U$ and non-zero $v\in T_{b_1}B^0$ we have 
$Im\langle dZ_{b_1}(v),\overline{dZ_{b_1}(v)}\rangle>0$.

Then the polarized integrable system $\pi: (\pi^{-1}(U),\omega^{2,0})\to U$ is isomorphic (as a polarized integrable system) to the ``canonical local model'' which is the polarized integrable system with the fiber over $b\in U$ given by $\Gamma\backslash (\Gamma\otimes \C)/(T_{Z(b)}Z(U))^\perp$ endowed with an obvious symplectic form and polarization (we are going discuss the symplectic form in a more general case below in Section 4.1.2).

Alternatively the local model is given by the quotient of $T^{\ast}U$ by the action of $\Gamma$ given by $(b,v)\mapsto (b,v+\alpha_b(\gamma)), b\in U,\gamma\in \Gamma, v\in T_b^{\ast}U$.

\begin{rmk} Locally, on the total space of the polarized integrable system one has an action of a real compact torus $\Gamma\backslash(\Gamma\otimes \R)$ by holomorphic symplectomorphisms preserving fibers, and each fiber is a torsor over this torus.

\end{rmk}

Notice that the local model for a polarized integrable system is endowed with a holomorphic Lagrangian section (zero section). The isomorphism between our integrable system and the local model is not unique. It is determined by a choice of holomorphic Lagrangian section over $U$. Gluing together local models we obtain a new polarized integrable system $\pi^{\prime}:X^{\prime}\to B^0$ with a Lagrangian section $B^0\to X^{\prime}$, which is canonically associated with the triple $(\GGamma,\langle\bullet,\bullet\rangle,\alpha)$ satisfying 1) and 2). For any $b\in B^0$ the fiber $(\pi^{\prime})^{-1}(b)$ has a structure of abelian group (zero is given by the Lagrangian section). The fiber $\pi^{-1}(b)$ is a torsor over $(\pi^{\prime})^{-1}(b)$. Isomorphism classes of polarized integrable systems with fixed $(\GGamma,\langle\bullet,\bullet\rangle,\alpha)$ are in one-to-one correspondence with elements of the group $H^1(B^0,\Omega^{1,cl}_{B^0}/\alpha(\GGamma))$. 

The data $(\GGamma,\langle\bullet,\bullet\rangle,\alpha)$ satisfying the conditions 1), 2) are equivalent to the following data:

i) A variation of polarized pure Hodge structure on $B^0$ of weight $-1$ given by $(R^1\pi_{\ast}\underline{\Z}_{X^0})^{\vee}$. The Hodge filtration is given by 
$$0=F^{-2}\subset F^{-1}\simeq (R^1\pi_{\ast}({\mathcal O}_{X^0}))^{\ast}\subset F^0=\GGamma\otimes {\mathcal O}_{B^0}.$$
ii) An isomorphism of vector bundles $\Psi:{T}_{B^0}\to (F^0/F^{-1})^{\ast}\simeq F^{-1}$. This isomorphism satisfies certain conditions which can be derived from the conditions 1), 2).

\begin{rmk}
 We can consider complex integrable systems with fibers which are compact complex tori without polarization. In this case the local model is determined by a submanifold $Z(U)\subset \Gamma^{\vee}\otimes \C$ such that $dim\,U={1\over{2}}rk\,\Gamma$ and for any $b\in U$ we have $T_{Z(b)}Z(U)\cap \Gamma^{\ast}_{\R}=0$, where $\Gamma$ is a lattice of even rank without skew-symmetric integer form.
\end{rmk}

\subsubsection{Case of mixed Hodge structure}

\begin{defn}
A semipolarized complex integrable system  is given by a holomorphic fibration of a complex analytic  symplectic manifold $\pi:X^0\to B^0$ where fibers are Lagrangian submanifolds which are semiabelian varieties with polarized abelian quotients.

\end{defn}

Our considerations in polarized case can be generalized to the semipolarized one.

Namely we have a local system of lattices $\GGamma\to B^0$ which is endowed with an integer skew-symmetric bilinear form $\langle\bullet,\bullet\rangle:\bigwedge^2\GGamma\to \underline{\Z}_{B^0}$, possibly degenerate.  Similarly to the pure case it is given by $\GGamma=(R^1\pi_{\ast}\underline{\Z}_{X^0})^{\vee}$.

This gives rise to an exact short sequence of local systems
$$0\to \GGamma_0\to \GGamma\to \GGamma^{symp}\to 0,$$
where $\GGamma_0$ is the kernel of the skew-symmetric bilinear form $\langle \bullet,\bullet\rangle$ and $\GGamma^{symp}$ is the symplectic quotient.

The local model is now given by a lattice $\Gamma$ endowed with a skew-symmetric form $\langle\bullet,\bullet\rangle:\bigwedge^2\Gamma\to \underline{\Z}_{B^0}$ and a local embedding $Z:U\to \Gamma^{\vee}\otimes \C$, where $U$ is a small neighborhood of a point $b\in B^0$.

Thus we have a local embedding $Z:U\to \Gamma^{\vee}\otimes \C$ such that the composition   $T_bU\stackrel{dZ}{\to}\Gamma^{\vee}\otimes \C\to \Gamma_0^{\vee}\otimes \C$ is surjection for any $b\in U$ and fibers of the corresponding submersion $U\to  \Gamma_0^{\vee}\otimes \C$ are  complex Lagrangian submanifolds 
of the symplectic leaves in the Poisson manifold $\Gamma^{\vee}\otimes \C$ (they are affine symplectic spaces parallel to the fibers of the natural map $\Gamma^{\vee}\otimes \C\to \Gamma_0^{\vee}\otimes \C$). Then   $Z(U)$ is a family of Lagrangian submanifolds over a domain in $\Gamma_0^{\vee}\otimes \C$. The positivity condition is also satisfied:
the restriction of the pseudo-hermitian form $i^{-1}\langle v,\overline{v}\rangle_{\ast}$ to $(dZ_b)(T_bU)\cap (\Gamma^{symp})^{\vee}\otimes \C, b\in U$ is positive.
Thus for any $Z_0\in \Gamma_0^{\vee}\otimes \C$ we have the corresponding Lagrangian submanifold $U_{|Z_0}$ which is endowed with a local system $\Gamma^{symp}\to U_{|Z_0}$ of symplectic lattices satisfying the positivity property. In other words, locally we have a family of (polarized) complex integrable systems parametrized by $\Gamma_0^{\vee}\otimes \C$.

The same double coset formula as in the polarized case describes fibers of the canonical local model of a semipolarized integrable system. The symplectic form on $\cup_{b\in U}\{b\}\times (\Gamma\backslash\Gamma\otimes \C/(T_{Z(b)}Z(U))^\perp$ is described such as follows. Its pull-back to $U\times (\Gamma\otimes \C)$ is the restriction of the canonical $2$-form on $(\Gamma^{\vee}\otimes \C)\times (\Gamma\otimes \C)$ obtained by the skew-symmetrization of the canonical pairing between $\Gamma$ and $\Gamma^{\vee}$.

Alternatively, we observe that $T^{\ast}U$ is a $\Gamma$-covering of the local model, hence the canonical symplectic structure on $T^{\ast}U$ descends to the local model giving the above symplectic structure. 

Similarly to the pure case we have locally a natural action of the connected abelian Lie group $\Gamma\backslash Ker(\Gamma\otimes \C\to \Gamma^{symp}\otimes \sqrt{-1}\R)$ (which is a product of a real torus and real vector space)  by holomorphic symplectomorphisms, such that fibers of the integrable system are torsors over this group.

Semipolarized integrable system gives rise to the following data:

1) A local system $\GGamma\to B^0$ of free abelian groups endowed with a skew-symmetric pairing 
$\langle\bullet,\bullet\rangle:\bigwedge^2\GGamma\to \underline{\Z}_{B^0}$. We denote by $\GGamma_0$ the kernel of this pairing.

2) A weight filtration $W_{\bullet}$
$$W_{-3}=0\subset W_{-2}\simeq \GGamma_0\subset W_{-1}=\GGamma$$
(notice that $W_{\bullet}$ is canonically determined by the pair $(\GGamma,\langle\bullet,\bullet\rangle)$).

3) A Hodge filtration
$$F^{-2}=0\subset F^{-1}\subset F^{0}=\GGamma\otimes {\mathcal O}_{B^0}.$$

4) An isomorphism of holomorphic vector bundles $\Psi:T_{B^0}\simeq (F^0/F^{-1})^{\ast}\subset (F^0)^{\ast}$.

The data 1)-4) are required to satisfy the following properties:

i) $gr_{-2}^{W_{\bullet}}(\GGamma)$ is a variation of pure Hodge structure of weight $-2$ concentrated in bidegree $(-1,-1)$.

ii) $gr_{-1}^{W_{\bullet}}(\GGamma)$ is concentrated in bidegrees $(-1,0)$ and $(0,-1)$.

iii) The pairing $\langle\bullet,\bullet\rangle$ induces a polarization on $gr_{-1}^{W_{\bullet}}(\GGamma)$.

iv) Locally and isomorphism $\Psi$ can be written as $\Psi=dZ$, where $Z$ is a local embedding of $B^0$ to $\Gamma^{\vee}\otimes \C$ as in the pure case.

Conversely, the data 1)-4) satisfying the properties i)-iv) define a semipolarized integrable system $X^{\prime}\to B^0$ endowed with a holomorphic Lagrangian section $B^0\to X^{\prime}$. An integrable system without Lagrangian section is determined by the associated integrable system with Lagrangian section (see Section 4.1.1) and a cohomology class in $H^1(B^0,\Omega^{1,cl}_{B^0}/\alpha(\GGamma))$.

Let now $\widehat{X}^0\to B^0$ be a fibration obtained from the initial semipolarized integrable system $X^0\to B^0$ by the fiberwise quotient by the natural action of the torus $\Gamma_0\otimes \C^{\ast}\simeq (\C^{\ast})^{r_0}, r_0=rk\,\Gamma_0$. Then $\widehat{X}^0$ is a Poisson manifold. Its symplectic leaves are total spaces of polarized integrable systems with bases which are submanifolds of $B^0$ obtained by fixing (locally) the value $Z_{|\Gamma_0}$ (e.g. in the case of  Hitchin systems with regular singularities we fix residues of the Higgs field at singularities). 

In other words we obtain a family of polarized  integrable systems (with canonical K\"ahler metrics on the base) which is parametrized (locally) by a domain in an affine space parallel to $\Gamma_0^{\vee}\otimes \C$.

\begin{rmk} Traditionally people speak about integrable systems as Poisson manifolds with symplectic leaves fibered by Lagrangian abelian varieties. The notion of semipolarized integrable system gives rise to such a structure (if we forget about polarization). But in a sense it is more precise. Namely, we have a variation of mixed Hodge structure (not visible in the traditional approach). Furthermore, the space of symplectic leaves carries locally a structure of an affine vector space. In practice the holonomy of the local system $\GGamma_0$ is finite (see Lemma 4.3.1 below). Also, under some mild assumptions (which are usually satisfied in practice) an integrable system in the traditional sense gives an integrable system in our sense by means of a simple topological construction, see Section 4.2.5 below.

\end{rmk}

\subsection{Integrable systems with central charge}

We are going to consider semipolarized integrable systems. Recall that the map $Z$ is defined locally up to a shift.

\begin{defn} A central charge for a semipolarized integrable system $\pi: X^0\to B^0$ is a holomorphic section $Z\in \Gamma(B^0,\GGamma^{\vee}\otimes {\mathcal O}_{B^0})$ such that the local isomorphism $\Psi: T_{B^0}\to (F^0/F^{-1})^{\ast}$ composed with the natural embedding $(F^0/F^{-1})^{\ast}\to (F^0)^{\ast}=\Gamma(B^0,\GGamma^{\vee}\otimes {\mathcal O}_{B^0})$ coincides with $dZ$ (cf. iv) in 4) in the previous subsection).

\end{defn}

In other words $Z$ is a homomorphism $\GGamma\to {\mathcal O}_{B^0}$. Clearly $dZ$ defines the only non-trivial $F^{-1}$-term of the Hodge filtration. For an integrable system with central charge we can locally embed $B^0$ as a submanifold in a {\it vector} space (not just an affine space as before). Not every integrable system has a central charge.

\begin{thm} If a semipolarized integrable system with holomorphic Lagrangian section has central charge then $[\omega^{2,0}]=0$.

\end{thm}

{\it Proof.} We keep the notation of Section 4.1. Locally we can identify (in $C^{\infty}$ sense) the total space $\pi^{-1}(U)$ with the product of $U$ and ${\mathbb T}_{\Gamma}:=Ker(\Gamma\otimes \C\to \Gamma^{symp}\otimes \sqrt{-1}\R)/\Gamma$. Then the holomorphic symplectic form $\omega^{2,0}$ is ${\mathbb T}_{\Gamma}$-invariant and its restriction to the tangent space at any point $(b,0)$ of the Lagrangian section is described such as follows. We have: $T_{(b,0)}(\pi^{-1}(U))\simeq Lie({\mathbb T}_{\Gamma})\oplus T_bU\subset (\Gamma\otimes \C)\oplus(\Gamma^{\vee}\otimes \C)$. The canonical skew-symmetric form on the latter space gives $\omega^{2,0}$ after restriction to $T_{(b,0)}(\pi^{-1}(U))$. Using the central charge $Z$ we define a complex-valued $C^{\infty}$ $1$-form $\beta$ on $\pi^{-1}(U)$ as the ${\mathbb T}_{\Gamma}$-invariant $1$-form, whose restriction to $T_{(b,0)}(\pi^{-1}(U))$ is the restriction of the $1$-form on $(\Gamma\otimes \C)\oplus(\Gamma^{\vee}\otimes \C)$ given by the pairing with $(0,Z(b))$. The direct calculation shows that $d\beta=\omega^{2,0}$. $\blacksquare$

The above Theorem gives an obstruction to the existence of central charge.

We remark that the condition in the above Theorem that the integrable system has Lagrangian section can be relaxed. Namely, let us recall that an integrable system without Lagrangian section is determined by the associated system with Lagrangian section and the ``twist", which is the cohomology class in  $H^1(B^0,\Omega^{1,cl}_{B^0}/\alpha(\GGamma))$. There is a morphism of sheaves of abelian groups ${\mathbb T}_{\GGamma}\to\Omega^{1,cl}_{B^0}/\alpha(\GGamma)$ over $B^0$ (here ${\mathbb T}_{\GGamma}=Ker(\GGamma\otimes \C\to (\GGamma^{symp}\otimes \sqrt{-1}\R)/\GGamma$)) given by
$$0\to {\mathbb T}_{\GGamma}\to (\GGamma\otimes \C)/\GGamma\stackrel{\alpha}{\to}\Omega^{1,cl}_{B^0}/\alpha(\GGamma).$$
Then one can easily generalize the above proof of the Theorem 4.2.2 to the case when the above twist belongs to the image of an element from $H^1(B^0,{\mathbb T}_{\GGamma})$. This generalization is useful for Hitchin integrable systems (see Section 8).

Let us now discuss the condition $[\omega^{2,0}]=0$ in several examples.

\subsubsection{$K3$ surfaces}

Let $\pi:X\to {\bf P}^1$ be an elliptic fibration of a $K3$ surface, and let $X^0\to B^0$ be the polarized integrable system obtained by throwing away singular fibers of $\pi$. Then $[\omega^{2,0}_{X^0}]\ne 0$, hence the integrable system does not have a central charge.
More generally complex integrable systems with the total space being a compact hyperk\"ahler manifold do not have central charge.

\subsubsection{Integrable systems from dimer models}
In [GoKen] the authors defined a class of integrable systems $X^0\to B^0$ for which $X^0$ is birationally symplectomorphic to the torus $(\C^{\ast})^{2n}$ endowed with the constant symplectic form $\sum_{i,j}\omega^{ij}dlog\,z_i\wedge dlog\,z_j$. Then $[\omega^{2,0}_{X^0}]\ne 0$, hence such integrable systems do not have  central charge.

\subsubsection{Systems birationally equivalent to those on cotangent spaces}

Let $\pi:X^0\to B^0$ be a semipolarized integrable system such $X^0$ is birationally symplectomorphic to a cotangent space $T^{\ast}M$ for some complex manifold $M$ (e.g. this is the case for Hitchin systems). We expect that under some mild conditions the central charge does exist. More precisely we can pull-back from $T^{\ast}M$ the Liouville form $pdq$, thus obtaining a meromorphic $1$-form $\lambda$ on $X^0$, such that $d\lambda=\omega^{2,0}$. The restriction of $\lambda$ to a semiabelian fiber $\pi^{-1}(b)$ is closed at generic points. In particular we can define residues of $\lambda_{|\pi^{-1}(b)}$ at  smooth components of the divisor of poles of $\lambda$. One can easily show that the residues are locally constant with respect to $b\in B^0$. Hence we obtain a finite collection of residues. If all of them are equal to zero then we can define the central charge $Z:\gamma^{\prime}\mapsto \int_{\gamma^{\prime}}\lambda$, where $\gamma^{\prime}$ is any loop in $\pi^{-1}(b)$ which sits in the complement to the divisor of poles and represents a class $\gamma\in H_1(\pi^{-1}(b),\Z)$. Since all residues are zero, the integral does not depend on the choice of representative $\gamma^{\prime}$.
This class of integrable systems contains so-called Seiberg-Witten integrable systems (see [Don]).

\subsubsection{Hitchin integrable systems} 

We will discuss this class of examples at length later in the paper. Let us just mention now that for a large class of $GL(n)$ Hitchin integrable systems with singularities (possibly irregular) one can define central charge. Hopefully it can be done for Hitchin systems associated with any reductive group.

This example can be put in the framework of log-families of Lagrangian submanifolds in non-compact Calabi-Yau $3$-folds (in the particular case of Hitchin systems we have log-families of spectral curves). We are going to discuss this class of examples later in  Sections 7,8.

\subsection{Families of integrable systems without central charge}

Suppose we are given  an analytic family of {\it full} complex integrable systems $\pi_t: (X_t,\omega^{2,0}_t)\to B_t$, where $t\in U$ and $U$ is a complex analytic manifold. We assume that there exist open dense subsets $X_t^0\subset X_t$ and $B_t^0\subset B_t$ such that the restriction of $\pi_t$ to $X_t^0$ gives rise to a polarized integrable system $\pi_t: (X_t^0,\omega^{2,0}_t)\to B_t^0$.
We assume that $\cup_{t\in U}X_t$ forms a locally trivial bundle over $U$ in the topological sense (we {\it do not} assume that   $\cup_{t\in U}X_t^0$ forms a locally trivial bundle over $U$).
We also assume that for every $t\in U$ we have: $H^1(X_t,\Q)=0$.

Let $B=\cup_{t\in U}B_t$ and $B^0\subset B$ be the open dense subset $\cup_{t\in U}B_t^0$.
We define a local system of lattices $\underline{\Gamma}\to B^0$ with fibers ${\GGamma}_{b}=H_2(X_t,\pi_t^{-1}(b),\Z), b\in B_t^0$. 
The long exact sequence of the pair $(X_t,\pi_t^{-1}(b))$ gives rise to a short exact sequence
$$0\to \underline{\Gamma}_{0,b}\to {\GGamma}_{b}\to {\GGamma}^{symp}_b\to 0,$$
where $\underline{\Gamma}_{0,b}=H_2(X_t,\Z)/Im(H_2(\pi_t^{-1}(b),\Z))$ and
 ${\GGamma}^{symp}_b=Ker(H_1(\pi_t^{-1}(b),\Z)\to H_1(X_t,\Z))$ is a sublattice of finite index in a symplectic lattice $H_1(\pi_t^{-1}(b),\Z)$, hence itself symplectic. The we see that ${\GGamma}$ carries a covariantly constant integer skew-symmetric pairing with the kernel ${\GGamma}_0$. Local system $\GGamma_0$ is constant along fibers of the projection $B^0\to U$.

Integration of $\omega_t^{2,0}$  gives a linear functional $Z: {\GGamma}\to \C$. The restriction of $Z$ to $\GGamma_0$ is constant along $B_t^0$ for any $t\in U$. Hence it defines (locally near $t\in U$) a map 
$$\chi_0: U\to {\GGamma}_{0,b}^{\vee}\otimes \C=Ker( H^2(X_t,\C)\to H^2(\pi_t^{-1}(b),\C)),\,\,\, t\mapsto [\omega^{2,0}_t]$$
for an arbitrary $b\in B_t^0$. 
We assume that it is in fact a local open embedding. Thus $U$ can be thought of as base of the universal family of integrable systems. This can be compared with Moser theorem which says that the universal family of real symplectic structures  on a symplectic manifold $(X,\omega)$ which are close to $\omega$ is parametrized by an open neighborhood of $[\omega]\in H^2(X,\R)$. 

The assumption that $\chi_0$ is a locally open embedding implies that $(B^0, \underline{\Gamma},\langle\bullet,\bullet\rangle, Z)$ defines a semipolarized integrable system with the base $B^0$. We warn the reader that the total space of this integrable system is bigger (it has even a bigger dimensions if $dim(U)>0$) than $\cup_{t\in U}X_t^0$, where $X_t^0=\pi_t^{-1}(B_t^0)$.

\subsection{Finiteness of the monodromy}

Finally we discuss the monodromy of the local system $\GGamma_0$, which is the kernel of the skew-symmetric form. We assume that the integrable system has central charge denoted by $Z$.
Then we obtain a locally well-defined map $Z:B^0\to \Gamma_0^{\vee}\otimes \C, Z\mapsto Z_{|\Gamma_0}$, where $\Gamma_0$ is a fiber of $\GGamma_0$. 
We will assume that our integrable system $\pi:X^0\to B^0$ is algebraic, i.e. in the definition of complex integrable system we have: $X^0$ is a smooth algebraic symplectic variety, $B^0$ is a smooth algebraic variety and $\pi$ is a regular map.

\begin{lmm} For an algebraic semipolarized integrable system the monodromy of the local system $\GGamma_0$ is finite.

\end{lmm}
{\it Proof.} As we discussed in the previous subsection the local system $\GGamma_0$ gives rise to a variation of pure Hodge structure of weight $-2$. Hence it admits a polarization. The monodromy preserves  positive quadratic form (polarization)  and is given by integer transformations. Hence its elements have finite order. $\blacksquare$

It follows from Lemma that we have a map $p_{B^0}:B^0\to (\Gamma_0^{\vee}\otimes \C)/G$, where $G$ is a finite group. Fibers of $p_{B^0}$ are smooth algebraic varieties endowed with K\"ahler metric.

\subsection{Local model near the discriminant }

We assume that our semipolarized integrable system $X^0\to B^0$ has central charge $Z$, and it is an open dense in a full complex integrable system $X\to B$.

Let us make the following

{\bf $A_1$-Singularity Assumption}

{\it We will assume that  $D=B-B^0$ is an analytic divisor. We will also assume that there exists an analytic divisor $D^1\subset D$ such that $dim\,D^1\le dim\,B^0-2$, the complement $D^0:=D-D^1$ is smooth, and such that our VMHS together with the central charge $Z$ (see Definition 4.2.1) has the following local model near $D^0$:

1) There exist local coordinates $(z_1,..,z_n, w_1,...,w_m)$ near a point of $D^0$ such that $z_1$ is small and $D^0=\{z_1=0\}$.

2) The map $Z: B^0\to \C^{2n+m}\simeq \Gamma^{\vee}\otimes \C$ is a multi-valued map given in coordinates by
$$(z_1,...,z_n,w_1,...,w_m)\mapsto (z_1,...,z_n,\partial_1F_0,...,\partial_nF_0,w_1,...,w_m),$$
where $\partial_i=\partial/\partial z_i$, and $F_0$ is given by the formula
$$F_0={1\over{2\pi i}}{z_1^2\over{2}}log\,z_1+G(z_1,..,z_n,w_1,...,w_m),$$
and $G$ is a holomorphic function. The Poisson structure on $\C^{2n+m}$ in the standard coordinates $(x_1,...,x_{2n+m})$ is given by the bivector $\sum_{1\le i\le n}\partial/\partial x_i\wedge \partial/\partial x_{i+n}$.

3) The function $F_0$ (called prepotential) satisfies also a positivity condition coming from the condition $i\langle dZ,\overline{dZ}\rangle>0$, which is satisfied for the restriction of $dZ$ to symplectic leaves $S_{c_1,...,c_m}:=\{(z_1,...,z_n,w_1,..., w_m)|w_i=c_i\}$. 

4) The monodromy of the local system $\GGamma$ about $D^0$ has the form $\mu\mapsto \mu+\langle \mu,\gamma\rangle \gamma$, where $\gamma$ is such that the pairing $\langle\gamma,\bullet\rangle\in \GGamma^{\vee}$ is a primitive covector. }

The map $Z$ is defined  up to a linear change of coordinates $$(z_1,...,z_n,z_{n+1},...,z_{2n},w_1,...,w_m)\mapsto (z_1,...,z_n,z_{n+1}+z_1,z_{n+2},...,z_{2n},w_1,...,w_m),$$
where $z_{n+i}=\partial_iF_0, 1\le i\le n$.
This can be interpreted as a local system of complex vector spaces endowed with a skew-symmetric form and a section. A choice of branch of $F_0$ allows us to identify the fibers of this local system with the standard $(\C^{2n+m})^{\ast}$ endowed the skew-symmetric form which is the the product of the standard symplectic form in $(\C^{2n})^{\ast}$ with the trivial form in $(\C^m)^{\ast}$. 

\begin{rmk} The $A_1$-Singularity Assumption allows only the simplest possible singularity at the discriminant $D$. There are other possibilities for such singularities. They are related to simply-laced Dynkin diagrams (Kodaira classification). We do not discuss more complicated singularities in this paper for two reasons:

i) in many examples (e.g. $GL(r)$ Hitchin systems discussed later) they do not appear;

ii) we do not understand fully the local geometry of other singularities.

\end{rmk}


\subsection{WCS for integrable systems with central charge}

We expect that for a large class of semipolarized integrable systems with central charge (including all algebraic ones) there exists a canonical WCS.

More precisely, suppose we are given a semipolarized integrable system $\underline{\Gamma}\to B^0$ with the central charge $Z$.  The local system $\GGamma$ gives rise to a canonical local system of torus (or quantum torus) Lie algebras $\underline{\g}$ (see Section 2.3, Examples 4, 5).

We will assume that the monodromy of 
$\underline{\Gamma}_0\to B^0$ is a finite group $G$ (as we have seen this is true in the case when $B^0$ is algebraic). Central charge $Z$ defines a submersion $p_{B^0}: B^0\to (\Gamma_0^{\vee}\otimes \C)/G$. Fix a  non-singular point $Z_0$ of the orbifold $(\Gamma_0^{\vee}\otimes \C)/G$ and let $M=B^0_{Z_0}=p^{-1}_{B^0}(Z_0)$. The restriction of the local system $\underline{\Gamma}_0$ to $M$ is trivial. Let us also fix $\theta\in \R$. Set $Y=Im(e^{-i\theta}Z)$. Then $Y$ defines a local embedding of every fiber $M$ of $p_{B^0}$ into $\Gamma_{\R}^{\ast}$ as an affine symplectic leaf which is parallel to $(\Gamma_{\R}^{symp})^{\ast}$. 

In order to construct WCS we would like to use the approach of Section 3. Recall that in Section 3 we gave a definition of the attractor tree bound by a family of cones $C^+$. Finiteness of the number of attractor trees was guaranteed by several assumptions, including the Mass Function Assumption. In this subsection we are going to reverse the logic. More precisely, in the hypothetic WCS we have an a priori idea of what are the relevant attractor trees. 

\begin{defn}
We call attractor tree good if it is a locally planar attractor tree in a fiber  $M$ of $p_{B^0}$ such that its tail edges hit transversally the discriminant $D$ at the locus $D^0$, and the velocity of any tail edge is proportional to the corresponding vector $\gamma$ (see part 4 of the $A_1$-Singularity  Assumption, Section 4.5). 

\end{defn}

Then one defines the family of  convex cones $C^+_{min}$ as the minimal closed subset of $(B^0)^{\prime}$ whose fibers under the natural projection to $B^0$ are closed convex cones and such that $C^+_{min}$ contains all velocities of the above-described good attractor trees. It is not a priori clear that those conic fibers are {\it strict} cones. 

We claim that the function $X=Re(e^{-i\theta}Z)$ plays the role of the mass function (see Section 3.3) and simultaneously gives a restriction on the velocities of good attractor trees.

\begin{lmm} Let us restrict $X$ to a fiber $M$ and identify the latter locally with a symplectic leaf in $\Gamma_{\R}^{\ast}$.
Let us also fix $v\in \Gamma_{\R}-\Gamma_{0,\R}$. Then
$${d\over{dt}}_{|t=0}X(m+t\iota(v),v)=-||\iota(v)||^2<0,$$
where $\iota: \Gamma_{\R}\to \Gamma_{\R}^{\ast}$ was defined in Section 3.1 and we understand the non-zero vector $\iota(v)$ as an element of $T_mM$ (recall that $M$ carries natural K\"ahler metric).

\end{lmm}
{\it Proof.} Follows from definitions. $\blacksquare$

\begin{cor} The function $X$ strictly decreases along the attractor flow. Moreover $X(m,v)>0$ at all inner points of a tail edge of a good attractor tree.

\end{cor}

{\it Proof.} It suffices to show that $X(m,v)$ approaches $0$ along a tail edge of a good attractor tree. $\blacksquare$

\begin{cor}  The function $X$ is positive on  edges of good attractor trees.

\end{cor}

{\it Proof.} By previous Corollary the result holds for tail edges. For other edges it follows by induction using the balancing conditions $\sum_i\gamma_{i}^{out}=\gamma^{in}$. $\blacksquare$

Consider now all functions $\widehat{X}(b,v)$ where $(b,v)\in (B^0)^{\prime}, Y(b)(v)=0$ which satisfy the properties:

i) $\widehat{X}(b,v)$ is linear in $v$ and strictly decreases along the attractor flow.

ii) $\widehat{X}(b,v)$ is strictly positive at inner points of tail edges of good attractor trees.

Every such function defines a closed subset in $(B^0)^{\prime}$ which is conic and convex in the direction of $\Gamma_{\R}$. Namely, we take $C_{\widehat{X}}^+=\{(b,v)|\widehat{X}(b,v)\ge 0\}$. We define $C^+:=\cap_{\widehat{X}}C_{\widehat{X}}^+$, where the intersection is taken over all such functions. The Corollaries 4.6.3 and 4.6.4 hold for the functions  $\widehat{X}(b,v)$. Therefore $C^+_{min}\subset C^+$.

\begin{conj}$C^+$ is a strict convex cone in the direction of $\Gamma_{\R}$.

\end{conj}

One can check that conjecture holds if $\underline{\Gamma}_0=0$.
We will discuss a motivation of a similar conjecture in the framework of Mirror Symmetry in Section 10.3.

\begin{defn} Canonical initial data associate with the tail of a good attractor tree with the velocity $k\gamma$ ($k\ge 1$ and $\gamma$ is primitive), the element ${1\over {k^2}}e_{k\gamma}$ of the torus Lie algebra.

\end{defn}

The canonical initial data are motivated by interpretation of WCS via DT-invariants (see Example 4, Section 2.3 and end of Section 3.3).

\begin{conj} Assuming Conjecture 4.6.5, Compactness Assumption, Tail Assumption and Mass Assumption, there is a unique WCS on $\R/\pi \Z\times B^0$ with the support in $C^+$ and the canonical initial data.

\end{conj}

In terms of DT-invariants for that WCS we have $\Omega_b(\gamma)=\Omega_b(-\gamma)$ for all generic $b\in B^0$.

In practical terms this means that starting with DT-invariants equal to $1$ which we assign to the smooth locus of the the discriminant divisor $B-B^0$ we can assign by induction  DT-invariants $\Omega_b(\gamma)\in \Z$ for any pair $(b,\theta)$ which does not belong to a wall and any $\gamma\in \Gamma_b$ such that $\theta=Arg(Z_b(\gamma))$. The collection $(\Omega_b(\gamma))$ satisfy the wall-crossing formulas from [KoSo1]. Algorithm for the construction follows from general considerations of Sections 3.2. Namely, for a fixed $(b,\gamma)$ we determine $\theta=Arg(Z_b(\gamma))$. Then we consider all good attractor trees on $\R/\pi \Z\times B^0$ such that
 their root vertex is $(b,\theta)$ and the root edge is $\gamma$.

Our assumptions imply that for generic $b$ there are finitely many such trees. They form a graph without oriented cycles (because $X=Re(e^{-i\theta}Z_b)$ is monotone along  attractor trajectories).
Hence we can order  edges of the graph is such a way that the lowest numbering receive vertices in $B-B^0$. Then we move from the lowest order vertices to the vertex $b$ using the WCF from [KoSo1] in order to calculate the DT-invariant for the outcoming edges of the graph. Using it last time for the root edge $\gamma$ we obtain $\Omega_b(\gamma)$.
Finally, varying $\theta$ and $Z_0$ we arrive to the WCS on $S^1_{\theta}\times B^0$.

\subsection{Metric on the base}

Recall the notation and assumptions of Sections 4.4, 4.5. We will assume that the monodromy of the local system $\GGamma_0$ is finite. Then we make the following

{\bf Completeness  Assumption} 

{\it The map $p_{B^0}$ extends uniquely to a complex analytic map $p_B:B\to (\Gamma_0^{\vee}\otimes \C)/G$. Fibers of $p_B$ are metric completions of the fibers of $p_{B^0}$. }

Under the Completeness Assumption we can work with  non-singular $B^0$ and then extend arising structures to the whole space $B$ uniquely. We expect that the Completeness Assumption holds in all realistic examples of full complex integrable systems.

For any $Z_0\in  (\Gamma_0^{\vee}\otimes \C)/G$ let us denote $p^{-1}_{B^0}(Z_0)$ by $B^0_{Z_0}$ and $p^{-1}_B(Z_0)$ by $B_{Z_0}$. The Completeness Assumption means that the K\"ahler metric on $B^0_{Z_0}$ extends to $B_{Z_0}$ making it into a length space (the metric is singular on $D\cap B_{Z_0}$).

Since $B_{Z_0}$ is a metric completion of $B^0_{Z_0}$, it can be canonically reconstructed from the latter as a topological space. Furthermore the complex structure on $B_{Z_0}$ can be reconstructed from the one on $B^0_{Z_0}$. Indeed we can extend it to  the set $B_{Z_0}-D^1$ using the local model described in the previous subsection. Then we extend the complex structure to the whole space $B_{Z_0}$ by the Hartogs principle,  taking the direct image of the sheaf on analytic function $\mathcal{O}_{B_{Z_0}-D^1}$.  These considerations explain that the WCS and the initial data can be canonically reconstructed from our semipolarized integrable system on $B^0$.

Next, for every $Z_0 \in (\Gamma_0^{\vee}\otimes \C)/G$ such that $Re(Z_0)=0$ (i.e. $Z_0\in (\Gamma_0^{\vee}\otimes i\R)/G$) we define a real $1$-form $\alpha_{Z_0}$ on
 $B^0_{Z_0}$ by the following formula
 $$\alpha_{Z_0}=\sum_{i,j}\omega^{ij} Re(z_i) d\,Im(z_j).$$
Here, locally on $B^0_{Z_0}\subset B^0$, we define functions $z_i(b):=Z_b(\gamma_i)$. In this definition $(\gamma_i)_{i=1,...,2n}, n=\,dim_\C B$ is a basis of a covariantly constant subspace $V^{sympl}\subset \GGamma_b\otimes \Q,\,\,b\in B^0_{Z_0}$ such that $V$ is complementary to $ \GGamma_{0,b}\otimes \Q$,
and $(\omega^{ij})_{1\le i,j\le 2n}$ is the inverse matrix to the symplectic pairing $\omega_{ij}=\langle \gamma_i,\gamma_j\rangle$.

 It is immediate to check that $\alpha_{Z_0}$ is well-defined and closed.

 In general, we expect that the following 
 holds.

{\bf  Potential Assumption.} {\it There is a smooth function $H_{Z_0}$ on $B^0_{Z_0}$ such that $dH_{Z_0}=\alpha_{Z_0}$.
 Moreover, the function $H_{Z_0}$ extends by continuity to $B_{Z_0}$, and gives a bounded below and proper map from $B_{Z_0}$ to $ \R$.}

Let us give some motivations for the Potential Assumption under the $A_1$-Singularity Assumption from the Section 4.5.

One can check that the integral of $\alpha_{Z_0}$ around a small loop around the divisor  $D^0\cap B_{Z_0}\subset B_{Z_0}$ vanishes, and that an antiderivative of $\alpha_{Z_0}$ (defined a priori on $B_{Z_0}^0$ near any point of $D^0\cap B_{Z_0}$) extends continuously to $D^0\cap B_{Z_0}$.

It is  natural to expect that $H^1(B_{Z_0},\C)=H^1(B_{Z_0}-D^1,\C)$
 (at least it holds if $B_{Z_0}$ is nonsingular analytic space because the removing of complex analytic subset of codimension at least 2 does not change the fundamental group).  There are good reasons to expect that $B_{Z_0}$ itself is simply connected (in fact  contractible, see below).
 Therefore, we conclude  that $\alpha_{Z_0}$ is exact, i.e. it can be written as $\alpha_{Z_0}=d H_{Z_0}$ for some real-valued function 
 $H_{Z_0}$ on $B^0_{Z_0}$. This function is strictly convex in the affine structure given by $Im(Z)$ because its tensor of second derivatives
 coincides with the metric tensor of the Riemannian metric $g_{B^0_{Z_0}}$ associated with the canonical K\"ahler metric on $B^0_{Z_0}$. By the $A_1$-Singularity Assumption the function $H_{Z_0}$ continuously extends to $D^0\cap B_{Z_0}$.

 The boundedness an properness of $H_{Z_0}$ can be checked in special cases. For example, let us consider   Seiberg-Witten integrable system, which is a polarized integrable system with central charge, fibers being  the elliptic curves $y+1/y-x^2=2u$ parametrized by $u\in B^0_{Z_0}=B^0=\C-\{-1,1\}$. Then the metric completion is $B_{Z_0}=B=\C$.  The function $H:=H_{Z_0}$ can be expressed as $H=K-2Im(F)$, where $$K=\sqrt{-1}\sum_{i,j}\omega^{ij}z_i\overline{z}_j$$
 is the potential of the K\"ahler metric on $B^0$ and $F$ is a holomorphic function on $B^0$ such that 
 $$dF=\sum_{i,j}\omega^{ij}z_iz_j.$$
One can check that $dF=const\cdot du$.

In this example $K$ is a transcendental function which is non-negative and grows as $u\to \infty$ as $C|u|\cdot log|u|$, where $C>0$. We see that the summand $K$ dominates $2Im(F)$, and hence $H$ is bounded from below and proper. We expect that the general case is similar to this example.

\begin{rmk}
In general, the convexity of the function $H_{Z_0}$ on $B^0_{Z_0}$ and its boundedness from below on $B_{Z_0}$ supports the idea that it has a unique global minimum $b^{min}_{Z_0}\in B_{Z_0}$ for any $Z_0$ satisfying the condition $Re(Z_0)=0$. 

Notice that the condition $Re(Z_0)=0$ implies that the map $b\mapsto Re(Z_b)\in (\GGamma_b/\GGamma_{0,b})^{\vee}\otimes \R$ identifies locally $B^0_{Z_0}$ with an open domain in a symplectic vector space.
Therefore we can define a canonical Euler vector field $Eu_{Z_0}$ on $B^0_{Z_0}$. We conjecture that the flow associated with  $Eu_{Z_0}$ extends to a continuous flow on   $B_{Z_0}$ which contracts as $t\to -\infty$ the space  $B_{Z_0}$ to the point  $b^{min}_{Z_0}$. In particular this conjecture implies that $B_{Z_0}$ is contractible.

If  $b^{min}_{Z_0}\in B_{Z_0}^0$ then the arising local geometry is related to the theory of cluster transformations. We will discuss it elsewhere.

In the example  of $SL_2$ Hitchin system with regular singularities on a smooth projective curve $C$, fixing purely imaginary $Z_0$ we obtain a complex integrable system over the base $B_{Z_0}$ which consists of quadratic differentials with fixed residues at their singularities, which are poles of second order. Then the point $b^{min}_{Z_0}$ can be identified with the unique Strebel differential on the curve $C$.

\end{rmk}

\section{Formal neighborhood of a wheel of projective lines and stability data}

Given a lattice $\Gamma$ endowed with an integer skew-symmetric form 
$\langle\bullet,\bullet\rangle$ , we can consider  stability data on the graded Lie algebra $\g=\oplus_{\gamma\in \Gamma-\Gamma_0}\Q\cdot e_{\gamma}$ where $[e_{\gamma_1},e_{\gamma_2}]=(-1)^{\langle\gamma_1,\gamma_2\rangle}\langle\gamma_1,\gamma_2\rangle e_{\gamma_1+\gamma_2}$ and $\Gamma_0=Ker\,\langle\bullet,\bullet\rangle$. The aim of this Section is to encode these data in terms of formal Poisson varieties endowed with some additional structures. Small variation of the central charge $Z:\Gamma\to \C$ corresponds to an isomorphic Poisson variety.

\subsection{Wheels of lines, wheels of cones and toric varieties}

Let $Y$ be a complex toric variety of dimension $n$. Then it is stratified by the orbits of the action of the torus $T\simeq (\C^{\ast})^n$. 

\begin{defn} Wheel of lines in $Y$ is a cyclically ordered collection  of   $1$-dimensional $T$-orbits $F_{i}, i\in \Z/m\Z, m\ge 3$  such that each 

a) $\overline{F}_i\simeq {\bf P}^1$ for any $i$.

b) $\overline{F}_{i-1}\cap \overline{F}_{i}=\{p_{i}\}$, where $p_i$ is a point;

c) $\overline{F}_i\cap \overline {F}_{j}=\emptyset$ if $|i-j|>1$;

\end{defn}

Let $\Gamma=Hom(T,\C^{\ast})$ be the group of characters.
Each intersection point $p_i$ is a $0$-dimensional $T$-invariant stratum, hence it defines  a closed strict rational convex cone $C_i\subset \Gamma^{\vee}\otimes \R$ of full dimension with interior points corresponding to $1$-parameter subgroups which attract  points of the open dense $T$-orbit of $Y$ to points $p_i$. Clearly the collection of cones $C_i$ forms a wheel of cones in the following sense.

\begin{defn} We will call  wheel of cones a cyclically ordered collection   of real closed polyhedral strict convex cones of dimension $n$, $C_i\subset \Gamma^{\vee}\otimes \R, i\in \Z/m\Z, m\ge 3$ such that:

a) $C_i\cap C_{i+1}$ is a face of codimension one in $C_i$ and $C_{i+1}$, $i\in \Z/m\Z$.
 
b) $int(C_i)\cap int(C_j)=\emptyset$ if $i\ne j$ (here $int$ means the interior).

\end{defn}

If $Y$ is smooth then all $C_i$ are isomorphic to octants, i.e. to $\R_{\ge 0}^n$ (modulo the action of $GL(n,\Z)$).
Let us denote $C_i\cap C_{i+1}, i\in \Z/m\Z$ by $C_{i,i+1}$. 

We will call the wheel of cones {\it admissible} if it satisfies the following two assumptions:

\vspace{2mm}

{\bf Connectedness Assumption} 

{\it  For any $i\in \Z/m\Z$ the cone $C_i$ is the convex hull of $C_{i-1,i}\cup C_{i,i+1}$ and for any $i,j, |i-j|>1$ 
the convex hull of the pair $C_{i,i+1},C_{j,j+1}$ contains all cones $C_{i+1},...,C_{j}$ or all cones $C_{j+1},...,C_{i}$.}

We call it the Connectedness Assumption because it implies the following property of the wheel of cones:

for any closed half-space $\alpha\subset \Gamma_{\R}$ such that $0\in \partial \alpha$ the set of indices $i$ for which $C_{i,i+1}$ belongs to $\alpha$, forms an interval in $\Z/m\Z$.

\vspace{2mm}

{\bf Non-degeneracy Assumption}

 {\it  $\cap_i C_i^{\vee}=\{0\}$, where $C_i^{\vee}=\{x\in \Gamma_{\R}|y(x)\ge 0, y\in C_i\}$ is the dual cone.}

\vspace{2mm}

One can check that  wheel of cones (and hence toric varieties) satisfying the above two assumptions do exist in any dimension $n\ge 3$. Indeed, let us consider a compact rational polyhedron $P\subset \R^{n-2}$ which contains the origin in the interior. Let $v_i, i\in \Z/3\Z$ be three cyclically ordered vectors in $\R^2$ which generate $\R^2$ and satisfy the condition $v_1+v_2+v_3=0$. Then we define the cone $C_i$ as the convex hull of the set $\{\R_{\ge 0}(v\oplus p)| p\in P, v= v_i \,\,or\,\, v=v_{i+1}\}$. One can check that in this way we obtain the wheel of cones.

We denote by $x^{\gamma}\in {\mathcal O}(T)$ the monomial corresponding to the vector $\gamma\in \Gamma$. Then $x^{\gamma_1}x^{\gamma_2}=x^{\gamma_1+\gamma_2}$. We identify the open orbit of $T$ with $T$ itself. 

For a wheel of lines $(F_i)_{i\in \Z/m\Z}$ we denote by $\widehat{Y}:=\widehat{Y}_{\cup_i\overline{F}_i}$ the completion of $Y$ along $\cup_i\overline{F}_i$.
The following result will be used in the next subsection.

\begin{thm} For an admissible wheel of cones one has:

a) $H^0(\widehat{Y},{\cal O}_{\widehat{Y}})=\C$.

b) $H^1(\widehat{Y},{\cal O}_{\widehat{Y}})=\prod_{\gamma \in \cup_{1\le i\le m}(C_{i,i+1})^{\vee}\cap \Gamma}\C \cdot x^{\gamma}.$

c) $H^i(\widehat{Y},{\cal O}_{\widehat{Y}})=0$ for $i\ge 2$.

\end{thm}

{\it Proof.} We are going to compute cohomology groups using a covering of $\widehat{Y}$, which consists of formal neighborhoods $U_{i}$ of  $F_i\cup \{p_i\}\cup F_{i+1}$. Then $U_{i,i+1}=U_i\cap U_{i+1}$ is the formal neighborhood of $F_{i+1}$ and $U_i\cap U_j=\emptyset$ for $|i-j|>1$. Notice that the schemes $F_i$ and $F_i\cup \{p_i\}\cup F_{i+1}, i\in \Z/m\Z$ are open affine subschemes of $\cup_i\overline{F}_i$.

The \v{C}ech cochain complex associated with this covering has the form

$$\bigoplus_{i\in \Z/m\Z}{\mathcal O}(U_i)\to \bigoplus_{i\in \Z/m\Z}{\mathcal O}(U_{i,i+1}).$$
This makes c)  clear, since there are no non-trivial \v{C}ech $i$-cochains where $i\ge 2$.

Next we observe that the algebra of functions on the formal neighborhood of $p_i$ is  $\widehat{\mathcal O}(U_i)\simeq\prod_{\gamma\in C_i^{\vee}\cap \Gamma}\C\cdot x^{\gamma}$. Similarly we define the completed vector space $\widehat{\mathcal O}(U_{i,i+1})=\prod_{\gamma\in C_{i,i+1}^{\vee}\cap \Gamma}\C\cdot x^{\gamma}$.
We will proceed by replacing the algebras of functions 
${\mathcal O}(U_i), {\mathcal O}(U_{i,i+1})$ by the corresponding completed vector spaces of formal series. After that we will explain how return to the actual algebras of functions. 

In order to compute $H^0$ we observe that by the Non-degeneracy Assumption we see that only monomial  $x^0=1$ appears in $H^0(\widehat{Y},{\cal O}_{\widehat{Y}})$, hence a) holds.
More generally, for any  $\gamma$ we can compute the input of $x^{\gamma}$ to $H^0$ and $H^1$. For that we draw a planar polygon with vertices corresponding to $p_i$ and edges corresponding to $F_i$.
We are interested in those points $p_i$ and lines $F_i$ for which the monomial $x^{\gamma}$ appears in the algebras ${\mathcal O}(U_{i})$ and ${\mathcal O}(U_{i,i+1})$ respectively. The union of the relevant vertices and edges is an open subset of the polygon. Hence it can be one of the following subsets:

i) empty subset;

ii) full polygon;

iii) an open interval which consists of a chain of $1\le k\le m$ consecutive open edges and $k-1$ vertices;

iv) a disjoint union of at least two open intervals from iii).

Case i) is clear since $x^{\gamma}$ does not appear in the cohomology at all. The case ii)  by the Non-degeneracy Assumption corresponds to $\gamma=0$, which gives  the input to both $H^0$ and $H^1$.

In the case iii) we have trivial input to $H^0$ and one-dimensional  input of $x^{\gamma}$ to $H^1$. Case iv) is impossible by the Connectedness Assumption. Hence we have proved b) by replacing the algebras of functions by their completions.
In order to finish the proof it suffices to check that the following complex of vector spaces is acyclic

$$\bigoplus_{i\in \Z/m\Z}\widehat{{\mathcal O}}(U_i)/{{\mathcal O}}(U_i)\to \bigoplus_{i\in \Z/m\Z}\widehat{{\mathcal O}}(U_{i,i+1})/{\mathcal O}(U_{i,i+1}).$$
It follows from definitions that ${\mathcal O}(U_i)$ consists of such series $f=\sum_{\gamma\in C_{i}^{\vee}\cap \Gamma}f_{\gamma}x^{\gamma}\in \widehat{{\mathcal O}}(U_i)$ that $Supp(f)$ is finite on rays parallel to the rays in $C_i^{\vee}$ corresponding to $(n-1)$-dimensional faces $C_{i-1,i}$ and $C_{i,i+1}$. Similarly ${\mathcal O}(U_{i,i+1})\subset \widehat{{\mathcal O}}(U_{i,i+1})$ consists of such series $g=\sum_{\gamma\in C_{i,i+1}^{\vee}\cap \Gamma}g_{\gamma}x^{\gamma}$ that $Supp(g)$ is finite on lines parallel to $C_{i,i+1}^\perp\subset C_{i,i+1}^{\vee}$.

In order to describe quotient spaces in the above complex let us introduce some notation.  Let $L$ be an oriented rational line in $\Gamma_{\R}$ with a positive primitive generator $l\in L\cap \Gamma$. Let $V\subset \Gamma_{\R}/L$ be a closed strict polyhedral cone. We denote by $\pi$ the natural projection $\Gamma_{\R}\to \Gamma_{\R}/L$. We denote by ${\mathcal F}(V,l)$ the vector space which is the quotient of the vector space of such series $f=\sum_{\gamma\in \pi^{-1}(V)\cap \Gamma}f_{\gamma}x^{\gamma}$ that $f_{\gamma-nl}=0$ for sufficiently large $n$ by the subspace of series for which for a given $\gamma$ we have $f_{\gamma\pm nl}=0$ for all sufficiently large $n\ge 0$. 
If we choose a splitting of the projection $\Gamma\to \Gamma\cap L$ then  
$${\mathcal F}(V,l)\simeq \{\sum_{\mu\in (\Gamma/\Gamma\cap L)\cap V}a_{\mu}x^{\mu}|a_{\mu}\in \C[[x^l]]/\C[x^l]\simeq \C((x^l))/\C[x^l,(x^l)^{-1}]\},$$
where $x^l$ is the monomial corresponding to the generator $l$.

Let us return to the proof. Let $l_i$ be a generator of $C_{i,i+1}^{\perp}\cap \Gamma\simeq \Z$ which is positive as a functional on $C_i$. Let $V_i$ be the image of $C_{i,i+1}^{\vee}$ under the projection $\Gamma_{\R}\to \Gamma_{\R}/\R\cdot l_i$. Then $\widehat{{\mathcal O}}(U_i)/{{\mathcal O}}(U_i)\simeq {\mathcal F}(V_{i-1},l_{i-1})\oplus {\mathcal F}(V_{i},-l_{i})$ and similarly
$\widehat{{\mathcal O}}(U_{i,i+1})/{\mathcal O}(U_{i,i+1})\simeq {\mathcal F}(V_{i},l_{i})\oplus {\mathcal F}(V_{i},-l_{i})$.
From this explicit description it is easy to see that the differential in the quotient complex is an isomorphism. This completes the proof. $\blacksquare$

\subsection{Deformations of   formal Poisson manifolds}

Let us choose an admissible wheel of cones $(C_i)_{i\in \Z/m\Z}$ .
Suppose that the character lattice $\Gamma$ is endowed with an integer skew-symmetric form $\langle \bullet,\bullet\rangle:\bigwedge^2\Gamma\to \Z$ with the kernel  $\Gamma_0\subset \Gamma$. Then the toric variety $Y$ and its completion $\widehat{Y}$ defined in the previous subsection carry $T$-invariant Poisson structures. Both schemes are stratified (by the closures of $T$-orbits in the case of $Y$ and by their completions in the case of $\widehat{Y}$).

Conversely, suppose we have a free abelian group $\Gamma$ endowed with a skew-symmetric integer form $\langle \bullet,\bullet\rangle$ and a wheel of cones satisfying the Connectedness and Non-degeneracy Assumptions. 
Recall that a toric variety is given by a $T$-torsor together with a choice of fan. 
In what follows we are going to use the canonical $T$-torsor $\TT_{can}:=\TT_{can}(\Gamma,\langle\bullet,\bullet\rangle)$ which is the spectrum of the algebra with $\C$-linear basis $e_{\gamma}, \gamma\in \Gamma$ and the multiplication rule
$e_{\gamma}e_{\mu}=(-1)^{\langle\gamma,\mu\rangle}e_{\gamma+\mu}$ (the choice of this torsor is motivated by applications to the theory of DT-invariants).  Let us choose the fan which consists of cones $C_i$ and their faces.  We will denote the corresponding toric variety by $Y_{can}$. Its completion along $\overline{F}_{can}=\cup_i\overline{F}_{i,can}$ will be called {\it canonical local toric model associated with the wheel of cones $(C_i)_{i\in \Z/m\Z}$} (or simply {\it local toric model})  and denoted by $\widehat{Y}_{can}$. Notice that all $1$-dimensional $T$-orbits ${F}_{i,can}\subset \widehat{Y}_{can}$ are endowed with distinguished coordinates given by $e_{\gamma_i}$, where $\gamma_i\in C_i^{\vee}\cap \Gamma$ is a  primitive vector such that $C_{i,i+1}\subset \gamma_i^{\perp}$.

In case of the local toric model we have a distinguished collection of rational functions which are central with respect to the Poisson bracket and which are parametrized by $\Gamma_0$. Namely, denote by $D_Y\subset Y$ the canonical toric divisor (i.e. the complement to the open $T$-orbit), and by $D_{\widehat{Y}}\subset \widehat{Y}$ its completion along $\overline{F}$.
Let ${\mathcal O}_{\widehat{Y}}(*D_{\widehat{Y}})$ be the sheaf  of rational functions having  poles at $D_{\widehat{Y}}$. An element $e_{\gamma}, \gamma\in \Gamma$ defines a section $s_{\gamma}$ of this sheaf. In particular the map $c:\Gamma_0\to \Gamma(\widehat{Y},{\mathcal O}_{\widehat{Y}}(*D_{\widehat{Y}})^{\times}), \gamma\mapsto s_{\gamma}$ defines a homomorphism of $\Gamma_0$ to the abelian group of invertible functions having poles at $D_{\widehat{Y}}$. The image of $c$ (equivalently the collection $(s_{\gamma})_{\gamma\in \Gamma_0}$) is by definition our distinguished collection.

Next we would like to describe stratified formal Poisson varieties which are locally isomorphic to the above local toric model and endowed with additional data called decoration. We need some preparations for that.

\begin{defn} A pair $(\ZZ,D)$ consisting of a (possibly formal) normal scheme $\ZZ$ with a reduced divisor $D$ such that all singularities of $\ZZ$ are contained in $D$ is called a local formal toric pair if the following condition is satisfied: for any closed point $x\in D$ the pair $(\widehat{\ZZ}_x,\widehat{D}_x)$ obtained by the completion at $x$ is isomorphic to a similar pair
$(\widehat{\ZZ}_x^{tor},\widehat{D}_x^{tor})$, where $\ZZ_x^{tor}$ is a toric variety and $D_x^{tor}$ is the canonical toric divisor (i.e. the complement to the open orbit).

\end{defn}

For any local formal toric pair the divisor $D$ is canonically stratified, where the stratification is induced by the canonical stratification of the corresponding canonical toric divisors  $D_x^{tor}$. All open strata are smooth.

With a $0$-dimensional stratum $\{x\}\subset D$ we associate a lattice $\Gamma_x$ and a torsor $\TT_x$ over the torus $Hom(\Gamma_x,\C^{\ast})$. Namely, $\Gamma_x$ is the lattice of characters of the torus of the corresponding local formal toric model. The torsor $\TT_x$ is defined such as follows.  Choose an isomorphism $\phi_x: (\widehat{\ZZ}_x,\widehat{D}_x)\simeq  (\widehat{\ZZ}_x^{tor},\widehat{D}_x^{tor})$. It induces an isomorphism of completed local algebras
$${\mathcal O}_{\widehat{\ZZ}_x}\simeq \prod_{\gamma\in \Gamma_x\cap C_x}\C\cdot x^{\gamma}\simeq {\mathcal O}_{\widehat{\ZZ}_x^{tor}},$$
where $C_x\subset \Gamma_x\otimes \R$ is a strict rational polyhedral cone.

Notice that there exists a natural homomorphism of groups 
$$G_x:=Aut({\mathcal O}_{\widehat{\ZZ}_x^{tor}},J_{\widehat{D}_x^{tor}})\to Hom(\Gamma_x,\C^{\ast}),$$
where $J_{\widehat{D}_x^{tor}}=\prod_{\gamma\in \Gamma_x\cap int(C_x)}\C\cdot x^{\gamma}$ is the ideal of $\widehat{D}_x^{tor}$.

It is easy to describe the above homomorphism at the level of Lie algebras. Namely,
$$Lie(G_x)=\prod_{\gamma\in \Gamma_x\cap C_x}(\Gamma_x^{\vee}\otimes \C)\cdot x^{\gamma},$$
where we interpret elements of $\Gamma_x^{\vee}\otimes \C$ as toric vector fields. Then the homomorphism is defined by taking the $x^0$-component of the vector field. It follows that 
$Hom(\Gamma_x,\C^{\ast})\simeq G_x/\overline{[G_x,G_x]}$ (quotient of the topological group $G_x$ by the closure of its commutant).
The torsor $\TT_x$ is the $Hom(\Gamma_x,\C^{\ast})$-torsor associated with the natural $G_x$-torsor consisting of isomorphisms $\phi_x$. 

Similarly, for any $1$-dimensional closed stratum $\overline{F}\simeq \C{\bf P}^1$ of $D$ which contains exactly two $0$-dimensional strata $x_0,x_{\infty}$ one can define a lattice $\Gamma_{\overline{F}}$ and a torsor $\TT_{\overline{F}}$ over $Hom(\Gamma_{\overline{F}},\C^{\ast})$. It is easy to see that there is a canonical isomorphism $(\Gamma_{\overline{F}},\TT_{\overline{F}})\simeq (\Gamma_x,\TT_x)$ where $x=x_0$ or $x=x_{\infty}$.

\begin{defn} Suppose we are given a lattice $\Gamma$ endowed with a skew-symmetric form $\langle\bullet,\bullet\rangle:\bigwedge^2\Gamma\to \Z$.

A decorated formal scheme is defined by the following data:

i) A formal   Poisson scheme $\widehat{X}$ of pure dimension $n=rk\,\Gamma$, endowed with a normal Poisson divisor $D_{\widehat{X}}\subset \widehat{X}$ such that the pair $(\widehat{X}, D_{\widehat{X}})$ is a local formal toric pair.

ii) The reduced (non-formal) scheme associated with $\widehat{X}$ is a wheel of parametrized projective lines 
$\overline{F}=\cup_{k\in \Z/m\Z}i_k({\bf P}^1)$ where each $i_k$ is an embedding of the projective line such that $p_k=i_k(\infty)=i_{k+1}(0)$ is the only intersection point of $\overline{F}_k=i_k({\bf P}^1)$ and $\overline{F}_{k+1}=i_{k+1}({\bf P}^1)$. Moreover all points $p_k$ and lines 
$\overline{F}_k$ are strata of the canonical stratification of $D_{\widehat{X}}$.

iii) A homomorphism $c_{\widehat{X}}:\Gamma_0\to \Gamma(\widehat{X},
{\mathcal O}_{\widehat{X}}(*D_{\widehat{X}})^{\times})$.

iv)  For any $k\in \Z/m\Z$  isomorphisms $\Gamma_{p_k}\simeq \Gamma$ and
$\TT_{p_k}\simeq \TT_{can}$.

We require that the data i)-iv) satisfy the following conditions:

a) for any $k\in \Z/m\Z$ there exists an isomorphism $\phi_k$ of the pair $(\widehat{X}_{p_k},\widehat{D}_{\widehat{X},{p_k}})$ (completions at the point $p_k$) with the completions at the corresponding point of the corresponding local formal toric pair (see Def. 5.2.1), which identifies the Poisson structures and compatible with the homomorphism $c_{\widehat{X}}$.

b) for any $k\in \Z/m\Z$ the composition of the isomorphisms 

$$(\Gamma_{p_k},\TT_{p_k})\simeq (\Gamma,\TT_{can})\simeq (\Gamma_{p_{k+1}},\TT_{p_{k+1}})$$

coincides with the composition of the isomorphisms

$$(\Gamma_{p_k},\TT_{p_k})\simeq (\Gamma_{\overline{F}_k},\TT_{\overline{F}_k})\simeq (\Gamma_{p_{k+1}},\TT_{p_{k+1}}).$$

\end{defn}

\begin{rmk} The local toric model $\widehat{Y}_{can}$ carries a natural structure of decorated formal Poisson scheme.

\end{rmk}

Next we would like to describe the deformation theory of decorated formal Poisson schemes. For that we need to study local symmetries of the local toric model preserving the structure of a decorated formal Poisson scheme. A toy-model of the result can be illustrated by the following Proposition.

\begin{prp} Let $(\Gamma, \langle\bullet,\bullet\rangle)$ be a lattice endowed with an integer skew-symmetric form, $C\subset \Gamma_{\R}$ be a closed rational strict convex cone such that $int\,C\ne \emptyset$.
 
Let $Y_C=Spec(\oplus_{\gamma\in C\cap \Gamma}\Q\cdot e_{\gamma})$, where $e_{\gamma}e_{\mu}=(-1)^{\langle\gamma,\mu\rangle}\langle\gamma,\mu\rangle e_{\gamma+\mu}$ be the corresponding toric variety, and let $\widehat{Y}_C$ be its completion at $0\in C$, i.e. $\widehat{Y}_C=Spf(\prod_{\gamma\in C\cap \Gamma}\Q\cdot e_{\gamma})$. Consider the group of such automorphisms of $\widehat{Y}_C$ that:

1) they preserve the completion of the toric stratification;

2) they preserve the Poisson structure induced by $\langle\bullet,\bullet\rangle$;

3) they preserve all elements $e_{\gamma}, \gamma\in \Gamma_0=Ker\, \langle\bullet,\bullet\rangle$
considered as rational functions on $\widehat{Y}_C$;

4) they are equal to $id$ on the torsor $\TT_{y_0}$, where $y_0$ is the the only $0$-dimensional toric stratum.

Then this group is a pronilpotent proalgebraic with the Lie algebra isomorphic to $\prod_{\gamma\in C\cap(\Gamma-\Gamma_0)}\Q\cdot e_{\gamma}$, which acts via $\{e_{\gamma},\bullet\}$.

\end{prp}

We are not going to prove the Proposition, since we are not going to use it. Let us explain informally its meaning. An automorphism of an affine Poisson variety or its completion can act non-trivially on the Poisson center. Inner Poisson derivations act trivially on the center. This explains the condition 3). The condition 4) allows us to exclude infinitesimal symmetries identical on the center but not inner which are given by $\{log\,e_{\gamma},\bullet\}$. 

In a similar way we conclude that the sheaf of Lie algebras of infinitesimal symmetries of  the local model $\widehat{Y}_{can}$ is naturally isomorphic to $\g_{can}={\mathcal O}_{\widehat{Y}_{can}}/{\mathcal O}_{\widehat{Y}_{can}}^{center}$ endowed with obvious Poisson bracket $\{\bullet,\bullet\}$. Here ${\mathcal O}_{\widehat{Y}_{can}}^{center}=Ker\, \{\bullet,\bullet\}$.

\begin{thm} We have:

1) $H^0(\widehat{Y}_{can}, \g_{can})=0$.

2) $H^1(\widehat{Y}_{can}, \g_{can})\simeq \prod_{\gamma\in \cup_iC_{i,i+1}^{\vee}\cap(\Gamma-\Gamma_0)}\Q\cdot e_{\gamma}$.

3) $H^{\ge 2}(\widehat{Y}_{can}, \g_{can})=0$.

\end{thm}

{\it Proof.} Analogous to the proof of Theorem 5.1.3. The only difference is that elements $e_{\gamma}, \gamma\in \Gamma_0$ are now excluded from considerations.
$\blacksquare$

The sheaf of pronilpotent Lie algebras $\g_{can}$ defines the  deformation functor $Def_{\g_{can}}$ from the category of commutative (not necessarily Arin) unital algebras over $\Q$ to the category of groupoids:
$ Def_{\g_{can}}(R)$ is the groupoid of $exp(\g_{can}\widehat{\otimes} R)$-torsors over $\widehat{Y}_{can}$. Notice that such torsors can be identified with decorated formal Poisson schemes $\widehat{X}$. The above theorem implies the following result.

\begin{prp} The deformation functor $Def_{\g_{can}}$ is representable by an affine scheme $\MM:=\MM_{\Gamma,\langle\bullet,\bullet\rangle,(C_i)_{i\in \Z/m\Z}}$ (with the trivial stacky structure) which is (non-canonically) isomorphic to the infinite-dimensional affine space over ${\mathbb A}^{\infty}=\varprojlim_{N}{\mathbb A}^{N}$ over $\Q$.

\end{prp}

There is an alternative description of $\MM$. Namely, let $G_i$ denotes the  proalgebraic group of automorphisms of the formal neighborhood $U_i$ and $G_{i,i+1}$ be a similar group for $U_{i,i+1}$ (we always assume that automorphisms are compatible with identifications i)-iv)).  We have a chain of embeddings:
$$\hookleftarrow G_1\hookrightarrow G_{1,2}\hookleftarrow G_2\hookrightarrow G_{2,3}\hookleftarrow G_3\hookrightarrow....$$
Then the product group $G=G_{1,2}\times G_{2,3}\times G_{3,4}\times...$ is endowed (as a scheme) with the free
 action of the group $H=G_1\times G_2\times G_3\times...\subset G\times G$ (namely the factors of $G_i\times G_{i+1}$ act on $G_{i,i+1}$ by left and right multiplication respectively).
By the above considerations with Lie algebras we conclude that the following holds.

\begin{prp} The scheme $\MM$ is the scheme of orbits of the above action of $H$ on $G$. It is also isomorphic to the double coset $\MM\simeq G_{diag}\backslash (G\times G)/H$.

\end{prp}

\begin{rmk} Suppose $(C_i^{\prime})_{i\in \Z/m^{\prime}\Z}$ be another admissible wheel of cones such that for any $j$ there exists $i$ such that $C_j^{\prime}\subset C_i$. Then there is natural embedding $\MM_{\Gamma,\langle\bullet,\bullet\rangle,(C_i)_{i\in \Z/m\Z}}\to \MM_{\Gamma,\langle\bullet,\bullet\rangle,(C_i^{\prime})_{i\in \Z/m^{\prime}\Z}}$. Furthermore if under this embedding $\widehat{X}$ is mapped into $\widehat{X}^{\prime}$ then $\widehat{X}^{\prime}$ is obtained from $\widehat{X}$  in the following way:

a) first we make a finite sequence of blow-ups of $\widehat{X}$ with centers at some strata (such blow-ups will be automatically stratified);

b) then in the resulting formal scheme we choose a wheel of lines each of which is a one-dimensional stratum, and take the completion along the wheel.

Notice that this procedure depends on purely combinatorial data.

\end{rmk}

\subsection{Relation  to stability data}

Let $(\Gamma,\langle\bullet,\bullet\rangle)$ be a  lattice endowed with an  integer skew-symmetric  form $\langle\bullet,\bullet\rangle$. Stability data on the torus Lie algebra $\g=\g_{\Gamma-\Gamma_0}$ are given by a central charge $Z:\Gamma\to \C$ and a collection of numerical DT-invariants $\Omega(\gamma)\in \Q, \gamma\in \Gamma-\Gamma_0$ (see [KoSo1] or Section 2.3, Example 1). The Support Property ensures that 
$\overline{\cup_{\gamma\in Supp(\Omega)}\R\cdot \gamma}\cap (Ker\,Z_{\R}-\{0\})=\emptyset$, where $Z_{\R}:\Gamma_{\R}\to \C$ is the $\R$-linear extension of $Z$. We will assume in this subsection that $rk\,Z_{\R}=2$. In this  case $Z_{\R}^{\ast}((\R^2)^{\ast})\subset \Gamma_{\R}^{\ast}$ is an oriented two-dimensional vector space.

\begin{defn} An admissible wheel of cones $(C_i)_{i\in \Z/m\Z}$ is compatible with the central charge $Z:\Gamma\to \C, rk\,Z_{\R}=2$ if

a) for any $i\in \Z/m\Z$  the intersection $l_i=(C_{i,i+1}-\partial C_{i,i+1})\cap Z_{\R}^{\ast}((\R^2)^{\ast})$ is an open ray;

b) rays $l_i, i\in \Z/m\Z$ go in the clockwise order with respect to the orientation in $Z_{\R}^{\ast}((\R^2)^{\ast})$.

\end{defn}

\begin{prp} In the above notation there exists an admissible wheel of cones $C_i, i\in \Z/m\Z$ compatible with $Z$ and such that $Supp(\Omega)\subset \cup_iC_{i,i+1}^{\vee}-\{0\}\subset \Gamma_{\R}-Ker\,Z_{\R}$.

\end{prp}

{\it Proof.} Let us choose an isomorphism $\Gamma_{\R}\simeq \R^n=\R^2\oplus \R^{n-2}$ in such a way that $Z_{\R}$ becomes a projection $(x_1,...,x_n)\mapsto x_1+ix_2$. Recall the example of the admissible wheel of cones from Section 5.1. In that example $\cup_iC_{i,i+1}^{\vee}-\{0\}$ contains an open neighborhood of $\R^2-\{0\}$ in $\Gamma_{\R}$ and is disjoint from $\R^{n-2}=Ker\,Z_{\R}$.
It follows from the Support Property that for sufficiently large $t>0$ the set $\delta_t(\cup_iC_{i,i+1}^{\vee})$ contains $Supp(\Omega)$, where $\delta_t(x_1,x_2,...,x_n)=(x_1,x_2,tx_3,...,tx_n)$. Then the cones $(\delta_t^{\ast})^{-1}(C_i), i\in \Z/m\Z$ obtained by application of the map which is inverse of conjugate to $\delta_t$ form an admissible wheel of (non-rational) cones. Then taking a small perturbation we obtain an admissible wheel of rational cones. This completes the proof. $\blacksquare$

We will need a stronger statement proof of which is analogous but lengthy and hence omitted.

\begin{prp} For stability data on $\g_{\Gamma-\Gamma_0}$ as above there exist an admissible wheel of cones $(C_i)_{i\in \Z/m\Z}$ as in the Proposition 5.3.2 as well as the following data:

1) a cyclic decomposition $\R^2=V_{1,1}\cup...\cup V_{1,k_1}\cup V_{2,1}\cup...\cup V_{2,k_2}\cup...\cup V_{m,1}\cup...\cup V_{m,k_m}, m\ge 3, k_i\ge 1$, where $V_{i,j}$ are closed strict sectors such that two consecutive sectors have a common edge, $Z^{-1}(\partial V_{i,j}-\{0\})\cap \Gamma=\emptyset$;

2) a cyclically ordered collection of closed strict convex cones $C(V_{i,j})\subset \Gamma_{\R}$ compatible with $Z$ and such that $Z(C(V_{i,j}))\subset V_{i,j}$, the set $Supp(\Omega)$ belongs to $\cup_{i,j}C(V_{i,j})$, and  for any $i,j$ the set $C(V_{i,j})-\{0\}$ belongs to $int(C_{i,i+1}^{\vee})$.
\end{prp}

Let us make a choice of sectors and cones as in the Proposition. Then our stability data on $\g_{\Gamma-\Gamma_0}$ give rise to the collection of elements 
$$g_{i,j}\in exp(\prod_{\gamma\in C(V_{i,j})\cap (\Gamma-\Gamma_0)}\Q\cdot e_{\gamma})\subset G_{i,i+1},$$
 where the latter group was defined in the end of the previous subsection. Let us associate with our stability data a point of $\MM$ represented by the coset of the element $(g_{1,1}g_{1,2}...g_{1,k_1},g_{2,1}...g_{2,k_2},...,g_{m,1}...g_{m,k_m})$.

\begin{thm} For given $\Gamma, \langle\bullet,\bullet\rangle, Z:\Gamma\to \C$ with $rk\,Z_{\R}=2$ the above map provides a bijection from the set of stability data on $\g_{\Gamma-\Gamma_0}$ with fixed central charge $Z$ to the set 
$\varinjlim \MM_{\Gamma,\langle\bullet,\bullet\rangle,(C_i)_{i\in \Z/m\Z}}$, where the inductive limit is taken with respect to subdivision maps described in the Remark 5.2.8 over all admissible wheel of cones $(C_i)$ compatible with $Z$.

\end{thm}

{\it Sketch of the proof.} In what follows all chains of cones $(C_i)$ will be compatible with the central charge $Z$.
The proof will consist of several steps.

Step 1. The moduli space $\MM_{\Gamma,\langle\bullet,\bullet\rangle,(C_i)_{i\in \Z/m\Z}}:=\MM_{(C_i)}$ can be defined for admissible wheels of {\it non-rational cones} via the double coset construction.

Step 2. For a special choice of cones $C_i$ we can identify the space of stability data on the Lie algebra $\g_{\Gamma-\Gamma_0}$ with the central charge $Z$ and $Supp(\Omega)\subset \cup_iC_{i,i+1}^{\vee}$ with $\MM_{(C_i)}$.  It can be done along the lines of the example with polyhedron in Section 5.1. More precisely, let us fix a convex polygon in $\R^2$ which contains the origin in the interior and has cyclically ordered vertices $v_1,v_2,...,v_m$. Let us fix a convex bounded closed polyhedron $P\subset \R^{n-2}$ which contains the origin in the interior. Let us define a decomposition $(\R^2)^{\ast}-\{0\}=\cup_{i\in \Z/m\Z}V_i$ into the union of semiclosed strict sectors defined by the condition $V_i=\{u\in (\R^2)^{\ast}|u(v_i)>u(v_{i-1}),u(v_i)\ge u(v_{i+1})\}$. We define  strict convex cones $C(V_i)\subset (v_i\oplus P)^{\vee}$ as $\{u\oplus w|u\in V_i, w\in P\}$. Then $$\cup_i(v_i\oplus P)^{\vee}-\{0\}=\sqcup_iC(V_i).$$

For such a choice  we have $C_{i,i+1}=\R_{\ge 0}(v_i\oplus P)$ and $C_i$ is the convex hull of $C_{i-1,i}, C_{i,i+1}$.

Below we will define a bijection between the space of stability data on the Lie algebra $\g_{\Gamma-\Gamma_0}$ with the central charge $Z$ and $Supp(\Omega)\subset \cup_iC_{i,i+1}^{\vee}$ with $\MM_{(C_i)}$.
 
First, let us introduce a collection of pronilpotent groups $G_i^{(1)}\subset G_i$ such that $Lie(G_i^{(1)})=
\prod_{\gamma\in C(V_i)\cap (\Gamma-\Gamma_0)}\Q\cdot e_{\gamma}$. Then following [KoSo1] we parametrize our stability data by the collection of elements $A_{V_i}\subset G_i^{(1)}, i\in \Z/m\Z$. Namely, we set $A_{V_i}=\prod^{\rightarrow}_{l\subset V_i}A_l$, where the product is taken in the clockwise order over the set of rays in $V_i$  with vertex at $0$, each factor is given by the formula
$A_l=\prod_{\gamma\in C(V_i)\cap (\Gamma-\Gamma_0), Z(\gamma)\in V_i}T_{\gamma}^{\Omega(\gamma)}$, and $T_{\gamma}:e_{\mu}\mapsto (1-e_{\gamma})^{\langle \gamma,\mu\rangle}$.

Then in the double coset description  of $\MM_{(C_i)}$ we take the orbit of the element $(id_G,A_{V_1},A_{V_2},...,A_{V_m})\in G\times G$. It is easy to see that this gives the desired bijection (cf. Prop. 5.2.5). In fact we have
$$\prod_iG_i^{(1)}\simeq G_{diag}\backslash(G\times G)/H\simeq \MM.$$

Step 3. Let us introduce a partial order on the set of wheels ${\mathcal C}=(C_i)$ compatible with $Z$. Namely we say that ${\mathcal C}^{\prime}=(C_j^{\prime})\le {\mathcal C}=(C_i)$ if for any $i$ there exists $j$ such that $C_{j,j+1}^{\prime}\subset C_{i,i+1}$ (equivalently, for any $j$ there exists $i$ such that $C_j^{\prime}\subset C_i$). The partial order $\le$ gives rise to the category with objects ${\mathcal C}=(C_i)$ such that $int(C_i)\cap \R^2$ is non-empty, and morphisms defined by the partial order (poset category).  Then we observe that if ${\mathcal C}^{\prime}\le {\mathcal C}$ then we have a natural embedding $\MM_{(C_i)}\to \MM_{(C_i^{\prime})}$ (notice that we do not need cones to be rational for all that).

Step 4. Assume that ${\mathcal C}^{\prime}\le {\mathcal C}$, and for any $i$ there exists $j$ such that $C_{j,j+1}^{\prime}=C_{i,i+1}$ and also we have $\cup_jC_j^{\prime}=\cup_iC_i$. Then the embedding from Step 3 is an isomorphism of affine schemes.

Step 5.  For an admissible wheel of cones ${\mathcal C}=(C_i)$ there exist ${\mathcal C}^{\prime}\le {\mathcal C}$ and ${\mathcal C}^{\prime\prime}\ge {\mathcal C}^{\prime}$ such that the conditions from Step 4 hold and ${\mathcal C}^{\prime\prime}$ is a wheel of cones from Step 2.

Let us comment on Step 5. In order to find ${\mathcal C}^{\prime\prime}$ one chooses the vertices $v_i, 1\le i\le m$ of the polygon in Step 2 in such a way that $v_i\notin C_{i,i+1}\cap \R^2$. Then one replaces the polyhedron $P$ from Step 2 by $\varepsilon P$, where 
$\varepsilon$ is a sufficiently small positive number.

Step 6. By previous steps an element from $\MM_{\mathcal C}$ gives an element from $\MM_{{\mathcal C}^{\prime\prime}}$, hence the stability data on $\g_{\Gamma-\Gamma_0}$ by Step 2.
This concludes the sketch of the proof.  $\blacksquare$

Assume we are given $\Gamma, \langle\bullet,\bullet\rangle$. Consider a continuous family of central charges $Z_x, x\in X$, where $X$ is a Hausdorff topological space and such that $rk\,Z_x=2$ for all $x\in X$. Consider a family $\sigma_x, x\in X$ (non-necessarily continuous) of stability data on $\g_{\Gamma-\Gamma_0}$ with central charges $Z_x$. Then we have the following result proof of which is omitted.

\begin{prp} The family $\sigma_x, x\in X$ is continuous if and only if there exists an open covering $X=\cup_{\alpha}U_{\alpha}$, collection of admissible wheels of cones ${\mathcal C}_{\alpha}=(C_{\alpha,i})$ and points $m_{\alpha}\in \MM_{{\mathcal C}_{\alpha}}$ such that for any $x\in U_{\alpha}$ the stability condition corresponding to $m_{\alpha}$ and having central charge $Z_x$ is identified by Theorem 5.3.3 with $\sigma_x$.

\end{prp}

\begin{rmk} In the case $rk\,Z_{\R}=1$ one can develop a similar theory by replacing the formal neighborhood of a wheel of lines by the one for chains of lines (in some interesting cases just one projective line is enough).

\end{rmk}

\subsection{Toric-like compactifications}

Let $\NN$ be a smooth algebraic variety over a field $k$ of characteristic zero, and $\NN_1$ be a normal scheme over $k$ which contains $\NN$ as an open subscheme.

\begin{defn} We say that $\NN_1$ is a toric-like compactification of $\NN$ if the pair $(\NN_1,D)$, where $D=\NN_1-\NN$ is a reduced divisor, is a local formal toric pair. 

\end{defn}

In what follows we  assume that $k=\C$, although this assumption can be relaxed.
Notice that $\NN_1$ does not have to be proper. It follows from the definition that $D$ is stratified and its stratification is compatible with the local toric picture.

With each $0$-dimensional stratum (point) $x\in D$ we can associate a free abelian group (lattice), which in the obvious notation can be written as $\Gamma_x=H^1((x\to \NN_1)^{\ast}(\NN\to \NN_1)_{\ast}\Z_{\NN})$ where $\Z_{\NN}$ is the constant sheaf and the arrows are natural embeddings. In the case $k=\C$ and analytic topology $\Gamma_x=H^1(B_x\cap \NN,\Z)\simeq \Z^{dim\,\NN}$, where $B_x$ is a small ball with the center in $x$.  

Let $\overline{F}$ be a closed $1$-dimensional stratum whose complement to the union of $0$-dimensional strata is isomorphic to $\C^{\ast}$. 
Then one has a similarly defined lattice $\Gamma_{\overline{F}}\simeq \Z^{dim\,\NN}$. In the complex analytic case one can define $\Gamma_{\overline{F}}=H^1(U_{\overline{F}}\cap \NN,\Z)$, where $U_{\overline{F}}$ is a tubular neighborhood of $\overline{F}$. If $x\in \overline{F}$ then by topological reasons we have a canonical isomorphism $\Gamma_x\simeq \Gamma_{\overline{F}}$.

\begin{defn} A chain of lines $\overline{F}_k, \overline{F}_k\ne \overline{F}_{k+1}, 1\le k\le m$ 
is given by a sequence of $1$-dimensional strata with parametrization $i_k:{\bf P}^1\simeq \overline{F}_k$ such that $i_{k+1}(0)=i_{k}(\infty), 0\le k\le m-1$, and $i_k(0), i_k(\infty)$ are the only $0$-dimensional strata of $\overline{F}_k$. 

\end{defn}
We will use the notation $p_k=i_k(0), k=1,...,m, p_{m+1}=i_{m}(\infty)$.

The above considerations give us a chain of canonical isomorphisms of lattices $\Gamma_{p_1}\simeq \Gamma_{\overline{F}_1}\simeq \Gamma_{p_2}\simeq \Gamma_{\overline{F}_2}\simeq ...\simeq\Gamma_{\overline{F}_m}\simeq \Gamma_{p_{m+1}}$. We denote the identified lattices by $\Gamma$. We assume that $\Gamma$ is endowed with a skew-symmetric pairing $\langle\bullet,\bullet\rangle:\bigwedge^2 \Gamma\to \Z$. We denote the kernel of this form by $\Gamma_0$.

 Suppose that we are given an automorphism $T:\NN\to \NN$. Suppose furthermore that there exist open subsets $\NN\subset U_1\subset \NN_1$ and $\NN\subset U_{m+1}\subset \NN_1$ such that $p_1\in U_1, p_{m+1}\in U_{m+1}$ and such that $T$ extends to an isomorphism $\overline{T}:U_1\to U_{m+1}$. Then $\overline{T}$ induces an isomorphism $\overline{T}_{1,m+1}:\Gamma_{p_1}\simeq \Gamma_{p_{m+1}}$.

Next we would like to formulate a list of assumptions under which we will construct a point of the moduli space of decorated formal Poisson schemes. Namely we assume that:

a) $\NN$ is endowed with a Poisson structure.

b) Automorphism $T$ preserves the Poisson structure.

c) The isomorphism $T_{1,m+1}$ coincides with the one obtained from the chain of isomorphisms of lattices.

d) We are given a homomorphism of abelian groups $c: \Gamma_0\to {\mathcal O}(\NN)^{\times}$ whose image belongs to the Poisson center.

e)  For any $1\le i\le m+1$ there exists an isomorphism $\phi_i$ of the pair 
$(\widehat{\NN}_{1,p_i},\widehat{D}_{p_i})$ (completions at $p_i$) with the corresponding pair in the local formal toric model. This is an isomorphism of formal Poisson schemes, where the local formal toric model is endowed with the Poisson structure given by the skew-symmetric form $\langle\bullet,\bullet\rangle_i:\bigwedge^2\Gamma_{p_i}\to \Z$ obtained from the skew-symmetric form $\langle\bullet,\bullet\rangle$ via the canonical isomorphism $\Gamma_{p_i}\simeq \Gamma$.

f) For any $1\le i\le m+1$ we are given an isomorphism $g_{p_i}:\TT_{p_i}\simeq \TT_{can}$ such that for any $1$-dimensional closed stratum $\overline{F}\simeq \C{\bf P}^1$ containing exactly two $0$-dimensional strata $\{x_0\},\{x_{\infty}\}$ the following two compositions coincide:

$$(\Gamma_{x_0},\TT_{x_0})\buildrel{(id,g_{x_0})}\over{\simeq}(\Gamma,\TT_{can})\buildrel{(id,g_{x_{\infty}}^{-1})}\over{\simeq}(\Gamma_{x_{\infty}},\TT_{x_{\infty}})$$

and

$$(\Gamma_{x_0},\TT_{x_0}){\simeq}
(\Gamma_{\overline{F}},\TT_{\overline{F}}){\simeq}(\Gamma_{x_{\infty}},\TT_{x_{\infty}}).$$

g)  For any $1\le i\le m+1$ there exists $\phi_i$ (see e)) such that the pull-back
$\phi_i^{\ast}(c(\gamma))$ is a function of weight $\gamma$ on the open toric stratum.

h) Let $C_i\subset \Gamma_{p_i}^{\vee}\otimes \R, 1\le i\le m$ be closed strict rational convex cones arising from the toric-like stratification. Then (after identification $T_{1,m}: C_1\simeq C_m$) we obtain an admissible wheel of cones.

i) The composition
$$\TT_{p_1}\buildrel{g_{p_1}}\over{\simeq}\TT_{can}
\buildrel{g_{p_{m+1}}^{-1}}\over{\simeq}\TT_{p_{m+1}}$$
coincides with the Poisson isomorphism induced by $\overline{T}$.

Under the above assumptions a)-i) let us consider the disjoint union of completions $\sqcup_{1\le i\le m}\widehat{\NN}_{1,\overline{F}_i}$ and then identify the formal neighborhood of $\infty_i\in \overline{F}_i$ with the one of $0_{i+1}\in \overline{F}_{i+1}$ for $1\le i\le m-1$ using the embeddings to $\NN$, and finally the formal neighborhood of $\infty_m\in \overline{F}_m$ with the one of $0_{1}\in \overline{F}_{1}$ using the isomorphism $\overline{T}_{1,m+1}$. In this way we obtain an admissible wheel of cones endowed with additional data giving us a point in the moduli space $\MM_{\Gamma,\langle\bullet,\bullet\rangle,(C_i)_{1\le i\le m+1}}$.

\section{WCS and Mirror Symmetry}

Considerations in this section will be mostly heuristic.  We are going to explain how the ideas of Mirror Symmetry in Strominger-Yau-Zaslow (SYZ for short) torus fibration picture can be combined with previous considerations of this paper. This will give us a WCS which conjecturally should coincide with the one constructed in Section 4.

\subsection{Reminder on Fukaya categories}

 Let $(X,\omega)$ be a compact smooth symplectic manifold of dimension $2n$ and ${\bf B}\in H^2(X,\R/2\pi \Z)\simeq Hom(H_2(X,\Z),\R/2\pi \Z)$ be the $B$-field. It is expected that for a sufficiently large $\lambda>0$ the triple $(X,\lambda\omega,\B$) gives rise to a $\Z/2\Z$-graded $A_{\infty}$-category ${\mathcal F}(X,\lambda\omega,\B)$ called the Fukaya category.  Some objects of the Fukaya category are pairs $(L,E)$ where $L\subset X$ is an oriented Lagrangian submanifold endowed with a spin structure such that $\B_{|L}=0$, and $E$ is a $U(1)$-local system on $L$. Space of morphisms $Hom((L_1,E_1),(L_2,E_2))$ is labeled by intersection points of $L_1$ and $L_2$. In order to define the $A_{\infty}$-structure one needs to 
choose an almost complex structure on $X$. Then higher composition maps are given by a properly defined count of pseudo-holomorphic discs $D$ ``weighted`` by $e^{-\int_{D}(-\lambda \omega+i\B)}$.
Not every Lagrangian submanifold $L$ can support an object of the Fukaya category. A sufficient condition for that is the absence of pseudo-holomorphic discs of Maslov index $2$ such that $\partial D\subset L$. More advanced picture which handles this problem was developed in [FOOO]. Convergence of series which defines higher composition maps is another big issue, which is not proved by this time. Typically people avoid this problem by working over Novikov ring of series in the above-mentioned weight.  In this case the approach can be made rigorous (see [FOOO]). Furthermore, in the presence of a top degree almost complex form one can define a $\Z$-graded version of the Fukaya category. In what follows we will assume that for ``sufficiently large $\omega$'' we can define the Fukaya category  ${\mathcal F}(X,\omega,\B)$ over the field of complex numbers. In some sense this category depends holomorphically on $[\omega]+i\B\in H^2(X,\C/2\pi i\Z)$. 

Under appropriate conditions one can define a version of ${\mathcal F}(X,\omega,\B)$ called {\it wrapped Fukaya category} in case when $X$ is non-compact and endowed with a proper map $H:X\to [0,+\infty)$ (see [AbSe]). In that case Lagrangian submanifolds supporting objects are non-compact, having ``good'' behavior ``at infinity''.  The space $Hom((L_1,E_1),(L_2,E_2))$ is defined by means of intersection points $L_1\cap exp(t\varphi_H)(L_2)$, where $\varphi_H=\{H,\bullet\}$ is the Hamiltonian vector field corresponding to $H$.

\subsection{SYZ picture and integrable systems}

Now suppose  that we have a real integrable system $\pi: (X,\omega)\to B$ with a Lagrangian section $s:B\to X$. Then over an open dense subset $B^0$ of a smooth $n$-dimensional manifold $B$ we have a Lagrangian torus fibration with marked zero points in the fibers. If $X$ is non-compact then we assume that $X$ is endowed with a proper function $H:X\to [0,+\infty)$ mentioned in the previous subsection, and that the function is a pull-back of a similar function on $B$. We will impose the condition $c_1(T_X)=0\in H^2(X,\Z)$, although it is not necessary  for some of our considerations.

The open submanifold $B^0$ carries a $\Z$-affine structure with local affine coordinates $x_i, 1\le i\le n$ in a neighborhood of $b\in B^0$ which are determined (up to a shift) by the condition $dx_i=\int_{\gamma_i}\omega$, where $\{\gamma_i\}_{1\le i\le n}$ is a basis in $\Gamma_b=H_1(\pi^{-1}(b),\Z)$. The vector fields $\partial/\partial x_i$ generate a covariantly constant lattice $T_{\Z}\subset TB^0$.
We will assume that $B$ carries a metric $g_B$ which is complete and which is Riemannian on $B^0$.  
Then it gives rise to an isomorphism $T^{\ast}B^0\simeq TB^0$ of the tangent and cotangent bundles.
In particular we have a lattice $T_{\Z}^{\vee}\subset T^{\ast}B^0$.
Notice that $\pi^{-1}(B^0)$ is canonically symplectomorphic to the total space of the Lagrangian torus fibration $(T^{\ast}B^0/T_{\Z}^{\ast},\omega_{T^{\ast}B^0/T_{\Z}^{\ast}}) $. The latter is in turn symplectomorphic to the total space of the ``rescaled'' torus bundle $(T^{\ast}B^0/\varepsilon T_{\Z}^{\ast}, \varepsilon^{-1}\omega_{T^{\ast}B^0/T_{\Z}^{\ast}})$. Using the metric we endow $\pi^{-1}(B^0)$ with an almost complex structure $J_{\varepsilon}$ compatible with the rescaled symplectic form $\varepsilon^{-1}\omega$ and the pull-back of the metric $g_B$. The structure $J_{\varepsilon}$ is singular on $X-\pi^{-1}(B^0)$. We can replace it by a non-singular almost complex structure $\overline{J}_{\varepsilon}$ compatible with $\varepsilon^{-1}\omega$, and such that $\overline{J}_{\varepsilon}=J_{\varepsilon}$ outside of a small $\delta=\delta(\varepsilon)$- neighborhood of $X-\pi^{-1}(B^0)$ such that $lim_{\varepsilon\to 0}\delta(\varepsilon)=0$.

Then as $\varepsilon\to 0$ the $\overline{J}_{\varepsilon}$-holomorphic curves converge to singular surfaces whose $\pi$-images are graphs in $B$ with edges which are $g_B$-gradient lines of affine functions (see a discussion of this result of Fukaya and Oh in [KoSo6]). At a vertex $v$ of the gradient graph the balancing condition is satisfied: $\sum_i\gamma_v^{i}=0$, where $\gamma^{i}_v$  denote adjacent  to $v$ edges which are identified with the corresponding integer affine functions.

In what follows we will assume that $dim_{\R}(B-B^0)\ge 2$ (cf. [KoSo6]). This condition is closely related to the condition $c_1(T_X)=0$. 
Then for a generic point $b\in B^0$ there is no gradient tree as above with the root at $b$ and external vertices at $B-B^0$. Such trees correspond to limits of $\overline{J}_{\varepsilon}$-holomorphic discs with boundaries on $\pi^{-1}(b)$. Indeed, the union of roots of such trees is a union of countably many hypersurfaces in $B$. We can call them ``walls".  The reader should  not mix them with walls in WCS for complex integrable systems (see Section 10 for discussion of these walls). Informally we can think that $B^0$ locally looks as a locally finite union of convex polyhedral domains separated by  walls and each polyhedral domain $P$ gives a family of objects in the Fukaya category parametrized by the tube domain $Log^{-1}(P)\subset (\C^{\ast})^n$, where $Log: (\C^{\ast})^n\to \R^n$ is the tropical map $(z_1,...,z_n)\mapsto (log|z_1|,...,log|z_n|)$. Then $Arg(z_i)$ correspond to $U(1)$-local systems on fibers of $\pi$. 
According to the Mirror Symmetry philosophy there is a complex variety $X^{\vee}$ (mirror dual to $X$) containing all  parametrized families of objects of  the Fukaya category as open subsets. Crossing a wall which separates polyhedral domains corresponds to a change of coordinates on $X^{\vee}$. If $g_B$ is locally given in affine coordinates by the Hessian matrix $(\partial^2H/\partial x_i\partial x_j)$ for some convex function $H$, then edges of gradient graphs will be $\Z$-affine segments in the dual affine structure on $B^0$.

\subsection{The case of complex integrable systems}

First, let us assume that we are given a polarized (full) integrable system $\pi:(X,\omega^{2,0})\to B$ endowed with a holomorphic Lagrangian section $s:B\to X$. Notice that $codim_{\R}(B-B^0)\ge 2$ automatically. Let us fix $\zeta\in \C^{\ast}$ and take $\omega_{\zeta}=Re(\zeta^{-1}\omega^{2,0})$ as the real symplectic form on $X$. As the $B$-field we take $\B_{\zeta}=Im(\zeta^{-1}\omega^{2,0})+\B_{can}$, where ${\bf B}_{can}\in H^2(X,\pi\Z/2\pi \Z)\simeq H^2(X,\Z/2\Z)$ is a ``canonical'' $B$-field defined in the Appendix.

\begin{rmk} Our choice of $\B_{can}$ is motivated by the appearance of the factor $(-1)^{\langle\gamma,\mu\rangle}$ in the theory of DT-invariants. 

\end{rmk}

Recall (see Section 4.7) that under some natural assumptions  there exists a proper continuous function $H_{\zeta}:B^0\to [0,+\infty)$ which gives the metric on the base.

Tube domains for mirror duals $X_{\zeta}^{\vee}$ in coordinate-free language belong to $Hom(\Gamma_b,\C^{\ast}), b\in B^0, \Gamma_b=H_1(\pi^{-1}(b),\Z)$. Hence they are endowed with a symplectic structure associated with the polarization. The above-discussed changes of coordinates preserve the symplectic structure. Hence each $X_{\zeta}^{\vee}, \zeta\in \C^{\ast}$ is a holomorphic symplectic manifold. One can explain the symplectic structure in a different way.  Indeed the polarization on fibers of ${X} \to B$ gives rise to a canonical holomorphic line bundle ${\mathcal L}$ on $X$ whose restriction on fibers ${X}\to B$ is ample (more precisely it is defined outside of the preimage of the discriminant locus, but we expect that it extends to the whole space). The cohomology class $c_1(\LL)\in H^2({X},\Z)$ can be interpreted as an element of the second Hochschild cohomology of the above wrapped Fukaya category  $\mathcal F(X, \omega_{\zeta}, {\bf B}_{\zeta})$. Equivalently it is a second Hochschild cohomology class of the derived category of coherent sheaves on the mirror dual $X_{\zeta}^{\vee}$. One can argue that this cohomology class is represented by a non-degenerate holomorphic Poisson bivector field. Its inverse is our holomorphic symplectic form.

Now we are ready to consider the semipolarized case. To simplify the exposition we fix $\zeta=1$ and omit the $B$-field ${\bf B}$ from the notation.

Namely, let  $\pi:(X,\omega^{2,0})\to B$ be a semipolarized integrable system with central charge $Z$ and holomorphic Lagrangian section. We assume that the monodromy of the local system $\underline{\Gamma}_0\to B^0$ is trivial so all its fibers can be identified with the fixed lattice $\Gamma_0$ (this can always be achieved by taking a finite cover, see Lemma 4.2.3). Then we obtain a holomorphic family of complex integrable systems $(X_{Z_0},\omega^{2,0}_{X_{Z_0}})\to B_{Z_0}$ parametrized by  $Z_0\in Hom(\Gamma_0,\C)$.  Our discussion of the Fukaya categories make  plausible  the proposal that the holomorphic family of the Fukaya categories  gives rise to a  holomorphic family  of mirror duals ${X}^{\vee}_{Z_0}:=(X_{Z_0}, Re(\omega^{2,0}_{X_{Z_0}}))^{\vee}$  parametrized by $Hom(\Gamma_0,\C)$. We argue that ${X}^{\vee}_{Z_0}$ carry a holomorphic symplectic form.

The total space of this family will be denoted by $X^{\vee}$. 

In fact it is a pull-back via the map $exp:Hom(\Gamma_0,\C)\to Hom(\Gamma_0,\C^{\ast})$ of an {\it algebraic} family of smooth complex symplectic varieties $X^{\vee,alg}\to Hom(\Gamma_0,\C^{\ast})$.

Here is an informal argument in favor of that. For each $Z_0=Z_{|\Gamma_0}$ we have the corresponding polarized integrable system $(X_{Z_0},\omega^{2,0}_{X_{Z_0}})\to B_{Z_0}$ as discussed before. Consider the real affine space $\alpha_{X_0}\subset Hom(\Gamma_0,\C)$ defined by the condition $Re(Z_0)=X_0$. 
This gives rise to a family of smooth symplectic manifolds $(X_{Z_0},Re(\omega^{2,0}_{X_{Z_0}}))$, where $Z_0\in \alpha_{X_0}$. The cohomology class $[Re(\omega^{2,0}_{X_{Z_0}})]$ is locally constant along $\alpha_{X_0}$. Then by a Moser-type theorem (suitably adopted to non-compact case) we conclude that all symplectic manifolds $(X_{Z_0},Re(\omega^{2,0}_{X_{Z_0}}))$ can be (non-canonically) identified up to symplectic isotopy. Notice that the $B$-field depends on $Im(Z_0)$. Therefore the corresponding Fukaya categories (and therefore their mirror duals) depend only on $Z$ modulo $Hom(\Gamma_0,2\pi i\Z)$.  The corresponding Fukaya categories (and hence their mirror duals) are periodic, with respect to the shifts, hence form a holomorphic family over the torus  $Hom(\Gamma_0,\C^{\ast})$. The algebraicity of this family is a conjecture, which we will discuss  below in Section 6.5. A priori the dual variety is just a complex manifold without an algebraic structure. The latter comes from additional considerations related to the {\it wrapped} Fukaya category.

From the above discussion we conclude that the total space of the family of mirror duals is a complex Poisson variety $X^{\vee,alg}$ (which we will call the mirror dual to our complex integrable system) endowed with a Poisson map to $Hom(\Gamma_0,\C^{\ast})$ (hence fibers of this map are symplectic leaves). We will use the term ``mirror dual'' being applied to the holomorphic family $X^{\vee}\to Hom(\Gamma_0,\C)$.

There is an alternative approach to the construction of $X^{\vee,alg}$ which we are going to explain below. Let $B_{\R}\subset B$ be the closure of the set $\{b\in B^0|Re\,(Z_b)_{|\Gamma_0}=0\}$. Let $X_{\R}=\pi^{-1}(B_{\R})\subset X$. Then $X_{\R}$ is a coisotropic submanifold. Therefore it carries a foliation with symplectic quotient which we will denote by $X_{\R}^{\prime}$. We have the corresponding real integrable system $((X_{\R}^{\prime},\omega_{X_{\R}^{\prime}})\to B_{\R}$. The fiber over $b\in B_{\R}^0:=B_{\R}\cap B^0$ is isomorphic to the compact real torus $(\Gamma_b\otimes \R)/\Gamma_b$. Notice that the cohomology class $[\omega_{X_{\R}^{\prime}}]=0$ (i.e. our real integrable system is exact). We can define the Fukaya category ${\mathcal F}(X_{\R}^{\prime},\omega_{X_{\R}^{\prime}},\B^{\prime}_{can})$, where $\B^{\prime}_{can}\in H^2(X_{\R}^{\prime},\pi \Z/2\pi \Z)$ is the canonical $B$-field associated with the integer skew-symmetric form on $\Gamma_b$.  Since $(X_{\R}^{\prime},\omega_{X_{\R}^{\prime}})$ is an exact symplectic manifold and the pairing of $exp(\B^{\prime}_{can})$ with the class of any pseudo-holomorphic curve belongs to $exp(\pi i\Z)=\{-1,+1\}\subset \Q$, we conclude that ${\mathcal F}(X_{\R}^{\prime},\omega_{X_{\R}^{\prime}},\B^{\prime}_{can})$ is defined over $\Q$ (again, we ignore here the convergence problem). The torus $(\Gamma_0\otimes \R)/\Gamma_0$ acts on $X_{\R}^{\prime}$ in the Hamiltonian way preserving Lagrangian torus fibers of the projection to $B_{\R}$. By Mirror Symmetry this corresponds to the holomorphic map $(X_{\R}^{\prime})^{\vee}\to Hom(\Gamma_0,\C^{\ast})$. 

\begin{conj}
The complex  Poisson manifold $(X_{\R}^{\prime})^{\vee}$ is algebraic and endowed with the map to $Hom(\Gamma_0,\C^{\ast})$.  It is canonically isomorphic to the complex algebraic Poisson variety $X^{\vee,alg}$ with its map to $Hom(\Gamma_0,\C^{\ast})$.
\end{conj}

So far we have been discussing semipolarized integrable systems with fixed holomorphic symplectic form. Let us consider the $\C^{\ast}$-family of holomorphic symplectic forms $\omega_{\zeta}^{2,0}=\omega^{2,0}/\zeta$ on $X$. Then the corresponding mirror dual Poisson varieties $X_{\zeta}^{\vee,alg}, \zeta\in \C^{\ast}$ form a local system of quasi-affine algebraic varieties over $\C^{\ast}$ (we will discuss it in Section 6.5). More precisely, recall that mirror duals were constructed first by fixing $Z_0=Z_{|\Gamma_0}$ and then by looking at the result as a family over either $Hom(\Gamma_0,\C)$ (this gives us a holomorphic family) or $Hom(\Gamma_0,\C^{\ast})$ (this gives us an algebraic family). After the rescaling the symplectic form to $\omega^{2,0}/\zeta$, the central charge $Z$ gets replaced by $Z/\zeta$, hence $Z_0=Z_{|\Gamma_0}$ gets replaced by $Z_0/\zeta$. Hence in the construction of mirror duals we fix $Z_0/\zeta$. Taking the union of mirror duals (for fixed $\zeta$) we obtain the mirror dual Poisson variety $X_{\zeta}^{\vee,alg}$. In fact they form a holomorphic local system of algebraic varieties over $\C^{\ast}$. 
This can be proved  by using of the Moser-type arguments (in case when we deal with polarized integrable systems having central charge  one can identify the fibers directly). In a similar way we obtain a holomorphic family over $\C^{\ast}$ of complex Poisson manifolds $X_{\zeta}^{\vee}$ each of which is endowed with a holomorphic Poisson morphism to $Hom(\Gamma_0,\C)$. Clearly the total space $X^{\vee}$ of the latter family is the universal cover of the total space of the former family $X^{\vee,alg}$.

Taking the fiber $X_1^{\vee,alg}:=X_{\zeta=1}^{\vee,alg}$ we obtain a Poisson variety endowed with a Poisson automorphism $T:X_1^{\vee,alg}\to X_1^{\vee,alg}$, which is equal to $id$ on the algebra of central functions (the latter is isomorphic to ${\mathcal O}(Hom(\Gamma_0,\C^{\ast}))$) (cf. with Section 5.4).

Finally, let us remark that if the monodromy of the local system $\underline{\Gamma}_0\to B^0$ is a finite group $G$ then the above considerations still work and give us the mirror dual Poisson variety $X^{\vee,alg}$ together with a Poisson morphism to $Hom(\Gamma_0,\C^{\ast})/G$ (here $\Gamma_0$ can be thought of as a fiber of the pull-back of $\underline{\Gamma}_0$ to the universal cover).

\subsection{Wall-crossing structure from the point of view of Mirror Symmetry}

Recall $g_B$-gradient trees from Section 6.2. In the case of semipolarized integrable systems with  central charge and holomorphic Lagrangian section we have K\"ahler metrics on the bases of the corresponding polarized integrable systems, hence edges of the gradient trees are straight segments in the dual $\Z$-affine structure (see Section 6.2).  In terms of the central charge, the dual affine structure for the symplectic form $Re(\omega^{2,0}/\zeta), |\zeta|=1$ is given by $Y_{\theta}:=Im(e^{-i\theta}Z)$ with fixed restriction of $Y_{\theta}$ to ${\Gamma_0}$, where $\zeta=e^{i\theta}\in \C^{\ast}$. As we briefly recalled in Section 6.2, the SYZ approach to Mirror Symmetry (see more on that in [KoSo6]) gives rise to an inductive procedure of constructing walls and changes of coordinates, starting with certain data assigned to generic points of the discriminant $B-B^0$. Namely, for a point $b\in B^0$ sufficiently close to a generic point of the discriminant, one counts limiting pseudo-holomorphic discs whose $\pi$-image on the base is a short gradient segment connecting the point $b$ with a point of $B-B^0$. 

The inductive procedure is a priori different from the one discussed above in Section 4. Nevertheless in this case one can prove by induction (moving along the oriented gradient tree from the discriminant to a given point) that the walls and the changes of coordinates in Mirror Symmetry story of Section 6.2 coincide with those in Section 4. In particular, the changes of coordinates preserve the Poisson structure on $X_{Z_0}^{\vee}, Z_0\in Hom(\Gamma_0,\C)$ and depend algebraically on the point of $Hom(\Gamma_0,\C^{\ast})$. They can be interpreted as Poisson transformations of $X^{\vee}$ identical on the Poisson center.

The above discussion gives an alternative approach to WCS constructed in Section 4. The initial data for which $\Omega(\gamma)=1$ (see Section 4) for $A_1$-singularities correspond to the count of pseudo-holomorphic discs in the standard $A_1$-singularity model (see e.g. [KoSo2], [ChLaLeu]).

\subsection{Algebraicity of the mirror dual}

For simplicity we are going to discuss the case of polarized integrable systems. We hope that our arguments can be extended to the semipolarized case. In particular, the basis described below should be a basis in the algebra over ${\mathcal O}(Hom(\Gamma_0,\C^{\ast}))$. One can speculate that it coincides with the canonical bases expected in the theory of cluster varieties (see [FoGo1]).

As we discussed before, the mirror dual  $X^{\vee}$ to an exact {\it real} integrable system $\pi:(X,\omega)\to B$ endowed with Lagrangian section $s:B\to X$ and the $B$-field which is a $2$-torsion, is an algebraic variety defined over $\Q$. More precisely, we expect that $X^{\vee}$ is a quasi-affine (maybe formal) scheme of finite type over $\Z$. In case if there exists a proper continuous function $H:B\to [0,+\infty)$ which is (strictly) convex with respect to the $\Z$-affine structure on $B^0$ and has ``good'' behavior at the discriminant $B-B^0$, we expect that $X^{\vee}$ will be a Zariski open in the spectrum of a finitely generated algebra $R$, which can be described such as follows.

For any $t\in \R$ let $L_t=exp(t\varphi_H(s(B)))\subset X$ be a Lagrangian submanifold obtained from the zero section $s(B)$ by the Hamiltonian shift along the vector field $\varphi_H=\{H,\bullet\}$. Morally all $L_t,t\in \R$ should correspond to isomorphic objects in the wrapped Fukaya category. More precisely, for $t_1<t_2$ let us consider the basis of the Floer complex $CF(L_{t_1},L_{t_2})$ given by intersection points $L_{t_1}\cap L_{t_2}$. We will assume that this intersection belongs to $\pi^{-1}(B^0)$. Convexity of $H$ implies that Maslov indices of all intersection points are zero. Hence the Floer differential is trivial. The composition
$$m_{t_1,t_2,t_3}: CF(L_{t_1},L_{t_2})\otimes CF(L_{t_2},L_{t_3})\to CF(L_{t_1},L_{t_3})$$
sends the tensor product of two basis elements to a {\it finite} $\Z$-linear combination of basis elements (this follows from the ``energy considerations'' with the function $H$). Hence we obtain a directed $\Z$-linear non-unital $A_{\infty}$-precategory (see [KoSo6]) with objects $L_t, t\in \R$ and $Hom(L_{t_1},L_{t_2})$ well-defined for $t_1<t_2$ only. 

Assume now that $H$ has a unique global minimum $b_{min}\in B^0$. It gives a common intersection point of all $L_t, t\in \R$, hence a canonical element $i_{t_1,t_2}\in Hom(L_{t_1},L_{t_2})$ which satisfies the property $i_{t_1,t_2}i_{t_2,t_3}=i_{t_1,t_3}$. In this way we identify all objects $L_t$ of our precategory. Then we define the algebra $R$ as the algebra of the endomorphisms of any of them.
The fact that $R$ is finitely generated is not entirely obvious. One can hope that it follows from more careful considerations with filtrations on $R$ coming from the function $H$.

By Mirror Symmetry the zero section $s(B)$ corresponds to the sheaf ${\mathcal O}_{X^{\vee}}$ (or maybe a line bundle on $X^{\vee}$) because $Hom((s(B),\C),(\pi^{-1}(b),\rho))\simeq \C$ for any $b\in B^0$ and any $U(1)$-local system $\rho$ on $\pi^{-1}(b)$. 

We conclude that $R\simeq {\mathcal O}(X^{\vee})$, and parameterizations of open subsets of $X^{\vee}$ by the tube domains give rise to embeddings of $R$ into the algebras of Laurent series obtained by completions of ${\mathcal O}(Hom(\Gamma_b,\C^{\ast}))$ with respect to closed strict convex cones. In the case if our real integrable system comes from a complex polarized one with the central charge, the algebraic symplectomorphism $T: X_1^{\vee}:=X^{\vee}\to X^{\vee}$ corresponds to an automorphism of $R$.
Being the mirror dual, the variety $X^{\vee}\subset Spec(R)$ carries an algebraic volume element $\Omega_{X^{\vee}}$. We expect that there is a $\Z PL$ map of $B^0$ to the skeleton $Sk(X^{\vee},\Omega_{X^{\vee}})$, where the skeleton is defined for the class of logartithmic Calabi-Yau manifolds  in a way slightly more general than in [KoSo2].

\begin{defn} (cf. [GroHaKe])
By a logarithmic Calabi-Yau manifold (log CY manifold for short) we will understand a complex non-compact algebraic  manifold ${\mathcal Y}^0$ endowed with a nowhere vanishing algebraic top degree form $\Omega_{{\mathcal Y}^0}$ which admits a compactification ${\mathcal Y}$ by a simple normal-crossing  divisor $D$ such that $\Omega_{{\mathcal Y}^0}$ has poles of order at most $1$ on $D=\cup_iD_i$ (i.e. $\Omega_{{\mathcal Y}^0}$ is a log-form) and there exists a point in ${\mathcal Y}-{\mathcal Y}^0$ and local coordinates $(z_1,...,z_n)$ such that $\Omega_{{\mathcal Y}^0}=\bigwedge_{1\le i\le n}{dz_i\over{z_i}}$ where $n=dim_{\C}{\mathcal Y}^0$.
\end{defn}

Having a log CY manifold ${\mathcal Y}^0$ one can assign to it a $\Z PL$ topological space $Sk({\mathcal Y}^0)$ of dimension $n$ with linear structure called the {\it skeleton of ${\mathcal Y}^0$}. The construction basically copies the one from [KoSo6]. Namely, for any compactification ${\mathcal Y}$ as in the above definition we define $Sk({\mathcal Y}^0,{\mathcal Y})$ as the complement to $\{0\}$ of the set of such  $\sum_i \lambda_iD_i, \lambda_i\in \R_{\ge 0}$ that if $\lambda_i>0$ then $\Omega_{{\mathcal Y}^0}$ has pole at $D_i$ and $\cap_{i|\lambda_i>0}D_i\neq \emptyset\}$.  Different choices of ${\mathcal Y}$ give rise to $\Z PL$-isomorphisms of $Sk({\mathcal Y}^0,{\mathcal Y})$. Hence we can use the notation $Sk({\mathcal Y}^0)$, and call it the skeleton of ${\mathcal Y}^0$. Integer points of $Sk({\mathcal Y}^0)$ correspond to certain divisorial valuations on the algebra of rational functions on ${\mathcal Y}^0$.

The notion of log CY manifold is analogous to the notion of maximally degenerate proper Calabi-Yau manifold over a non-archimedean field (see [KoSo2], [KoSo6]). More precisely, for a proper Calabi-Yau manifold over a non-archimedean field we defined in [KoSo2], Sect. 6.6 the notion of its skeleton. It is a compact $\Z PL$ topological space of dimension less or equal then the dimension of the Calabi-Yau manifold. The dimensions are equal if and only if the Calabi-Yau manifold is maximally degenerate.

\begin{defn} A log CY ${\mathcal Y}^0$ of complex dimension $n$ is called good if its skeleton $Sk({\mathcal Y}^0)$ coincides with the closure of the set of points of  $Sk({\mathcal Y}^0)$ each of which has a neighborhood 
homeomorphic to a real $n$-dimensional ball.

\end{defn}

In the language of compactifications by simple normal crossing divisors (s.n.c. divisors for short) this means that for a subset of divisors $D_{i_1},...,D_{i_k}, k\le n-1$ of the compactifying divisor $D$ such that $\Omega_{{\mathcal Y}^0}$ has poles of degree $1$ at each of them and such that  $D_{i_1}\cap D_{i_2}\cap...\cap D_{i_k}\ne \emptyset$, there is another subset $D_{i_{k+1}},...,D_{i_n}$ of divisors  such that $D_{i_j}\subset D$ and $\Omega_{{\mathcal Y}^0}$ has poles of degree $1$ at each of them, which altogether satisfy the property  $D_{i_1}\cap D_{i_2}\cap...\cap D_{i_n}=\{pt\}$.

Also, if in the above notation $k=n-1$ the intersection $D_{i_1}\cap D_{i_2}\cap...\cap D_{i_{n-1}}$ is isomorphic to $\C{\bf P}^1$ and the iterated residue of $\Omega_{{\mathcal Y}^0}$ at the intersection coincides with the form $dz/z$ in the chosen coordinate $z$ on $\C{\bf P}^1$. This implies that 
$Sk({\mathcal Y}^0)$ is an oriented topological pseudomanifold of dimension $n$ with possible singularities in codimension $\ge 2$. 

We remark that typically a log CY is good except of rather pathological examples.

\begin{defn} For a given good log CY ${\mathcal Y}^0$ its s.n.c. compactification ${\mathcal Y}$ is very good if each intersection $D_{i_1}\cap D_{i_2}\cap...\cap D_{i_{n-1}}$ as above contains exactly two $0$-dimensional strata (points $z=0,\infty\in \C{\bf P}^1$).
\end{defn}

Then we claim that any very good compactification ${\mathcal Y}$ defines a natural $\Z$-linear structure on $Sk({\mathcal Y}^0)$ with singularities in codim $\ge 2$.
For this we use the isomorphisms of $\Gamma_x\simeq \Gamma_{\overline{F}}$, where $x\in \{0,\infty\}\subset \overline{F}\simeq \C{\bf P}^1$ (see Section 5.2). These isomorphisms allow us to identify canonical  $\Z$-linear structures on $n$-dimensional octants corresponding to the strata $0$ and $\infty$. The definition of very good s.n.c. compactification and the construction of the corresponding $\Z$-linear structure on $Sk({\mathcal Y}^0)$ can be generalized in a straightforward way  to the case of toric-like compactifications.

Assume that we are given a proper toric-like compactification  ${\mathcal Y}$ of a good log CY ${\mathcal Y}^0$ such that the form $\Omega_{{\mathcal Y}^0}$ has poles of degree $1$ at all components of the divisor ${\mathcal Y}-{\mathcal Y}^0$. This compactification is very good in the above sense. One can describe the corresponding singular $\Z$-linear structure  using the language of non-archimedean analytic geometry as in [KoSo6]. Namely, let us extend scalars and consider ${\mathcal Y}^0$  as an algebraic variety over the non-archimedean field $\C((t))$. 
Then we have a continuous map from the Berkovich spectrum of the $\C((t))$-analytic space  $({\mathcal Y}^0)^{an}$ to $Sk({\mathcal Y}^0)$. It is a non-archimedean $n$-dimensional  torus fibration in the sense of [KoSo2], Section 4 (see also [KoSo6]) outside of the codimension $\ge 2$ subset of $Sk({\mathcal Y}^0)$. Hence it defines a $\Z$-affine structure outside of this  subset (see [KoSo2], Sections 4.1, 6.6 for more details). Since our variety was in fact defined over $\C$ (i.e. the corresponding family is constant with respect to the parameter $t$), we see that the $\Z$-affine structure is in fact linear (i.e. it is a vector structure).

Let us return to our integrable system.
Now for a point $b\in B^0$ which does not belong to walls we define a valuation $\nu_b$ on $R$ by taking the pull-back of the valuation on Laurent series on the completed algebra of functions on the torus given by the formula $val(\sum_{\gamma}c_{\gamma}e_{\gamma})=min\{Y_b(\gamma)|c_{\gamma}\neq 0\}$. We expect that the corresponding map $\nu: B^0-\{walls\}\to Sk(X^{\vee},\Omega_{X^{\vee}})$ extends to walls giving the desired $\Z PL$ map of $B$ to the skeleton. We keep the same notation $\nu$ for the extended map. One can hope that it is a homeomorphism.

\subsection{Relation to chains of lines}

Suppose we have a polarized complex integrable systems with central charge and a holomorphic Lagrangian section. 
Adding the parameter $\zeta \in \C^{\ast}$ to the considerations of the previous subsection we obtain a local system $X_{\zeta}^{\vee}$ of symplectic algebraic quasi-affine varieties endowed with algebraic volume forms $\Omega_{X_{\zeta}^{\vee}}$ (notice that since $\Gamma_0=\{0\}$ there is no difference between $X_{\zeta}^{\vee}$ and $X_{\zeta}^{\vee,alg}$).  By functoriality of the skeleton we obtain a local system of skeleta $Sk(X_{\zeta}^{\vee}):=Sk(X_{\zeta}^{\vee},\Omega_{X_{\zeta}^{\vee}})$. The symplectomorphism $T: X_1^{\vee}\to X_1^{\vee}$ gives rise to a $\Z PL$ map $T:Sk(X_{1}^{\vee})\to Sk(X_{1}^{\vee})$.

For any point $b\in B^0$ we have a map $\psi_b:\widetilde{\R^2-\{0\}}=\widetilde{\C}^{\ast}\to Sk(X_{1}^{\vee})-\{0\}$ such that the deck transformation $log\,\zeta\mapsto log\,\zeta+2\pi i$ of the universal covering  $\widetilde{\C}^{\ast}$ goes to the automorphism $T$. Namely, $\psi_b(log\,\zeta)$ is defined as the image of the point $b\in B^0$ under the above-defined map $\nu$ for the polarized complex integrable system $(X,{{\omega^{2,0}}\over{\zeta}})\to B$ (we identify $Sk(X_{\zeta}^{\vee})$ with  $Sk(X_{1}^{\vee})$ using the holonomy of the local system of skeleta over $\C^{\ast}$).

This map is piecewise linear with respect to the natural affine structure on $\widetilde{\R^2-\{0\}}$.  Let us assume we have chosen a very good compactification of $X_1^{\vee}$. Then it defines a $\Z$-linear structure on $Sk(X_{1}^{\vee})$ with conical singularities in codimension $\ge 2$. It follows that for {\it generic } $b\in B^0$ the image of $\psi_b$ is disjoint from the locus of singularities of the $\Z$-linear structure.  We will assume that this is the case.
Moreover we are going to assume that the point $\psi_b(0)=\psi_b(log\, 1)$ does  not belong to the locus of nonlinearity of the $\Z PL$ map $T:Sk(X_{1}^{\vee})\to Sk(X_{1}^{\vee})$ in the above $\Z$-linear structure (it suffices to assume that 
$\psi_b(0)$ does not belong to any rational hyperplane).

Consider the image under the map  $\psi_b$ of the set $\{log\,\zeta|0\le |Im\,log\,\zeta|\le 2\pi \}$. This is the fundamental domain for the natural $\Z$-action on the universal covering $\widetilde{\C}^{\ast}$. We will make the following assumption.

{\bf Monodromy Assumption} {\it  Let $f(t), t\in [0,1]$ be a path $t\mapsto log\,\zeta=2\pi i t$. Then the holonomy of the $\Z$-linear structure from $f(0)$ to $f(1)$ coincides with the map on the tangent spaces induced by $T$.}

Notice that the tangent map $d\psi_b$ can be thought of as an $\R$-linear map $\C\to T^{\Z}_{\psi_b(log\,\zeta),Sk(X_{1}^{\vee})}\otimes \R:=\Gamma_{\psi_b(log\,\zeta)}^{\vee}\otimes \R$ in the obvious notation. Dualizing we obtain an $\R$-linear map $Z_{\psi_b(log\,\zeta)}: \Gamma_{\psi_b(log\,\zeta)}\to \C$.  The Monodromy Assumption implies that after the identification $\Gamma_{\psi_b(0)}=\Gamma_{\psi_b(2\pi i)}$ induced by the holonomy of the $\Z$-linear structure, the transformation $T$ identifies the ``central charges" $Z_{\psi_b(0)}$ and $Z_{\psi_b(2\pi i)}$.

The image of the restriction of $\psi_b$ to the fundamental domain is compact modulo the natural global action of the group $\R_{>0}$ on $Sk(X_{1}^{\vee})-\{0\}$. Hence we can cover an open neighborhood of the image by a chain of rational convex polyhedral cones, which are disjoint from the locus of singularities of the $\Z$-linear structure. Furthermore, we can shrink the cones and sequence of cones $C_1,...,C_{m+1}, C_{m+1}=T(C_1)$ such that after the identification of $C_1$ and $C_{m+1}$ by $T$  we obtain an admissible wheel of cones (cf. Section 5.4, condition h)). Arguments here are similar to those in the proof of Proposition 5.3.2. 
This change of the collection of cones can be interpreted in terms of a sequence of blow-ups and blow-downs at toric strata of a toric-like compactification of $X_{1}^{\vee}$. Hence in a certain toric-like compactification of $X_{1}^{\vee}$ we obtain a chain of lines as in Section 5.4. We also expect that in our situation there is a natural choice of the additional data needed for obtaining a decorated formal Poisson scheme (see Definition 5.2.2).

The above considerations can be generalized to the semipolarized case.

Recall that construction of Section 5.4 gives rise to certain stability data of algebro-geometric origin on the graded Lie algebra $\g_{\Gamma-\Gamma_0}$. The above considerations relate that construction to Mirror Symmetry.

\subsection{Extension to $\zeta=0$}

Suppose we are given a polarized complex integrable system $\pi:(X,\omega^{2,0})\to B$ endowed with a holomorphic Lagrangian section $s:B\to X$.  Let $\underline{\Gamma}\to B^0$ be the corresponding local system of lattices over the complement to the discriminant. Suppose we are given a class $\beta\in H^1(B^0,\underline{\Gamma}^{\vee}\otimes (\R/2\pi \Z))$ which comes from the class $\beta_X\in H^2(X,\R/2\pi \Z)$ which vanishes on $s(B)$ and on fibers of $\pi$. Let us consider the holomorphic family of the Fukaya categories ${\mathcal F}(X,Re(\omega^{2,0}/\zeta), \B=Im(\omega^{2,0}/\zeta)+\beta_X)$. Then mirror duals $X_{\zeta}^{\vee}, \zeta\in \C^{\ast}$ form a holomorphic family of symplectic algebraic varieties over $\C$. We denote the holomorphic symplectic form on $X_{\zeta}^{\vee}$ by $\omega_{X_{\zeta}^{\vee}}$.

\begin{defn} Dual integrable system is a complex integrable system $Y\to B$ such that its restriction to $B^0$ is obtained by:

a) taking dual abelian varieties to fibers of $\pi$;

b) replacing a) by the torsor corresponding to $\beta$.

\end{defn}

Notice that the definition of the dual integrable system does not depend on a choice of polarization. The dual to a polarized integrable system is not polarized in the proper sense. It is ``fractionally polarized", i.e. the corresponding cohomology class is rational, and its positive integer multiple gives a polarization.
The above definition is not quite satisfactory since we ignore the discriminant locus $B-B^0$. 

Notice that there is a holomorphic Lagrangian section $B^0\to Y$ of the dual integrable system. Let us assume that it extends to the section $B\to Y$.

\begin{conj} The above family $X_{\zeta}^{\vee}$ of mirror duals endowed with holomorphic  symplectic forms $\zeta\omega_{X_{\zeta}^{\vee}}$ extends to $\zeta=0$ holomorphically in such a way that the fiber at $\zeta=0$ is holomorphically symplectomorphic  to the total space of the dual integrable system.
\end{conj}

In a similar way we can consider the case of semipolarized integrable system $\pi:X\to B$ with central charge, holomorphic Lagrangian section (in order to relate our considerations to DT-invariants we can take $\beta=\B_{can}$). In that case we start with mirror duals to the integrable systems $X_{Z_0/\zeta}\to B_{Z_0/\zeta}$, where $Z_0/\zeta\in Hom(\Gamma_0,\C)$ is fixed, and then combine them into a local system of complex Poisson manifolds $X_{\zeta}^{\vee}$ endowed with holomorphic Poisson maps to $Hom(\Gamma_0,\C)$.  Recall that by Conjecture 6.3.2 the manifolds $X_{\zeta}^{\vee}$ are expected to be pull-backs via the exponential map $Hom(\Gamma_0,\C)\to Hom(\Gamma_0,\C^{\ast})$ of Poisson algebraic varieties defined over $\Z$ and fibered over $Hom(\Gamma_0,\C^{\ast})$. Then the above conjecture is formulated such as follows.

\begin{conj}
 The local system of Poisson manifolds  $X_{\zeta}^{\vee}$ over $\C^{\ast}$ extends holomorphically (after rescaling the Poisson structure by $1/\zeta$)  to $\zeta=0$. This extension is compatible  with the Poisson morphism to $Hom(\Gamma_0,\C)$. Furthermore, the fiber at $\zeta=0$ is the total space of a (fractionally) semipolarized integrable system $X^{dual}\to B$ whose restriction to $B^0$ has semiabelian Lagrangian fibers with abelian quotients which are dual abelian varieties to the corresponding abelian quotients of fibers of $\pi$.
\end{conj}

This conjecture gives in particular a mathematical interpretation of the picture of twistor family for the total space of the Hitchin system proposed in [GaMoNe2].

Two conjectures below are formulated for simplicity in the polarized case.
There are versions of them in the case of semipolarized integrable systems with central charge and holomorphic Lagrangian section. 

\begin{conj} Let us fix a point $b\in B^0$ in the base of a complex integrable system $\pi: X^0\to B^0$ with abelian fibers endowed with a complex Lagrangian section $B^0\to X^0$. Let us fix a point $e^{i \theta}\in S^1$ such that the pair $(e^{i\theta},b)$ does not belong to the wall in $M=S^1\times B^0$.
Then the  constant family of complex symplectic manifolds $X^{\vee}_{te^{i \theta}}$  over an open ray $l_{\theta}=\R_{>0}e^{i\theta}$  can be extended to a $C^{\infty}$ family over the closed ray  $\R_{\ge 0}e^{i\theta}$ in such a way that the fiber at $t=0$ is a real integrable system over
$Sk_{\theta}$. Here $Sk_{\theta}$ is the skeleton of $(X^0)^{\vee}_{\zeta=e^{i \theta}}$. 
\end{conj}

\begin{conj} For any $e^{i \theta}\in S^1$ the corresponding $Sk_{\theta}$ is  $\Z PL$-manifold isomorphic to $B$
which is endowed with the affine structure derived from the symplectic form 
$Re(e^{-i \theta}\omega^{2,0})$ on $X^0$.
\end{conj}

Let us discuss their relationship of the Conjectures 6.7.4, 6.7.5 with the Conjecture 6.7.3 which predicts holomorphic extension at $\zeta=0$ of the local system of holomorphic Poisson manifolds $X_{\zeta}^{\vee}$. Assuming such an extension let us consider a holomorphic section $\zeta\mapsto f(\zeta)$ of the extended family over a small disc $|\zeta|<\varepsilon$. Let us assume that $f(0)\in \pi^{-1}(B^0)$. Then the restriction 
of $f$ to $l_{\theta}$ gives rise to a real analytic path $f_{\theta}(t)$ in $X^{\vee}_{e^{i \theta}}$. Recall that there is an analytic cover map 
$X_{\zeta}^{\vee}\to X_{\zeta}^{\vee,alg}$, where the latter is (again conjecturally) a quasi-affine algebraic variety over $\Q$.
\begin{conj} For generic $\theta$ the analytic path $f_{\theta}$ define a valuation $val_{\theta}$ on the algebra ${\mathcal O}(X^{\vee,alg}_{e^{i \theta}})$. This valuation depends only on the projection of $f(0)$ to $B^0$ and defines a point in $Sk_{\theta}$. After a continuous extension to $B$ we obtain in this way  a homeomorphism  $B\simeq Sk_{\theta}$ for any $\theta$. 
\end{conj}

The monodromy of the local system of skeleta $Sk_{\theta}$
around $S^1$ is given by the $\Z PL$-automorphism $T$ discussed before in Section 6.6. In terms of chains  of lines this means that
there is a finite decomposition $\cup_{1\le k\le m+1}V_k=S^1$ in the union of strict semiclosed sectors such that two consecutive ones have a common ray, and such that the limiting points at $t=0$ of the above analytic paths $f_{\theta}(t)$ are the same  as long as $e^{i\theta}\in V_k$. Furthermore, on the intersection $V_{1}\cap V_{m+1}$ we have an identification of the corresponding skeleta given by the automorphism $\overline{T}_{1,m+1}$.

\section{Wall-crossing structures and DT-invariants for non-compact Calabi-Yau $3$-folds}

There is a class of non-compact Calabi-Yau $3$-folds which gives rise to  complex integrable systems with central charge (and hence to the corresponding wall-crossing structures). Such Calabi-Yau manifolds admit compactifications by a normal crossing divisor where the holomorphic volume form has poles of degree at least one (such a variety is not a log CY since we allow poles of order strictly bigger than one). This class of ``good'' Calabi-Yau threefolds include those of the type $\{uv+P(x,y)=0\}$ where $P(x,y)$ is a polynomial such that the equation $P(x,y)=0$ defines a smooth affine curve. Presumably this class includes non-compact Calabi-Yau $3$-folds associated with Hitchin systems  (possibly with irregular singularities) for all gauge groups generalizing [DiDoPa].

\subsection{Deformation theory of non-compact Calabi-Yau 3-folds}

Let $\overline{X}$ be a complex projective variety of complex dimension $3$, and $D=\cup_iD_i\subset \overline{X} $ be a normal crossing divisor such that algebraic variety $X:=\overline{X}-D$ has a nowhere vanishing top degree form $\Omega^{3,0}_X$ which extends to $\overline{X}$ with poles of order $n_i\ge 1$  on $D_i$. We will also fix an ample line  bundle $\LL$ on $\overline{X}$ which defines a projective embedding of $\overline{X}$.  In this subsection we are going to discuss the deformation theories related to the pair $(\overline{X},D)$ or the triple $(\overline{X},D,\LL)$ (later we will see that these deformation theories are equivalent, see Proposition 7.1.3). We impose the following assumptions which will later guarantee that the global moduli stack  $\MM_{(\overline{X},D,\LL)}$ of deformations will be a smooth orbifold :

\vspace{2mm}

{\bf A1} $H^{2,0}(\overline{X})=0$.

\vspace{2mm}

{\bf A2} There exists a component $D_{i_0}$ such that $ord_{D_{i_0}}\Omega^{3,0}_X\ge 2$ and such that the restriction homomorphism $H^1(\overline{X},\Z)\to H^1(D_{i_0},\Z)$ is an isomorphism.

\vspace{2mm}

Loosely speaking we are going to consider deformations  of $\overline{X}$ which preserve the isomorphism class of the restriction of the ample line bundle $\LL_{|D_{i_0}}$ (see Assumption $A2$) as well as some ``decoration" of $D$, which consists of jets of order $n_i-1$ at each $D_i$ for which $n_i=ord_{D_i}\Omega^{3,0}_X>1$ (we skip a more precise formulation of the latter condition).

We will see that the Assumption $A2$ implies that the deformed varieties carry  $(3,0)$-forms with same orders of poles at the deformed smooth components $D_i$. The condition $K_{\overline{X}}=-\sum_in_iD_i$ is preserved under deformations.

Let us now describe precisely the corresponding moduli problems. We will work in analytic topology.
Let $T_{\overline{X},D}:=T_{\overline{X},D,\Omega^{3,0}_{X}}$ be the sheaf of holomorphic  vector fields on $\overline{X}$  satisfying the property that the contraction of such a vector field with $\Omega^{3,0}_X$ has poles of order $1$ on each $D_i$ (i.e. the contraction is a logarithmic form). Then $T_{\overline{X},D}$ is a sheaf of Lie subalgebras of the sheaf of graded Lie algebras of polyvector fields $\Lambda_{\overline{X},D}$ (here $k$ vector fields are placed in degree $1-k$ for $k=0,1,2,3$) which  satisfy the property that their contractions with $\Omega^{3,0}_X$ are logarithmic forms. Recall that a sheaf of $L_{\infty}$-algebras over a field of characteristic zero (e.g. Lie algebras or DGLAs) gives rise to the corresponding formal deformation theory (see e.g. [Ko2], [KoSo8]). We will denote the formal moduli space associated with an $L_{\infty}$-algebra $\g$ by $\MM_{\g}$.

Let us consider the following deformation theories:    
    
a) The one associated with the DGLA $\g_0=R\Gamma(\overline{X},T_{\overline{X},D})$. These are deformations of 
the pair $(\overline{X},D)$ preserving certain decoration on $D$.

b) The one associated with the differential graded Lie algebra (DGLA for short)
$\g_1=R\Gamma(\overline{X},\Lambda_{\overline{X},D},div_{\Omega^{3,0}_X})$, where $div:=div_{\Omega^{3,0}_X}$ is the divergence operator associated with the holomorphic volume form $\Omega^{3,0}_X$.
This deformation theory does not have a straightforward  geometric interpretation.

c) The one associated with the DGLA subalgebra 

$$\g_2=R\Gamma(\overline{X},T_{\overline{X},D}\stackrel{ div_{\Omega^{3,0}_X}}{\longrightarrow} {\cal O}_{\overline{X},D}).$$ 
These are deformations of the pair  $(\overline{X},D)$ preserving the same decorations as in a), but also this time preserving the section $(\Omega^{3,0}_X)^{-1}$ of the anticanonical bundle 
$-K_{\overline{X}}$.

Here ${\cal O}_{\overline{X},D}$ is the sheaf of functions on $\overline{X}$ such that being multiplied by $\Omega^{3,0}_X$ they have pole of order at most one at $D$ (i.e. they are degree zero polyvector fields from $\Lambda_{\overline{X},D}$).

\begin{prp} The moduli space $\MM_{\g_1}$ is naturally isomorphic to a formal submanifold of the formal neighborhood of $0\in H^3(X,\C)$. The moduli space $\MM_{\g_2}$ is a formal submanifold of $\MM_{\g_1}$.

\end{prp}

{\it Proof.}  DGLA $\g_1$ maps quasi-isomorphically to a DGLA $$\g_1^{\prime}:=R\Gamma(\overline{X},(i_{X\to \overline{X}})_{\ast}({Sym}^{\bullet}(T_X[1])[-1],div_{\Omega^{3,0}_X}))$$
$$\simeq R\Gamma({X},({Sym}^{\bullet}(T_X[1])[-1], div_{\Omega^{3,0}_X})).$$  
This follows form the classical result of Grothendieck that the complex of differential log-forms on $\overline{X}$ endowed with de Rham differential embeds quasi-isomorphically to $(i_{X\to \overline{X}})_{\ast}(\Omega^{\bullet}_X,d_{dR})$. Similarly the homomorphism of DGLAs 
$$ R\Gamma({X},({Sym}^{\bullet}(T_X[1])[-1], div_{\Omega^{3,0}_X}))\to  R\Gamma({X}_{an},({Sym}^{\bullet}(T_{X_{an}}[1])[-1], div_{\Omega^{3,0}_X}))$$
relating Zariski and analytic topologies, is a quasi-isomorphism.
The embedding of the abelian DGLA $\underline{\C}_{X}[2]$ to the sheaf of DGLAs$ ({Sym}^{\bullet}(T_{X_{an}}[1])[-1],div_{\Omega^{3,0}_X})$ 
which maps $1_{X}$ to $(\Omega^{3,0}_X)^{-1}$ is also a quasi-isomorphism.  This implies that the corresponding moduli space is a formal neighborhood of a point in the affine space $H^3(X,\C)$ .

There is an obvious embedding of the complex of sheaves corresponding to $\g_2$ into the one corresponding to $\g_1$. We need to check that it induces an embedding at the level of hypercohomology. Contracting both complexes with $\Omega^{3,0}_X$ we convert  polyvector fields into logarithmic forms. Then $\g_2$ is quasi-isomorphic to
${\mathbb H}^{\bullet}(\Omega^2_{\overline{X}}(log\,D)\to \Omega^3_{\overline{X}}(log\,D),d)$. By Hodge theory this is embedded  into the hypercohomology
${\mathbb H}^{\bullet}(\Omega^{\bullet}_{\overline{X}}(log\,D),d)$. Since $\g_1$ is quasi-isomorphic to an abelian DGLA, the same is true for $\g_2$, and moreover  $\MM_{\g_2}$ is a formal submanifold of $\MM_{\g_1}$ (see e.g. [KaKoPa], Proposition 4.11 (ii)). 
$\blacksquare$

Consider the natural $L_{\infty}$-morphism $\g_2\to \g_0$.

\begin{prp} Under the Assumption $A2$ this morphism induces an isomorphism of the moduli space $\MM_{\g_2}\to \MM_{\g_0}$.

\end{prp}

{\it Proof.} Hodge theory implies  that the morphism $\g_2\to \g_0$ induces an epimorphism on cohomology. Then it is easy to show that $\g_0$ is quasi-isomorphic to an abelian DGLA (see e.g. [KaKoPa], Proposition 4.11 (iii)). Thus we have a surjection $\MM_{\g_2}\to \MM_{\g_0}$ which is a smooth fibration. The tangent space to a fiber is isomorphic to $H^0(\overline{X},{\mathcal O}_{\overline{X},D})$. By Assumption $A2$ it is trivial. This proves the Proposition. $\blacksquare$

Now we would like to discuss the formal deformation theory which takes into account the ample line bundle $\LL$.


Let $At(\LL)$ denotes the sheaf of Lie algebras of infinitesimal symmetries of the pair $(\overline{X},\LL)$ (Atiyah algebra of $\LL$). It fits into an exact short sequence
$$0\to {\mathcal O}_{\overline{X}}\to At(\LL)\to T_{\overline{X}}\to 0.$$

Let us denote by $At_{i_0}(\LL)$ the subsheaf of $At(\LL)$ of infinitesimal symmetries vanishing at the divisor $D_{i_0}$ from the Assumption $A2$. Then we have a short exact sequence
$$0\to {\mathcal O}_{\overline{X}}(-D_{i_0})\to At_{i_0}(\LL)\to T_{\overline{X},D_{i_0}}\to 0,$$
where $T_{\overline{X},D_{i_0}}\subset T_{\overline{X}}$ is a sheaf of vector fields vanishing identically on $D_{i_0}$. One can check that $T_{\overline{X},D}$ is a subsheaf of $T_{\overline{X},D_{i_0}}$. Let us define a sheaf of DGLAs $\g_3$ as the fiber product of the morphisms

$$At_{i_0}(\LL)\to T_{\overline{X},D_{i_0}}\leftarrow T_{\overline{X},D} $$
over $T_{\overline{X},D_{i_0}}$. Then we have an exact sequence of sheaves
$$0\to {\mathcal O}_{\overline{X}}(-D_{i_0})\to \g_3\to T_{\overline{X},D}\to 0.$$

The sheaf of DGLAs $\g_3$ controls the deformation theory of $\overline{X}$ (together with decorations on $D$) endowed with a line bundle  such that the restriction of the line bundle to $D_{i_0}$ is identified with $\LL_{D_{i_0}}$. 
\begin{prp} The natural map $\MM_{\g_3}\to \MM_{\g_0}$ is an isomorphism.

\end{prp}
{\it Proof.} There are natural maps to the tangent sheaf $T_{\overline{X}}$ of both $At(\LL)$  and $T_{\overline{X},D}$. The fiber product of these two maps is  a sheaf of Lie algebras which we will denote by $At_{\overline{X},D}(\LL)$. Let $\g_4=R\Gamma(\overline{X},At_{\overline{X},D}(\LL))$ be the corresponding DGLA. The moduli space $\MM_{\g_4}$ is smooth because $\MM_{\g_0}$ is smooth and there is a formal bundle $q:\MM_{\g_4}\to \MM_{\g_0}$ with fibers which are smooth with tangent spaces isomorphic to $H^1(\overline{X},{\mathcal O}_{\overline{X}})$ (for that we need the condition $h^{2,0}(\overline{X})=0$ which is the Assumption $A1$).
Restricting $At_{\overline{X},D}(\LL)$ to $D_{i_0}$ we obtain the map $p:\MM_{\g_4}\to \widehat{Pic(D_{i_0})}_{\LL}$ which is the formal neighborhood of $\LL$ in the Picard group thought of as the moduli space of line bundles. This is an epimorphism at the level of tangent spaces since the map $H^1(\overline{X},{\mathcal O}_{\overline{X}})\to H^1(D_{i_0},{\mathcal O}_{D_{i_0}})$ is an isomorphism by Assumption $A2$.
The morphisms $p$ and $q$ give rise to an isomorphism $\MM_{\g_4}\simeq \MM_{\g_0}\times  \widehat{Pic(D_{i_0})}_{\LL}$.
Notice that the fiber of the map $p$ over $\LL_{|D_{i_0}}$ is isomorphic to the moduli space $\MM_{\g_3}$.  This proves the desired isomorphism 
 $\MM_{\g_3}\simeq \MM_{\g_0}$. $\blacksquare$

%

The above propositions imply that there are three canonically isomorphic formal moduli spaces: $\MM_{\g_3}\to \MM_{\g_0}\simeq \MM_{\g_2}$. The deformation theory associated with $\g_3$ is convenient for the construction below of the global moduli space.

Assume $\LL$ is an ample bundle. Let us choose $N\ge 1$ such $\LL^{\otimes N}$ gives a projective embedding of $\overline{X}$.
Then we can consider {\it non-formal} deformation theory corresponding to the above formal deformation theory associated with ${\g_3}$. More precisely, consider the scheme $\MM^{\prime}$ which parametrizes the following data:

1) Smooth projective subvarieties $\overline{X}^{\prime}\subset {\bf P}^{m-1}, m=rk\,H^0(\overline{X},\LL^{\otimes N})$ such that $rk\,H^0(\overline{X}^{\prime},{\mathcal O}(k))=rk\,H^0(\overline{X},\LL^{\otimes kN})$ for all sufficiently large $k>0$. We also assume that $\overline{X}^{\prime}$ satisfy $A1,A2$.

2) Normal crossing divisors $D^{\prime}=\cup_{i}D_i^{\prime}\subset \overline{X}^{\prime}$ together with a bijection between the set of irreducible components of $D^{\prime}$ and $D$.

3) A chosen holomorphic volume element $\Omega^{3,0}_{{X}^{\prime}}$ on $X^{\prime}:=\overline{X}^{\prime}-D^{\prime}$ such that $n_i^{\prime}:=ord_{D_i^{\prime}}\Omega^{3,0}_{{X}^{\prime}}=n_i$.

4) For components $D_i$ with $n_i\ge 2$  chosen isomorphisms $D_i\simeq D_i^{\prime}$ preserving stratifications by other divisors as well as an isomorphism $\cup_{n_i\ge 2}D_i\simeq\cup_{n_i^{\prime}\ge 2}D_i^{\prime}$ which is required to induce the above isomorphisms of individual divisors.

5) Chosen isomorphisms ${\mathcal O}(1)_{|D_{i_0}}\simeq \LL^{\otimes N}_{|D_{i_0}}$.

6) For any $i$ such that $n_i\ge 2$ and any $x\in D_i$ a chosen isomorphism of the formal neighborhood of $x$ with the formal neighborhood of the corresponding (see 4)) $x^{\prime}\in D_i^{\prime}$, preserving stratifications by other divisors as well as holomorphic volume elements, and defined up to the action of the formal completion of the Lie subalgebra 
$T_{\overline{X},D}^{div}\subset T_{\overline{X},D}$ of vector fields with zero divergence. Moreover, we require that the above formal isomorphisms can be chosen in such a way that they depend locally algebraically in Zariski topology on $x\in \cup_{n_i\ge 2}D_i$.

\begin{rmk} The last condition can be formulated in terms of jet spaces.

\end{rmk}

The scheme $\MM^{\prime}$ is smooth by the above-discussed formal deformation theory. The group $GL(m)$ acts on $\MM^{\prime}$ with finite stabilizers because $\Gamma({X}^{\prime},T_{\overline{X}^{\prime},D})=0$. Hence the quotient of $\MM^{\prime}$ by this action is a smooth Deligne-Mumford stack (orbifold).
Let us denote it by $\MM$. This will be the base of our complex integrable system.

\begin{rmk} Locally in analytic topology a neighborhood of a point $m\in \MM$ corresponding to $(\overline{X}^{\prime},D^{\prime})$ is embedded as a maximal isotropic submanifold in the Poisson manifold $H^3({X}^{\prime})$ by taking the cohomology class $[\Omega^{3,0}_{{X}^{\prime}}]$ (the period map).  The Poisson structure on the third cohomology comes from the observation that it is dual to $Hom(\Gamma^{\prime},\C)$, where $\Gamma^{\prime}=H_3({X}^{\prime},\Z)/tors$ carries an integer skew-symmetric intersection form.

Also we remark that one can generalize the above considerations by allowing $\Omega^{3,0}_X$ to extend to some components $D_i$ without zeros and poles. 

\end{rmk}

\subsection{WCS and integrable systems associated with the moduli space}

By Remark 7.1.5 we see that $\MM$ carries a local system $\underline{\Gamma}$ with the fiber over $u\in \MM$ given by $H_3(X^{\prime},\Z)/tors$, where $X^{\prime}$ is the corresponding non-compact Calabi-Yau $3$-fold. The intersection form endows $\underline{\Gamma}$ with an integer skew-symmetric form $\langle\bullet,\bullet\rangle$ while  the period map can be interpreted as a central charge
$Z:\gamma\mapsto \int_{\gamma}\Omega^{3,0}_X$. It gives rise to a holomorphic family of homomorphisms $Z_u:\Gamma_u\to \C, u\in \MM$, so we have a local embedding of $\MM$ into $\Gamma^{\vee}_u\otimes \C$ such that the image of $Z$ is a family of Lagrangian submanifolds in symplectic leaves of the Poisson structure on $\Gamma^{\vee}_u\otimes \C$  dual to the $2$-form $\langle\bullet,\bullet\rangle$. As we discussed previously, this family of Lagrangian manifolds is parametrized by the kernel of 
$\langle\bullet,\bullet\rangle$ and each Lagrangian manifold is the base of a complex integrable system.

\begin{prp} The mixed Hodge structure on $H^3(X,\C)$ can have non-trivial components in $H^{1,2},H^{2,1},H^{2,2}$ only.

\end{prp}

{\it Proof.} Since $X$ is smooth we have $W_3H^3(X)=H^3(X)$. Hence it suffices to show that $F^3H^3(X)=0$. Recall that the latter space can be defined as
$H^3(R\Gamma(\overline{X},0\to 0\to 0\to \Omega^3_{\overline{X}}(log\,D))$ which is equal to 
$H^0(\overline{X},\Omega^3_{\overline{X}}(log\,D))\simeq H^0(\overline{X},{\cal O}_{\overline{X},D})=0$ since $n_{i_0}>1$ by Assumption $A2$.
$\blacksquare$

After twisting by the Tate motive $\Z(1)$ we obtain a variation of mixed Hodge structure. It satisfies all the conditions from Section 4.1.2 except of (possibly) the condition 3) iii). By general reasons explained there it gives us a complex integrable system with fibers being semiabelian varieties, a central charge and holomorphic Lagrangian section. If the positivity condition 3) iii) is satisfied, then our complex integrable system is semipolarized.

Let us now discuss the positivity condition more precisely.
The tangent space to the base of each of the integrable systems is isomorphic to the image of $H^1(\overline{X},\Omega^2_{\overline{X}})$ in  $H^1(\overline{X},\Omega^2_{\overline{X}}(log\,D))$. Then we need positivity of the restriction of the pseudo-hermitian form on the latter space to the image of the former one. It is convenient to dualize the above embedding. The dual space can be identified with the image of the composition $$H^3(\overline{X},D)\to H^3(\overline{X})\to {\mathcal H}^3(\overline{X}, \Omega^0_{\overline{X}}\to \Omega^1_{\overline{X}}),$$
 where $H^3(\overline{X},D)\simeq H^3(X)^{\ast}$ is the cohomology of pair with complex coefficients and  ${\mathcal H}^3(\overline{X}, \Omega^0_{\overline{X}})=H^2(\overline{X}, \Omega^1_{\overline{X}})\simeq (H^1(\overline{X},  \Omega^2_{\overline{X}}))^{\ast}$ since $h^{3,0}(\overline{X})=0$ . One can identify the image of the map $H^3(\overline{X},D)\to H^3(\overline{X})$ with
$Ker(H^3(\overline{X})\to \oplus_iH^3(D_i))$ using the long exact sequence of the cohomology of pair.  Using Hodge decompositions of $H^3(\overline{X})$ and $H^3(D_i)$ we conclude that the positivity condition is equivalent to the following

\vspace{2mm}

{\bf Positivity Assumption A3} 
{\it For any non-zero differential form $\alpha$ representing an element of 
$Ker\left(H^1(\overline{X},\Omega^2_{\overline{X}})\to \oplus_i H^1(D_i,\Omega^2_{{D_i}})\right)$ we have $\int \alpha\wedge \overline{\alpha}>0$.}

\vspace{2mm} 
\begin{rmk}
 
Assumption $A3$ holds e.g. for the compactification of the total space of Hitchin system on ${\bf P}^1$ or in the case of Hitchin systems related to $ALE$ spaces as in [DiDoPa].

\end{rmk}

We recall that having a K\"ahler metric we can enlarge $B^0$ to the ``full" base $B$ defining the latter as the completion of $B^0$ with respect to the metric. We do not know how to define the integrable system with the base $B$, but this is not necessary for the construction of attractor trees and WCS.

Let us discuss the Assumption $A3$. 
The vector space $\bigwedge^2H^3(\overline{X},\Q)$ contains an element $\overline{\delta}$ which is the K\"unneth component of the diagonal. Then the vector space
$\bigwedge^2(W_3H^3(X,\Q))= Im\left( \bigwedge^2H^3(\overline{X},\Q)\to  \bigwedge^2H^3({X},\Q)))\right)$
contains the image of $\overline{\delta}$ which we will denote by $\delta$.  Then $A3$ is essentially equivalent to the claim that $\overline{\delta}$ is non-degenerate and together with the Lagrangian vector
subspace $Im\left(H^{2,1}(\overline{X})\to W_3H^3(X,\Q)\otimes \C\right)$ it defines the hermitian metric  with non-degenerate skew-symmetric part on the ambient vector space 
$$Im\left(H^3(\overline{X},\C)\to H^3({X},\C)\right).$$
Since the latter is the tangent space to the base of our integrable system, we conclude that it carries the K\"ahler metric.

In order to apply the algorithm of construction of WCS from Section 4.6 we also need a family of cones. It is more natural to discuss that piece of data in the framework of the Support Property for the DT-invariants of the Fukaya category of $X$. But the very existence of the Fukaya category and an appropriate stability condition is a non-trivial question. We are going to discuss it below.

\subsection{Fukaya categories for non-compact Calabi-Yau 3-folds and  stability conditions}

We start with general remarks. Our definition of WCS is motivated  by the theory of Donaldson-Thomas invariants for the Fukaya category. Namely, for a ``good''  non-compact Calabi-Yau threefold $X$ one should have a well-defined Fukaya category endowed with a stability condition, whose central charge is the period map of the holomorphic $3$-form. As we have already discussed in this section, this should give us a family of polarized integrable systems  whose bases sweep the moduli space of deformations of $X$ (equivalently, a semipolarized integrable system). The base of each polarized integrable system is endowed with a K\"ahler metric. Furthermore, the theory of DT-invariants from [KoSo1,5] says that  tangent spaces of points of the base should carry strict convex cones (this is equivalent to the so-called Support Property from [KoSo1]). We propose some sufficient conditions for the above picture to be true. In particular we impose the Assumptions $A1-A3$ which give part of the desired structures.  In this subsection we discuss the conditions under which the Fukaya category and a stability condition do exist. This gives additional to Assumptions $A1-A3$ sufficient conditions for $X$ to belong to a ``good'' class. On the other hand we expect that {\it for any Calabi-Yau $3$-fold $X$, compact or non-compact, one should be able to define collections of integers $\Omega_u(\gamma)$ parametrized by the open subspace in the moduli space of deformations of $X$, and which coincide with the DT-invariants of the Fukaya category of $X$ endowed with an appropriate stability condition in the case when the latter can be defined}.

The content of this subsection is purely motivational. 

1) In order to have a well-defined Fukaya category ${\mathcal F}(X)$ we need to ensure that holomorphic discs cannot touch the divisor $D$. 

2) In order to have a stability condition on ${\mathcal F}(X)$ we need to ensure compactness of the space of special Lagrangian submanifolds (SLAGs) in a fixed homology class. 

Having 1) and 2) we can speak about virtual Euler characteristic of the moduli spaces of SLAGs, hence to define DT-invariants. They should satisfy the wall-crossing formulas from [KoSo1]. For that we need to ensure the Support Property from [KoSo1].

Suppose we have ensured 2). We claim that the Support Property is satisfied by the Assumption $A2$. 
Indeed we can choose logarithmic forms $\alpha_i, 1\le i\le n$ on $X$ which are representatives of a basis in $H^3(X,\R)$. Then for points $x\in X$ which are sufficiently close to $D$ we have
$$||\alpha_i(x)||\le C(\alpha_i)\sqrt{|\Omega^{3,0}(x)|^2\over{|\omega^{1,1}(x)|^3}},$$
for any $1\le i\le n$. Here $\omega^{1,1}$ is a chosen K\"ahler form on $\overline{X}$, and we take any norm $||\bullet ||$ of the functional $\alpha_i(x)$. The inequality follows from the fact that the form $\Omega^{3,0}$ has poles of order at least one at any component $D_i$ of $D$, and there exists a component $D_{i_0}$ where it grows faster. 
Then considerations from Remark 1 of [KoSo1] can be applied in the non-compact case in the same way as in the compact one. This gives the Support Property and strict convex cones in WCS constructed in the previous subsection.

Now we turn to a discussion of 1). Let $(X,\omega)$ be a real symplectic manifold, possibly non-compact. We fix an almost complex structure $J$ which is integrable outside of a compact subset and compatible with $\omega$. In other words, $X$ is a complex manifold ``at infinity'' and $\omega(v,J(v)), v\in T_xX$ defines an almost K\"ahler form. If we want the Fukaya category to be $\Z$-graded we ask for a differential form $\Omega_X$ such that $\Omega_X$ has type $(n,0)$ in the complex structure defined by $J$ and does not have zeros on $X$.  

In order to ensure that no pseudo-holomorphic discs ``go to infinity'' one can impose one of the following sufficient conditions:

a) There exists smooth proper function $f:X\to [0,+\infty)$ with the property that for sufficiently large $c>0$ the hypersurface $f=c$ is smooth (i.e. it does not contain critical points of $f$), and the Levy form of this hypersurface is non-negative
(it suffices to require that $\overline{\partial}\partial(f)\ge 0$ outside of a compact).
The desired property of pseudo-holomorphic discs follows from the maximum principle.

b) There is a compact manifold $\overline{X}$ containing $X$ and such that outside of a compact it is an embedding of complex manifolds, and such that $D=\cup_iD_i:=\overline{X}-X$ is a normal crossing divisor which satisfies the following positivity condition:

{\it if a  rational curve $C\subset \overline{X}$ contains a smooth component belonging to some $D_i$ then its intersection index with $D$ is non-negative.}

If this condition is satisfied then a family of pseudo-holomorphic discs $S_t$ in $X$ cannot converge as $t\to \infty$ to a degenerate disc $S$ such that $S\cap D\ne \emptyset$. Indeed the intersection index of $S_t$ with $D$ is zero. The same should be true for $S$ since the intersection index is a homological invariant. But $S$ must have a component intersecting some $D_i$ with strictly positive index, and possibly other components which are rational curves belonging to $D$. Therefore the intersection index of $S$ and $D$ is strictly positive. This contradiction shows that $S$ cannot exist.

Condition b) looks weaker than the condition a) since it deals with rational curves only.
There is a class of examples where both a) and b) are satisfied. For that class the function $f$
is an extension to $X$ of the function
$$f(x)=\sum_in_ilog\left({1\over{dist(x,D_i)}}\right)+r(x),$$
where $D=\sum_in_iD_i$ as a divisor, $x\mapsto dist(x,Z)$ is the distance function to a set $Z$, and $r(x)$ is a smooth function on $\overline{X}$.
Then $\overline{\partial}\partial(f)=\alpha-\sum_in_i\delta_{D_i}$, where $\alpha$ is a smooth $(1,1)$-form on $\overline{X}$ which is non-negative outside of $D$, not equal to zero on $D$, and $\delta_{D_i}$ is the delta-distribution for the component $D_i$. Since the integral of the LHS is zero for any holomorphic curve in a neighborhood of $D$, we conclude that the intersection index of $D=\sum_in_iD_i$ with such a curve is non-negative. Notice also that such a curve does not have to be rational.

Summarizing, if a) or b) are satisfied then one can hope that there exists well-defined Fukaya category ${\mathcal F}(X,\omega)$.

Now we turn to a discussion of the condition 2), which should lead to a well-defined count of SLAGs.
These considerations are purely heuristic. The idea is to consider a flow on Lagrangian submanifolds of $X$ defined by the differential one form $\mu(L):=dArg(\Omega_X)_{|L}$. The direct computation shows that the function $vol(L)=\int_L(\Omega_X)_{|L}$ decreases along the flow (it is the gradient flow with respect to the $L^2$-metric on differential $1$-forms).
Stable points of the flow are SLAGs. Thus in order to achieve compactness of the space of SLAGs it is sufficient to ensure that there exists a real-valued function $H$ with the following property:

{\it if a compact Lagrangian $L$ belongs to the set $H\le c$ where $c$ is sufficiently big, then the trajectory of $L$ along the above flow also belongs to the same set.}

Reason for that: if $H$ has a local maximum at $x\in L$ then a small shift of $L$ along the flow makes the maximum smaller. We discuss below sufficient  conditions for $H$ to be a desired function. We expect that that under some conditions on the volume form there exists $H$ with the desired property.

Since the above condition is local, one can assume that $L$ is a real Lagrangian manifold in a flat K\"ahler  vector space with coordinates $(z_i,\overline{z}_i), 1\le i\le n$ given by
$\overline{z}_i=z_i+\sqrt{-1}P(z_1,...,z_n)$, where $P$ is a real formal power series which starts with terms of degree greater or equal than $3$. 

The top degree form is given by $\Omega=dz_1\wedge...\wedge dz_n(1+\sum_{1\le i\le n}(a_i+\sqrt{-1}b_i)z_i+O(z^2))$, where $a_i,b_i\in \R$.

We are looking for a function $H$ which can be written as
$$H=c+\sum_i(H_iz_i+\overline{H}_i\overline{z}_i)+\sum_{i,j}(H_{ij}z_iz_j+\overline{H}_{ij}\overline{z}_i\overline{z}_j)+\sum_{ij}r_{ij}z_i\overline{z}_j+O(z^4),$$ where $H_{ij}=H_{ji}, r_{ij}=r_{ji}$.

Local expression for $Arg\,\Omega_{|L}={\Omega_{|L}\over{\overline{\Omega}_{|L}}}$ is $\sum_ib_iz_i-{1\over{2}}\sum_{ij}\partial_{iij}(P)z_j$, where $\partial_{ijk}$  denotes the partial derivative with respect to $z_i,z_j,z_k$. 

Normal vector field to $L$ (i.e. the vector field of the flow) is given by
$$n=(\sqrt{-1})^{-1}\sum_ib_i(\partial/\partial z_i-\partial/\partial \overline{z}_i)-$$
$${1\over{2}}\sum_{ij}\partial_{iij}(P)(\partial/\partial z_i-\partial/\partial \overline{z}_i).$$
One can check that the partial derivative $\partial H/\partial n=(dH,n)$ is non-negative if
$$\sum_iIm(H_i)+2\sum_iRe(H_{ii})+\sum_iRe(r_{ii})\le 0.$$

This inequality can be written in a more invariant way. Notice that the function
$log({\Omega\over{dz_1\wedge...\wedge dz_n}})$ is well-define up to a constant. Therefore its Poisson bracket (with respect to the Poisson structure given by the K\"ahler form) is well-defined.
Let us denote by $(c_i)_{1\le i\le n}$  the (non-negative) spectrum of the quadratic form $H^{(2)}$
defined by the quadratic part of $H$ (with respect to the Hermitian metric given by the K\"ahler metric). Let $\Delta=\sum_{ij}\partial/\partial{z}_i\partial/\partial{\overline{z}}_i$ be the Laplace operator defined by the K\"ahler metric. Then the sufficient condition looks such as follows:

$$\Delta(H)+Re\{log\left({\Omega\over{dz_1\wedge...\wedge dz_n}}\right),H\}-\sum_ic_i\ge 0.$$
Moreover,if the inequality is strict for sufficiently large values of $H$, then all SLAGs belong to a compact subset of $X$.
One can check that the strictness of the inequality is not always achieved. A counter example is given by $X=\C^{\ast}$ and $\Omega=dz/z$. Then we have infinitely many SLAGs (circles) which ``go to infinity'' on the corresponding cylinder.  In this case the above inequality becomes an equality. One can hope that if poles of $\Omega$ at $D$ have order greater or equal than $2$, then the above inequality is strict. In that case one can hope to have a stability condition on the Fukaya category of $X$ and well-defined count of SLAGs (hence the corresponding DT-invariants in the case when $X$ is a $3CY$ manifold).

\section{Hitchin integrable systems for $GL(r)$}

In this section $C$ will denote a connected smooth projective irreducible curve over $\C$.
Although we are going to discuss Hitchin systems with the gauge group $GL(r)$, we hope that our constructions admit generalizations to arbitrary reductive groups.

\subsection{Reminder on non-singular case}
Let us recall some basics on Hitchin systems with non-singular Higgs fields assuming that the genus of the curve $C$ is bigger than $1$. Recall that $GL(r)$-Higgs bundle is a rank $r$ vector bundle over $C$ endowed with a morphism $\phi:E\to E\otimes K_C$ (equivalently, a morphism $T_C\to End(E)$) called the Higgs field. Here $K_C=T^{\ast}_C$ is the canonical  sheaf of $C$. The moduli space $\MM_{Higgs}(r,d)$ of stable Higgs bundles of rank $r$ and degree $d$ on $C$ is the total space of a polarized complex integrable system. The base $B=\prod_{1\le i\le r}\Gamma(C,K_C^{\otimes\,i})$ carries a structure of vector space. The projection map $\pi: \MM_{Higgs}(r,d)\to B$ assigns to a pair $(E,\phi)$ a collection $(Tr\,\phi, Tr\,\phi^2,...,Tr\,\phi^r)\in B$. This map has the following geometric meaning. Higgs bundle is the same as a  coherent sheaf $E_1$ on $T^{\ast}C$, which is pure and supported on a compact curve $S\subset T^{\ast}C$, called the spectral curve of $(E,\phi)$. The purity means that $E_1$ has no non-trivial subsheaves with zero-dimensional support. The direct image of $E_1$ under the canonical projection $T^{\ast}C\to C$ is isomorphic to $E$. Generically $S$ is smooth and $E_1$ is the direct image of a line bundle on $S$. For given $x\in C$ points of the intersection $S\cap T_x^{\ast}C$ correspond to eigenvalues of the linear map $\phi_x$ understood as an endomorphism of $E_x$ associated with a non-zero tangent vector to $C$. The point $\pi(E,\phi)\in B$ can be thought of as a collection of coefficients of the characteristic equation 
$p(x,y):=det(\phi_x-y\cdot id_E)=0$ which defines the spectral curve $S$.

We denote by $B^0\subset B$ the locus of smooth connected spectral curves. The fiber $\pi^{-1}(b), b\in B^0$ is a torsor over the Jacobian $Jac(S_b)$ of the corresponding spectral curve $S_b$. It consists of line bundles on $S_b$ of a certain degree. Therefore $B^0$ can thought of as a space of smooth connected projective curves $S$ in an open complex symplectic variety $(T^{\ast}C,\omega_{T^{\ast}C})$ with the homology class $r[C]\in H_2(T^{\ast}C,\Z)$.
Our integrable system depends on degree $d$, but the associated integrable system over $B^0$ with holomorphic Lagrangian section (see Section 4.1.1) does not depend on $d$ and has as fibers Jacobians $Jac(S_b)$.

The above construction can be generalized. Namely, instead of $T^{\ast}C$ we can consider an arbitrary complex symplectic surface $Y$ and smooth compact curves $S\subset Y$. More generally, one can replace $Y$ by a higher-dimensional complex quasiprojective symplectic manifold (more generally, K\"ahler manifold) and consider smooth compact complex Lagrangian submanifolds in it. Then $B^0$ is the analytic space of such Lagrangian submanifolds. 
One can show that $B^0$ is smooth (see next subsection). The fiber of the integrable system over the Lagrangian submanifold $L\in B^0$ is the Albanese variety $Pic_0(L)^{\ast}$.  Later we are going to generalize the above picture to the case of singular (possibly irregular)  Higgs fields $\phi$.

\subsection{Smoothness of the moduli space of deformations of complex Lagrangian submanifolds}

The reader can skip this subsection, since its results will not be used in the paper. We also remark that more general smoothness results in the dg-setting were obtained in [BanMan]. Nevertheless we present a different approach which has some benefits on its own.

In order to demonstrate smoothness we need the following result.

\begin{prp} Let $L\subset X$ be a compact complex Lagrangian submanifold in a quasiprojective symplectic manifold $X$. Then the moduli space of deformations of $L$ is smooth (i.e. the deformation theory is unobstructed).

\end{prp}

{\it Proof.} The proof will consist of several steps.

{\it Step 1. } We start with general remarks about deformations in the case of characteristic zero. If we study the formal deformation theory which is controlled by an $L_{\infty}$-algebra with (possibly) non-trivial cohomology in strictly positive degrees, then the corresponding deformation functor from the category of Artin algebras to sets is represented by a pro-Artin scheme, say, $Y$.

Let us now fix $k>0$ and consider the formal deformation theory with the deformation functor  on Artin algebras given by $R\mapsto Hom(Spec(\C[t]/(t^k)\otimes R),Y))$. The corresponding moduli space is the formal neighborhood of the map to the basis point in the formal scheme $Hom(Spec(\C[t]/(t^k)),Y)$. More generally we can consider any map $f\in Hom(Spec(\C[t]/(t^k)),Y)$ and study its formal deformation theory, thus getting a formal scheme $Y_f$.

Let us now recall the following result due to Z.Ran (see [Ra]): {\it $Y$ smooth if and only if for any $k,f$ the tangent space to $Y_f$ at the point $f$ is a free $\C[t]/(t^k)$-module.}

\begin{rmk} This tangent space can be identified with space of such maps $Spec(\C[t]/(t^k)\otimes \C[s]/(s^2))\to Y$ that their restriction to $Spec(\C[t]/(t^k)$ coincides with $f$. The structure of $\C[t]/(t^k)$-module on the tangent space then comes from the natural action of the monoid $\C[t]/(t^k)$ endowed with the operation of multiplication. More precisely an element $a(t)$ of the monoid acts on $\C[t]/(t^k)\otimes \C[s]/(s^2)$ as $t\mapsto t, s\mapsto a(t)s$. The proof of Ran's result in one direction is straightforward: if $Y$ is smooth then the tangent space  to the scheme $Hom(Spec(\C[t]/(t^k)),Y)$ is a free $\C[t]/(t^k)$-module for any $k\ge 1$.

\end{rmk}

{\it Step 2.} Another general remark is that for any finite-dimensional Artin algebra $R$ one can speak about smooth projective varieties over $Spec(R)$. The degeneration of Hodge-to-de Rham spectral sequence holds for such varieties. Indeed the de Rham cohomology forms a free $R$-module. Hodge cohomology coincides with de Rham cohomology at the marked point of $Spec(R)$. Then at the generic point of $Spec(R)$ Hodge cohomology can only drop. But this is impossible, because the spectral sequence implies that the rank of the de Rham cohomology must also drop, but it is constant.

{\it Step 3.}

\begin{lmm} Deformation theory of $L$ as a complex Lagrangian submanifold coincides with its deformation theory as a complex submanifold.
\end{lmm}

{\it Proof of Lemma.} Suppose $R$ is a finite-dimensional Artin algebra. Consider a family $L_s, s\in Spec(R)$ of complex  submanifolds of $X$. Then the restriction of the holomorphic symplectic form $\omega^{2,0}_X$ to each $L_s$ is a closed $2$-form. Let us assume  that it is non-trivial. Then we have a non-trivial family of de Rham  cohomology classes, which is equal to zero at the marked point of $Spec(R)$.  By Step 2 we arrive to the contradiction, since the family of such cohomology classes must be constant. 

{\it Step 4.}
Let us take $R=Spec(\C[t]/(t^k))$ and consider a family of {complex \it Lagrangian} submanifolds $L_s$ over $Spec(R)$ which coincides with the given $L$ at the marked point. Let us consider their first order  infinitesimal deformations as {\it submanifolds}, forgetting the Lagrangian structures. We are allowed to do that by Step 3. But the first order deformations of a manifold are given by sections of the normal bundles. Since $L$ is Lagrangian the normal bundle can be identified with the space of $1$-forms on $L$. Hence the tangent space can be identified with the space of (global) $1$-forms on $L\times Spec(\C[t]/(t^k))$. By Step 2 this space does not jump. This concludes the proof of the Proposition. $\blacksquare$

\subsection{Hitchin systems with irregular singularities}
 
Typically Hitchin systems on $C$ are studied for at most logarithmic singularities of the Higgs fields. The irregular case is less developed. 

For any point $x_0\in C$ we denote by $K_{x_0}$ the field of Laurent series at $x_0$, i.e. $K_{x_0}\simeq \C((t))$, where $t$ is a coordinate on the formal disc centered at $x_0$. The algebraic closure $\overline{K}_{x_0}$ is the field of Puiseux series: $\overline{K}_{x_0}\simeq \cup_{N\ge 1}\C((t^{1/N}))$.
The Galois group $Gal(\overline{K}_{x_0}/K_{x_0})$ is a topological group isomorphic to  $\widehat{\Z}$. Its topological generator acts on $\overline{K}_{x_0}$ as $t^{1/N}\mapsto e^{{2\pi i}\over {N}}t^{1/N}$.
Denote by ${\mathcal O}_{\overline{K}_{x_0}}\subset \overline{K}_{x_0}$ the ring of integers, ${\mathcal O}_{\overline{K}_{x_0}}\simeq \cup_{N\ge 1}\C[[t^{1/N}]]$.

\begin{defn}

A singular term  at the point $x_0\in C$ is an orbit of the Galois action of $\widehat{\Z}$ on the vector space $\overline{K}_{x_0}/{\mathcal O}_{\overline{K}_{x_0}}$.
\end{defn}

In a local coordinate $t=x-x_0$ one can represent a singular term as 
a finite  sum $c=\sum_{\lambda\in \Q_{<0}}c_{\lambda}(x-x_0)^{\lambda}$ 
considered modulo the action on coefficients  of the finite cyclic group ${\Z}/N\Z$ given by $c_{\lambda}\mapsto c_{\lambda}e^{2\pi i \lambda}$,
where the number $N$ called {\it  the ramification index of the singular term}  is defined as the minimal $N\ge 1$ such that all exponents $\lambda$ with $c_{\lambda}\ne 0$ belong to ${1\over{N}}\Z$.

\begin{defn} An irregular data on a smooth projective curve $C$ is given by a tuple $\{x_i\}_{i\in I}, \{r_{i,\alpha}\},\{N_{i,\alpha}\},(c_i^{\alpha})_{\alpha\in Q_i},r)$ where:

a) $x_i\in C, i\in I$ is a finite collection of distinct points ;

b) $r\ge 1$ is an integer number  called rank;

c) a finite collection of distinct singular terms $c^{\alpha}_i, \alpha\in Q_i$ at every point $x_i$; to each term we assign  the multiplicity $r_{i,\alpha}\ge 1$, and we require that the singular term  has the ramification index $N_{i,\alpha}$.

We require that for each $x_i$ the following identity holds
$\sum_{\alpha}r_{i,\alpha}N_{i,\alpha}=r$.

\end{defn}

\begin{rmk} In the case of  $GL(r)$-Hitchin system without singularities we can add dummy marked points $x_i$ and set all $c^{\alpha}_i=0$ and $r_{i,\alpha}=r$ for all marked points $x_i$.

\end{rmk}

Let us compactify $T^{\ast}C$ to $\overline{T^{\ast}C}=T^{\ast}C\cup C_{\infty}$ where $C_{\infty}\simeq C$ is the divisor at infinity. The canonical holomorphic symplectic form $\omega_{T^{\ast}C}$ has pole of order $2$ at $C_{\infty}$. 
Let us consider various smooth projective surfaces $W$ together with regular maps $f:W\to\overline{T^{\ast}C}$ which are birational. We also demand that the pull-back $\omega_W=f^{\ast}(\omega_{T^{\ast}C})$ has only poles but does not have zeros, and such that the complement to the set of poles of $\omega_W$ is isomorphic (via the morphism $f$) to $T^{\ast}C$. Equivalently, such $W$ is obtained by a sequence of blow-ups
$...\to W_2\to W_1\to W_0=\overline{T^{\ast}C}$, where $W_i=Bl_{p_i}(W_{i-1})$, and each point $p_i\in W_{i-1}$ is either a smooth point of a divisor in $W_{i-1}$ where the form $\omega_{W_{i-1}}$ has pole of order at least $2$, or the point $p_i$ is the intersection of two divisors where the form has poles of order at least $1$ (notice that in our case the situation when the orders of both poles are equal to $1$ is impossible).
Such surfaces naturally form a projective system. The set $Irr_1(f^{-1}(C_{\infty}))$ of irreducible components at which the symplectic forms $\omega_W$ have poles of order $1$ forms an inductive system of sets. It is easy to see that the inductive limit of this system can be identified with the set of singular terms. More precisely, to a {\it non-zero} singular term $c$ represented by a formal germ $\widehat{c}$ at $x_0\in C$ we associate a unique divisor $D\in Irr_1(f^{-1}(C_{\infty}))$ on an appropriate blow-up, such that the graph of $d\widehat{c}$ intersects  the divisor $D$ at smooth point as long as $x\to x_0$. If $c=0$ the corresponding  divisor is the exceptional divisor of the blow-up of  $\overline{T^{\ast}C}$ at the point $(x_0)_{\infty}\in C_{\infty}\simeq C$ corresponding to $x_0$.  

Given an irregular data $(\{x_i\},\{r_{i,\alpha}\},\{N_{i,\alpha}\},\{c_i^{\alpha}\},r)$ we consider a minimal surface $W$ as above such that the divisors corresponding to all singular terms are contained in $W$.  It is easy to see by induction in the number of blow-ups that all elements of 
$Irr_1(f^{-1}(C_{\infty}))$ are disjoint smooth rational curves each of which contains only one double point of the divisor $f^{-1}(C_{\infty})$. The complement to this double point is an affine line. Thus to a singular term $\sigma=(x_0,c)$ we have associated an affine line, which we will denote by 
${\mathbb A}^1(\sigma)$.

\begin{lmm} This affine line carries a naturally defined coordinate (i.e. it is canonically identified with ${\mathbb A}^1$).

\end{lmm}

{\it Proof.} Consider the $1$-form, which is the pull-back of the Liouville $1$-form $ydx$ from $T^{\ast}C$. Take any rational curve transversal to ${\mathbb A}^1(\sigma)$. Then the residue of the $1$-form at the intersection point does not change if we vary the transversal curve. Indeed, by Stokes theorem the comparison of residues reduces to the computation  of the integral of the symplectic form over a  $2$-dimensional chain. The latter can be made arbitrarily small. Hence the residue is well-defined and gives the desired coordinate. $\blacksquare$

\begin{defn} An additive refined irregular data on $C$ consist of:

a) an irregular data $(\{x_i\},\{r_{i,\alpha}\},\{N_{i,\alpha}\},\{c_i^{\alpha}\},r)$;

b) for any pair of indices $i,\alpha$ a finite subset 
$\Sigma_{i,\alpha}\subset {\mathbb A}^1(x_i,c_i^{\alpha})\simeq {\mathbb A}^1$ such that 
$|\Sigma_{i,\alpha}|\ge 1$;

c) a map $\Lambda_{i,\alpha}:\Sigma_{i,\alpha}\to \{Partitions\}-\{0\}$ such that $\sum_{z\in \Sigma_{i,\alpha}}|\Lambda_{i,\alpha}(z)|=r_{i,\alpha}$, where for a partition $\lambda=(\lambda_1,...,\lambda_k)$ we use the notation $|\lambda|=\lambda_1+2\lambda_2+3\lambda_3+...$.

\end{defn}

We will often skip the adjective ``additive" and will speak simply about refined irregular data. In Section 8.4 we are going introduce multiplicative refined irregular data.

The pair $(\Sigma_{i,\alpha},\Lambda_{i,\alpha})$ can be thought of as a conjugacy class in $gl(r_{i,\alpha},\C)$. More precisely, $\Sigma_{i,\alpha}$ corresponds to the set of eigenvalues, and the partition $\Lambda_{i,\alpha}(z), z\in \Sigma_{i,\alpha}$ describes multiplicities of the corresponding Jordan blocks.

In the case of Hitchin system with regular singularities all $c_i^{\alpha}=0$, all $r_{i,\alpha}=r$, and the above conjugacy classes in $gl(r,\C)$ can be thought of as the conjugacy classes of the residues of the Higgs field with has poles of order $1$ at the marked points.

With a refined irregular data we are going to associate a pair consisting of a complex symplectic surface $Y$ which contains $T^{\ast}C$ as a Zariski dense open subset and homology class $\beta\in H_2(Y,\Z)$ (the fundamental class of the spectral curve).

Let us describe the construction. Let $W$ be the above-described smooth projective surface constructed for a given irregular data. Now we would like to take the refinement into account. For each pair of indices $i,\alpha$ and any $z\in \Sigma_{i,\alpha}\subset {\mathbb A}^1:={\mathbb A}^1(x_i,c_i^{\alpha})\subset W$ we will make a sequence of blow-ups of $W$ starting with the blow-up at $z$. More precisely, let $D_{i,\alpha}$ be the closure of ${\mathbb A}^1$. We consider a sequence of blow-ups such that the center of each blow-up is a point on the strict transform of 
$D_{i,\alpha}$ corresponding to the point $z$. The number of elements in the sequence of blow-ups is $max\{k|(\Lambda_{i,\alpha}(z))_k\ne 0\}$, where $(\lambda)_k$ denotes the $k$-th component of the partition $\lambda$.

{\it We denote by $\overline{Y}$ the resulting projective surface}. The complement to the divisor of poles of the holomorphic symplectic form is denoted by $Y$. If the irregular data is non-empty then $Y$ is strictly bigger than $T^{\ast}C$, since it contains a chain $E_{1,i,\alpha,z}, E_{2,i,\alpha,z},...$ of rational curves associated with triples $(i,\alpha,z), z\in \Sigma_{i,\alpha}$. Here $E_{1,i,\alpha,z}$ are exceptional divisors of the blow-ups at which the symplectic form is regular. The numeration is chosen in such a way that  $E_{1,i,\alpha,z}$ appears at the first blow-up in case if at least one $c_i^{\alpha}$ is non-zero. If all $c_i^{\alpha}=0$ and $z$ is the intersection of the closure of vertical fiber $T_{x_i}^{\ast}C$ with an appropriate component of the divisor $W-T^{\ast}C$ we add to the chain the closure $\overline{T_{x_i}^{\ast}C}$ and call it $E_{1,i,\alpha,z}$.

We illustrate the above discussion by the figure below.
Notice that the numbers $0, 1, 2, 3$ on the figure refer to the order of poles of the symplectic $2$-form. The figure contains divisors (and their intersection points) obtained at all steps of our sequence of blow-ups. In order to see the final diagram of them, one should keep all the divisors on the figure as well as only those intersection points which can be reached from the ``southwest corner" without crossing lines.

\vspace{2mm}

\centerline{\epsfbox{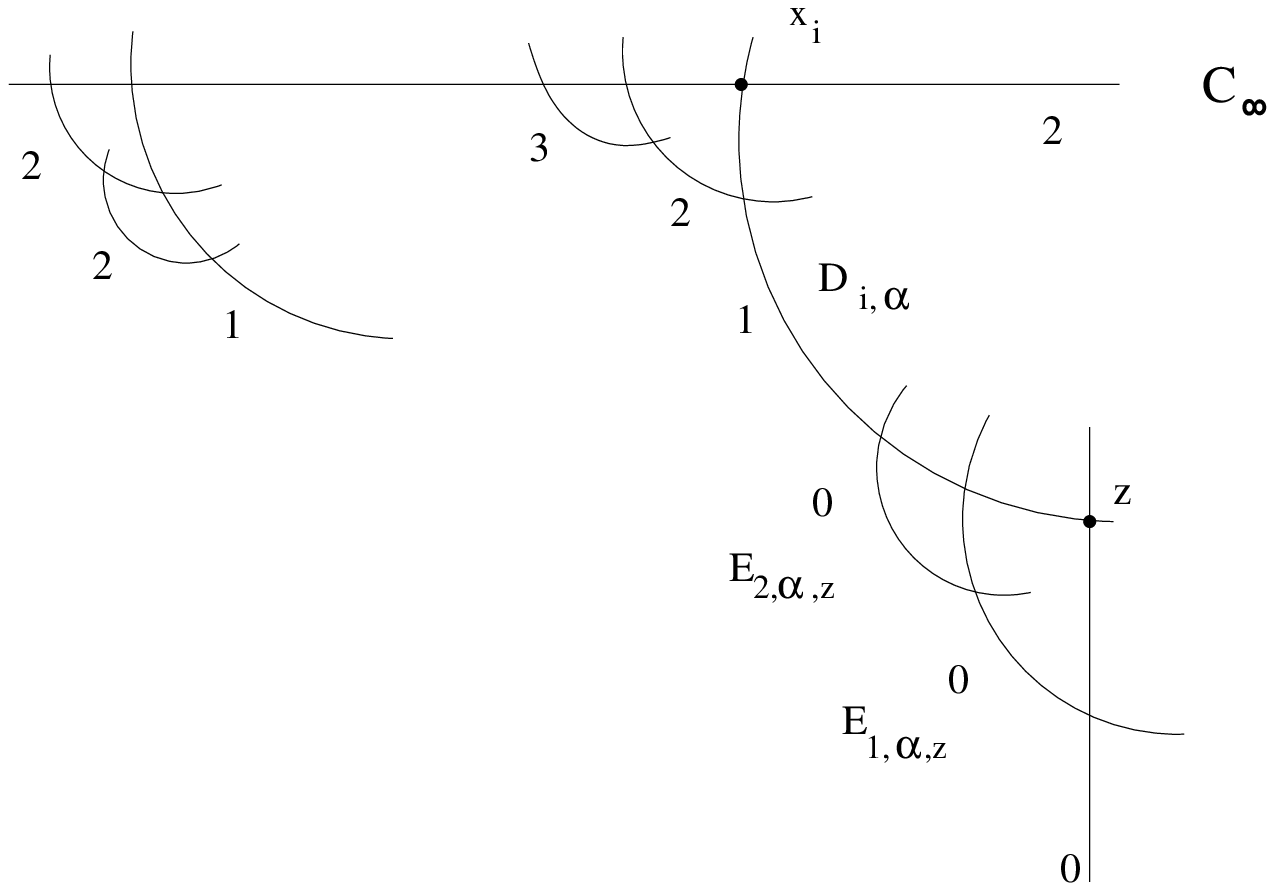}}

\vspace{2mm}

The homology class $\beta\in H_2(Y,\Z)$ is uniquely determined by the following intersection indices:

i) $\beta\cdot [E_{k,i,\alpha,z}]=(\Lambda_{i,\alpha}(z))_k\in \Z_{\ge 1}$;

ii) for generic $x\in C$ we require $\beta\cdot [\overline{T_{x_i}^{\ast}C}]=r$.

Now, having $Y$ and $\beta$ we construct the integrable system in the natural way explained in the end of Section 8.1. The base $B$ consists of compact effective divisors in $Y$ of class $\beta$. We call those divisors  {\it irregular spectral curves} (cf. [Bo1] in the non-refined case). We will assume that the subspace $B^0\subset B$ of smooth connected curves is nonempty. This assumption seems to be satisfied  almost always, e.g. in case when the genus $g(C)>0$ or in case if  $g(C)=0$ and Hitchin system with regular singularities, provided the additive Deligne-Simpson problem (see e.g. [Kos]) has a solution.

The fiber over $b\in B^0$ is the Jacobian of the spectral curve $S_b$. The corresponding complex integrable system is polarized and has a holomorphic Lagrangian section (zero section). One can also consider a version of this construction, when one takes as fibers the torsors over the Jacobians parametrizing line bundles of a given degree $d\in \Z$. In this case the corresponding integrable system does not have in general a Lagrangian section.

In the case when all partitions are of the type $(1)$ (i.e. $|\lambda_{i,\alpha}(z)|=1$ for all $i,\alpha$, which is an analog of Hitchin system with regular singularities and semisimple monodromy) we extend the above integrable system to a semipolarized one by varying points $z\in \Sigma_{i,\alpha}$. The projection of any smooth spectral curve to $W$ is a smooth curve intersecting transversally the divisor of $1$-st order poles of $\omega_W$. The lattice $\Gamma$ can be identified with the integer first homology of the punctured spectral curve (punctures are  intersection points with the preimage of the curve $C_{\infty}$). This is an example of the integrable system associated with a log-family of Lagrangian subvarieties (in this case curves) in the cotangent bundle (see [Ko1]). In Section 8.4. below we are going to describe an analog of the above construction in the case when the restriction  $|\lambda_{i,\alpha}(z)|=1$ on the partitions is dropped.

\begin{prp} The base of  the integrable system associated with refined irregular data is an affine space.

\end{prp}

{\it Proof.}  Consider a unique line bundle ${\mathcal L}\to \overline{Y}$ which is trivialized on $C_{\infty}$ and such that $c_1({\mathcal L})=\beta$. More precisely we require that $P.D.(c_1({\mathcal L}))=(i_{Y\to \overline{Y}})_{\ast}(\beta)$ where  $i_{Y\to \overline{Y}}$ is the natural embedding and  $P.D.$ denotes the Poincar\'e dual. Then the restriction of ${\mathcal L}$ to $\overline{Y}-Y$ is trivial and trivialized. Consider the space of sections $s\in \Gamma(\overline{Y}, {\mathcal L})$ such that
$s_{|\overline{Y}-Y}=1$. This is an affine space. On the other hand it can be identified with the base $B$ by taking the divisor of zeros of $s$. This proves that $B$ is an affine space. $\blacksquare$.

\begin{exa} In case of $GL(r)$-Hitchin system with logarithmic singularities and semisimple monodromy the above construction  can be also described such as follows. 

The base of the corresponding polarized  integrable system is obtained via the above-described procedure. We make blow-ups at all $(x_i)_{\infty}\in C_{\infty}$ corresponding to $x_i\in C\simeq C_{\infty}$. For each $1\le i\le n$ we fix $r$ distinct points on the exceptional divisor $D_i\simeq {\mathbb A}^1$ corresponding to the eigenvalues of the residues of the Higgs field at $x_i$. 

  Let us consider curves $\overline{\Sigma}$ in $\overline{Y}$ which satisfy the following properties:

a) $\overline{\Sigma}$ intersects each $D_i$ with the multiplicity $1$ at each of the chosen $r$ marked points.

b) $\overline{\Sigma}$ intersects each vertical fiber of $T^{\ast}C$ with intersection index $r$.

c) $\overline{\Sigma}$ does not intersect the preimage of $C_{\infty}$ under the blow-up.
 
The space of such curves forms an affine space which is the base of our Hitchin system.
The total space of the latter is birational to the twisted cotangent bundle  to  the moduli space of vector bundles $Bun_{GL(r),x_1,...,x_n}$ on  $C$ endowed with  a choice of a full flag at each point $x_i$. The parameters of the twist correspond to eigenvalues of the residues of Higgs fields at the marked points. We consider only the stable locus in the moduli space of Higgs bundles with logarithmic singularities.

In case when we do not mark points on $D_i, 1\le i\le n$ (i.e. we do not fix the eigenvalues) we obtain a log family of spectral curves as in [Ko1] which forms the base of a semipolarized integrable system.

\end{exa}

\subsection{The Betti version}

Recall that Zariski open part $B^0\subset B$ can be identified with the set of possible irregular data for which the corresponding spectral curve is smooth. Recall that an irregular data includes points $z\in \Sigma_{i,\alpha}$ representing eigenvalues of the residue of the Higgs field at marked points $x_i$. Let us vary each point $z$ along the affine line where it naturally belongs  in such a way that different points do not coincide. Morally this variation corresponds to adding the local system $\underline{\Gamma}_0$ to the picture (indeed fixing the eigenvalues of the monodromy as well as singular terms at marked points
corresponds to a choice of symplectic leaf). In this subsection we are going to discuss the Betti version of the story (cf. also [Bo1]).

Let us recall the formal classification of irregular connections over the formal punctured disc or, equivalently, over the field $K:=\C((t))$. Let $E$ be an $r$-dimensional $K$-vector space endowed with a connection $\nabla: E\to E\otimes \Omega^1_{K}$. Then as a $D$-module $(E,\nabla)$ admits a canonical decomposition 
$$(E,\nabla)\simeq\oplus_{\alpha}(E_{c^{\alpha}},\nabla^{c^{\alpha}})$$
over a finite collection of singular terms $(c^{\alpha})$. We denote the ramification index of $c^{\alpha}$ by $N_{\alpha}$. More precisely the above decomposition can be described such as follows.
For each $\alpha$ let us choose a representative $\overline{c}^{\alpha}\in \C((t^{1/N_{\alpha}}))$ of the singular term $c^{\alpha}$. Then each $D$-module $(E_{c^{\alpha}},\nabla^{c^{\alpha}})$ is isomorphic to the direct image of the canonical $N_{\alpha}$-covering $Spec (\C((t^{1/N_{\alpha}})))\to Spec(\C((t)))$ of the $D$-module which is the tensor product of a vector bundle of rank $r_{\alpha}$  endowed with a connection with regular singularities and a rank one $D$-module $M_{c^{\alpha}}$, which is also a vector bundle with a connection $\nabla^{\alpha}$, such that the generator $m_{\alpha}\in M_{c^{\alpha}}$ satisfies the condition $\nabla^{\alpha}(m)=m\otimes d(\overline{c}^{\alpha})$. We also have $r=\sum_{\alpha}r_{\alpha}N_{\alpha}$.

Let $C$ be a smooth projective curve with marked points $x_i, 1\le i\le n$ and $(E,\nabla)$ be an algebraic vector bundle with connection on the punctured curve $C-\{x_i\}_{1\le i\le n}$. Let us choose a formal coordinate at each point $x_i$. Then the formal expansions of $(E,\nabla)$ at the marked points give rise to a collection of singular terms 
$(c_i^{\alpha})$ as well as to a collection of positive integers $r_{i,\alpha}, N_{i,\alpha}$ derived from the above formal classification. In this way we obtain the {\it canonical irregular data associated with $(E,\nabla)$} (it is canonical in the sense that it does not depend on the choice of formal coordinates at the marked points). 

{\it Irregular Riemann-Hilbert correspondence} gives a topological description (``Betti side") of the complex analytic stack of algebraic vector  bundles with connection with prescribed  irregular data (``de Rham side"). Namely with the pair $(E,\nabla), rk\, E=r$ we associate a local system of rank $r$ (i.e. a locally constant sheaf in analytic topology) on $C-\{x_i\}_{1\le i\le n}$ endowed with the so-called Stokes structure. The local system consists of analytic germs of flat sections of $(E,\nabla)$. We are going to recall the notion of Stokes structure on $\C$ following [KaKoPa]. It describes  the local picture on a general  curve $C$.

 For any marked point $x_i$ and a generic ray $x_i+\varepsilon e^{i\varphi}$ for sufficiently small $\varepsilon >0$ we have a filtration of the  local system of flat sections by the exponential growth of a section restricted to the ray. Terms of the filtration can be identified with intersection points of our ray with the closed real analytic curve 
$$\theta\mapsto exp(Re(c^{\alpha}_i(x_i+\varepsilon e^{i\theta})))e^{i\theta},$$
i.e. when $\theta=\varphi$.

Those filtrations are subject to the conditions described in [KaKoPa]. They  are called a {\it Stokes structure at $x_i$}. In particular, normalization of each curve is a circle $S_{x_i,\alpha}^1$ with the winding number about $x_i$ equals to $N_{i,\alpha}$. On the union $\cup_{\alpha}S_{x_i,\alpha}^1$ we have the associated graded local system. The rank of this local system is $r_{i,\alpha}$.
Collection of Stokes structures for all points $x_i, 1\le i\le n$ is called the {\it Stokes structure for $(E, \nabla)$}.

Let us fix an irregular data. Then the result of Malgrange [Mal] says that there is  a one-to-one correspondence between algebraic connections on $C-\{x_i\}_{1\le i\le n}$ producing our irregular data and local systems on $C-\{x_i\}_{1\le i\le n}$  endowed with the Stokes structure with the  singular terms and discrete parameters $r_{i,\alpha}, N_{i,\alpha}$ derived from the irregular data.

Recall that in Section 8.2 we also defined refined irregular data for the moduli space of Higgs bundles (``Dolbeault side"). Now we would like to introduce a similar notion on the Betti side.

\begin{defn} Multiplicative refined irregular data are defined exactly in the same way as additive irregular data with the only change that $\Sigma_{i,\alpha}\subset \C^{\ast}\subset \C={\mathbb A}^1(\C)$.

\end{defn}
A local system endowed with a Stokes structure defines a 
multiplicative refined irregular data. Namely, the set $\Sigma_{i,\alpha}$ is defined as the set of eigenvalues of the auxiliary local systems on circles $S^1_{x_i,\alpha}$ and the partitions correspond to the sizes of Jordan blocks.

\begin{defn} Let us fix a multiplicative refined irregular data $\sigma$ on $C$. Denote by $M_{Betti}(\sigma)$  the Artin stack over $\C$ of local systems of rank $r$ on $C-\{x_i\}_{1\le i\le n}$ endowed with Stokes structure and such that the corresponding  multiplicative refined irregular data coincides with $\sigma$
(in particular the rank of the local system on $S_{x_i,\alpha}^1$ is equal to the number $r_{i,\alpha}$ from $\sigma$). 
\end{defn}

We denote by  
$M_{Betti}^{simp}(\sigma)$ the algebraic space over $\C$ of isomorphism classes of objects of $M_{Betti}(\sigma)$ which are simple as objects of the abelian category of local systems endowed with Stokes structure.
Equivalently, the corresponding holonomic $D$-module on $C-\{x_i\}_{1\le i\le n}$ is simple. The space $M_{Betti}^{simp}(\sigma)$ is smooth. By analogy with the case of regular singularities we expect that it carries a symplectic structure. Moreover we expect that it is a quasi-affine scheme. In order to explain the latter point it is convenient to consider a larger Artin stack 
$M_{Betti}^{\prime}(\sigma)$ obtained by weakening of some conditions in the definition of $M_{Betti}(\sigma)$. Recall that in the definition of 
$M_{Betti}(\sigma)$ we required  that the monodromies of local systems along the circles  $S_{x_i,\alpha}^1$ belong to certain conjugacy classes $C_{i,\alpha}:=C_{i,\alpha}(\sigma)\subset GL(r_{i,\alpha},\C)$. In the definition of $M_{Betti}^{\prime}(\sigma)$ we relax this condition and say that the monodromies belong to the closures $\overline{C}_{i,\alpha}$.


Let us denote by $M_{Betti}^{coarse}(\sigma)$ the affine scheme  $Spec( {\mathcal O}(M_{Betti}^{\prime}(\sigma)))$.

\begin{que} Is $M_{Betti}^{simp}(\sigma)$ an open subscheme of $M_{Betti}^{coarse}(\sigma)$?

\end{que}

Positive answer will imply that  $M_{Betti}^{simp}(\sigma)$ is a quasi-affine scheme.

Here are some arguments in favor the positive answer to the Question 8.4.3 in the case of regular singularities. First we observe that there are many functions on $M_{Betti}^{\prime}(\sigma)$ given by traces of holonomies along closed loops. We claim that the sizes of Jordan blocks of monodromies around points $x_i$ can be also detected. 
Let us illustrate the claim by an example. Let $y\in C-\{x_i\}_{1\le i\le n}$ be a base point. Then our local system gives rise to an $r$-dimensional representation $\rho$ of the fundamental group $\pi_1(C-\{x_i\}_{1\le i\le n},y)$. Let us assume for simplicity that for some point $x_{i_0}$ the monodromy $\rho(l_{i_0})$ is unipotent of order $k\ge 1$, where $l_{i_0}$ is a based loop which is freely homotopic to a small loop surrounding $x_{i_0}$. Then $(\rho(l_{i_0})-id)^k=0$. Hence $Tr\left(\rho(l_{i_0})-id)^k \rho(l)\right)=0$ for any $l\in \pi_1(C-\{x_i\}_{1\le i\le n},y)$. This equation gives identities between traces of monodromies. Similar considerations can be applied to $\bigwedge^i\rho, 1\le j\le k$. In this way we recover information about sizes of Jordan blocks. We expect that similar arguments work in general case.

The conclusion is that conjecturally for each multiplicative refined irregular data $\sigma$ we have a smooth symplectic quasi-affine variety $M_{Betti}^{simp}(\sigma)$. This variety depends on continuous parameters which are eigenvalues of monodromies on  $S_{x_i,\alpha}^1$. Let us allow the eigenvalues to vary in such a way that they do not coincide. The total space should be a Poisson manifold. We can also enumerate the eigenvalues in the following way.

\begin{defn} A combinatorially refined irregular data is an irregular data endowed with a collection of integers $s_{i,\alpha}\ge 1$ and a  collection of maps $\Psi_{i,\alpha}: \{1,...,s_{i,\alpha}\}\to \{Partitions\}-\{0\}$ such  that $\sum_{1\le j\le s_{i,\alpha}}|\Psi_{i,\alpha}({j})|=r_{i,\alpha}$.

\end{defn}

For a given  combinatorially refined irregular data $\tau$ we  construct a larger moduli space   $M_{Betti}^{simp,en}(\tau)$ (enumerated version of $M_{Betti}^{simp}(\sigma)$)  whose set of $\C$-points is the disjoint union of $M_{Betti}^{simp}(\sigma)(\C)$, over the set pairs $(\sigma, f), f=(f_{i,\alpha})$, where each $f_{i,\alpha}$ is a bijection $\Sigma_{i,\alpha}\simeq \{1,...,s_{i,\alpha}\}$. Here we assume that the map $\Lambda_{i,\alpha}\circ f_{i,\alpha}^{-1}=\Psi_{i,\alpha}$  where $\Psi_{i,\alpha}$ are the maps from the definition of the combinatorially refined irregular data $\tau$. 

It is easy to see that the space $M_{Betti}^{simp,en}(\tau)$ is fibered over the hypersurface $\mathcal{H}_{Betti}(\tau)$ in
 $\prod_{i,\alpha}((\C^{\ast})^{s_{i,\alpha}}-Diag)$ given by the equation 

$$\prod_{i,\alpha}\prod_{1\le j\le s_{i,\alpha}}\lambda_{i,\alpha,j}^{r_{i,\alpha,j}}=(-1)^{\sum_{i,\alpha}(N_{i,\alpha}-1)r_{i,\alpha}}.$$
This equation comes from the fact that the product of determinants of monodromies around marked points $x_i$ is equal to $1$. The appearance of the factor of $-1$ in the RHS is due to the presence of coverings in the definition of Stokes structure.

\begin{que} Is $M_{Betti}^{simp,en}(\tau)$ a smooth Poisson quasi-affine scheme, with the space of symplectic leaves identified with $\mathcal{H}_{Betti}(\tau)$?
\end{que}

We have not discussed above the origin of the symplectic structure on $M_{Betti}^{simp}(\sigma)$ (hence the Poisson structure on $M_{Betti}^{simp,en}(\tau)$). The case of semisimple monodromy and Laurent series  (i.e. no fractional powers appear) 
was studied in [Bo1].

One can also define an affine scheme  $M_{Betti}^{coarse,en}(\tau)\supset M_{Betti}^{simp,en}(\tau)$  endowed with a surjective map   to the complex torus $\overline{\mathcal H}_{Betti}(\tau)$ which is the closure of  ${\mathcal H}_{Betti}(\tau)$ in the ambient torus  (in other words it is a shifted subtorus in 
$\prod_{i,\alpha}(\C^{\ast})^{s_{i,\alpha}}$ given by the above equation).

In order to define $M_{Betti}^{coarse,en}(\tau)$ let us consider the  moduli stack  $M_{Betti}^{\prime,en}(\tau)$ of the following structures: local systems on $C-\{x_i\}_{1\le i\le n}$ endowed with the Stokes structure and decompositions of the associated local systems on $S^1_{x_i,\alpha}$ into the direct sums labeled by $j, 1\le j\le s_{i,\alpha}$. For each direct summand  the monodromy has only one eigenvalue $\lambda_{i,\alpha,j}\in \C^{\ast}$. The unipotent part of the monodromy is dominated by the partition $\Psi_{i,\alpha}(j)$. The latter means that the conjugacy class of the unipotent part of the monodromy belongs to the closure of the unipotent conjugacy class corresponding to $\Psi_{i,\alpha}(j)$. Finally we define $M_{Betti}^{coarse,en}(\tau):= Spec({\mathcal O}(M_{Betti}^{\prime,en}(\tau)))$. The main point of this definition is to allow the eigenvalues $\lambda_{i,\alpha,j}$ to coincide for different values of $j$.

\begin{rmk}
1) The space $M_{Betti}^{coarse,en}(\tau)$ seems to be an analog of $X$-variety in the theory of cluster varieties. Presumably one can also define an analog of $A$-variety (see [FoGo1]).



2) Rescaling $c^{\alpha}_i\mapsto c^{\alpha}_i/\zeta,\, \zeta\in \C^{\ast}$  gives rise to  non-linear local systems over $\C^{\ast}$ of all versions of $M_{Betti}$. Taking the fiber over $\zeta=1$ we see that it is endowed with an automorphism given by the monodromy. This automorphism corresponds to the Coxeter automorphism in the theory of cluster algebras.

\end{rmk}

\subsection{Semipolarized irregular systems}

First we would like to describe an additive analog $M_{Dol}^{sm,irred,en}(\tau)$  of the space $M_{Betti}^{sm,en}(\tau)$. Here the notation $sm, irred$ means  {\it smooth, irreducible} correspondingly and refers to spectral curves. Let us comment on the notation. Suppose that the genus $g(C)>1$, and we are dealing with regular Hitchin system (i.e. there are no marked points and the irregular data is empty). Then by Corlette-Simpson result there is a one-to-one correspondence between simple local systems on $C$ and stable Higgs bundles of degree zero. Recall that Higgs bundles can be identified with coherent sheaves on $T^{\ast}C$ with pure $1$-dimensional compact support. Under this identification line bundles of a certain degree give a Zariski open subset in the moduli space of stable Higgs bundles of degree zero. This observation motivates our notation in the general case.

Recall (see Section 8.2) that with an (additive) refined irregular data $\sigma$ we can associate a complex integrable system with central charge and zero section. Strictly speaking we consider only an open part $B^0(\sigma)$ of the base  consisting of smooth irreducible spectral curves in the surface associated with $\sigma$. Fibers are Jacobians of spectral curves. Part of the refined irregular data consists of sets of points $\Sigma_{i,\alpha}\subset {\mathbb A}^1(x_i,c_i^{\alpha})$. In the case of Hitchin systems with regular singularities they are eigenvalues of the residues of the Higgs field at marked points. 

Let us fix a combinatorially refined irregular data $\tau$. Similarly to the case of 
$M_{Betti}^{simp,en}(\tau)$ we define   the moduli space $M_{Dol}^{sm, irred, en}(\tau)$ by allowing enumerated points  $z=z_{i,\alpha,j}\in \Sigma_{i,\alpha}$ to vary along the corresponding affine lines ${\mathbb A}^1(\C)=\C$ in such a way that they do not collide. 
The  resulting space  $M_{Dol}^{sm, irred, en}(\tau)$ is the total space of a family of polarized integrable systems parametrized by the hypersurface $\mathcal{H}_{Dol}:=\mathcal{H}_{Dol}(\tau)$ in  $\prod_{i,\alpha}(\C^{s_{i,\alpha}}-Diag)$ singled out by the equation $$\sum_{i,\alpha}\sum_{1\le j\le s_{i,\alpha}} z_{i,\alpha,j}r_{i,\alpha,j}=0.$$ 
This equation follows from the condition that sum of the residues of the Liouville form restricted to the spectral curve vanishes.

Denote by $B^0=B^0(\tau)$ the total space of the fibration over $\mathcal{H}_{Dol}$ whose fibers are bases of the above polarized integrable systems.  The projection $p: B^0\to \mathcal{H}_{Dol}$ is a smooth morphism of smooth algebraic varieties.
Fiber $p^{-1}(h), h\in\mathcal{ H}_{Dol}$ can be identified by the previous considerations with the moduli space of smooth irreducible spectral curves in the appropriate surface $\overline{Y}=\overline{Y}(h)$. We will construct a local system of lattices $\underline{\Gamma}\to B^0$ endowed with a covariantly constant integer skew-symmetric pairing 
$\langle\bullet,\bullet\rangle:\bigwedge^2 \underline{\Gamma}\to \underline{\Z}_{B^0}$ and a central charge $Z: \underline{\Gamma}\to {\mathcal O}^{an}_{B^0}$ in such a way that we obtain a semipolarized integrable system with central charge. Moreover the local system $\underline{\Gamma}_0:=Ker\, \langle\bullet,\bullet\rangle$ will be trivial, i.e.   $\underline{\Gamma}_0\simeq \Gamma_0\otimes \underline{\Z}_{B^0}$, where $\Gamma_0$ is a fixed lattice. The restriction $Z_{|\underline{\Gamma}_0}$ will be identified with the composition $B^0\to\mathcal{ H}_{Dol}\hookrightarrow Hom(\Gamma_0,\C)$.

Let us explain how to define the dual local system $\underline{\Gamma}^{\vee}$ and the central charge. Recall that an irregular  spectral curve $S$ contains pairwise disjoint effective divisors 
$D_{i,\alpha,z_{i,\alpha,j}}$, where $deg\,D_{i,\alpha,z_{i,\alpha,j}}=|\Lambda_{i,\alpha}(z_{i,\alpha,j})|=|\Psi_{i,\alpha}(j)|$. 
These divisors  are intersections of $S$ with the chain of rational curves $E_{k, i,\alpha, z=z_{i,\alpha,j}}$ defined in Section 8.2.

Then we define the fiber of $\underline{\Gamma}^{\vee}$ over $S$ as the abelian group of $1$-chains on $S$ whose boundaries are $\Z$-linear combinations of divisors $D_{i,\alpha,z_{i,\alpha,j}}$ considered modulo the boundary of $2$-chains.  
Clearly the abelian groups depend continuously on parameters and hence define a local system. Also for any $b\in B^0$ we have a short exact sequence

$$0\to H_1(S_b,\Z)\to \underline{\Gamma}_b^{\vee}\to {\Gamma}_0^{\vee}\to 0,$$
where
$${\Gamma}_0^{\vee}=Ker(f:\oplus_{i,\alpha}\Z^{s_{i,\alpha}}\to \Z),$$
where $f:( n_{i,\alpha,j})_{1\le j\le s_{i,\alpha}}\mapsto \sum_{i,\alpha,j}n_{i,\alpha,j}r_{i,\alpha,j}$.

Dualizing we obtain local systems $\underline{\Gamma}$ and its trivial local subsystem $\underline{\Gamma}_0$. The skew-symmetric pairing on $\underline{\Gamma}$ is the pull-back of the symplectic structure on  $H_1(S,\Z)$.

Next we are going to describe the central charge.
One observes that the fiber over $b\in B^0$ of the  local system $\underline{\Gamma}$ can be identified (up to torsion) with the quotient  $H_1(S_b-D,\Z)/P_b$, where $D=\sqcup_{i,\alpha} D_{i,\alpha}$ is the union of the exceptional irreducible divisors in the surface $\overline{Y}$, where the pull-back of the symplectic form $\omega_{T^{\ast}C}$ to   $\overline{Y}\supset T^{\ast}C$  does not have poles, and $P_b$ is the subgroup generated by $1$-chains $\gamma_{i,\alpha}$ sitting in a small neighborhood in $S_b$ of the intersection $S_b\cap  D_{i,\alpha}$ and such that the linking number of $\gamma_{i,\alpha}$ and  $D_{i,\alpha}$ is equal to zero.
One observes that  integrals of the pull-back of the Liouville form $ydx$ over the elements of $H_1(S_b-D,\Z)/P_b$ are well-defined. This gives us the central charge $Z$.

 The group of automorphisms of the tuple $(B^0,\underline{\Gamma},\langle\bullet,\bullet\rangle,Z)$ contains a finite subgroup $\prod_mSym_{k_m}$. The latter is the the product of symmetric groups with each factor acting on points $z_{i,\alpha,j}$ which belong to the same affine line indexed by $(i,\alpha)$ in such a way that it permutes those points $z_{i,\alpha,j}$ which are endowed with the same partition.  The quotient by the group of automorphisms will be a semipolarized integrable system with central charge and non-trivial local system of lattices $\underline{\Gamma}_0$.

\subsection{Conjectures about mirror duals for Hitchin systems}

Fix a combinatorially refined irregular data $\tau$. The above considerations give rise to a family of  polarized integrable systems parametrized by a variety ${\mathcal H}_{Dol}(\tau)$ which is an open dense subset in 
$Hom(\Gamma_0,\C)$, as well as a semipolarized integrable system which we denote by $(X^0(\tau),\omega^{2,0}(\tau))\to B^0(\tau)$ (all endowed with holomorphic Lagrangian sections). For an individual polarized integrable system corresponding to an irregular data $\sigma$ we have defined the full base $B(\sigma)$ in terms of irregular spectral curves (see Section 8.2), but our integrable systems so far have been defined over the locus $B^0(\sigma)$ of smooth irregular spectral curves.

\begin{conj} 1) There exists a full semipolarized complex integrable system $X(\tau)\to B(\tau)$ containing $X^0(\tau)\to B^0(\tau)$ as an open dense complex integrable subsystem.

2) The corresponding individual polarized integrable systems $X(\tau)_{Z_0}\to B(\tau)_{Z_0}$ in the notation of Section 4.7  have full bases $B(\tau)_{Z_0}=B(\sigma)\simeq \C^{dim_{\C}\,B(\sigma)}$, where $\sigma$ is determined by $\tau$ and $Z_0\in Hom(\Gamma_0,\C)$.

3) The mirror dual $X(\tau)^{\vee,alg}$ to the integrable system $X(\tau)\to B(\tau)$ in the sense of Section 6.3  is an affine scheme which contains $M^{simp,en}_{Betti}(\tau)$ as Zariski open subset.
\end{conj}

\begin{rmk} Strictly speaking part 3) of the above conjecture should be corrected. Recall that the algebraic mirror dual in Section 6.3 was fibered over  the algebraic torus $Hom(\Gamma_0,\C^{\ast})$. The variety $M^{simp,en}_{Betti}(\tau)$ is fibered over the  ${\mathcal H}_{Betti}(\tau)$ which is an open subset in a torsor over $Hom(\Gamma_0,\C^{\ast})$. This discrepancy should be corrected by some twist.
Probably this twist is related to the choice of the canonical $B$-field $B_{can}$ which is $2$-torsion.
 
\end{rmk}

\begin{que} Does the mirror dual coincide with the $Spec({\mathcal O}(M_{Betti}^{coarse,en}(\tau)))$?

\end{que}

Now we are going to discuss the family over $\C^{\ast}$ related to the rescaling $\omega^{2,0}\mapsto \omega^{2,0}/\zeta, \zeta\in \C^{\ast}$. In the case of irregular Hitchin systems this rescaling corresponds to the rescaling $c_i^{\alpha}\mapsto c_i^{\alpha}/\zeta$ of the singular terms. Assuming the above conjecture, we obtain by taking the mirror dual, a holomorphic family of Poisson varieties $X^{\vee,alg}_{\zeta}(\tau)$ over $\C^{\ast}$ containing $M^{simp,en}_{Betti,\zeta}(\tau)$ as open subvarieties. They are locally constant in analytic topology on $\C^{\ast}$, hence we get a local system of algebraic varieties.  

Recall from Section 6.3 that we also have a complex analytic mirror dual $X(\tau)^{\vee}\to Hom(\Gamma_0,\C)$ which is obtained from $X(\tau)^{\vee,alg}\to Hom(\Gamma_0,\C^{\ast})$ via the exponential map. Introducing the parameter $\zeta$ we obtain a complex analytic family $X_{\zeta}(\tau)^{\vee}, \zeta\in \C^{\ast}$. It contains the pull-back $\widetilde{M}^{simp,en}_{Betti,\zeta}(\tau)$ of the space $M^{simp,en}_{Betti,\zeta}(\tau)$ as an open dense subset. We can also define a larger space 
$\widetilde{M}^{coarse,en}_{Betti,\zeta}(\tau)\supset \widetilde{M}^{simp,en}_{Betti,\zeta}(\tau)$ (see Remark 8.3.6).

Recall Conjecture 6.7.3 (extension to $\zeta=0$).  In our case it says that the analytic family  of Poisson varieties $X_{\zeta}(\tau)^{\vee}$ admits analytic extension to $\zeta=0$ with the fiber at $\zeta=0$ isomorphic to $X^{dual}(\tau)$ (see Conjecture 6.7.3). In our case it is reasonable to expect that  $X^{dual}(\tau)\simeq X(\tau)$ since Jacobians of spectral curves are principally polarized abelian varieties.  

\begin{conj} The local system of Poisson varieties  $\widetilde{M}^{coarse,en}_{Betti,\zeta}(\tau)$ admits an analytic extension to $\zeta=0$ with the fiber at zero isomorphic to $M^{sm,irred,en}_{Dol}(\tau)$.
\end{conj}
The total space of the extended to $\zeta=0$ family should be a complex algebraic variety if we endow  fibers $\widetilde{M}^{coarse,en}_{Betti,\zeta}(\tau)$   with algebraic structures coming from the de Rham description of fibers via inverse Riemann-Hilbert correspondence as algebraic vector bundles endowed with irregular $\zeta$-connections. Furthermore, the de Rham description makes the above conjecture almost evident similarly to the well-known case if Hitchin systems without singularities.

Let us describe explicitly the geometric meaning of a germ of holomorphic section of the above analytic family at $\zeta=0$.
Let us fix an additive refined irregular data $\sigma$ on the curve $C$. In particular it gives us a positive integer $r$ (the rank). Consider a holomorphic vector bundle $E$ of rank $r$ over $C\times D_{\varepsilon}$, where $D_{\varepsilon}=\{\zeta\in \C| |\zeta|<\varepsilon\}$ and $\varepsilon>0$ is sufficiently small. We will think of it as a family $E_{\zeta}$ of holomorphic vector bundles. Consider a relative along $C$ meromorphic connection $\nabla$ such that:

a) $\nabla$ has finite order poles at the marked points $x_i, 1\le i\le n$;

b) $\nabla$ has the pole of order $1$ along $C\times \{\zeta=0\}$.

In particular for all $\zeta\in \C^{\ast}$ we have a meromorphic connection $\nabla_{\zeta}$ on the vector bundle $E_{\zeta}$.

Using the formal classification of irregular connections we assign to each $\nabla_{\zeta}$ irregular data with singular terms at the points $x_i$. We require that they coincide with the given irregular data after rescaling $c_i^{\alpha}\mapsto c_i^{\alpha}/\zeta$ of singular terms (here we forget about the refinement). 

Let us choose an additive refined irregular data compatible with $\tau$. In other words we choose points points $z_{i,\alpha,j}\subset {\mathbb A}^1(x_i,c_i^{\alpha})\simeq \C$ such that $z_{i,\alpha,j_1}\ne z_{i,\alpha,j_2}$ for $j_1\ne j_2$.
Moreover we assume that eigenvalues $\lambda_{i,\alpha,j}$ of the monodromy of $\nabla_{\zeta}$ along the circle $S^1_{x_i,\alpha}$ have the form
$$\lambda_{i,\alpha,j}=exp({z_{i,\alpha,j}\over{\zeta}}).$$

Furthermore we assume that for any $i,\alpha,j$ the sizes of the Jordan blocks with eigenvalues $\lambda_{i,\alpha,j}$ as above form the partition $\Psi_{i,\alpha}(j)$ coming from the combinatorially refined irregular data $\tau$ as
as long as  
$${z_{i,\alpha,j_1}-z_{i,\alpha,j_2}\over{\zeta}}\notin 2\pi \sqrt{-1}\Z $$
for $j_1\ne j_2$. We call the latter condition {\it non-exceptionality condition}.
It depends on the formal type of the irregular connection $\nabla_{\zeta}$ at $x_i$  and does not depend on the Stokes structure.


\begin{rmk} One can generalize the above story by allowing the curve $C$ and irregular data to depend analytically on $\zeta$. Then we will require that
$$\lambda_{i,\alpha,j}=exp({z_{i,\alpha,j}\over{\zeta}}+O(1)).$$ 
\end{rmk}

Under the above assumptions the limit $lim_{\zeta\to 0}\zeta\nabla_{\zeta}$ does exist and defines a meromorphic Higgs field $\varphi$ (with poles at the points $x_i$) on the vector bundle $E_0:=E_{|\zeta=0}$. It defines the refined irregular data coinciding with the given one.
The closure $S\subset \overline{Y}$ of the non-compact spectral curve $S^0$ given by $det(\phi-yId)=0\subset T^{\ast}(C-\{x_i\}_{1\le i\le n})\subset Y$ is an irregular spectral curve in our sense.

Assume that $S:=S_b$ is smooth and irreducible (i.e. it corresponds to a point $b\in B^0$). There is a natural line bundle $\LL\to S^0$ corresponding to 
$(E_0)_{|C-\{x_i\}_{1\le i\le n}}$. Extending it to $S$ (the ambiguity for such an extension is a lattice of finite rank) we obtain a line bundle $\LL\to S$, hence a point in $Jac^d(S_b)=Jac^d(S)$ for some degree $d$ belonging to a fiber of the integrable system at $\zeta=0$.
Suppose that for our refined irregular data we have chosen $\zeta\in \C^{\ast}$ such that non-exceptionality condition holds for $\zeta, \{z_{i,\alpha,j}\}$.  Then our connection  $\nabla_{\zeta}$ defines a point $f(\zeta)$ in the covering $\widetilde{M}_{Betti,\zeta}^{coarse,en}(\tau)$ via taking the Stokes structure of $(E_{\zeta}, \nabla_{\zeta})$. The point $f(\zeta)$ depends holomorphically on $\zeta$ such that 
$\zeta\ne 0$ and
$\zeta\ne ({z_{i,\alpha,j_1}-z_{i,\alpha,j_2}})/2\pi \sqrt{-1}k, k\in \Z-\{0\}.$

\begin{conj} 1) Let us assume the Conjecture 8.5.4. Then the map $\zeta\mapsto f(\zeta)$ extends to a germ of an analytic curve at $\zeta=0$. 

2) Moreover  the value $f(0)$ is the point  of the space $M^{sm,irred,en}_{Dol}(\tau)$ corresponding to the line bundle $\LL\to S$.

\end{conj}

Our definition of $\widetilde{M}_{Betti}^{coarse,en}(\tau)$ is a bit artificial. Probably the following version of it will  behave  better. Fix an  irregular 
data $\eta$ (no refinement is chosen). Consider an Artin stack 
${M}_{Betti}(\eta)$ parametrizing local systems on $C-\{x_i\}_{1\le i\le n}$ with Stokes structures compatible with irregular data. We denote by 
${M}_{Betti}^{an}(\eta)$ the corresponding analytic stack.
Let us define an analytic stack $\widetilde{M}_{Betti}^{an}(\eta)$ by adding the following additional data: for each pair $i,\alpha$ a choice of covariantly  constant endomorphism $L_{i,\alpha}$ of the associated local system ${\mathcal E}_{i,\alpha}$ on the circle $S^1_{x_i,\alpha}$ such  that $exp(L_{i,\alpha})$ is the monodromy automorphism for ${\mathcal E}_{i,\alpha}$. In the case when all monodromies are semisimple and have distinct eigenvalues  the set of choices of ``logarithms" $L_{i,\alpha}$ for all $i,\alpha$  is a torsor over the lattice $\Z^{\{(i,\alpha,j)\}}$. This stack can be thought of as a replacement of $\widetilde{M}_{Betti}^{coarse,en}(\tau)$ (the enumeration of eigenvalues is lost). We claim that there exists another Artin stack ${M}_{DR}(\eta)$ such that the (irregular) Riemann-Hilbert correspondence gives an isomorphism of analytic stacks $\widetilde{M}_{Betti}^{an}(\eta)\simeq 
{M}_{DR}^{an}(\eta)$. In the example of connections with regular singularities the stack
${M}_{DR}(\eta)$ is the moduli stack of vector bundles on $C$ endowed with meromorphic connections which have first order pole at the marked points.
Introducing the parameter $\zeta\in \C^{\ast}$ as above by rescaling the singular terms of the irregular data we obtain an algebraic family 
${M}_{DR,\zeta}(\eta)$ of Artin stacks. We believe that it can be extended as an algebraic family of Artin stacks over $\C$ with the fiber at $\zeta=0$ being the stack $M_{Dol}(\eta)$ of semistable generalized Higgs bundles of type $\eta$. The latter are roughly certain coherent sheaves on the Poisson surface $\overline{Y}:=\overline{Y}(\eta)$ (see Section 8.2) with pure one-dimensional support.



Let us  discuss the analog of a combinatorial refinement in this setting. Assume that for each pair $i,\alpha$ we are given a finite unordered collection of partitions $\lambda^{i,\alpha,j}$ (possibly with repetitions) such the sum of weights 
$\sum_j|\lambda^{i,\alpha,j}|$ is equal to the rank $r_{i,\alpha}$ of 
${\mathcal E}_{i,\alpha}$. Then we pose the following condition: {\it  the conjugacy class of the linear operator $L_{i,\alpha}$ (considered as  a linear endomorphism  of the fixed fiber of ${\mathcal E}_{i,\alpha}$) belongs to the closure of the set of such linear operators that for each of its eigenvalues $\mu_j$ the set of  Jordan blocks with the eigenvalue $\mu_j$ defines a collection of partitions which coincides with the given collection of partitions
$\lambda^{i,\alpha,j}$.}
This defines a closed substack of ${M}_{DR}(\eta)$. Presumably one can define a similar substack of $M_{Dol}(\eta)$. Introducing the parameter $\zeta$ we obtain as above a family of Artin substacks over $\C$. Finally we remark that Mirror Symmetry naturally gives us an {\it analytic} family of analytic stacks ${M}_{DR,\zeta}^{an}(\eta)\simeq \widetilde{M}_{Betti,\zeta}^{an}(\eta) $  over $\zeta\in \C^{\ast}$ and  its limit at $\zeta=0$ given by $M_{Dol}^{an}(\eta)$. Also the fiber over $\zeta\ne 0$ have some ``remnants" of the algebraic structure on ${M}_{Betti,\zeta}(\eta)$. The algebraic structures on ${M}_{DR}(\eta)$ and $M_{Dol}(\eta)$ familiar in Geometric Langlands Correspondence seem to play no role in the case of Mirror Symmetry considered before.

\subsection{Remarks about $SL(r)$ case}

In the case of $SL(r)$ Hitchin systems the above considerations have to be modified. Namely, in the definition of the irregular data we impose an additional condition: {\it sum of all branches of the singular terms at each marked point $x_i$ is equal to zero (modulo series which are regular at $x_i$)}.

This condition can be reformulated such as follows. For with each singular term $c_i^{\alpha}=\sum_{\lambda\in \Q_{<0}}a_{\lambda,i}^{\alpha}(x-x_i)^{\lambda}$ we associate its trace 
$$Tr(c_i^{\alpha})=N_{i,\alpha}\cdot \sum_{\lambda\in \Z_{<0}}a_{\lambda,i}^{\alpha}(x-x_i)^{\lambda}\in \C[(x-x_i)^{-1}].$$

Then the above condition says that for each $1\le i\le n$ we have $\sum_{\alpha}Tr(c_i^{\alpha})=0$.

We can impose a similar condition for spectral curves. A spectral curve can be thought of as a graph of a multivalued  closed $1$-form on $C-\{x_i\}_{1\le i\le n}$. We demand that the sum of all branches of the $1$-form vanishes identically.
Fix combinatorially refined irregular data $\tau$. We denote by $B_{SL(r)}(\tau)\subset B(\tau)$ the subspace of spectral curves $S$ which satisfy the above condition. 

\begin{prp} One has $dim\,B_{SL(r)}(\tau)=dim\,B(\tau)-g(C)$, where $g(C)$ is the genus of $C$.

\end{prp}
{\it Proof.} For any spectral curve in $B(\tau)$ the sum of branches of the corresponding multivalued $1$-form is a holomorphic $1$-form on $C$. Also the space of $1$-forms $\Omega^1(C)$ acts on $B(\tau)$ by adding the graph of the $1$-form. This gives an isomorphism $B(\tau)\simeq B_{SL(r)}(\tau)\times \Omega^1(C)$. The result follows. $\blacksquare$

\begin{rmk}
In the case of Hitchin systems with regular singularities the above condition means that the sum of eigenvalues of the singular part of the Higgs field at each $x_i$ is equal to zero.
\end{rmk}

Considering the corresponding local system of {\it symplectic} lattices $\underline{\Gamma}^{symp}$  (which is the quotient of the bigger local system $\underline{\Gamma}$) we see that the fiber of $\underline{\Gamma}^{symp}_S$ is $Prym(S):=Ker(H_1(S,Z)\to H_1(C,\Z))$. The fibers are polarized but not principally polarized. Another choice would be $Coker(H^1(C,\Z)\to H^1(S,\Z))$. The mirror duals to the integrable systems corresponding to \begin{large}                                                                                                                                                                                                                                                                                                                                                                                                                                                           \end{large}these two choices are different and should correspond to $M_{Betti}^{SL(r)}(\tau)$ and $M_{Betti}^{PGL(r)}(\tau)$.

\subsection{Relation to non-compact Calabi-Yau 3-folds}

Having a spectral curve $S\subset Y\subset \overline{Y}$ as above and line bundles $\mathcal{L}_i\to Y, i=1,2$ such that the restrictions of $\mathcal{L}_i, i=1,2$ to $Y-S$ are trivialized and such that $\mathcal{L}_1\otimes \mathcal{L}_2\simeq {\mathcal O}_Y(S)$ we can construct a non-compact Calabi-Yau threefold (total space of the conic bundle over $Y$). Namely, let us fix a section $t\in \Gamma(Y, \mathcal{L}_1\otimes \mathcal{L}_2)$ such that $t_{|S}=0$. Then we consider a subvariety $X$ of the total space $tot(\mathcal{L}_1\oplus \mathcal{L}_2)$ which consists of pairs $(l_1,l_2)\in \mathcal{L}_{1,y}\oplus \mathcal{L}_{2,y}, y\in Y$ such that $l_1\cdot l_2=t(y)$. Writing it in local coordinates as $x_1x_2=f(y_1,y_2)$ we see that $X$ carries a well-defined nowhere vanishing holomorphic volume form $\Omega_X^{3,0}$ locally given by ${dx_1\over{x_1}}dy_1dy_2$. Then such $X$ satisfies the assumptions $A1-A3$ from Section 7. Furthermore the corresponding moduli space of deformations of $X$ discussed in Section 7 is essentially the same as the space $B$ of spectral curves discussed in this section. The fiber of the local system $\underline{\Gamma}$ is isomorphic to
$H_3(X,\Z)$ modulo torsion. An easy calculation with exact sequences of fiber bundles shows that if $\mathcal{L}_1$ is trivial then 
 $H_3(X,\Z)\simeq H_2(Y,S,\Z)$. Then the symplectic lattice $\underline{\Gamma}^{symp}$ is  smaller than $H_1(S,\Z)$, while the kernel lattice $\underline{\Gamma}_0$ is bigger than the one for the lattice described in Section 8.5.

In the case when the monodromies along $S^1_{x_i,\alpha}$ are  semisimple  with non-coinciding eigenvalues one can spell the above discussion in terms of log-families of spectral curves.
In particular there is an open dense part in the moduli space of deformations of $X$  which is isomorphic to the space of log-families of smooth curves in  $T^{\ast}C$.  In order to be in agreement with $3$-dimensional story we need to go from $GL(r)$ Hitchin integrable systems to $SL(r)$ Hitchin integrable systems. 
Then, as we discussed above, for each spectral curve $S_b, b\in B^0_{SL(r)}(\tau)$ the lattice $\Gamma_b^{symp}$ is isomorphic to $H_2(Y,S,\Z)^{symp}\simeq Prym(S_b)\simeq H_3(X,\Z)/Ker(\langle\bullet,\bullet\rangle)$. Periods of the restriction $ydx_{|S_b}$ of the canonical $1$-form  can be identified with periods of the holomorphic volume form $\Omega^{3,0}_X$.

\begin{rmk} The reader remembers that when discussing singular Hitchin systems we fixed the essentially irregular part of the Higgs field. Coefficients of those fixed Puiseux series as well as the conformal structure on $(C,\{x_i\}_{1\le i\le n})$ can be thought of ``external'' parameters for the integrable systems in question. In terms of Calabi-Yau $3$-fold $X$ this means that we have extra parameters arising from the full moduli space of deformations of $X$.

\end{rmk}

\section{Wall-crossing structures for compact Calabi-Yau 3-folds and split attractor flow}

Let X be a compact complex Calabi-Yau 3-fold endowed with an ample line
bundle (polarization). We make a simplifying assumption that $H^1(X,\Q) = 0$
(otherwise the considerations below should be changed slightly).
The moduli stack $\MM :=\MM_X$ of complex structures on $X$ is a smooth
Deligne-Mumford stack (orbifold). The moduli stack $\LL := \LL_X$ of pairs $(X_{\tau},\Omega^{3,0}_{X_{\tau}})$
parametrizing pairs (complex structure $\tau$, holomorphic volume form) is a $\C^{\ast}$-bundle $p : \LL_X \to \MM_X$. In what follows we ignore those points of the stacks which have non-trivial stabilizers. Thus we will often abuse the terminology and speak about moduli spaces, not stacks.
Locally $\LL_X$ is embedded into $H^3(X,\C)$ via the period map $(X_{\tau},\Omega^{3,0}_{X_{\tau}})\mapsto [\Omega^{3,0}_{X_{\tau}}]$.

It is known (see [Don], [DonMar]) that $\LL_X$ is the base of non-polarized complex
integrable system with fibers given by intermediate Jacobians of the underlying
Calabi-Yau 3-folds. Since $X$ is compact the fibers in general are nonalgebraic.
Instead the fibers carry natural pseudo-K\"ahler metrics of signature $(1, n)$, where $n = {1\over{2}} rk H^3(X)-1$. This integrable system can be considered as a special case
of the one mentioned in Remark 4.1.2. The corresponding local system of lattices
has as fibers $\Gamma_b = H^3(X_{\tau} ,\Z), b\in \LL_X$, where 
$p(b) =\tau$. It also has the central
charge $Z_b(\gamma) =\int_{\gamma}\Omega^{3,0}_{X_{\tau}}$.

The only difference with the case of polarized integrable systems is that
now we do not have the positivity constraint on the skew-symmetric bilinear form
on $\Gamma_b$.

On the other hand in this case one has a well-defined local system of the
Fukaya categories ${\mathcal F}(X_{\tau} )$ over $\MM_X$. It is expected that a choice of point $b\in \LL_X$
defines a stability condition on ${\mathcal F}(X_{\tau}) , \tau = p(b)$ with the central charge $Z_ b:\Gamma_b\to \C$ as above and for which semistable objects are SLAGs endowed with local systems. Hence
we can speak (cf. Section 7.3) about DT-invariants $\Omega_b(\gamma), b\in \LL_X,\gamma\in \Gamma_b-\{0\}$.
We conclude that there is a corresponding WCS (see Section 2.3, Example
6). 

Moreover, using this WCS we can construct a non-archimedean symplectic
orbifold ${\mathcal X}$ along the lines of Section 4.6. More precisely, the construction of
Section 4.6 gives rise to a family ${\mathcal X}_{\zeta}, \zeta\in \C^{\ast}$ of such orbifolds, but all of them
are canonically isomorphic due to the natural $\C^{\ast}$-action on $\LL_X$.

We do not expect that $X$ is isomorphic to an open domain in an algebraic
orbifold. It is not even clear whether it is isomorphic to an open domain in
a complex analytic orbifold (the problem arises because of the expected overexponential
growth of  $\Omega_b(\gamma)$ as $|\gamma|\to \infty$. 
As a result, the second of the three
approaches to DT invariants discussed in the Introduction (namely the one with
the wheels of lines, see Section 6.6) cannot be applied.

Also the first named approach (via attractor flow and trees, see Section 3)
should be modified. More precisely, the initial WCS should in addition to what
was discussed in Section 4.6 (which are the values $1$ for DT-invariants at generic
conifold points) depend on infinitely many integer parameters, which are values
of DT-invariants at so-called attractor points (see e.g. [Den]). More precisely, in
the compact case the volume of $X$ is finite, so it can be used to normalize the
central charge. The normalized function (considered as a function on the total space of the local system $\GGamma\to \MM_X$) has countably many  minimal points. Their lifts to $\LL_X$ are called attractor points (see Sections 9.2, 9.3 about the details).

At this time we do not know how to extract those additional data directly
from the geometry of the above-described integrable system over $\LL_X$. In a
sense the additional data live ``at the infinity” of the moduli space $\LL_X$. Since
we are lacking the second approach to DT-invariants in the compact case, we
have modify the Conjecture 1.2.2 and only claim that the DT-invariants coming
from the Fukaya categories can be canonically reconstructed from the values of
the ``tropical” DT-invariants at the attractor points.

Also, the string theory suggests (see e.g. [AlManPerPi]) that there exists a complex
analytic contact orbifold ${\mathcal Y}$ with $dim\,{\mathcal Y} = dim\,{\mathcal X} + 1$ (which is called in
physics the twistor space for the quaternion-K\"ahler moduli space of hypermultipltes).
The orbifold 
${\mathcal X}$ is a formal germ of a “divisor at infinity” of ${\mathcal Y}$. The
structure of ${\mathcal Y}$ is still a mystery, but one can hope that 
${\mathcal Y}$ can be used for the
description of DT-invariants along the lines of the second approach (via the
wheels of lines).

\subsection{Split attractor flow and black holes: motivation from supergravity}
As a physics motivation for our considerations we will briefly explain the concept of the split attractor flow from the theory of supersymmetric black holes.

Let ${\cal M}_{CFT}$ be the ``moduli space" of  unitary $N=2$ superconformal field theories. It is believed that in case if there are no chiral fields of dimension $(2,0)$ then ${\cal M}_{CFT}\simeq {\cal M}_A\times {\cal M}_B$, where for CFTs associated with a $3CY$ manifold $X$ the moduli space  ${\cal M}_A$ is the space of complexified K\"ahler structures on $X$ while {\cal ${\cal M}_B$ is the moduli space of complex structures on $X$.

Recall critical superstring theory in $10$ dimensions    can degenerate to a family of  superconformal field theories with central charge $\widehat{c}=6$ over a $4$-dimensional flat space-time. The latter is $\R^4$ endowed with a singular metric satisfying Einstein equation with matter. The metric has singularities at black holes. Assuming time invariance we obtain a metric $g$ on $\R^3\setminus \{x_1,...,x_n\}$, where $x_i$ are positions of stationary black holes.
This family can be interpreted as a  map $h: \R^3\setminus \{x_1,...,x_n\}\to {\cal M}_{CFT}$ which satisfies together with the metric $g$ a complicated system of equations.

Let us assume that our CFT is of geometric origin and comes from a $3CY$ manifold $X$ such that $H^{1,0}(X)=H^{2,0}(X)=0$. Assume that the K\"ahler component of $h$ is constant. Then according to [Den], [DenGrRa] (see also [DenMo]) the set of pairs $(h,g)$ is in one-to-one correspondence with the set of   maps

$$\phi:\R^3\setminus \{x_1,...,x_n\}\to {\cal M}_X$$
(here ${\cal M}_X:=\MM_B(X)$ is the moduli space of complex structures on $X$)
coming from the following ansatz.
Namely, the map $\phi$ is obtained by the projectivization of the map
$\hat{\phi}:\R^3\setminus \{x_1,...,x_n\}\to \LL$, where $\R^3\setminus \{x_1,...,x_n\}$ is endowed with the standard flat Euclidean metric (which is different from the metric $g$) and $\LL:=\LL_X$ is the Lagrangian cone of the moduli space of deformations of $X$ endowed with a holomorphic volume form (it is locally embedded into $H^3(X,\C)$ via the period map). The cone $\LL$ is the total space of a $\C^{\ast}$-bundle over $\MM_X$. We endow $\LL$ with an integral affine structure via the local homeomorphism $Im:\LL\to H^3(X,\R), (\tau, \Omega^{3,0}_{\tau})\mapsto Im([\Omega^{3,0}_{\tau}])$, where $\tau\in {\cal M}_X$ is a complex structure on $X$, and $\Omega^{3,0}_{\tau}$ is the corresponding holomorphic volume form.
Then the ansatz comes from harmonic maps $\hat{\phi}$ which are locally of the form
$Im\circ \hat{\phi}(x)=\sum_{1\le i\le n}{\gamma_i\over{|x-x_i|}}+v_{\infty}$, where $\gamma_i, 1\le i\le n$ are elements of the charge lattice $\Gamma=H_3(X,\Z)\simeq H^3(X,\Z)$ (their meanings are the charges of black holes) and $v_{\infty}$ is the boundary condition ``at infinity" satisfying the constraint $\sum_{1\le i\le n}\langle \gamma_i,v_{\infty}\rangle =0$ (see [Den]). This gives us $\phi$.

The image of $\phi$ is an ``amoeba-shaped" $3$-dimensional domain in ${\cal M}_X$.  Hypothetically, connected components of the moduli space of maps $\phi$ with given $v_{\infty}, \gamma_i, 1\le i\le n$ have some cusps which are in one-to-one correspondence with split attractor trees (see [Den]). When we approach to such a cusp the $3d$ amoeba degenerates to a split attractor tree. This is somehow similar to the conventional ``tropical'' story, when holomorphic rational curves in Gromov-Witten theory degenerate at cusps to the gradient trees
on the base of SYZ torus fibration.
In fact edges of the split attractor tree are the gradient trajectories of the multivalued function $|F_{\gamma}|^2=|\int_{\gamma}\Omega^{(3,0)}|^2/|\int_X\Omega^{(3,0)}\wedge \overline{\Omega^{(3,0)}}|$ considered as a function on $\MM_X$.  Any edge is locally a projection of an affine line in $\LL$ with the slope $\gamma\in \Gamma$. 
If the split attractor flow (lifted from  ${\cal M}_X$ to $\LL$) starting at $v_{\infty}$ in the direction $\gamma$ hits the wall of marginal stability  where $\gamma=\gamma_1+\gamma_2+...+\gamma_k , Arg(\int_{\gamma_i}\Omega^{(3,0)})=Arg(\int_{\gamma}\Omega^{(3,0)}), 1\le i\le k$ then all $\gamma_1,...,\gamma_k$ belong (generically) to a two-dimensional plane.

We are  going to explain below that using our wall-crossing formulas it is possible to find all $\Omega(\gamma):=\Omega(b,\gamma), b\in \LL,\gamma\in \Gamma, \langle Im\,b,\gamma\rangle=0$ starting with a collection of integers $\Omega(b_{\gamma},\gamma)$ at the ``generalized attractor points" given by  conifold points and points $b_{\gamma}\in \Lambda$  defined by the equation $Im\,b_{\gamma}=\gamma$. The points $\C^{\ast}b_{\gamma}\in {\cal M}_X$ are external vertices of the split attractor trees. The wall-crossing formulas are used at the internal vertices of the trees for the computation of $\Omega(b,\gamma)$. The numbers $\Omega(b_{\gamma},\gamma)$ can be arbitrary.

\subsection{Affine structure on the Lagrangian cone}

Let us fix $\gamma\in \Gamma$. The wall of second kind associated with $\gamma$ (see [KoSo1]) is given explicitly by the set $$\LL_{\gamma}= \{(\tau,\Omega^{3,0}_{\tau})\in \LL|\langle Im([\Omega^{3,0}_{\tau}]),\gamma\rangle=0,\langle Re([\Omega^{3,0}_{\tau}]),\gamma\rangle>0\}.$$
(Notice that the condition
$\langle Re([\Omega^{3,0}_{\tau}]),Im([\Omega^{3,0}_{\tau}])\rangle>0$ holds on $\LL$). In what follows we will locally identify $\LL$ with the cone in $H^3(X,\C)$ and denote the point corresponding to $(\tau,\Omega^{3,0}_{\tau})$ simply by $\Omega^{3,0}$. We endow $\LL$ with an integer affine structure given locally $\Omega^{3,0}\mapsto Im(\Omega^{3,0})\in H^3(X,\R)$.

We define a multivalued function $F_{\gamma}:\MM_X\to \R_{\ge 0}$  by the formula:
$$F_{\gamma}(\Omega^{3,0})={|\langle \Omega^{3,0},\gamma\rangle|\over{\sqrt{\langle Re(\Omega^{3,0}),Im(\Omega^{3,0})\rangle}}}.$$
(The RHS does not depend on a choice of the lift  to $\LL$).

We define the {\it volume function} $v:\LL\to \R_{>0}$  from the equality
$$\langle Re(\Omega^{3,0}),Im(\Omega^{3,0})\rangle ={-1\over{2i}}\langle \Omega^{3,0},\overline{\Omega^{3,0}}\rangle=v(\Omega^{3,0})^2.$$

The moduli space ${\MM}_X$ carries the Weil-Petersson metric.  Let us recall its definition.
 Let us fix a form $\Omega^{3,0}_{X_0}$ such that $v(\Omega^{3,0}_{X_0})=1$.
The tangent space to $\LL$ at a point $\Omega^{3,0}_{X_0}$ can be identified with the term $F^2(H^3(X_0,\C))$ of the Hodge filtration, which  is decomposed into the direct sum $ H^{3,0}(X_0)\oplus H^{2,1}(X_0)$ by Hodge theory. Hence the tangent space to $\MM$ at the point $[X_0]$ is identified with
$H^{2,1}(X_0)$.  The latter space carries a natural Hermitean norm. This gives the metric on the tangent space $T_{[X_0]}\MM_X$.

\begin{thm}
Let us fix a non-zero $\gamma\in \Gamma$. 
One can  lift  the gradient flow of $|F_{\gamma}|^2$  to a flow 
on the wall $\LL_{\gamma}$ whose trajectories are straight lines with the slope $\gamma$ in the affine structure given by $Im(\Omega^{3,0})$. More precisely, the integral curve $\dot{x}=grad\,|F_{\gamma}|^2(x)$ near the point $x(0)=x_0\in \MM_X$ coincides as unparametrized curve with the image  of the  straight line $Im(\Omega_t^{3,0})=Im(\Omega_0^{3,0})+t\gamma$, where $\langle Im(\Omega_0^{3,0}),\gamma\rangle=0,
\langle Re(\Omega_0^{3,0}),\gamma\rangle>0$, and $\Omega_0^{3,0}$ belongs to a $\C^{\ast}$-fiber over $x_0$ of the bundle $\LL\to \MM_X$.
\end{thm}

{\it Proof.}  For any point $x_0:=[X_0]\in \MM$ such that $F_\gamma(x_0)\ne 0$
there is a unique   lift $\Omega_0^{3,0}\in \LL_{\gamma}$  with $v(\Omega_0^{3,0})=1$.
The tangent space $T_{x_0}\MM$ is identified with variations $\Omega_0^{3,0}\mapsto \Omega_0^{3,0}+\delta \Omega^{3,0}$
such that $\delta \Omega^{3,0}\in H^{2,1}$.

Let us  compute the variation of the function $log\,|F_{\gamma}(\Omega^{3,0})|^2=log|\langle \Omega^{3,0},\gamma \rangle|^2-log\,\langle Re(\Omega^{3,0}),Im(\Omega^{3,0})\rangle$. We obtain:
$$\delta log\,\langle \Omega^{3,0},\gamma\rangle+\delta log\,\langle \overline{\Omega^{3,0}},\gamma\rangle-\delta log\,\langle Re(\Omega^{3,0}),Im(\Omega^{3,0})\rangle=$$
$$2{Re(\langle \delta \Omega^{3,0},\gamma \rangle)\over{\langle \Omega_0^{3,0},\gamma \rangle }}-2i{Im(\langle \delta \Omega^{3,0},\overline{\Omega_0^{3,0} }\rangle)\over{\langle \Omega_0^{3,0},\overline{\Omega_0^{3,0}} \rangle }}.$$

The last term vanishes because $\delta \Omega^{3,0}$ belongs to $H^{2,1}$.
Since $\langle \Ome_0,\gamma \rangle= \langle Re(\Ome_0),\gamma \rangle>0$  and $\delta log|F_{\gamma}|^2=\delta|F_{\gamma}|^2/|F_{\gamma}|^2$
we have:
$$\delta |F_{\gamma}|^2=a Re(\langle \delta \Ome,\gamma \rangle)$$
 where $a\in \R$ is a real constant depending on $\Ome_0,\gamma$.

Next we observe that by the Hodge decomposition the element $\gamma\in \Gamma=H^3(X,\Z)\subset H^3(X,\R)$ can be written as $\gamma=\gamma^{3,0}+\overline{{\gamma}^{3,0}}+\gamma^{2,1}+\overline{{\gamma}^{2,1}}$, where  the upper index denotes the $(p,q)$-Hodge component. By orthogonality condition we conclude that
$\delta |F_{\gamma}|^2=a Re(\langle \delta \Ome,\gamma^{1,2} \rangle)=-a Im(\langle \delta \Ome,\overline{i \gamma^{2,1}}\rangle)$.

The RHS is by definition the pairing of two tangent vectors in $T_{x_0}\MM_X\simeq H^{2,1}$ with respect to the Weil-Petersson metric.
Hence  $grad\,|F_{\gamma}|$ at $x_0$ is proportional to $i \gamma^{2,1}$. But $i\gamma^{2,1}$ is the projection of the tangent vector
$i(\gamma^{3,0}+\gamma^{2,1})$ at $\Omega_0^{3,0}\in \LL$ whose imaginary part is
$\gamma/2= Im( i\gamma^{3,0}+i\gamma^{2,1})$ . This concludes the proof. 
$\blacksquare$


In what follows we will need to know the behavior of the volume function $v(\Ome)$ along the gradient trajectory of the function $|F_{\gamma}|^2$. We choose a parametrization of the trajectory such that $Im(\dot{\Omega}^{3,0})=\gamma$ or equivalently that $\dot{\Omega}^{3,0}=2i(\gamma^{3,0}+\gamma^{2,1})$.

Along the gradient trajectory we have:
$${d\over {dt}}(log\,v)={1\over{2}}{d\over {dt}}(log\,v^2)=Re\left({\langle \Ome,\overline{\dot{\Omega}^{3,0}}\rangle}\over {\langle \Ome,\overline{\Ome}\rangle}\right)=$$
$$Re\left(\langle \Ome, {\overline{\dot{\Omega}^{3,0}}}-\dot{\Omega}^{3,0}\rangle)\over {\langle \Ome,\overline{\Ome}\rangle}\right)=Re\left(\langle \Ome, -2iIm({\dot{\Omega}^{3,0}})\rangle\over {-2iv^2}\right)=$$
$${Re(\langle \Ome,\gamma\rangle)\over {v^2}}={\langle \Ome,\gamma\rangle\over {v^2}}.$$
(In the course of the computations we use the Lagrangian property of $\LL$ which gives  $\langle \Ome,\dot{\Omega}^{3,0}\rangle=0$, as well as the equality $Im(\dot{\Omega}^{3,0})=\gamma$).

Therefore ${d\over {dt}}v={\langle \Ome,\gamma\rangle\over {v}}=\pm F_{\gamma}$.
Notice that if $\langle \Ome,\gamma\rangle\ne 0$ then similarly to the proof of Theorem 9.2.1 we have
$${d\over {dt}}(|F_{\gamma}|^2)(\Ome)=|F_{\gamma}|^2(\Ome)\left(2{Re(\langle \dot{\Omega}^{3,0},\gamma \rangle)\over{\langle \Ome,\gamma \rangle }}-2{Re(\langle \dot{ \Omega}^{3,0},i\overline{\Ome }\rangle)\over{\langle \Ome,i\overline{\Ome} \rangle }}\right).$$

Now we can use the formula $\dot{\Omega}^{3,0}=2i(\gamma^{3,0}+\gamma^{2,1})$.

Then $2i\gamma^{3,0}=c\Ome$, where $c\in \C$. Hence the input of the summand $2i\gamma^{3,0}$ to the RHS of the above formula is
$$|F_{\gamma}|^2(\Ome)\left(2{Re(c\langle {\Ome},\gamma \rangle)\over{\langle \Ome,\gamma \rangle }}-2{Re(c\langle{ \Ome},i\overline{\Ome }\rangle)\over{\langle \Ome,i\overline{\Ome} \rangle }}\right)=0.$$
(Notice that $\langle {\Ome},\gamma \rangle>0$ and $\langle{ \Ome},i\overline{\Ome }\rangle>0$ by our assumptions, hence the expression in the big brackets simplifies to $2Re(c)-2Re(c)=0$).

Therefore in the RHS we have the contribution of the summand $2i\gamma^{2,1}$ only.
This gives us
$${d\over {dt}}(|F_{\gamma}|^2)=|F_{\gamma}|^2 2{Re(\langle 2i\gamma^{2,1},\overline{\gamma^{2,1}\rangle)}\over{\langle \Ome,\gamma \rangle }}=8|F_{\gamma}|^2
{\langle Re(\gamma^{2,1}),Im({\gamma^{2,1})\rangle}\over{\langle \Ome,\gamma \rangle }}.$$

Using the formula ${d\over {dt}}v={\langle \Ome,\gamma\rangle\over {v}}$ we conclude that
$${d^2v\over {dt}^2}={4\over v}\langle Re(\gamma^{2,1}),Im(\gamma^{2,1})\rangle.$$
Notice that properties of polarized Hodge structures imply that
$$\langle Re({\gamma}^{2,1}),Im({\gamma}^{2,1})\rangle\le 0\,.$$

Recall that the gradient lines are projections of  unparametrized straight lines (see Theorem 9.2.1). Then our computations imply  the following statement.

\begin{prp} The volume function $v$ is concave in parameter $t$
on the 
straight line  $Im(\Omega_t^{3,0})=Im(\Omega_0^{3,0})+t\gamma$, where
$\langle Im\,\Omega^{3,0}_0,\gamma\rangle=0$. 

\end{prp}

\begin{rmk} Notice that in the above computations we could replace $\gamma$ by any vector in $\Gamma_{\R}$. Recall also that we can identify locally $\LL$ with the real symplectic vector space $H^3(X,\R)$ via the map $[\Omega^{3,0}]\mapsto Im\,[\Omega^{3,0}]$ which introduces the affine structure on $\LL$. The real $2$-dimensional subspace spanned by $Im([\Omega_0^{3,0}])$ and $\gamma$ is isotropic. 
Then the Proposition 9.2.2 can be reformulated in geometric terms as a statement that
the  function $v$ gives rise to a (positive) concave homogeneous function of degree $+1$ on any real isotropic plane in $H^3(X,\R)$ . 

\end{rmk}

The above considerations motivate the following definitions.

\begin{defn} Let $\gamma \in \Gamma=H^3(X,\Z)$. We call $\Ome_{\gamma}\in \LL$ a $\gamma$-attractor point if $Im(\Ome_{\gamma})= \gamma$.
\end{defn}

Since a lift of the gradient trajectory of $|F_{\gamma}|^2$ is a straight line in $\LL$, the critical points of $F_{\gamma}$ can appear only in the limit $t\to \infty$. Hence the limiting point in $\MM_X$ is the projection of a $\gamma$-attractor point.
Thus we see that the projections of attractor points  are local minima of the multivalued functions $F_{\gamma}$ (in physics literature these projections are called attractor points).

Moreover it is easy to see that critical points of $F_{\gamma}$ are all local minima and are either projection of attractor points or belong to the locus $F_{\gamma}^{-1}(0)$.
The equation $F_{\gamma}=0$ defines a complex hypersurface in the (universal cover) of the space $\MM_X$. Points of this hypersurface are absolute minima of $F_{\gamma}$, and moreover they form a set of points where the function $F_{\gamma}$ is not differentiable. 

We expect that for a generic gradient line of $F_{\gamma}$ on  $\MM_X$ there are three possibilities:

1) The gradient line hits the projection of a $\gamma$-attractor point.

2) The gradient line reaches in finite time a point in the boundary of the metric completion of  $\MM_X$ with respect to the Weil-Petersson metric. This point is called {\it conifold point}. We will assume that the conifold   points form an analytic  divisor in the above completion (which is expected to be a complex analytic space). 

3) The gradient line reaches in finite time a point in the locus $F_{\gamma}^{-1}(0)$.

Then for given $\gamma$ the universal covering of $\MM_X$ splits into a disjoint union of three open domains corresponding to these three possibilities and a closed subset of measure zero.

The divisor of conifold points has (roughly) the following structure which we explain in the framework of complex integrable systems. Consider a polarized integrable system with central charge endowed with conical structure, i.e. with a $\C^{\ast}$-action which rescales the central charge and preserves the discriminant. Taking the quotient  by the $\C^{\ast}$-action we obtain a local model for the divisor of conifold points in the completion of $\MM_X$.  In particular we expect the typical singularity will be of the type $A_1$. In terms of $X$ this means that we approach a point in the completion of $\MM_X$ where $X$ develops an ordinary  double point with local equation $\sum_{1\le i\le 4}x_i^2=0$.

We can restate these three possibilities in the language of straight lines on $\LL_X$. Let us fix $\gamma$ and consider in the universal cover of $\LL_X$ the  ray (or interval) $Im(\Omega^{3,0}_t)=Im(\Omega_0^{3,0})+t\gamma, t\in [0,t_0)$, where $\langle Im(\Omega_0^{3,0}),\gamma\rangle=0$, and $t_0\in (0,+\infty]$ is the maximal possible value of $t$ for which the map $t\mapsto \Omega^{3,0}_t$ is well-defined.

The case 1) means that $t_0=+\infty$. The restriction of the function $v$ on the ray has strictly positive derivative and the limit of the derivative as $t\to +\infty$ is non-zero: $lim_{t\to +\infty}{dv\over{dt}}>0$.

In the case 2) we have $t_0<+\infty$ and $lim_{t\to t_0}{dv\over{dt}}=0$.

We claim that in the case 3) that $t_0<+\infty$ but the limit of the derivative of the function $v$ is strictly negative as $t\to t_0$. Indeed, the picture in $\MM_X$ means that there exists finite $t_1$ such that the derivative of $v$ at $t_1$ is equal to zero. It is easy to see that although the gradient line of $F_{\gamma}$ stops at such  point, we can continue the corresponding ray in $\LL$ to some $t>t_1$. In terms of the gradient trajectories this means that we consider another gradient trajectory of the function $F_{\gamma}$  and move along it in the opposite direction (i.e. in the direction of increasing values of $F_{\gamma}$) for $t>t_1$. Therefore for $t>t_1$ the derivative ${dv\over{dt}}$ becomes negative. By concavity of $v$  we conclude that we cannot extend the ray indefinitely, hence $t_0<+\infty$. We expect that in this case the image of the corresponding interval in $\LL_X/\R_{>0}$ is everywhere dense. This property distinguishes the case 3) from the case 2) purely in terms of affine geometry of $\LL_X$ (without use of function $v$).

\subsection{Trees and generalized attractor points}

Let us discuss  abstract attractor trees. Basically, it is the same as the tropical trees discussed previously. The difference is in the notion of attractor point.

Suppose that $\LL$ is a smooth $C^{\infty}$-manifold which admits an open covering $\LL=\cup_{i\in I}U_i$ with transition functions belonging to the group $Aut(\Gamma,\langle\bullet,\bullet\rangle)$, where $\Gamma\simeq \Z^{2n}$ is a lattice endowed with the integer non-degenerate skew-symmetric form $\langle\bullet,\bullet\rangle$. We assume that  each $U_i$ endowed with the induced $\Z$-linear structure is isomorphic to an open  cone in $\R^{2n}$.
Notice that because we have a $\Z$-linear structure on $\LL$ we can speak about integer points in $\LL$. Also the conical structure implies that $\R_{>0}$ acts on $\LL$.

Let us define a $(2n-1)$-dimensional manifold $\LL^{\prime}_{\Z}$ as the set of pairs $(u,\gamma)$, where $u\in \LL, \gamma\in T_u\LL-\{0\}$ such that in the above $\Z$-linear local coordinates $\gamma$ is integer and $\langle u,\gamma\rangle=0$. This definition is similar to the definition of $M^{\prime}_{\Z}$ from Section 3.1. We define the {\it attractor flow} on $\LL^{\prime}_{\Z}$ by the formula $\dot{u}=\gamma, \dot{\gamma}=0$.

\begin{defn} We define $\LL^{\prime,attr}_{\Z}\subset \LL^{\prime}_{\Z}$ as a set of points $(u,\gamma)$ such that the trajectory of the attractor flow  starting at $(u,\gamma)$ exists for all $t\in [0,+\infty)$.

\end{defn}
In local coordinates the projection to $\LL$ of such a trajectory can be written as $t\mapsto u_t:=u+t\gamma$.

\begin{defn} A generalized attractor point is a connected component of the interior $Int(\LL^{\prime,attr}_{\Z}) \subset \LL^{\prime}_{\Z}$.
\end{defn}

Below we provide some explanations.

With any $(u,\gamma)\in \LL^{\prime,attr}_{\Z}$ we associate a map $f_{(u,\gamma)}: (0,+\infty)\to \LL$ given by $t\mapsto t^{-1}u_t$. 

There are two possibilities:

a) The limit $lim_{t\to +\infty}f_{(u,\gamma)}(t)$ does exist. Then in local coordinates this limit is equal to $\gamma$. This condition is open in $\LL^{\prime}_{\Z}$. Generalized attractor points of this type can be identified with  $\gamma$-attractor points from the previous subsection.

b) The limit $lim_{t\to +\infty}f_{(u,\gamma)}(t)$ does not exist.

In the case a) the limit is an integer point in $\LL$.  It is easy to see that for any integer point  $u_0\in \LL$ the set of pairs 
$(u,\gamma)\in  \LL^{\prime,attr}_{\Z}$ such that $lim_{t\to +\infty}f_{(u,\gamma)}(t)=u_0$ is a non-empty open connected subset in  $\LL^{\prime}_{\Z}$ (it is a star-shaped domain). Hence we conclude that integer points in $\LL$ give generalized attractor points. The inclusion is not a bijection. The complement to the image corresponds to the interior of the domain described in case b). The latter can be thought of as the set of integer points in the ``boundary" of $\LL$.
Such ``integer boundary points" do appear in practice. For example take $\LL=\LL_X$ and points $(u,\gamma)$ where $u$ corresponds to the point in the moduli space $\MM_X$ close to the cusp and $\gamma$ belongs to a Lagrangian sublattice in $H_3(X,\Z)$ invariant under the monodromy. In the mirror dual picture such classes $\gamma$ correspond to Chern classes of coherent sheaves on the dual Calabi-Yau  with at most $1$-dimensional support ($D0-D2$ branes in the language of physics).

Assume that we are given an open subset $\LL^{\prime,conif}_{\Z}\subset 
\LL^{\prime}_{\Z}$ which is preserved by the attractor flow for $t\ge 0$ and is disjoint from  $\LL^{\prime,attr}_{\Z}$. 
For example in the situation when $\LL=\LL_X$ described in the previous subsection we define $\LL^{\prime,conif}_{\Z}$ as the interior of the set of points described in the case 2) there. As we mentioned in Section 9.2 this probably means that the projection of the corresponding trajectory of the attractor flow to $\LL/\R_{>0}$ is not everywhere dense.


We will be talking about metrized rooted trees below. As in [KoSo4] those are trees with lengths assigned to edges. There are internal edges and tail edges. Internal edges have finite (positive) length and tail edges have possibly infinite length.
Also, the root vertex is adjacent to exactly one edge.

\begin{defn} A tropical tree in $\LL$ is given by

a) A metrized rooted  tree $T$ with edges oriented toward  tails.

b) A continuous map $\phi:T\to \LL$, smooth outside of vertices, with the following properties:

outside of vertices the map $t\mapsto (\phi(t),\phi^{\prime}(t))$ is a trajectory of the attractor flow on  $\LL^{\prime}_{\Z}$;

at each internal vertex $v$ the following balancing condition is satisfied:

$$\sum_{e\in {v^{out}}}\phi^{\prime}(e)=\phi^{\prime}(e^{in}(v)),$$
where $v^{out}$  denote the set of outcoming  from $v$ edges and $e^{in}(v)$ is the only edge incoming to $v$
(the derivative $\phi^{\prime}(e)$ is constant along the edge $e$);

each tail edge the map $\phi$ is a trajectory of the attractor flow belonging either to  $\LL^{\prime,conif}_{\Z}$ or to $\LL^{\prime,attr}_{\Z}$;

for any vertex $v$ directions $\phi^{\prime}$ of all edges in the set $v^{out}$ are different and belong to an open half-space in a rank $2$ symplectic sublattice  in the tangent space $T_{\phi(v)}\LL$;

\end{defn}
Having a tropical tree we can (and sometime will) interpret its edges as trajectories in $\LL^{\prime}_{\Z}$.

Let $\phi: T\to \LL$ be a rooted tropical tree in $\LL$. Abusing the notation we will simply denote it by $T$.
Suppose we are given a volume function $v:\LL\to \R_{>0}$ which satisfies the  Proposition 9.2.2 (i.e. it is concave along  edges). Then the following result holds.

\begin{prp} If the volume function increases along  tail edges of $T$ then it increases along every edge (we consider orientation of the tree toward tails).

\end{prp}

Let us define the function $F: \LL^{\prime}_{\Z}\to \R$ as the derivative along the attractor flow of the pull-back of $v$ under the natural projection. This function is an analog of the multivalued function $F_{\gamma}$ from Section 9.2. The function $F$ is invariant under the $\R_{>0}$-action on $\LL^{\prime}_{\Z}$. Proposition 9.3.4 means that edges of any tropical tree  belong to the domain $F>0$. Moreover, the concavity of the function $v$ on edges and the balancing condition imply that the value of $F$ at the root vertex is strictly bigger than the sum of limiting values of $F$ on tail edges (cf. [DenMo]). Such limiting values are strictly positive for edges of the tree hitting attractor points and equal to zero for conifold and ``integer boundary  points".

Thus we see that the function $F$ imposes the ``energy-like" restrictions on tropical trees.

Our considerations with the function $F$ motivates the following {\it Finiteness Assumption} (cf. assumptions in Section 3.2):

{\it For each point $(u,\gamma)\in \LL^{\prime}_{\Z}$ outside of a set of measure zero the number of tropical trees rooted at $(u,\gamma)$ is finite.}

By analogy with Proposition 3.2.6 one can design a procedure which as we expect produces the WCS on $\LL$ starting with ``initial data" given by integer numbers $\Omega(u,\gamma)$ assigned to generalized attractor points and irreducible components of the divisor of conifold points. We assign arbitrary integers to generalized attractor points and assign integers equal to $1$ to conifold points with $A_1$ singularities. Similarly to Section 3.2 this WCS can be understood as an integer-valued function on $\LL^{\prime}_{\Z}$ which discontinuous at polyhedral walls and satisfies the Support Property and WCF from [KoSo1].


The procedure is similar to those described in Sections 3.2 and 4.6. 

By the Finiteness Assumption there are finitely many tropical trees in $\LL$ rooted at $(u,\gamma)$. The union of all such trees is a finite directed graph $G:=G(u,\gamma)$ without oriented cycles because of monotonicity of the  function $F$. Each inner vertex of $G$ belongs to a symplectic plane. Then we can start with a generalized attractor point or conifold point with $A_1$ singularity which are tails of $T$ and move backward to $(u,\gamma)$. We  will use the wall-crossing formulas.  Recall that they have the following form
$$\prod^{\longleftarrow}T_{\gamma^{out}}^{\Omega(P,\gamma_{out})}=\prod^{\longrightarrow}T_{\gamma^{in}}^{\Omega(P,\gamma_{in})},$$
where $P$ is a vertex of $G$ and $T_{\gamma^{out}}^{\Omega(P, \gamma_{out})}$ are symplectomorphisms of the $2$-dimensional symplectic subspace in the tangent space at $P$ corresponding to the edges of $G$ outcoming from $P$ (similarly for incoming edges).
Since we know by induction the numbers $\Omega(P,\gamma_{out})$ for outcoming edges, we can calculate
$\Omega(P,\gamma_{in})$ from the wall-crossing formula and proceed further toward $P$. Finally, it gives us the desired number $\Omega(u,\gamma)$. 

\begin{rmk} As we already mentioned in Remark 3.2.7, there is no guarantee that the result of the application of the above procedure is indeed a WCS.
 At some strata of codimension 2 the cocycle condition can fail. The geometric structure of walls on $\LL$ is very involved, and we do not understand it completely.
At the moment we have the following (maybe too optimistic) picture: it is sufficient (and maybe even necessary) to put the constraint $\Omega=1$ at conifold points (assuming that all conifold points have $A_1$ singularities).
The integer values of $\Omega$ at all generalized attractor points can be chosen arbitrarily. Then we obtain a WCS.

\end{rmk}

\subsection{Remarks on the support of DT-invariants}

Recall that in the definition of WCS we required an existence of a strict convex cone. In the case of Calabi-Yau $3$-folds this property is called Support Property (see [KoSo1]), since it gives a bound on the support of the function $\Omega(u,\gamma)$ (numerical DT-invariants). Heuristic arguments in favor of that given in the Remark 1 [KoSo1] were based on the following simple geometric fact. Let $\eta$ be a closed $3$-form on $X$. Then there exist $C:=C_{\eta}>0$ such that for any SLAG $L$ we have $|\int_L\eta|\le C|\int_L\Omega^{3,0}_X|$. Equivalently, we have a constraint on the homology class $\gamma=[L]\in H_3(X,\Z)$. The constant $C$ depends in an essential way on the metric on $X$. In this subsection we propose an alternative approach to the Support Property based entirely on the affine geometry of $\LL=\LL_X$. Namely, the volume function $v$ gives us the following constraint on the pair $(u,\gamma)$ such that $\Omega(u,\gamma)\ne 0$: the derivative of $v$ at $u$ in the direction $\gamma$ is positive. Recall that this property was deduced from three facts:

i) the function $v$ is a strictly positive function on $\LL$ of homogeneity degree $1$ with respect to the $\R_{>0}$ action;

ii) the function $v$ is concave on germs of $2$-dimensional isotropic subspaces in $\LL$;

iii) the derivative of $v$ along a trajectory in $\LL^{\prime,conif}_{\Z}$ is positive.

We claim that there are infinitely many perturbations of $v$ which still obey i)-iii). Namely, let us consider any smooth function $\delta v$ on $\LL$ which is homogeneous of degree $1$ and such that $Supp(\delta v)/\R_{>0}\subset \LL/\R_{>0}$ is {\it compact}. Consider the function $v_{\varepsilon}:=v+\varepsilon \delta v$. The compactness of $Supp(\delta v)/\R_{>0}$ implies that the properties i) and ii) are satisfied for sufficiently small $\varepsilon$. The property iii) follows from ii). More generally, any function $v^{\prime}$ which satisfies i) and ii) and coincides with $v$ outside of a compact (modulo the action of $\R_{>0}$) satisfies also iii). Here is a reason for that: the condition of monotonicity of such a function is sufficient to check on parts of the trajectories which are close to conifold points, where the function coincides with $v$. For any function $v^{\prime}$ satisfying i)-iii) let us consider the set $C_{v^{\prime}}\subset tot(T\LL)$ consisting of pairs $(u,\dot{u})$ such that $\langle u,\dot{u}\rangle=0$ and $dv_{|T_u\LL}(\dot{u})\ge 0$. Let us define $C_{univ}$ as the intersection of $C_{v^{\prime}}$ over all $v^{\prime}$ satisfying i)-iii). It is easy to see that for any $u\in \LL$ the intersection $C_{univ}\cap T_u\LL$ is a {\it strict closed convex cone} in the hyperplane $u^{\perp}:=Ker\,\langle u,\bullet\rangle$. The above inductive construction implies that $Supp(\Omega(u,\gamma))\subset C_{univ}$.

\section{Analog of WCS in Mirror symmetry}

\subsection{Pair of lattices, volume preserving transformations and WCS in a vector space}

In the case of SYZ picture of Mirror Symmetry the construction of mirror dual involves transformations which 
locally preserve the volume form rather than a Poisson structure. In this case $\g$ is the Lie algebra of divergence-free vector fields on $Hom(\Gamma, \C^{\ast})$ and there is no distinguished skew-symmetric form on $\Gamma$. In notation of section 6.2 the lattice 
$\Gamma$ is $\Gamma_b$, the first homology group of a fiber of a {\it real } integrable system at a given point $b\in B^0$. Differently from the  case of symplectomorphisms when the dimension of the graded component is equal to $1$ (see Section 2.3, Example 4), we now have 
 $dim\,\g_{\gamma}=n-1$, where $n=rk\,\Gamma$ for $\gamma\ne 0$. Explicitly, the Lie algebra of vector fields on the algebraic torus
 $Hom(\Gamma, \C^{\ast})$
  is spanned by elements $x^\gamma\partial_\mu$ where $\gamma \in \Gamma,\mu \in \Gamma^\vee$ satisfying the linear relations
   $$x^\gamma\partial_{\mu_1}+x^\gamma\partial_{\mu_2}=x^\gamma\partial_{\mu_1+\mu_2}\,\,.$$
 Derivation $\partial_\mu$ is a constant vector field in logarithmic coordinates.
 The commutator rule is given by
 $$[x^{\gamma_1}\partial_{\mu_1},x^{\gamma_2}\partial_{\mu_2}]=x^{\gamma_1+\gamma_2}\bigl( (\mu_1,\gamma_2)\partial_{\mu_2}
 - (\mu_2,\gamma_1)\partial_{\mu_1} \bigr)\,\,.$$
 The subalgebra $\g$ of divergence-free vector fields is spanned by elements $x^\gamma\partial_\mu$ with $(\mu,\gamma)=0$.
  It is obviously graded by the lattice $\Gamma$. Similarly to the symplectic (and also Poisson) case (see the beginning of Section 3.1),  
  the 
   graded complement to $\g_0$ is a Lie subalgebra
   $\g'=\oplus_{\gamma\ne 0}\g_\gamma$ in $\g$
   (notice that an analogous property  {\it does not } hold for the Lie algebra of all vector fields).
   
   One can generalize the above considerations to the following situation. Suppose 
   we are given two lattices $\Gamma_1,\Gamma_2$ and an integer  pairing between them $(\bullet,\bullet): \Gamma_2\otimes \Gamma_1\to \Z$. We do not assume that the pairing is non-degenerate. We denote by $\Gamma_{1,0}\subset \Gamma_1$ and $\Gamma_{2,0}\subset \Gamma_2$ the corresponding kernels of the pairing.

  Then we consider the Lie algebra $\g:=\g_{\Gamma_1,\Gamma_2,(\bullet,\bullet)}$
   spanned by elements $x^\gamma\partial_\mu$ where $\gamma\in \Gamma_1$ and $\mu \in \Gamma_2$ such that $(\mu,\gamma)=0$,
   satisfying the same relations as above. It contains the Lie subalgebra $$\g^{\prime}:=\oplus_{\gamma\in \Gamma_1-\Gamma_{1,0}}\g_{\gamma}.$$ The previous special case corresponds to  
   $\Gamma_1=\Gamma,\,\Gamma_2=\Gamma^\vee$.
   In general,  $\g$ can be thought as the Lie algebra of divergence-free vector fields on a torus, preserving a collection of coordinates 
    and commuting with a subtorus action. Explicitly, if we omit the condition $(\mu,\gamma)=0$ of being divergence-free,
     in some coordinates $(x_1,\dots,x_{a+b+c})$ for $a,b,c\in\Z_{\ge 0}$ we get vector fields of the form
     $$\prod_{i=1}^{a+b}x_i^{k_i}\cdot x_j\partial/\partial x_j,\,\,\,(k_i)_{1\le i\le a+b}\in\Z^{a+b},\,a+1\le j\le a+b+c\,.$$

From the point of view of SYZ picture of Mirror Symmetry we have a real integrable system $X\to B$ with the dimension of the total space $X$ equal to $2b$. Moreover we have chosen a Lagrangian zero section as well as $c$ other Lagrangian sections (more precisely, a homomorphism from $\Z^c$ to the abelian group of Lagrangian sections). Also we assume that there is an $a$-dimensional space of deformations of the above structure which is an $a$-dimensional vector subspace  in $H^2(X,\R)$ defined over $\Q$.
The mirror dual $X^{\vee}$ is the complex manifold of complex dimension $b$ depending on $a$ holomorphic parameters and carrying $c$ line bundles. The Mirror Symmetry preserves the parameter $b$ and exchanges $a$ and $c$.

  Now we are ready to describe the analog of a WCS for Lie algebra $\g^{\prime}$.
      
      The main difference with the formalism from 2.1 is that now walls are hyperplanes in $\Gamma_{2,\R}^*:=\Gamma_2^\vee\otimes\R$ (and not in the dual space to the grading lattice $\Gamma_1$). We define a  wall as a hyperplane in $\Gamma_{2,\R}^*$ given by $\mu^\perp$, where $\mu\in \Gamma_2-\Gamma_{2,0}$ (one may assume that $\mu$ is primitive).
      With any wall $H\subset\Gamma_{2,\R}^*$ we associate a graded Lie subalgebra 
$$\g_H:=\bigoplus_{\gamma\in \Gamma_1}\g_{H,\gamma}\subset \g^{\prime}$$
 spanned by $x^\gamma\partial_\mu$
 such that $(\mu,\gamma)=0$ and $\gamma\in\Gamma_1-\Gamma_{1,0}$. As in the Poisson case, this Lie algebra is abelian. It is convenient to associate with any $\gamma$ as above a nonzero constant vector field on the hyperplane $H$ equal to $\iota( \gamma ):=(\bullet,\gamma)\in\Gamma_2^\vee\subset\Gamma_{2,\R}^*$. In SYZ picture the trajectories of this vector field are (possible) parts of tropical trees corresponding
to analytic discs with the boundary on a small Lagrangian torus, the fiber of SYZ fibration.

 Also with any $\Q$-vector subspace $V \subset \Gamma_{2,\R}^*$ which is the intersection of two walls we associate a graded Lie algebra
 $\g_V$ 
(which is {\it not} a subalgebra of $\g^{\prime}$) such as follows. As a $\Gamma_1$-graded vector space $\g_V$ will be equal to the direct sum $\oplus_{H\supset V}\g_H$ over all walls containing $V$. The Lie bracket on $\g_V$ is defined as follows.
Let  $(x^{\gamma}\partial_{\mu})_{H}$ where $\gamma\in \Gamma_1-\Gamma_{1,0}, \mu\in \Gamma_2-\Gamma_{2,0}$ denotes the element
$x^{\gamma}\partial_{\mu}\in\g_H$ considered as an element of $\g_{H}\subset\g_V$, where    $H=\mu^{\perp}$ is a wall containing $V$. 
Then we define the Lie bracket by the formula:
$$[(x^{\gamma_1}\partial_{\mu_1})_{H_1},(x^{\gamma_2}\partial_{\mu_2})_{H_2}]=
(x^{\gamma_3}\partial_{\mu_3})_{H_3},$$
in case if $H_{i}=\mu_i^{\perp}, i=1,2,3$, $\gamma_3=\gamma_1+\gamma_2$, $\mu_3=(\mu_1,\gamma_2){\mu_2}-(\mu_2,\gamma_1){\mu_1}$ and $\mu_3\notin \Gamma_{2,0}$. Otherwise, i.e. if $\mu_3\in \Gamma_{2,0}$ (and as one can easily see $\mu_3=0$), we define the commutator to be equal to zero.

  As in Section 2.1.3, we consider the pronilpotent case by choosing a strict convex cone $C\subset \Gamma_1\otimes\R$, and working with 
      $\g_C:=\prod_{\gamma\in \Gamma\cap C-\Gamma_{1,0}}\g_\gamma$.       
 Then for a given functional $\phi:\Gamma_1\to \Z$ which is nonnegative and proper on  the closure of $C$, we consider finite-dimensional nilpotent quotients 
      $$\g_{C,\phi}^{(k)}=\oplus_{\gamma \in \Gamma_1-\Gamma_{1,0}|\phi(\gamma)\le k}\g_{C,\gamma}=\g_C/m_{C,\phi}^{(k)}$$
      where $m_{C,\phi}^{(k)}:=\prod_{\gamma \in \Gamma_1-\Gamma_{1,0}:\phi(\gamma)>k}\g_{C,\gamma}$ is an ideal in $\g_C$.
 
Similarly we define the Lie algebras  $\g_{H,C,\phi}^{(k)}$.
 and 
$\g_{V, C,\phi}^{(k)}$.

 Let us fix finitely many walls $H_i, i\in I$. We define the set $WCS_k(\{H_i\}_{i\in I},C,\phi)$ of wall-crossing structures for $\g_{C,\phi}^{(k)}$ which are supported on the union $\cup_{i\in I}H_i$ in the following way. 
 First we observe that the walls $H_i,i\in I$ give rise to the natural stratification of $\Gamma_{2,\R}^*$. 
 Then an element  $WCS_k(\{H_i\}_{i\in I}, C,\phi)$ is a map which associates  an element $g_\tau$ of the group 
 $\exp(\g_{H_i,C,\phi}^{(k)}), i\in I$, where $\tau \subset H_i$  is a  co-oriented stratum 
of codimension one in $\Gamma_{2,\R}^*$ (notice that $\tau$ is an open subset of $H_i$).


The only condition on this map says that  for any generic closed loop $f:\R/\Z\to \Gamma_{2,\R}^*$ surrounding a codimension two stratum $\rho\subset V, codim_{\R}V=2$, the product of images  of the corresponding  elements $exp(g_{\tau_{t_i}})$ in 
$exp(\g_{V,C,\phi}^{(k)})$ over the finite sequence of intersection points $f(t_i)$ of the loop with  walls $H_i$ is equal to the identity. 

Now we take the inductive limit of the sets $WCS_k(\{H_i\}_{i\in I}, C,\phi) $ over all finite collections $\{H_i\}_{i\in I}$ of walls and after that we take the projective limit over $k$. 
The resulting set $WCS_{\g,C}$ is our analog of WCS relevant to the Mirror Symmetry. Analogously to Section 2 we can generalize our considerations and define $WCS_{\g,C}$ as a sheaf of sets on $\Gamma_{2,\R}^{\ast}$ and generalize even further assuming that $\Gamma_1$ and $\Gamma_2$ are local systems of lattices on a topological space, say, $M$. Instead of the central charge we now have a morphism of sheaves of abelian groups $\Gamma_2\to \underline{Cont}_M$. It is not clear a priori why the sheaf  $WCS_{\g,C}$ is non-trivial, and how to introduce ``coordinates" on its stalk at zero. Although the structure we have defined is different from $WCS$ discussed in Section 2, we will abuse the language and still call it the wall-crossing structure.

Next we would like to discuss an analog of the initial data.

\begin{prp} Let $V=H_1\cap H_2, H_i=\mu_i^{\perp},i=1,2$ be an intersection of two walls and $\gamma\in\Gamma_1-\Gamma_{1,0}$ be a vector
such that $(\mu_i,\gamma)=0,i=1,2$.
Then the natural map  $$\bigoplus_{H:V\subset H}\g_{H,\gamma}\to \bigl(\g_V/[\g_V,\g_V]\bigr)_\gamma$$ is injective.

\end{prp}

{\it Proof.} Follows immediately from the  formula for the bracket. $\blacksquare$

Then we can define the analog of the initial data in the following two ways depending on a choice of a sign.
For any $\gamma\in \Gamma_1-\Gamma_{1,0}$ and any $H=\mu^{\perp}$ such that $ (\mu,\gamma)=0$ we will construct  elements $a_H^{(k),\pm}(\gamma)\in \g_{H,C,\phi,\gamma}^{(k)}$.  Namely, it follows from the above Proposition that for any stratum $\tau\subset H$ as above  such $\iota(\gamma)$ belongs to the closure $\overline{\tau}$ of the stratum $\tau$, the $\gamma$-component of $log(g_{\tau})$ does not depend on $\tau$. We denote it by $a_H^{(k),+}(\gamma)$. Similarly we define $a_H^{(k),-}(\gamma)$ using the strata $\tau$ such that $-\gamma\in \overline{\tau}$. Next we define elements $a_H^{(k),\pm}\in \g_{H,C,\phi}^{(k)}$ as $\sum_{\gamma}a_H^{(k),\pm}(\gamma)$. After that we define elements $a^{(k),\pm}\in \oplus_H \g_{H,C,\phi}^{(k)}$ where the sum is taken over the set of walls (notice that the sum is finite).

\begin{conj} The set $WCS_k(C,\phi)$ is identified via passing to initial data  $a^{(k),+}$    with the set $\oplus_H \g_{H,C,\phi}^{(k)}$, where the sum is taken over all walls. Similar statement is true for $a^{(k),-}$.

\end{conj}

Finally, taking the projective limit over $k$ we define the elements $a^{\pm}$. These elements play a role of the initial data for the  sheaf $WCS_{\g,C}$  in a vector space.

\begin{rmk} Notice that for a fixed $V$ there is a homomorphism $\g_V\to \g^{\prime}$ given by the natural inclusions $\g_H\to \g^{\prime}$ for all $H\supset V$. Hence for a small open subset $U$ in $\Gamma_1\otimes \R$ we have a cocycle with values in $exp(\g_C)$. In plain terms, it is given by an element $g_{x_1,x_2}\in exp(\g_C)$ defined for  two points $x_1,x_2\in U$ which do not belong to any wall. Cocycle condition means that $g_{x_1,x_2}g_{x_2,x_3}=g_{x_1,x_3}$. We see that we have a picture similar to the one from Section 2. Hence we can use the above transformations in order to glue a Calabi-Yau manifold (possibly over a non-archimedean field) from  open coordinate charts. 

In the passage from $WCS$ to the glued manifold we lose some data. This point is clear when we look at the initial data. Indeed if assume the above Conjecture, we see that the direct sum of all $\g_H$ is ``bigger"  than $\g^{\prime}$. One can speculate that the whole $WCS$ contains the information sufficient for reconstruction of both mirror dual Calabi-Yau manifolds. The initial data can be thought of as association of a rational number to any triple $(\gamma,\mu,k)$, where $\gamma\in \Gamma_1-\Gamma_{1,0}, \mu\in \Gamma_2-\Gamma_{2,0}, k\in \Z_{\ge 1}$ and $\gamma,\mu$ are primitive vectors such that $(\gamma,\mu)=0$. The rational number is the coefficient of $x^{k\gamma}\partial_{\mu}\in \g_{H}, H=\mu^{\perp}$ in the initial data. Notice that the above conditions are symmetric with respect to the exchange of  $\Gamma_1$ and $\Gamma_2$.

\end{rmk}

\begin{rmk}
Notice that the Lie algebra $\g_{\Gamma_1,\Gamma_2, (\bullet,\bullet)}$ does not change if we replace $\Gamma_2$ by a sublattice of finite index.  Then at first sight it looks like that our WCS depends on $\Gamma_2\otimes \Q$ only. But one can recover a finer ``integer" structure of the WCS. At the level of the corresponding pronilpotent groups we consider subgroups generated by the elements $T_{\gamma,\mu}:=exp\left(\sum_{n\ge 1}{x^{n\gamma}\partial_{\mu}\over{n}}\right)$.  In the case $\Gamma_1=\Gamma_2$ and
 $(\bullet,\bullet)=\langle\bullet,\bullet\rangle$ being skew-symmetric, as in the Poisson case considered in Section 4, we have similar transformations $T_{\gamma}:=T_{\gamma,\gamma}=exp\left(\{Li_2(x^{\gamma}),\bullet\}\right)$ (cf. Remark 2.3.1).

\end{rmk}

\subsection{Pair of local systems of lattices from non-archimedean point of view}

Recall the discussion in  Section 6.2 of  the geometry of the base of the real integrable system which appears in SYZ picture of Mirror Symmetry.

Here we would like to recall the origin of the integrable system following [KoSo6] and [KoSo2].
It can be approached  either in the framework  of Gromov-Hausdorff collapse or using the language of non-archimedean geometry of Berkovich.
In the former approach we have a family $X_t, t\to 0$ of maximally degenerate polarized complex Calabi-Yau manifolds. It was conjectured in [KoSo6] that for sufficiently small $t$ the manifold $X_t$ contains an open subset $X_t^{\prime}$ which is in Gromov-Hausdorff metric close to the total space of the real integrable system $\pi_t: X_t^{\prime}\to B^0$  over an open smooth manifold $B^0\subset B$ of some metric space $B$. The latter is the Gromov-Hausdorff limit of $X_t, t\to 0$ (see [KoSo6] for details). The restriction of the metric to $B^0$ is a smooth Riemannian metric, which  locally given by the matrix of second derivatives of a convex function $\Phi$ on $B^0$ which satisfies the real Monge-Amp\`ere equation $det(\partial^2 \Phi/\partial x_i \partial x_j)=const$.

For any point  $b\in B^0$ we can define the lattice $\Gamma_{1,b}:=H_2(X_t, \pi_t^{-1}(b),\Z)$ (the latter stabilizes as $t\to 0$ along a ray). This family of lattice gives rise to a local system $\underline{\Gamma}_1\to B^0$.  It can be extended to a local system on $S^1_t\times B^0$.

It was conjectured in [KoSo6] (and hence was assumed in Section 6.2) that $codim (B-B^0)\ge 2$. Typically $B$ is a topological manifold homeomorphic to a sphere, complex projective space or a torus.

The approach via non-archimedean geometry
gives rise to the non-archimedean integrable system $\pi: {\mathcal X}^{an}\to B$ (see [KoSo2]). Here ${\mathcal X}^{an}$ is the compact polarized  analytic manifold (in the framework of Berkovich theory) over a non-archimedean field $K$ (typically $K=\C((t))$) corresponding to the collapsing family $X_t$, and $B$ is a PL space called skeleton of  ${\mathcal X}^{an}$. \footnote{In a recent paper [MusNic] the notion of the skeleton was generalized to all $K$, including the case of mixed and positive characteristic.} 
Furthermore,  $n:=dim_{\R}B=dim_{K}{\mathcal X}^{an}$, and the projection $\pi$ locally over $B^0$ looks as the map $(K^{\ast})^n\to \R^n$ given by $(z_1,...,z_n)\mapsto (log|z_1|,...,log|z_n|)$. Conjecturally the skeleton coincides with the Gromov-Hausdorff limit, hence the same notation.
Thus  $B$ carries a (singular) $\Z$-affine structure, which is non-singular on the open smooth manifold $B^0$. For the ``discriminant locus" $B^{sing}=B-B^0$ we have (conjecturally) the condition $dim(B^{sing})\le n-2$. 

We are going to assume that the singularity of the singular integral affine structure on $B$ is of the $A_k, k\ge 1$ type  at $B^{sing}$. Here the discriminant $B^{sing}$ is a closed subset of $B$, which contains an open dense topological submanifold $B_2^{sing}$ such that $dim(B_2^{sing})=n-2$ and $dim(B^{sing}-B_2^{sing})\le n-3$ and $B$ is a PL manifold near any point of $B_2^{sing}$. 
Furthermore $B_2^{sing}$ locally looks  as a topological submanifold in $\R^n$ given in the standard coordinates $(x_1,...,x_n)$ by the equations
$$x_1=f(x_3,...,x_n), x_2=0,$$
where $f(x_3,...,x_n)$ is a continuous function.

The integral affine structure on $B-B_2^{sing}$ in this local model coincides with the standard one on the open set $U_+$ which is the complement to the closed set  $x_1\ge f(x_3,...,x_n), x_2=0$.
On the open subset $U_-$ which is the complement to the closed subset 
 $x_1\le f(x_3,...,x_n), x_2=0$ the integral affine structure is the standard one in the coordinates $x_1^{\prime}=x_1+k\cdot max\{0,x_2\}, x_2^{\prime}=x_2,...,x_n^{\prime}=x_n$.

One can see that $B_2^{sing}$ belongs to a canonical germ of a hypersurface which is (in the affine structure)  an integer $(n-1)$-dimensional hyperplane  (outside of $B^{sing}$) endowed with a locally constant integer vector field. In the local picture the hyperplane is given by the equation $x_2=0$, and the vector field is given by $k\cdot sign(f(x_3,...,x_n)-x_1)\partial/\partial x_1$.

Recall that the approach  with collapse gives a local system of lattices $\underline{\Gamma}_1\to B^0$. In the non-archimedean picture we have $\underline{\Gamma}_{1,b}=H_{2,Betti}({\mathcal X}^{an},\pi^{-1}(b),\Z), b\in B^0$, where $H_{i,Betti}$ denotes properly defined Betti homology of the analytic space ${\mathcal X}^{an}$ over the field $K=\C((t))$.

There is a natural projection $p: \underline{\Gamma}_1\to T^{\Z}$, where $T^{\Z}:=T^{\Z}_{B^0}\subset T_{B^0}$ is the locally covariant lattice which defines the $\Z$-affine structure on $B^0$. Denoting $\underline{\Gamma}_{1,0}=Ker(p)$ we obtain an exact sequence of lattices 
$$0\to \underline{\Gamma}_{1,0}\to \underline{\Gamma}_1\to T^{\Z}.$$

In what follows we will assume that the local system $\underline{\Gamma}_{1,0}$ is trivial (this condition is automatically satisfied in most of examples).

\begin{rmk}

The above geometry (including the integrable systems $\pi_t: X_t^{\prime}\to B^0$) arises also in the situation when in the non-archimedean integrable system  $\pi: {\mathcal X}^{an}\to B$ the total space ${\mathcal X}^{an}$ is non-compact but $\pi$ is proper.  In terms of Gromov-Hausdorff collapse one can think of a family $X_t, t\to 0$ of non-compact complex manifolds endowed with complete K\"ahler metrics. The limiting metric on $B^0$ is given by the second derivatives of a convex function but this time not necessarily obeying the real Monge-Amp\`ere equation.

\end{rmk}

\begin{rmk} Let us assume that in the previous Remark the K\"ahler forms $\omega_t$ on $X_t$ have  homology classes $(log|t|)^{-1}\beta$, where $\beta$ belong to the image of $H^2(X_t,\Z)$ in $H^2(X_t,\R)$ (this image does not depend on $t$ when $t$ is sufficiently small. Let us assume that for any $b\in B^0$ the homomorphism  $H_1(\pi_t^{-1}(b),\Z)\to H_1(X_t,\Z)$ is equal to zero (e.g. it is sufficient to assume that  $H_1(X_t,\Z)=0$). Then we can define a local system $\underline{\Gamma}_1\to B^0$ in the following way. Its fiber over $b\in B^0$ is the set of pairs $(\gamma,v)$, where $\gamma\in T^{\Z}_b=H_1(\pi_t^{-1}(b),\Z)\simeq \Z^n, n=dim_{\C}X_t=dim_{\R}B^0$, and $v\in \R$ is an element of $\Z$-torsor $\Z+lim_{t\to 0}\int_{\delta_t}(log|t|)^{-1}\omega_t$. Here $\delta_t$ is a  $2$-chain in $X_t$ with the boundary in $\pi_t^{-1}(b)$ such that the boundary $\partial\delta_t$ is $C^1$-close to a closed geodesic in the torus $\pi_t^{-1}(b)$ representing homology class $\gamma$.
The existence of the limit follows from the condition that $\omega_t$ is close to a semiflat metric, for $|t|\ll 1$.

Then we have a short exact sequence
$$0\to \underline{\Gamma}_{1,0}=\underline{\Z}_{B^0}\to \underline{\Gamma}_1\to T^{\Z}\to 0.$$

\end{rmk}

Now we discuss the origin of the local system $\GGamma_2$. Let us first assume that the non-archimedean field $K$ carries a discrete valuation $val(K^{\ast})=\Z\subset \R$. In this case the base $B^0$ carries a sheaf of affine functions with {\it integer coefficients}. Then we define $\GGamma_2$ to be this sheaf. We have a short exact sequence
$$0\to \underline{\Gamma}_{2,0}=\underline{\Z}_{B^0}\to \underline{\Gamma}_2\to (T^{\Z})^{\ast}\to 0.$$
Then we have the natural pairing $\GGamma_1\otimes \GGamma_2\to   T^{\Z}\otimes (T^{\Z})^{\ast}\to \underline{\Z}_{B^0}$. In general one should consider  Calabi-Yau manifolds over the field $\C((t_1,....,t_m)), m=rk\, \GGamma_{2,0}$. This can be thought of as a family of $n$-dimensional bases $B^0$ endowed with $\Z$-affine structure and parametrized (locally) by a domain in $\Gamma_{2,0}^{\vee}\otimes \R$. The total space of this family of bases can be identified locally with a domain in $\Gamma_{2}^{\vee}\otimes \R=\Gamma_{2,\R}^{\ast}$ (cf. Section 4.6). Therefore we can speak about the WCS on the total space of this family.

Recall that according to the general philosophy recalled in Section 6.2 in order to construct the mirror dual to the Fukaya category of the Calabi-Yau manifold near the cusp, we should count holomorphic discs with boundaries on fibers of the SYZ fibration. From the non-archimedean point of view such discs become tropical trees in $n$-dimensional bases $B^0$ depending on parameters in $\Gamma_{2,0}^{\vee}\otimes \R$. In the next subsection we are going to discuss such trees for a fixed value of the parameter in $\Gamma_{2,0}^{\vee}\otimes \R$.

\subsection{Tropical trees, Finiteness Assumption and WCS}

We keep the notation of the previous subsection. In particular $B$ denotes the base of a non-archimedean integrable system. It is a PL space. We do not assume that the valuation on the non-archimedean field $K$ is integer. This means that we fix a $\Z$-affine structure on $B^0$ corresponding to an arbitrary (not necessarily integer) point in $\Gamma_{2,0}^{\vee}\otimes \R$. In what follows we will not utilize $\GGamma_2$. The reader should keep in mind that our objects depend on parameters from $\Gamma_{2,0}^{\vee}\otimes \R$.

\begin{defn} A tropical tree in $B$ is an oriented toward the root metrized tree $T$ (see Section 3.2) endowed with a continuous map $f:T-\{Tail\,vertices\}\to B^0$ together with a continuous lift $f^{\prime}: T-\{Vertices\}\to tot(\underline{\Gamma}_1)-tot(\underline{\Gamma}_{1,0})$ such that each edge lifts to a piece of the trajectory of the attractor flow
$$\dot{b}=\iota(\gamma), \dot{\gamma}=0, (b,\gamma)\in tot(\underline{\Gamma}_1)-tot(\underline{\Gamma}_{1,0})\subset tot(\underline{\Gamma}_1\otimes \R).$$

We assume that at each internal vertex $v$ we have the balancing condition $\sum_i\gamma^{out}=\gamma^{in}$ (cf. Definition 3.2.1), and all  $\gamma_i^{out}$ are pairwise distinct and there are $i_1, i_2$ such that $\gamma_{i_1}^{out}$ is not parallel to $\gamma_{i_2}^{out}$.

Furthermore we assume that the tail vertices belong to $B_2^{sing}$ and the germs of edges near tail vertices belong to the above-described canonical hypersurface and the speeds of the $f^{\prime}$-lifts of the tail edges are proportional (with minus sign) to the  canonical locally constant vector field described above.

\end{defn}

 Next we are going to discuss the analogs of the Finiteness Assumption (cf. Section 9.3)   as well as the existence of strict convex cones (cf. Conjecture 4.6.5). The idea is to use ``tropical metrics" on $B^0$ which are non-singular as well as limits of such, which can also contain $\delta$-functions supported on ``tropical effective divisors". We warn the reader that the conditions discussed below are still not sufficient for the finiteness of the number of trees (which is a tropical analog of Gromov compactness for pseudoholomorphic curves). Our conditions guarantee the existence of strict convex cones, boundedness of lengths and finiteness of the number of tails of  tropical trees with the given generic root and the velocity of the root edge. One has to put some extra constraints on the behavior of the affine structure or the ``tropical metric" at singularities of codimension $\ge 3$ in order to achieve the finiteness. Hopefully it can be done. We expect that such conditions are implicit in [GroSie1,2], where the  procedure for construction of the mirror dual Calabi-Yau gives  a WCS in our sense. 

Recall that in complex geometry it is natural to consider non-negative $(1,1)$-currents as limits of sequences of K\"ahler metrics on a given manifold. There is an analog of this notion in the non-archimedean geometry. 
More precisely, there a sheaf of monoids ${\mathcal Psh}\subset \underline{Cont}$ on ${\mathcal X}^{an}$ of continuous  plurisubharmonic functions (see [BouFaJo1]). The sheaf of abelian groups ${\mathcal Ph}$ is defined as  ${\mathcal Psh}\cap (-{\mathcal Psh})$. The quotient ${\mathcal Psh}/{\mathcal Ph}$ is by definition the sheaf of non-negative $(1,1)$-currents. One can define similar sheaves on $B$ in the following way. For example, a continuous of plurisubharmonic function on an open subset of $B$ is a  continuous function on this subset such that its pull-back to ${\mathcal X}^{an}$ is a plurisubharmonic functions.
A germ of a non-negative $(1,1)$-current at a point $b\in B^0\subset B$ is the same as  a germ of a convex (in the affine structure) function modulo a germ of affine function.

Having a non-negative $(1,1)$-current $\phi$ on $B$ and a tropical tree $T$ on $B$ we can define the integral $\int_T\phi\in \R_{\ge 0}$, assuming that $\phi$ is sufficiently regular (e.g. smooth) near the root $b_0$ of $T$. The reader can think of the integral as the limit of integrals of non-negative $(1,1)$-currents over holomorphic discs. Furthermore one can show that this integral depends only on the velocity of $T$ at $b_0$ and hence gives rise to a linear functional on $\underline{\Gamma}_{1,b_0}$.
Therefore $\phi$ defines a half-space where the velocities of the tropical trees at roots must belong (at least at the points where $\phi$ is smooth and strictly convex). Moreover, for any such point $b_0$ we can consider a small variation  $\phi_{\varepsilon}$, e.g. by taking 
$\phi+\varepsilon \eta$, where $\eta$ is a an arbitrary smooth function supported in a small neighborhood of $b_0$. Similarly to Section 9.4  we can intersect the corresponding half-spaces for all $\varepsilon>0$ and obtain a strict convex cone in $T^{\Z}_{b_0}\otimes \R$. But we need a strict convex cone in $\underline{\Gamma}_{1,b_0}\otimes \R$. This can be achieved by considering variations of $\phi$ of more general type which change its cohomology class $[\phi]\in H^2_{Betti}({\mathcal X}^{an})\otimes \R$.

Therefore if there exists  a non-negative $(1,1)$-current $\phi$ which is smooth and strictly convex at  any point of $B^0$, then there exists a family of closed strict convex cones $C_b\subset \underline{\Gamma}_{1,b}\otimes \R, b\in B^0$ such that the velocities of tropical trees rooted at $b$ belong to $C_b$ (cf. Conjecture 4.5.3). We call {\it Support Property} the existence of such family of cones. 

Assuming that such $\phi$ does exist  one can show that the total length of a tropical tree with given root end the velocity of the root edge is bounded. Here the length is measured with the respect to some auxiliary  Riemannian metric on $B^0$ obtained from $\phi$ by local considerations. 
This argument is not sufficient to guarantee the analog of Finiteness Assumption, i.e. the finiteness of the number of tropical trees with given root $b\in B^0$ and the velocity of the root edge $\gamma\in \underline{\Gamma}_{1,b}-\underline{\Gamma}_{1,0,b}$. The finiteness can fail if there exists an infinite sequence of tropical trees with given $(b,\gamma)$ and increasing numbers of tail edges. Hence we need to find a restriction which guarantees the boundedness of the number of tail edges. In order to achieve that we will need a smaller class of non-negative $(1,1)$-currents. Namely they will be smooth and strictly convex on $B^0$ and satisfy the property that  $\int_l\phi\ge 1$ for any arbitrarily small piece $l$ of the trajectory  of the canonical locally-constant vector field, such that $l$ hits the discriminant. \footnote{In order to illustrate the latter condition consider delta-currents corresponding to compact curves sitting at the preimage of $B_2^{sing}$ in the total space of the integrable system. The integral of such a current over a holomorphic disc is bounded from below by the intersection index. Our assumption is a ``tropical" version of this fact.} For such currents $\int_T\phi$ gives an upper bound for the number of tail edges of any tropical tree $T$. Together with the above-discussed upper bound for the length of $T$ it will imply the finiteness of the number of tropical trees with fixed $(b,\gamma)$.

First let us make a simplifying assumption that $B$ carries  a $\Z PL$-structure which is compatible with $\Z$-affine structure on $B^0$ and such that $B^{sing}=B-B^0$ is a $\Z PL$ subset of $B$. In terms of the local model near $A_k$-singularities discussed above this is equivalent to the fact that the function $f(x_3,...,x_n)$ is piecewise-linear with rational slopes. In such a situation we suggest to take as $\phi$ the $(1,1)$-current associated with the Gromov-Hausdorff limit $S$ of a family of ample effective divisors  $S_t\subset X_t, t\to 0$ (we call $S$ tropical effective divisor). Let us require $S$ contains $B^{sing}$, does not contain germs of the canonical hyperplanes $x_2=0$ (see the description of the local model), and the intersection number of  $S$ with any germ of a trajectory of the canonical locally-constant vector field near a generic point point of $B_2^{sing}$ is greater or equal than one. More precisely such a germ can be understood as the projection to $B$ of a non-archimedean analytic disc $D$ in ${\mathcal X}^{an}$. In the collapse picture it is represented by a family of complex holomorphic discs $D_t\subset X_t$. Then the above intersection number is defined as the usual intersection number $S_t\cdot  D_t$ for 
$|t|\ll1$.

\begin{exa} Let $n=2, k=1$. The tropical effective divisor $S$ is the union of three rays: $S_1=\{x_1=0, x_2\ge 0\}$, $S_2=\{x_1=0, x_2\le 0\}$ and $S_3=\{x_1+x_2=0, x_2\ge 0\}$.
The rays $S_1,S_3$ are taken with the multiplicity $+1$, while the ray $S_2$ is taken with the multiplicity $+2$. Here we use the focus-focus $\Z$-affine structure on $\R^2-\{(0,0)\}$ which is the standard affine structure on $\R^2-\{x_2=0,x_1\ge 0\}$ and has as local affine coordinates
$(x_1+max(0,x_2),x_2)$ near the ray $x_2=0,x_1\ge 0$.

Notice that   $S_1\cup S_2$ (union as sets) is the limit as $\varepsilon \to 0$ of the family of straight lines $x_1=-\varepsilon, x_2\in \R$ in $B^0=\R^2-\{(0,0)\}$. Similarly $S_2\cup S_3$ is the limit of the family of straight lines (in the $\Z$-affine structure) given by $(S_2+(0,\varepsilon))\cup (S_3+(0,\varepsilon))$. There are two types of germs of tropical discs: $D_{\pm}:=D_{\pm,\delta}$.
Here $D_-$ is given by $\{-\delta<x_1\le 0, x_2=0\}$, and $D_+$ is given by $\{0\le x_1<\delta, x_2=0\}$, where $0<\delta \ll 1$. Then one can easily check that $D_-\cdot (S_1+S_2)=1, D_-\cdot (S_2+S_3)=0$. Similarly, $D_+\cdot (S_1+S_2)=0, D_+\cdot (S_2+S_3)=1$. Therefore $D_-\cdot S=D_+\cdot S=1$, where $S=S_1+2S_2+S_3$.

\end{exa}

The following figure illustrates the Example.

\vspace{2mm}

\centerline{\epsfbox{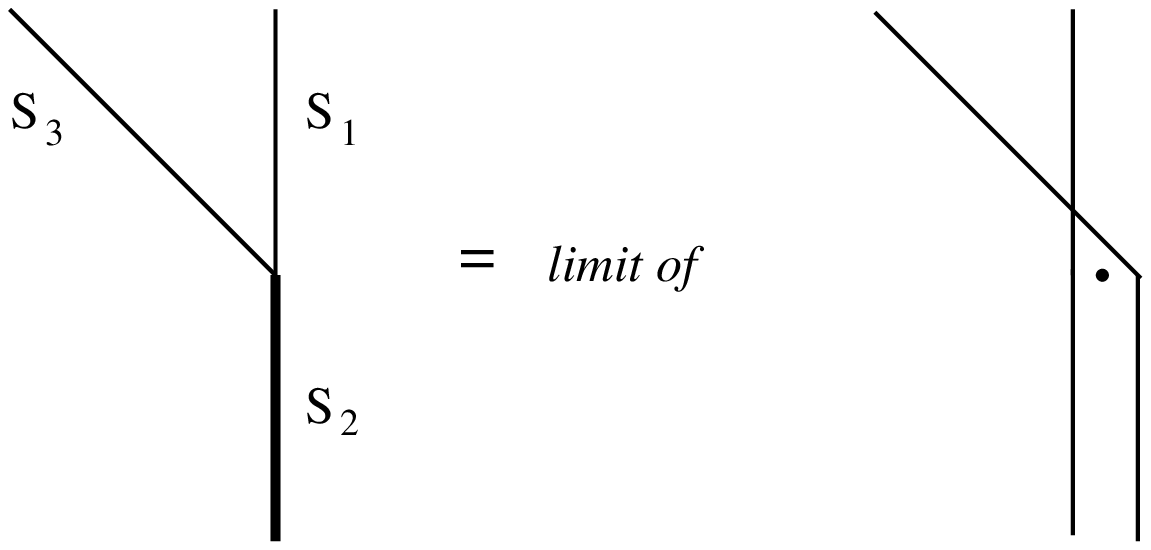}}

\vspace{2mm}

If we have a tropical effective divisor $S$ the above constraints on the intersection numbers with germs of trajectories and $\phi$ is the corresponding non-negative $(1,1)$-current then the integral $\int_T\phi$ is well-defined for any tropical tree $T$ such the root of $T$ does not belong to $S$. Furthermore this integral gives an upper bound on the number of tails of $T$. In order to have an upper bound for {\it all} tropical trees it is sufficient to find a collection $(S_{\alpha})$ of effective tropical divisors satisfying the above constraint and such that $\cap_{\alpha}S_{\alpha}\subset B^{sing}$. It is easy to find such a collection using the ampleness of the corresponding divisors in $X_t, t\to 0$.

We conclude that we can achieve (under the appropriate conditions) the Finiteness Assumption as well as the Support Property. Then using the analog of the procedure described in the Section 4 we construct the WCS. 
After that, using a non-negative $(1,1)$-current which is smooth and strictly positive on $B^0$, we can endow $B^0$ with the dual $\Z$-affine structure. It can understood as the base of the canonical non-archimedean integrable system which is glued from tube domains in the non-archimedean torus $(\C((t))^{\ast})^n$ (see [KoSo2]). The walls of the WCS become curved in the new affine structure. The transformations corresponding to walls can be used in order to modify the total space of the above-mentioned canonical integrable system. As the result, we obtain a non-archimedean integrable system which can be extended to the integrable system with the base $B$. Its total space can be thought of as the analytic space corresponding to the mirror dual family $X_t^{\vee}, t\to 0$.

\begin{rmk} The approach of  [GroSie1,2] gives an example of WCS discussed above. In their case $B^{sing}$ as a $\Z PL$ subset of $B$. Their notion of ``slab" corresponds to the notion of tropical effective divisor discussed above. What we call WCS corresponds to the notion of ``scattering diagram" in the loc.cit.

\end{rmk}


\section{Appendix}

\subsection{Canonical $B$-field}

Let us consider a fibration $\pi:X^0\to B^0$ whose fibers are compact tori, endowed with a section. Denote by $\underline \Gamma$ the local 
system $\left(R^1\pi_*(\mathbb{Z})\right)^\vee$ of first homology groups of fibers. We assume that  $\underline \Gamma$ is endowed with a skew-symmetric pairing
$\langle\bullet,\bullet\rangle$, possibly degenerate. The goal of this subsection is to define  a canonical cohomology class in $H^2(X^0,\Z/2\Z)$ naturally associated with the pairing.

First, let us consider an individual fiber $\Gamma$.
Skew-symmetric pairing on $\Gamma$ gives a {\it symmetric} pairing on $\Gamma\otimes \Z/2\Z$. Hence we can consider the group $V$ of polynomials $P$ of degree at most $2$ on
$\Z/2\Z$-vector space $\Gamma\otimes \Z/2\Z$ such that $P(0)=0$ and the bilinear form $(x,y)\mapsto P(x+y)-P(x)-P(y)$ is proportional (with the factor in
$\Z/2\Z$) to the form $(x,y)\mapsto \langle x,y\rangle\mod 2$.

We have a short exact sequence
$$0\to Hom(\Gamma,\Z/2\Z)\to V\to \Z/2\Z\to 0,$$
hence a dual sequence
$$0\to \Z/2\Z\to V^{\vee}\to \Gamma\otimes \Z/2\Z\to 0.$$
We have the natural map $\Gamma\to \Gamma\otimes \Z/2\Z$. Let $W$ be the fiber product 
$$W=lim(\Gamma\to \Z/2\Z\leftarrow V^{\vee}).$$

Then we have a short exact sequence
$$0\to \Z/2\Z\to W\to \Gamma\to 0$$

Passing to the classifying spaces we obtain a fibration
over the torus $K(\Gamma,1)$ with the fiber being the Eilenberg-MacLane space $K(\Z/2\Z,1)$.
Going from the local model to the global picture we obtain a fibration with the fiber $K(\Z/2\Z,1)$ over $X^0$. Its characteristic class is 
the desired class in 
$H^2(X^0,\Z/2\Z)$.  

\vspace{3mm}

{\bf References}

\vspace{2mm}

[AbSe] M. Abouzaid, P. Seidel, An open string analogue of Viterbo functoriality, Geom. Topol.
14 (2010), no. 2, 627–718, arXiv:0712.3177.
\vspace{2mm}

[AlManPerPi] S. Alexandrov, J. Manschot, D. Persson, B. Pioline, Quantum hypermultiplet moduli spaces in N=2 string vacua: a review,   arXiv:1304.0766.
\vspace{2mm}

\vspace{2mm} [BanMan] R. Bandiera, M. Manetti, 
On coisotropic deformations of holomorphic submanifolds, arXiv:1301.6000.

\vspace{2mm}

[Bo1] P. Boalch, Geometry and braiding of Stokes data; fission and wild character varieties, arXiv:1111.6228.

\vspace{2mm}

[BouFaJo1] S. Boucksom, C. Favre, M. Jonsson, Singular semipositive metrics in non-archimedean geometry, arXiv:1201.0187.

\vspace{2mm}

[BouFaJo2] S. Boucksom, C. Favre, M. Jonsson, Solution to a non-archimedean Monge-Amp\`ere equation, arXiv:1201.0188.

\vspace{2mm}
[ChLaLeu] K. Chan, S. Lau, N. Leung, SYZ mirror symmetry for toric Calabi-Yau manifolds,
arXiv:1006.3830.

\vspace{2mm}

[DelPh] E. Delabaere, F. Pham, Resurgent methods in semi-classical analysis, Annales IHP, Section A, 71:1, 1999, 1-94.
\vspace{2mm}

[Den] F. Denef, Supergravity flows and D-brane stability, arXiv:hep-th/0005049.

\vspace{2mm}

[DenGrRa] F. Denef, B. Greene, M. Raugas, Split attractor flows and the spectrum of BPS D-branes on the quintic, arXiv:hep-th/0101135.

\vspace{2mm}

[DenMo] F. Denef, G. Moore, Split states, Entropy enigmas, holes and halos, arXiv:hep-th/0702146.

\vspace{2mm}

[DiDoPa]  D. Diaconescu, R. Donagi, T. Pantev, Intermediate Jacobians and ADE Hitchin Systems, arXiv:hep-th/0607159. 

\vspace{2mm}
[Don] R. Donagi, Seiberg-Witten integrable systems,  arXiv:alg-geom/9705010.
\vspace{2mm}

[DonMar] R. Donagi, E. Markman, Cubics, integrable systems, and Calabi-Yau threefolds, arXiv:alg-geom/9408004.
\vspace{2mm}

[EynOr] B. Eynard, N. Orantin, Invariants of algebraic curves and topological expansions, arXiv: math-ph/0702045.

\vspace{2mm}

[FOOO] K. Fukaya, Y-G. Oh, H.Ohta, K. Ono, Lagrangian intersection Floer theory,Studies in Advanced Mathematics, 2010.

\vspace{2mm}

[FoGo1] V. Fock, A. Goncharov, Cluster ensembles, quantization and the dilogarithm, Invent.
Math. 175 (2009), no. 2, 223-286, see also  arXiv:math/0311245.

\vspace{2mm}

[GaMoNe1] D. Gaiotto, G. Moore, A. Neitzke,  Four-dimensional wall-crossing via three-dimensional field theory, arXiv:0807.4723.

\vspace{2mm}
[GaMoNe2] D. Gaiotto, G. Moore, A. Neitzke,  Wall-crossing, Hitchin Systems, and the WKB Approximation, arXiv:0907.3987.

\vspace{2mm}

[GaMoNe3] D. Gaiotto, G. Moore, A. Neitzke,  Wall-crossing in coupled 2d-4d systems, arXiv:1103.2598.
\vspace{2mm}

[GaMoNe4] D. Gaiotto, G. Moore, A. Neitzke, Framed BPS states, arXiv:1006.0146.

\vspace{2mm} 

[GoKen] A. Goncharov, R. Kenyon, Dimers and cluster integrable systems, arXiv:1107.5588.

\vspace{2mm} 

[GroHaKe] M. Gross, P. Hacking, S. Keel, Mirror symmetry for log Calabi-Yau surfaces I, arXiv:1106.4977.

\vspace{2mm}

[GroSie1] M. Gross, B. Siebert, Mirror Symmetry via Logarithmic Degeneration Data I,
arXiv:math/0309070.

\vspace{2mm}

[GroSie2] M. Gross, B. Siebert, From real affine geometry to complex geometry, arXiv:math/0703822, published in Ann. Math., 174:3, p. 1301-1428, 2011.

\vspace{2mm}

[GroSie3] M. Gross, B. Siebert, Theta functions and mirror symmetry,
arXiv:1204.1991.

\vspace{2mm}

[KaKoPa] L. Katzarkov, M. Kontsevich, T. Pantev, Hodge theoretic aspects of mirror symmetry, arXiv:0806.0107.
\vspace{2mm}

[Ko1] M. Kontsevich, Holonomic $D$-modules and positive characteristic, arXiv:1010.2908.

\vspace{2mm}

[Ko2] M. Kontsevich, Deformation quantization of Poisson manifolds, I, arXiv:q-alg/9709040. 

\vspace{2mm}

[KoSo1] M. Kontsevich, Y. Soibelman, Stability structures, motivic Donaldson-Thomas invariants and cluster transformations, arXiv:0811.2435.

\vspace{2mm}

[KoSo2] M. Kontsevich, Y. Soibelman, Affine structures and non-archimedean analytic spaces, math.AG/0406564.

\vspace{2mm}
[KoSo3] M. Kontsevich, Y. Soibelman, Notes on A-infinity algebras, A-infinity categories and non-commutative geometry. I, math.RA/0606241.

\vspace{2mm}

[KoSo4] M. Kontsevich, Y. Soibelman,Deformations of algebras over operads and Deligne's conjecture, arXiv:math/0001151.
\vspace{2mm}

[KoSo5] M. Kontsevich, Y. Soibelman, Cohomological Hall algebra, exponential Hodge structures and motivic Donaldson-Thomas invariants, arXiv:1006.2706.

\vspace{2mm}

[KoSo6] M. Kontsevich, Y. Soibelman,
Homological mirror symmetry and torus fibrations, arXiv:math/0011041.

\vspace{2mm}

[KoSo7] M. Kontsevich, Y. Soibelman, Motivic Donaldson-Thomas invariants: summary of results, arXiv:0910.4315. 

\vspace{2mm}
[KoSo8] M. Kontsevich, Y. Soibelman, Deformation theory I, draft of the book, available at www.math.ksu.edu/$\sim$ soibel.

\vspace{2mm}
[Kos] V. Kostov, On the Deligne-Simpson problem, ArXiv:math/0011013.

\vspace{2mm}
[Mal] B. Malgrange, \'Equations differenti\'elles \`a coefficients polynomiaux, 1991, Birkhauser. 

\vspace{2mm}
[Mor1] D. Morrison, Geometric Aspects of Mirror Symmetry, arXiv:math/0007090.

\vspace{2mm}

[MusNic] M. Mustata, J. Nicaise, Weight function on non-archimedean analytic spaces and the Kontsevich-Soibelman skeleton, arXiv:1212.6328.

\vspace{2mm}

[Ra] Z. Ran, Lifting of cohomology and unobstructedness of certain holomorphic maps, arXiv:math/9201267 [pdf, ps, other]
 
 \vspace{2mm}
 
[Vo] A. Voros, The return of the quartic oscillator (the complex WKB method), Annales Institut H. Poincar\'e, 29:3, 1983, 211-338.

\vspace{2mm}

Addresses:

M.K.: IHES, 35 route de Chartres, F-91440, France, {maxim@ihes.fr}

Y.S.: Department of Mathematics, KSU, Manhattan, KS 66506, USA, {soibel@math.ksu.edu}

\end{document}